\documentclass[a4paper,11pt]{amsart}
\usepackage{amssymb,amsmath,amsthm,stmaryrd,bm,enumitem,hyperref, xcolor} %graphicx 

\def\res{\hbox{ {\vrule height .3cm}{\leaders\hrule\hskip.3cm}}\hskip5.0\mu}

\newcommand\beqn{\begin{equation}}
\newcommand\eeqn{\end{equation}}
\newcommand\beqny{\begin{eqnarray}}
\newcommand\eeqny{\end{eqnarray}}
\newcommand\beqnyn{\begin{eqnarray*}}
\newcommand\eeqnyn{\end{eqnarray*}}

\usepackage{amsmath,amsfonts,amssymb,amsthm}
\newtheorem{theorem}{Theorem}[section]
\newtheorem{lemma}[theorem]{Lemma}

\newtheorem{corollary}[theorem]{Corollary}
\newtheorem{definition}[theorem]{Definition}
\newtheorem{remark}[theorem]{Remark}

\def\r{\rho}

\numberwithin{equation}{section}

\newcommand{\op}[1]{\operatorname{\text{\rm #1}}}  

\title[Analysis of singularities of area-minimizing currents, Part I]{Analysis~of~singularities~of~area-minimizing~currents, Part~I: planar frequency, branch points of rapid decay, and weak locally uniform approximation}

\author{Brian Krummel and Neshan Wickramasekera}

\begin{document}

\begin{abstract}
This paper is the first in a five-part series developing a new geometric framework for $n$-dimensional locally area-minimizing rectifiable currents, initially in ${\mathbb R}^{n+m}$ ($m \geq 2$) and extending to general Riemannian ambient spaces in Part~V. The foundational 1983 work of Almgren---later clarified by De~Lellis and Spadaro---established the optimal upper bound $n-2$ for the Hausdorff dimension of the singular set. Our program shifts the focus towards the local structure of the current. The work provides a generalization of the classical two-dimensional regularity theory established by White, Chang, and Micallef--White, and yields a geometrically more direct proof of Almgren’s theorem together with new asymptotic and local topological conclusions.

In this paper we introduce the \emph{planar frequency}, an intrinsic frequency function for the current defined relative to a fixed plane and base point. In contrast with the frequency functions used in the classical Almgren program, the planar frequency is defined directly through geometric integrals with respect to the current's weight measure; 
it does not require 
iterative construction of center manifolds (a foundational idea used in the classical approach), nor a multi-valued parameterization. 
Planar frequency 
is $\leq 1$ on stationary cones,  and we prove that it is approximately monotone 
subject to a decay condition on the current's $L^{2}$-height relative to 
the chosen plane and base point.

Planar frequency forces a geometric dichotomy at singularities. 
The current either satisfies at 
small scales a weak 
non-planar approximation property where it is significantly better  approximated by a (scale-dependent) non-planar cone than by any plane, or this fails and the current decays rapidly to a plane.
Planar-frequency monotonicity prevents 
oscillation between these behaviors under scaling, decomposing fixed-integer-density singularities into  
disjoint sets ${\mathcal S}$ and ${\mathcal B}$:
\begin{itemize}
\item ${\mathcal S}$ is  
where the weak non-planar approximation holds in a locally uniform, quantified form.
\item points in ${\mathcal B}$ are rapid-decay branch points, where the current converges to a unique tangent plane with rate 
$o(\rho^{1 + \alpha})$ as the scale $\rho \to 0,$ for some locally uniform constant $\alpha >0$.
\end{itemize}
This division---a quantitative alternative to the classical qualitative branch/non-branch point split---forms the basis of our structural analysis. Additionally, it yields the singular set Hausdorff dimension more directly:  
for ${\mathcal S},$ 
the built-in quantitative non-degeneracy yields dimension control via the same reasoning as for non-branch points---even though ${\mathcal S}$ may contain branch points in addition to all (integer-density) non-branch points;
for ${\mathcal B}$, decay and planar frequency monotonicity streamline the 
use of the center manifold, ensuring this device is only needed in its canonical setting. 

The subsequent papers complete our structural study. Part II employs 
the weak non-planar approximation to prove  
${\mathcal H}^{n-2}$-a.e.\ tangent-cone uniqueness, ${\mathcal H}^{n-2}$-nullity of the branch set in ${\mathcal S}$   
and local $(n-2)$-rectifiability 
of fixed-density subsets of ${\mathcal S}.$
Parts~III~and~IV analyze the current 
near ${\mathcal B}$ by exploiting   
its uniform decay. Part~III treats points of planar frequency~$\neq 2$ without using center manifolds, relying instead on 
planar-frequency monotonicity. Part~IV addresses frequency-2 points 
using the center manifold, where we capitalize on a geometrically direct canonical relationship: 
fixed-density frequency-2 points 
locally lie in the critical nodal set of a \emph{single} approximating normal 
map over a single center manifold.
This long-known feature---not exhibited by classical iterative center-manifold constructions for arbitrary branch points---is central to our work in Part~IV.
Together, Parts~III~and~IV yield higher-order expansions as the primary analytic result, with geometric and topological corollaries: 
rectifiability conclusions for ${\mathcal B}$ and a frequency criterion 
for the current   
near a branch point 
to be a topological $n$-disk admitting a $C^{1, \mu}$ parameterization.

Broadly, this framework highlights geometric mechanisms underlying 
the higher-codimension regularity theory for area-minimizing currents.
By establishing intrinsic uniform a priori estimates as the primary organizing principle, the work aligns the higher-codimension theory more closely with the spirit of codimension 1 regularity and classical PDE. This shift facilitates a unified approach: 
recovering classical results---such as the dimension bound and weak stratification for the singular set---and yielding new   
structural conclusions for the current and its singularities.
\end{abstract}

\maketitle

\tableofcontents

\section{Introduction}
\subsection{Historical context and the main result of the present article} Let $T$ be an $n$-dimensional locally area minimising rectifiable current in 
an open subset of ${\mathbb R}^{n+m}$ 
(or more generally, in a sufficiently smooth $n+m$ dimensional Riemannian manifold).
The monumental 1983 work of Almgren, published posthumously as~\cite{Almgren}, established that ${\rm sing} \, T$, 
the interior singular set of $T$, has Hausdorff dimension at most $n-2$.  This dimension bound is optimal whenever the codimension $m \geq 2,$ as illustrated by  multiplicity 1 currents induced by certain holomorphic  
subvarieties of ${\mathbb C}^{n}$ (which are always locally area minimising). For instance, each of $T_{1} = \{(z, w) \, : \, z^{2} = w^{3}\} \subset {\mathbb C} \times {\mathbb C} \approx {\mathbb R}^{4}$ and 
$T_{2} = \{(z, w) \, : \, zw = 0\} \subset {\mathbb C} \times {\mathbb C} \approx {\mathbb R}^{4}$ is the support of a two (real) dimensional locally area minimising current in ${\mathbb R}^{4}$ with an isolated singularity at the origin. On the other hand if the codimension of $T$ is $1$, then ${\rm sing} \, T$ is empty if $n \leq 6$, and has Hausdorff dimension $\leq n-7$ if $n \geq 7$. This dimension bound, which is again sharp,  was the final outcome of a series of works, spanning approximately the decade 1960--1970, due to De~Giorgi (\cite{DG}),  Fleming (\cite{FW}), Almgren (\cite{A1}), Simons (\cite{Simons}) and Federer (\cite{Fed70}).  The major difference between the case of codimension $\geq 2$ and that of codimension 1, as far as singularities are concerned, is that  in codimension $\geq 2$ the current may have \emph{branch points}, i.e.\ non-immersed points of the current where there is at least one tangent cone supported on an $n$-dimensional plane (as in example $T_{1}$ above, which has an isolated branch point at $(0, 0)$). In codimension 1, De~Giorgi's work (\cite{DG}) implies that points of the current with planar tangent cones are embedded points.

Almgren's strategy in \cite{Almgren} for bounding the dimension of the singular set was to establish first $(n-2)$-dimensionality of the set of non-branch-point singularities (such as the point $(0, 0)$ in example $T_{2}$ above). It is easily seen that at a non-branch-point singularity, each tangent cone is a singular cone whose spine---i.e.\ the maximal subspace of translation invariance---is of dimension $\leq n-2$. Based on this fact, Almgren developed an elementary, very general argument (\cite[Corollary~2.27]{Almgren}) to show that the Hauadorff dimension of the set of non-branch-point singularities is at most $n-2$. Once this result was in place, it remained to bound the size of the branch set. This was by far the deepest and most involved part of \cite{Almgren}, for which Almgren developed a powerful set of techniques and ideas that have since been highly influential in a wide variety of geometric and PDE theoretic problems.  
More recently,  De~Lellis and Spadaro published Almgren's theory in more modern and concise language  (\cite{DeLSpa1},\cite{DeLSpa2},\cite{DeLSpa3}), providing technical streamlining of certain parts of the original version \cite{Almgren}, and generating a renewed interest in this profound work by making it more accessible.

In this and the subsequent papers \cite{KrumWicb}--\cite{KrumWice} we develop a new geometric framework for the analysis of interior singularities of locally area minimizing rectifiable currents of codimension $\geq 2$. This work is aided by analytic ingredients from our earlier work on multi-valued Dirichlet energy minimizing functions (\cite{KrumWic2}). Our program shifts the focus from studying the singular set itself to analyzing more broadly the local structure of the current near singularities. It achieves a three-fold outcome:
 \begin{itemize}
 \item first, it provides a generalization of the classical structure theory for  locally 
 area-minimizing two-dimensional rectifiable currents. This theory was established by the combined work of White (\cite{Whi83}), Chang (\cite{Chang}) and Micallef--White (\cite{MicWhi95}), which built upon an extended Almgren's theory in dimension 2; 
 \item second, it yields local structure results for the singular set: this includes rectifiability conclusions in general, and $C^{1, \alpha}$ regularity for the branch set under further conditions; 
 \item third, some of the new ingredients developed for this structural study yield Almgren's singular-set dimension bound in a more geometrically direct way.
 \end{itemize}
 
To maximize expository clarity and highlight the main new ideas, we shall assume in the present paper and  in \cite{KrumWicb}--\cite{KrumWicd} that the ambient space is an open subset 
of the Euclidean space ${\mathbb R}^{n+m}.$ Extensions to general Riemannian ambient spaces, which require technical adjustments to the Euclidean arguments, are discussed in \cite{KrumWice}.

 In the framework we develop, in place of the classical qualitative division of singularities (based on the tangent cone type) as branch and non-branch points, 
 the first step---accomplished in the present article---is to establish a quantitative decomposition of the singular set into two disjoint pieces. 
 
 In the classical Almgren program, substantial complications arise in the analysis of the branch set from the lack of an estimate at branch points giving decay towards a unique tangent plane. 
 Bypassing this uniqueness question, the classical approach proceeds by treating  all branch points equally and utilizing a technically demanding iterative construction of ``center manifolds'' for an arbitrary branch point. These center manifolds facilitate a blow-up procedure with sheet separation for the blow-up, ensuring persistence of singularities in passing from the current to a blow-up at a branch point; this, in turn, allows branch-set (Hausdorff) dimension control for the current via singular-set dimension bounds for the blow-ups. By contrast, for structural analysis the uniqueness-of-tangent-planes question is a central starting issue---and one whose resolution does not require sheet-separation; in fact, rapid merging of sheets is a favorable scenario for this. Thus one expects to be able to address this uniqueness question without involving center manifolds.
 
 In view of these considerations, it is natural to consider decomposing the singular set---specifically, singularities with a fixed (integer) density $q$---as the disjoint union of a set 
${\mathcal B}_{q}$ of branch points $Z$ where the scaled current $\eta_{Z, \rho \, \#} \, T$ \emph{decays rapidly} (in the $L^{2}$ sense) to a unique tangent plane  as $\rho \to 0$, and the complementary set ${\mathcal S}_{q}.$  It turns out that this is a considerably more effective way to proceed; once the rapid-decay requirement for ${\mathcal B}_{q}$ is chosen correctly, $T$ can be shown to satisfy a certain locally uniform weak approximation-by-non-planar-cones property at each point in ${\mathcal S}_{q}$ and each sufficiently small scale. This property lends itself well not only to bounding the size of ${\mathcal S}_{q}$ with very little effort---via essentially the same method as Almgren's argument for the set of non-branch-points---but also to performing an asymptotic analysis of the current at generic (i.e., ${\mathcal H}^{n-2}$ a.e.) points in ${\mathcal S}_{q}$. 
For ${\mathcal B}_{q}$, the existence of a unique tangent plane at each point, together with the accompanying \emph{locally uniform} decay estimate
built into its definition, offers two key advantages: it enables a more direct argument for efficient dimension control and provides a crucial foundation for {higher-order asymptotic analysis of the current at typical (i.e., ${\mathcal H}^{n-2}$-almost all)\ branch points.

This decomposition result (Theorem~\ref{branch and cones thm-intro} and Corollary~\ref{immediate-cor} below) is our main theorem here. It asserts, roughly speaking, the following: \emph{given an integer $q \geq 2$ and a numbers $\varepsilon, \beta \in (0, 1)$, there is a fixed $\alpha \in (0, 1),$ depending only on $n$, $m$, $q$, $\varepsilon$ and $\beta$ such that if ${\rm sing}_{q} \, T$ denotes the set of density $q$ singularities of $T$, then: 
\begin{itemize}
\item[{\rm (i)}] there is an ambient open set $V_{q, \varepsilon, \beta}$ with 
${\rm sing}_{q} \, T \subset V_{q, \varepsilon, \beta}$ and a relatively closed set ${\mathcal B}_{q, \varepsilon, \beta} \subset V_{q, \varepsilon, \beta}$ of branch points in ${\rm sing}_{q} \, T$ such that for each $Z \in {\mathcal B}_{q, \varepsilon, \beta}$ the scaled current $\eta_{Z, \rho \, \#} \, T$ converges to a multiplicity $q$ plane at a rate $O(\rho^{\alpha})$, and 
\item[{\rm (ii)}] for each 
$Z \in {\mathcal S}_{q,\varepsilon, \beta} = {\rm sing}_{q} \, T \setminus {\mathcal B}_{q, \varepsilon, \beta}$, the current $T$ satisfies the following property: for every point $Z^{\prime}$ with density $\geq q$ and sufficiently close to $Z,$ and every sufficiently small scale $\rho^{\prime}$ (depending on $Z$ but independent of $Z^{\prime}$), the scaled current $\eta_{Z^{\prime}, \rho^{\prime} \, \#} \, T$ is close to some non-planar cone ${\mathbf C}_{Z^{\prime}, \rho^{\prime}}$ and is significantly closer to ${\mathbf C}_{Z^{\prime}, \rho^{\prime}}$ than to any (single) plane, in a precise sense determined by $\varepsilon$ and $\beta$.
\end{itemize}} 

The set ${\mathcal B}_{q, \varepsilon, \beta}$ entirely consists of branch points by definition,  but note that the theorem allows that ${\mathcal S}_{q, \varepsilon, \beta}$ may also contain branch points. 
If a point $Z \in {\mathcal S}_{q, \varepsilon, \beta}$ is not a branch point, it is not difficult to see that for all sufficiently small $\rho$ (depending on $Z$), the scaled current $\eta_{Z, \rho \, \#} \, T$ is much closer to a non-planar cone than to any plane. The significance of the theorem is that this non-degeneracy property holds \emph{even if 
$Z \in {\mathcal S}_{q, \varepsilon, \beta}$ is a branch point}; 
additionally, there is a \emph{locally uniform} choice for how small the scale needs to be for the property to hold.

In proving this result, the main difficulty one has to overcome is ruling out the possibility that there may be a point $Z \in {\mathcal S}_{q, \varepsilon, \beta}$ such that 
$\eta_{Z, \sigma \, \#} \, T$ indefinitely oscillates, as $\sigma \to 0$, between rapidly decaying to a plane and being much closer to a non-planar cone than to a plane (and decaying slowly, or not at all, to a plane). To rule out this behaviour (and for a number of other key purposes subsequently), here we introduce an intrinsic frequency 
function for $T$, which we call the \emph{planar frequency function}, which satisfies, among other things, a certain monotonicity property. The planar frequency function $\rho \mapsto N_{T, P, Z}(\rho)$ (for $\rho \in (0, \rho_{0})$ for some $\rho_{0}>0$) is defined  relative to a given fixed $n$-dimensional plane $P$ and a given base point $Z \in {\rm spt}  \, T$ and in terms of 
geometric quantities integrated over the current $T$ (see Definition~\ref{freq defn})  
whenever $T$ has no boundary in the cylinder ${\mathbf C}_{\rho_{0}}(Z, P) = \{X \, : \pi_{P}(X - Z) < \rho_{0}\}$  (where $\pi_{P} \, : \, {\mathbb R}^{n+m} \to P$ is the orthogonal projection) and all points of $T$ in ${\mathbf C}_{\rho_{0}}(Z, P)$ are within a bounded  distance to the affine plane $Z + P.$

Before proceeding to a precise statement of the decomposition theorem, let us digress briefly to contrast  the planar frequency function and its utility with the role of the \emph{Almgren frequency function} in \cite{Almgren}. 

In \cite{Almgren}, Almgren introduced a frequency function  for area minimizing currents $T$ based on his frequency function for the linearized setting (i.e.\ for Dirichlet energy minimizing multi-valued functions), also introduced in \cite{Almgren}. This frequency function is associated to the normal height of $T$ relative to  a carefully constructed, sufficiently smooth (curved) center manifold that well-approximates the average height of $T$ off a plane to which $T$ is close. Thus this frequency function is neither ``planar'' in our sense, nor is it intrinsic; it is defined in terms of a choice of a Lipschitz multi-valued normal map on the chosen center manifold, whose graph well-approximates $T$ over an appropriate interval of scales.  It is established in \cite{Almgren} that corresponding to any given branch point, a sequence of such center manifolds and normal maps can be chosen iteratively. These correspond to intervals of scales at which $T$ is sufficiently close to a plane. Such intervals exist by the definition of branch point, but there may not be a single such interval with end point $0$ since tangent-cone uniqueness is not known and a non-planar tangent cone may exist at the branch point.  (Moreover, even if there is a unique tangent plane with a decay estimate, the center manifold might not pass through the branch point if the rate of decay is too slow.) It is then shown in \cite{Almgren} that the corresponding ``center manifold frequency functions'' satisfy an approximate monotonicity property. This is then   
used to rule out the possibility that on approach to a branch point the normal maps decay to $0$ infinitely rapidly, or in other words, the ``sheets'' of $T$ decay towards each other infinitely rapidly. This makes it possible to deduce the dimension bound for the branch set from the corresponding bound for the ``linear'' setting, i.e.\ for multi-valued Dirichlet energy minimizing functions arising as blow-ups of the normal maps. This way of proceeding requires a more extensive technical framework than via the planar frequency function. The significant technical challenges of this part of \cite{Almgren} are however, arguably, unavoidable given the basic premise: to \emph{directly} bound the size of ${\rm sing} \, T$ without addressing the question of uniqueness of the tangent planes at branch points. In our approach, by contrast, a center manifold is only needed in the entire program (including in the study of the local structure of ${\rm sing} \, T$ and of the current $T$ on approach to typical singular points) to analyze branch points where there is sufficiently fast  decay towards the tangent plane. This is a technically considerably simpler and canonical situation: for each of these special fast-decay branch points $Z$, a \emph{single} center manifold arises which passes through $Z$ and all nearby fast-decay branch points, and moreover, the corresponding normal map contains fast-decay branch points in its critical nodal set. We shall discuss this point in more detail towards the end of the introduction.

In our approach, initially, we are not concerned with ruling out the infinite order decay of sheets of $T$ towards each other. Our first use of the planar frequency function is to rule out the oscillatory behaviour of the current described above (and thereby to reach the decomposition theorem, Theorem~\ref{branch and cones thm-intro}). The key properties of the planar frequency function that enable us to do this are: 
\begin{itemize}
\item[(a)] $N_{T, P, Z}(\rho)$ is approximately monotonically non-decreasing on any interval $I$ of scales $\rho$ on which the rescaled current 
$\eta_{Z, \rho \, \#} \, T$ is close to $P$ and decays towards $P$ at any fixed rate $o(\rho^{\alpha})$ for some positive $\alpha$ (Theorem~\ref{mono freq thm} below); 
\item[(b)] $N_{T, P, 0}(\cdot)$ takes values $\leq 1$ whenever $T$ is any cone (having vertex $0$) with ${\rm spt} \, T \neq P$ and 
${\rm spt} \, T \cap P^{\perp} = \{0\}$ (Lemma~\ref{freq of cones lemma}). 
\end{itemize}

We can now give a precise statement of the decomposition theorem, which is proved employing the planar frequency function together with relatively elementary ingredients of Almgren's theory (specifically, the theory of Dirichlet energy minimising multi-valued functions,  strong Lipschitz approximation theorem and  convergence results for blow-up sequences of area 
minimisers relative to a plane). As mentioned above, this theorem is the first crucial  step in our asymptotic analysis of area minimisers near typical singularities. 
Here we let ${\mathcal P}$ denote the set of $n$-dimensional planes (passing through the origin) in ${\mathbb R}^{n+m},$
and use the notation $E(T, P, {\mathbf B}_{\rho}(Z)) = \left(\rho^{-n-2} \int_{\mathbf{B}_{\rho}(Z)} \op{dist}^2(X,Z +P) \,d\|T\|(X)\right)^{1/2}$ where 
$P \in {\mathcal P}$. 

\begin{theorem}[Theorem~\ref{branch and cones thm}]\label{branch and cones thm-intro}
For every integer $q \geq 2$ and $\varepsilon, \beta  \in (0, 1)$ there exist $R = R(n,m,q,\varepsilon,\beta) > 10$, $\delta = \delta(n,m,q,\varepsilon,\beta) \in (0, 1)$ and $\alpha = \alpha(n,m,q,\varepsilon,\beta) \in (0, 1)$ such that if $T$ is an $n$-dimensional locally area-minimizing rectifiable current of $\mathbf{B}_R(0)$ with  
\begin{equation}\label{branch and cones hyp1-intro} 
	(\partial T) \llcorner \mathbf{B}_R(0) = 0 \quad \mbox{and} \quad \|T\|(\mathbf{B}_R(0)) \leq (q + \delta) \,\omega_n R^n
\end{equation}
then 
\begin{equation*}
	\mathbf{B}_1(0) \cap \{ X : \Theta(T,X) \geq q \} = \mathcal{S} \cup \mathcal{B}
\end{equation*}
where $\mathcal{S}$  and $\mathcal{B}$ (which depend on $q$,$\varepsilon$ and $\beta$) are disjoint sets for which the following hold:
\begin{enumerate}[itemsep=3mm,topsep=0mm]
	\item[{\rm (I)}]  ${\mathcal S}$ is locally compact, and satisfies the following \emph{weak locally uniform approximation-by-non-planar-cones} property: for each $Z_0 \in \mathcal{S}$ there exists $\rho_0 \in (0,1/16]$ (depending on $Z_0$) such that for every $Z \in \op{spt} T \cap \mathbf{B}_{\rho_0}(Z_0)$ with $\Theta(T,Z) \geq q$ and every $\rho \in (0,\rho_0],$ one of the following assertions {\rm (i)} or {\rm (ii)} holds true:  
	\begin{enumerate}[itemsep=3mm,topsep=3mm]
	\item[{\rm (i)}]  there exists an integral cone $\mathbf{C} = \mathbf{C}_{Z, \rho}$ supported on a union of $n$-dimensional planes meeting along a common $(n-2)$-dimensional subspace such that $T$ is close to $\mathbf{C}$ in $\mathbf{B}_{\rho}(Z)$ and $T$ is significantly closer to $\mathbf{C}$ than to any plane $P$ in $\mathbf{B}_{\rho}(Z)$ in the sense that 
	\begin{align}
		\label{branch and cones concl2-intro} & \hspace{9mm} \rho^{-n-2} \int_{\mathbf{B}_{\rho}(Z)} \op{dist}^2(X,Z+\op{spt}\mathbf{C}) 
			\,d\|T\|(X) < \varepsilon^2\; \mbox{and}\\
		\label{branch and cones concl3-intro} &\int_{\mathbf{B}_{\rho}(Z)} \op{dist}^2(X,Z+\op{spt}\mathbf{C}) \,d\|T\|(X) 
			\\&\hspace{.4in}+ \int_{\mathbf{B}_{\rho/2}(Z) \cap \{\op{dist}(X,\op{spine} \mathbf{C}) \geq \rho/16\}} \op{dist}^2(Z+X,\op{spt}T) \,d\|\mathbf{C}\|(X) 
				\nonumber 
			\\& \hspace{1in}\leq\, \beta^2 \inf_{P \in {\mathcal P}} \int_{\mathbf{B}_{\rho}(Z)} \op{dist}^2(X,Z+P) \,d\|T\|(X); \nonumber\\
	\label{branch and cones concl3-1-intro} & \{ Y \in \mathbf{B}_{\rho}(Z) : \Theta(T,Y) \geq q \} \subset 
			\{ Y \in \mathbf{B}_{\rho}(Z) : \op{dist}(Y,Z+\op{spine}\mathbf{C}) < \varepsilon \rho \}; 
			\end{align}
	\item[{\rm (ii)}]  there is an $(n-3)$-dimensional linear subspace $L$ such that 
	\begin{equation}\label{branch and cones concl4-intro}
		\{ Y \in \mathbf{B}_{\rho}(Z) : \Theta(T,Y) \geq q \} \subset \{ Y \in \mathbf{B}_{\rho}(Z) : \op{dist}(Y,Z+L) < \varepsilon \rho \} . 
	\end{equation}
	\end{enumerate}
	
	\item[{\rm (II)}]  $\mathcal{B}$ is relatively closed in ${\mathbf B}_{1}(0),$ and: 
	\begin{enumerate}[itemsep=3mm,topsep=3mm]
	\item[{\rm (i)}] if $Z_0 \in \mathcal{B}$ then $\Theta(T, Z_{0}) = q$ and there is a unique $n$-dimensional plane $P_{Z_0}$ such that 
	\begin{equation}\label{branch and cones concl1 intro} 
	E(T,P_{Z_0},\mathbf{B}_{\sigma}(Z_0)) \leq C_0 \Big(\frac{\sigma}{\rho}\Big)^{\alpha} E(T, P_{Z_0}, \mathbf{B}_{\rho}(Z_{0}))
	\end{equation}
	for all $0 < \sigma \leq \rho \leq 2$, where $C_{0} = C_{0}(n, m, q, \varepsilon, \beta) \in (0, \infty)$ is a constant;
	 \item[{\rm (ii)}]  if $Z_{0} \in {\mathcal B}$, then $P_{Z_0}$ (as in \eqref{branch and cones concl1 intro}) taken with multiplicity $q$ and oriented appropriately is the unique tangent cone to $T$ at $Z_0$; 
	
	\item[{\rm (iii)}]  if $Z_{0} \in {\mathcal B}$ then the planar frequency ${\mathcal N}_{T, {\rm Pl}}(Z_{0}) = \lim_{\rho \rightarrow 0^+} N_{\widetilde{T},P_{Z_0},Z_0}(\rho)$ exists and ${\mathcal N}_{T,{\rm Pl}}(Z_{0}) \geq 1+\alpha$ 
	(where $\widetilde{T} = T \res {\mathbf B}_{15/8}(Z_{0})$ and $N_{\widetilde{T},P_{Z_0},Z_0}(\rho)$ is as in Definition~\ref{freq defn} taken with $\widetilde{T}$ in place of $T$); 
	 \item[{\rm (iv)}] for any two points $Z_{1}, Z_{2} \in {\mathcal B} \, \cap\, {\mathbf B}_{1/2}(0)$, we have that $$\|\pi_{P_{Z_{1}}} - \pi_{P_{Z_{2}}}\| \leq C|Z_{1} - Z_{2}|^{\alpha}\left(E(T, P_{Z_{1}}, {\mathbf B}_{2}(Z_{1})) +E(T, P_{Z_{2}}, {\mathbf B}_{2}(Z_{2})\right)),$$ where $C = C(n, m, q, \epsilon, \beta)$ is a constant and for each $P \in {\mathcal P}$, 
	$\pi_P : \mathbb{R}^{n+m} \rightarrow P$ is the orthogonal projection onto $P$.
	\end{enumerate}
\end{enumerate}
\end{theorem}

Theorem~\ref{branch and cones thm-intro} is, in part, inspired by arguments in~\cite{Wic08},\cite{Wic14} for analysis of stable codimension 1 integral varifolds.  As was observed in~\cite{Wic08},\cite{Wic14}, using blow-up arguments and the Hardt-Simon inequality~\cite{HardtSimon} (see Lemma~\ref{hardt simon lemma} below), one can regard density $q$ singular points as having planar frequency $\geq 1$ (if the latter exists).  In the setting of \cite{Wic08},\cite{Wic14} however, unlike here, the validity of conditions 
(\ref{branch and cones concl2-intro})  and (\ref{branch and cones concl3-intro}) (for appropriately chosen, fixed $\epsilon$ and $\beta$) at \emph{some} scale and for some cone ${\mathbf C}$ (which, in that setting, is made up of half-hyperplanes meeting along a common $(n-1)$-dimensional boundary) automatically implies decay of the varifold to a unique non-planar cone of the same type. (See the more recent work \cite[Theorem~3.1]{MW} for a much more definitive result for stable hypersurfaces in this vein.) Moreover, a frequency function was not necessary (nor used) in that setting.

The proof of Theorem~\ref{branch and cones thm-intro} proceeds by  
setting $\mathcal{S}$ to be the density $ \geq q$ singular points $Z_0$ in ${\mathbf B}_{1}(0)$ such that for some $\sigma > 0$, the decay (by a fixed factor) of height excess of $T$ relative to an optimal plane from scale $\sigma$ to $\sigma/2$ fails. It follows that the planar frequency function $N_{T,P,Z_0}(\rho)$ is close to 1 for any plane $P$ and $\rho \in (\sigma/4, \sigma/2]$ (Lemma~\ref{not decaying lemma1} and Lemma~\ref{not flat lemma1}).  Thus by monotonicity of the planar frequency function (Theorem~\ref{mono freq thm}), if $T$ again decays towards a plane $\widetilde{P}$ from scale $\sigma/2$ down to a smaller scale $\sigma_{1} < \sigma/2$, the planar frequency function of $T$ relative to $\widetilde{P}$ must remain close to 1 at scales $\in (\sigma_{1}, \sigma/2]$.  Hence $T$ must be significantly closer to a non-planar cone than to any multiplicity $q$ plane in $\mathbf{B}_{\rho}(Z_0)$ for all $\rho \in (\sigma_{1},\sigma]$ (Lemma~\ref{freq2cone lemma}). 
More precisely, either $T$ is weakly close to an area-minimizing cone but $T$ is not close to any multiplicity $q$ plane in $\mathbf{B}_{\rho}(Z_0)$ (Lemma~\ref{not decaying lemma2}), or $T$ is close to a multiplicity $q$ plane in $\mathbf{B}_{\rho}(Z_0)$ and $T$ is significantly closer to the graph of a homogeneous degree one Dirichlet energy minimizing $q$-valued function over $P$ (Lemma~\ref{freq2cone lemma} and Lemma~\ref{not flat lemma1}).  
Since the non-planar cone in either case is either a union of planes meeting along an $(n-2)$-dimensional subspace or has an axis (i.e.~spine) of dimension $\leq n-3$, the assertion of Theorem~\ref{branch and cones thm-intro}(I) follows.  Setting $\mathcal{B}$ to be the complement of ${\mathcal S}$ in the set of density $ \geq q$ singular points in ${\mathbf B}_{1}(0)$, it follows (from the definition of ${\mathcal S}$) that at every point of $Z \in {\mathcal B}$ the current decays to a unique plane at the rate $|X - Z|^{\alpha}$ or faster for some fixed $\alpha \in (0, 1)$. 

Let $q \geq 1$, $\epsilon, \beta \in (0, 1),$ and let $R = R(n, m, q, \epsilon, \beta),$ $\delta = \delta(n, m, q, \epsilon, \beta)$ be as in Theorem~\ref{branch and cones thm-intro}. Let $T$ be a locally area minimizing rectifiable current in some open set in ${\mathbb R}^{n+m}$. Since for every point $Z \in {\rm sing}_{q} \, T$ (the set of density $q$ singularities of $T$) there is a scale 
$\sigma_{Z}>0$ such that hypothesis \eqref{branch and cones hyp1-intro} of Theorem~\ref{branch and cones thm-intro} is satisfied with $T_{Z} = \eta_{Z, \sigma_{Z} \, \#} \, T$ in place of $T$, by setting 
$V = \bigcup_{Z \in {\rm sing}_{q} \, T} \, {\mathbf B}_{\sigma_{Z}}(Z)$, ${\mathcal S}_{q, \varepsilon, \beta} = \bigcup_{Z \in {\rm sing}_{q} \, T} \,\eta_{Z,\sigma_{Z}}^{-1} {\mathcal S}_{Z}$ where ${\mathcal S}_{Z}$ is the set ${\mathcal S}$ given by Theorem~\ref{branch and cones thm-intro} with $T_Z$ in place of $T$, and 
\begin{equation*}
    {\mathcal B}_{q, \varepsilon, \beta} = {\rm sing} \, T \cap V \cap \{Z \, : \, \Theta(T, Z) \geq q\} \setminus {\mathcal S}_{q, \varepsilon, \beta} , 
\end{equation*}
we have the following immediate consequence of Theorem~\ref{branch and cones thm-intro}:

\begin{corollary}\label{immediate-cor}
Let $T$ be an $n$-dimensional locally area minimizing rectifiable current in an open set in ${\mathbb R}^{n+m}$. 
For every integer $q \geq 1$ and $\varepsilon, \beta  \in (0, 1)$ there exists an open set $V = V_{q, \epsilon, \beta}$ in ${\mathbb R}^{n+m}$ with ${\rm sing}_{q} \, T \subset V$ such that 
\begin{equation*}
	\op{sing} T \cap V \cap \{ X : \Theta(T,X) \geq q \} = \mathcal{S}_{q, \varepsilon, \beta} \cup \mathcal{B}_{q, \varepsilon, \beta}
\end{equation*}
where $\mathcal{S}_{q,\varepsilon,\beta}$, $\mathcal{B}_{q,\varepsilon,\beta}$ are disjoint, locally compact sets having the following properties:

\begin{itemize}
\item[{\rm (a)}] conclusion {\rm (I)} of Theorem~\ref{branch and cones thm-intro} holds with ${\mathcal S}_{q, \varepsilon, \beta}$ in place of ${\mathcal S};$
\item[{\rm (b)}] $\mathcal{B}_{q, \varepsilon, \beta}$ is relatively closed in $V$ and  if $Z_{0} \in \mathcal{B}_{q, \varepsilon, \beta}$ then: 
	\begin{enumerate}[itemsep=3mm,topsep=3mm]
	\item[{\rm (i)}] $\Theta(T, Z_{0}) = q,$ and there is a number $\sigma_{Z_{0}}>0$ such that for every $Z \in {\mathcal B}_{q, \varepsilon, \beta} \cap {\mathbf B}_{\sigma_{Z_{0}}}(Z_{0})$ there is a unique $n$-dimensional plane $P_{Z}$ such that the plane $P_{Z}$ with multiplicity $q$ and appropriate orientation is the unique tangent cone to $T$ at $Z$ 
	\begin{equation}\label{branch and cones concl1 intro-cor} 
		E(T, P_{Z}, {\mathbf B}_{\sigma}(Z)) \leq C_0 \Big(\frac{\sigma}{\rho}\Big)^{\alpha} E(T, P_{Z, 0}, {\mathbf B}_{\rho}(Z)) 
	\end{equation}
	for all $0 < \sigma \leq \rho \leq \sigma_{Z_0},$ where $C_{0} = C_{0}(n,m, q, \varepsilon, \beta) \in (0,\infty)$ and $\alpha = \alpha(n, m, q, \varepsilon, \beta) \in (0, 1)$ are the constants as in Theorem~\ref{branch and cones thm-intro}{\rm (II)(i)};
    \item[{\rm (ii)}] the planar frequency ${\mathcal N}_{T, {\rm Pl}}(Z_{0}) = \lim_{\rho \rightarrow 0^+} N_{\widetilde{T},P_{Z_0},Z_0}(\rho)$ exists and ${\mathcal N}_{T,{\rm Pl}}(Z_{0}) \geq 1+\alpha$ 
	(where $\alpha$ is as in (\ref{branch and cones concl1 intro-cor}) and $\widetilde{T} = T \res {\mathbf B}_{\rho_{0}}(Z_{0})$ for appropriately small $\r_{0}>0,$ and $N_{\widetilde{T},P_{Z_0},Z_0}(\rho)$ is as in Definition~\ref{freq defn} taken with $\widetilde{T}$ in place of $T$).
    \end{enumerate}
\end{itemize}
\end{corollary}

\subsection{Hausdorff dimension bound for the singular set: simplifications over Almgren's argument} Once Theorem~\ref{branch and cones thm-intro} is established, our primary focus will be on using it to analyse the asymptotic behaviour of the current on approach to a typical singular point, establishing uniqueness of tangent cones, existence and uniqueness of non-zero multi-valued tangent functions (blow ups) at branch points as well as 
structural properties of the singular set itself including its $(n-2)$-rectifiability and local-finiteness-of-measure properties, and the topological nature of the current near branch points satisfying additional conditions. We carry out this asymptotic  analysis in \cite{KrumWicb}, \cite{KrumWicc} and \cite{KrumWicd}. As a first consequence of Theorem~\ref{branch and cones thm-intro} though we obtain that the Hausdorff dimension of ${\rm sing} \, T$ is $\leq n-2$, which is Almgren's main theorem in \cite{Almgren}. 

\begin{corollary}\label{Alm big reg thm}
Let $T$ be an $n$-dimensional locally area-minimizing rectifiable current in an open subset of $\mathbb{R}^{n+m}$.  The Hausdorff dimension of ${\rm sing} \, T$ is at most $n-2$.
\end{corollary}

\begin{proof}[Proof outline] 
Theorem~\ref{branch and cones thm-intro} provides a way to reach this conclusion avoiding a considerable part of the technical complexity of the argument in \cite{Almgren}. We here sketch how the simplifications arise.
Let ${\mathcal S}_{q,\varepsilon,\beta}$ and ${\mathcal B}_{q,\varepsilon,\beta}$ be as in Corollary~\ref{immediate-cor}.  First, by \eqref{branch and cones concl3-1-intro} and \eqref{branch and cones concl4-intro},  there is $\gamma(\varepsilon, \beta) >0$ with $\gamma(\varepsilon, \beta) \to 0$ as $(\varepsilon, \beta) \to (0,0)$
such that ${\mathcal S}_{q, \varepsilon, \beta}$ satisfies ${\mathcal H}^{n-2+ \gamma(\varepsilon, \beta)}({\mathcal S}_{q, \varepsilon, \beta}) = 0$. Having established this, to make 
a similar conclusion for ${\mathcal B}_{q, \epsilon, \beta}$,  we can capitalize on the validity of the 
decay estimate \eqref{branch and cones concl1 intro} and the approximate monotonicity of the planar frequency function. Write ${\mathcal B}_{q, \varepsilon, \beta} = {\mathcal B}_{q, \varepsilon, \beta}^{(\neq 2)} \cup 
{\mathcal B}_{q, \varepsilon, \beta}^{(=2)}$, where $${\mathcal B}_{q, \varepsilon, \beta}^{(\neq 2)} = \{Z \in {\mathcal B}_{q, \varepsilon, \beta} \, : \, 
{\mathcal N}_{T, {\rm Pl}}(Z) \neq 2\} \;\; \mbox{and}$$ 
$${\mathcal B}_{q, \varepsilon, \beta}^{(=2)} = \{Z  \in {\mathcal B}_{q, \varepsilon, \beta} \, : \, {\mathcal N}_{T, {\rm Pl}}(Z)  = 2\}.$$ Since a locally uniform estimate for $T$ holds at every point $Z \in {\mathcal B}_{q, \varepsilon, \beta},$ giving decay of $T$ towards a unique tangent plane $P_{Z}$, we obtain, as a straightforward consequence of the monotonicity formula for the planar frequency function, that any tangent function (blow-up) of $T$ at a point $Z \in {\mathcal B}_{q, \epsilon, \beta}$ relative to $P_{Z}$ is 
a \emph{non-zero}, homogeneous (of degree $\mu_{Z} = {\mathcal N}_{T, {\rm Pl}}(Z) \geq 1+ \alpha$), $q$-valued Dirichlet energy minimizing function $\varphi_{Z}  \, : \, P_{Z} \to {\mathcal A}_{q}(P_{Z}^{\perp})$ (which, in case $\mu_{Z}$ is an integer $\geq 2$, may well be $q$ copies of a single-valued harmonic function). Here ${\mathcal A}_{q}(P_{Z}^{\perp})$ denotes the space of 
``unordered $q$-tuples of points in $P_{Z}^{\perp},$'' or more precisely, the space $\{\sum_{k=1}^{q} \llbracket a_{j} \rrbracket \, : \, a_{j} \in P_{z}^{\perp} \;\; \forall j = 1, \ldots, q\}$ where $\llbracket a\rrbracket$ denotes the Dirac mass at $a \in P_{Z}^{\perp}$.
Letting $\varphi_{Z, \, a}(x)$ be the average of the the $q$ values of $\varphi_{Z}(x)$ for $x \in P_{Z}$ (so that either $\mu_{Z}$ is an integer and $\varphi_{Z, \, a}$ is a non-zero single-valued homogeneous harmonic polynomial of degree $\mu_{Z}$, or $\varphi_{Z, \, a} \equiv 0$), and  writing $\ell_{Z, \, \xi}(x) = \varphi_{Z, \, a}(\xi) + D\varphi_{Z, \, a}(\xi) \cdot (x- \xi)$ for any fixed $\xi \in P_{Z}$ and all $x \in P_{Z}$, it can readily be checked that if $Z \in {\mathcal B}_{q, \varepsilon, \beta}^{(\neq 2)}$ then 
the set $\widetilde{S}(\varphi_{Z}) \equiv \{\xi \in P_{Z} \, : \, {\mathcal N}_{\varphi_{Z} - \ell_{Z, \, \xi}}(\xi) \geq {\mathcal N}_{\varphi_{Z}}(0)\}$ is a linear subspace of $P_{Z}$ with ${\rm dim} \, S(\varphi_{Z}) \leq (n-2)$. Here 
${\mathcal N}_{\psi}(y)$ is the (Almgren) frequency of $\psi$ at $y$. Indeed, if $\xi \in \widetilde{S}(\varphi_{Z})$ then $\varphi_{Z, \, a}(\xi) = 0$ and $D\varphi_{Z, \, a}(\xi) = 0$, which 
means that $\widetilde{S}(\varphi_{Z}) = S(\varphi_{Z}) \equiv \{\xi \in P_{Z} \, : \, {\mathcal N}_{\varphi_{Z}}(\xi) = {\mathcal N}_{\varphi_{Z}}(0)\},$ and 
by a standard argument $S(\varphi_{Z})$ is a linear subspace of $P_{Z}$ of dimension $\leq n-2$ (and is the subspace along which $\varphi_{Z}$ is invariant under translation). This readily implies  that the Hausdorff dimension of ${\mathcal B}_{q, \varepsilon, \beta}^{(\neq 2)}$ is $\leq n-2.$ (See \cite{KrumWicc} for details). 

Note that  if $\varphi_{Z}$ is $q$ copies of a single-valued \emph{quadratic} harmonic polynomial then $\widetilde{S}(\varphi_{Z}) = P_{Z}$; since this may be the case if $Z \in {\mathcal B}_{q, \varepsilon, \beta}^{(=2)}$, the analysis of 
${\mathcal B}_{q, \varepsilon, \beta}^{(=2)}$ needs to proceed differently. 
To bound the dimension of ${\mathcal B}_{q, \varepsilon, \beta}^{(=2)}$, we proceed as in \cite{Almgren} by utilising a center manifold (following its construction given in \cite{DeLSpa2})  with considerable added simplifications arising from the fact that the current at every point $Z \in {\mathcal B}_{q, \varepsilon, \beta}^{(=2)}$ satisfies a locally uniform estimate giving  decay to a unique plane at a \emph{quadratic} rate in the scale. (For instance, complications necessitating having to consider ``intervals of flattening'' and correspondingly infinite sequences of center manifolds, as in \cite{Almgren},\cite{DeLSpa2}, \cite{DeLSpa3}, are all removed by the uniform decay estimates.)
 This quadratic decay estimate is an immediate consequence of the fact that the planar frequency at every point in ${\mathcal B}_{q, \varepsilon, \beta}^{(=2)}$ is $2$, and it means that the center manifold in this setting is canonical: about any point in ${\mathcal B}_{q, \varepsilon, \beta}^{(=2)}$, there is a a center manifold (which is unique up to an additive 
term of order $o(|X- Z|^{2})$) that contains that point and 
\emph{all} nearby points in ${\mathcal B}_{q, \varepsilon, \beta}^{(=2)}$; 
moreover, these points are all contained in the critical nodal set of the $q$-valued normal map $N$ over the center manifold (whose graph provides a quantitative approximation for the current),  the Almgren frequency function associated with $N$ is approximately monotone, and $N$ satisfies a \emph{uniform} quadratic $L^{2}$-decay estimate about each point of ${\mathcal B}^{(=2)}_{q, \epsilon, \beta}$. This allows us to use a standard argument to produce  non-trivial, average-free homogeneous $q$-valued blow-ups of $N$ at any point in ${\mathcal B}^{(=2)}_{q, \epsilon, \beta},$ which are locally Dirichlet energy minimizing and hence have translation invariance along a subspace of dimension no greater than $n-2$. This together with the uniform $L^{2}$-decay estimate satisfied by $N$ at points in ${\mathcal B}^{(=2)}_{q, \epsilon, \beta}$ (which passes to the blow-up) leads directly to the conclusion that ${\rm dim}_{\mathcal H}({\mathcal B}^{(=2)}_{q, \epsilon, \beta}) \leq n-2$.  
(See \cite{KrumWicd} for details; in fact for any density $q$ branch point with decay rate to the tangent plane \emph{quadratic or higher}, there is a single center manifold that locally contains that point and all nearby such branch points, and thus the same argument carries over to provide the dimension bound for branch points with quadratic or faster decay). 

Since ${\rm sing}_{q} \, T \equiv \{Z \in {\rm sing} \, T\, : \, \Theta(T, z) = q\} \subset {\mathcal S}_{q, \varepsilon, \beta} \cup {\mathcal B}_{q, \varepsilon, \beta},$ combining these results yields that ${\rm dim}_{\mathcal H} \, ({\rm sing}_{q} \,  T) \leq n-2$ for every integer $q \geq 2$. Since by a standard tangent cone analysis the set of singularities of non-integer density has Hausdorff dimensnion $\leq n-3$, the corollary follows. 
\end{proof}

\subsection{Uniqueness of tangent cones and the structure of singularities away from rapid-decay branch points: choosing $\epsilon$ and $\beta$} \label{non-branch-structure} 
Since every point in ${\mathcal B}_{q, \varepsilon, \beta}$ has a unique planar tangent cone, it follows that 
 the set ${\mathcal S}_{q, n-2}$ of all density $q$ non-branch point singularities (where each tangent cone has spine dimension $\leq n-2$) is contained in ${\mathcal S}_{q, \varepsilon, \beta}.$ As mentioned above,   
Theorem~\ref{branch and cones thm-intro} however allows that ${\mathcal S}_{q, \varepsilon, \beta}$ may also contain branch points, and in fact it allows a priori that ${\mathcal S}_{q, \varepsilon, \beta} \setminus {\mathcal S}_{q, n-2}$ may be a set (consisting of branch points) of positive $(n-2)$-dimensional Hauadorff measure.
In \cite{KrumWicb}, we rule out this latter possibility. In fact in \cite{KrumWicb} we shall take the above locally uniform weak approximation property of the current at points in ${\mathcal S}_{q, \varepsilon, \beta}$ (to be precise, conclusions~(\ref{branch and cones concl2-intro}), (\ref{branch and cones concl3-intro}), (\ref{branch and cones concl4-intro}) in Theorem~\ref{branch and cones thm-intro}) as a starting point and prove the following: \emph{there is a choice of $\varepsilon = \varepsilon(n, m, q) \in (0, 1)$ and $\beta = \beta(n, m, q) \in (0, 1)$ such that for ${\mathcal H}^{n-2}$ a.e.\ point $Z \in {\mathcal S}_{q}  \equiv {\mathcal S}_{q, \varepsilon, \beta}$, the current $T$ has a unique tangent cone ${\mathbf C}_{Z}$ supported on a union of (at least two) planes intersecting along an $(n-2)$-dimensional subspace, and moreover, $T$ satisfies a decay estimate giving convergence of $\eta_{Z, \rho \, \#} \, T$ to ${\mathbf C}_{Z}$ as $\rho \to 0^{+}$ at a rate $o(\rho^{\mu})$ for some fixed $\mu = \mu(n, m, q) \in (0, 1)$ depending only on $n$, $m$ and $q$}. From this it follows that ${\mathcal H}^{n-2}({\mathcal S}_{q} \setminus {\mathcal S}_{q, n-2}) = 0$ and that \emph{${\mathcal S}_{q}$ is  
locally $(n-2)$-rectifiable (with locally finite $(n-2)$-dimensional Hausdorff measure)}. Moreover, since  the tangent cone at every point in 
${\mathcal B}_{q, \varepsilon, \beta}$ is unique, it then follows that the \emph{current has a unique tangent cone (equal to the sum of a finite number of planes) at ${\mathcal H}^{n-2}$ a.e.\ point}.

\begin{remark}\label{top-stratum-rem} {\rm
The work of Naber--Valtorta~\cite{NV15} implies that 
${\mathcal S}_{q, n-2}$  is countably $(n-2)$-rectifiable.   
However, one cannot apply the results of~\cite{NV15} to ${\mathcal S}_{q, \varepsilon,\beta}$ prior to ruling out (as done in \cite{KrumWicb}, once $\varepsilon,$ $\beta$ are chosen appropriately depending only on $n$, $m$ and $q$) the possibility that ${\mathcal S}_{q, \varepsilon, \beta}$ contains a set of positive $(n-2)$-dimensional Hausdorff measure consisting of branch points of $T$. Our argument ruling out this possibility yields simultaneously the uniqueness of tangent cones and rectifiability conclusions, independently of the rectifiability theorem of \cite{NV15}.  
}\end{remark}

\subsection{Higher order asymptotics at rapid-decay branch points: structure of the branch set and of the full singular set} \label{higher-order-asym} Let us now discuss briefly the higher-order behaviour of $T$ near typical branch points, analysed in \cite{KrumWicc}, \cite{KrumWicd}. Let $T$ and ${\mathcal B}_{q, \varepsilon, \beta}$ be as in Corollary~\ref{immediate-cor}.
Fix $\varepsilon = \varepsilon(n, m, q)$, $\beta = \beta(n, m, q)$ in the manner just described (in Section~\ref{non-branch-structure}). With this choice of $\varepsilon, \beta$, write 
$${\mathcal B}_{q} = {\mathcal B}_{q, \varepsilon, \beta},$$
$${\mathcal B}_{q}^{(\neq 2)} = \{Z \in {\mathcal B}_{q} \, : \, 
{\mathcal N}_{T, {\rm Pl}}(Z) \neq 2\} \;\; \mbox{and}$$ 
$${\mathcal B}_{q}^{(=2)} = \{Z  \in {\mathcal B}_{q} \, : \, {\mathcal N}_{T, {\rm Pl}}(Z)  = 2\}.$$

As pointed out in the proof-outline of Corollary~\ref{Alm big reg thm} above,  
every tangent function of $T$ at any point $Z \in {\mathcal B}_{q}$ is a non-zero locally Dirichlet energy minimizing $q$-valued function 
$\varphi_{Z} \, : \, P_{Z} \to {\mathcal A}_{q}(P_{Z}^{\perp})$ which is homogeneous of degree ${\mathcal N}_{T, {\rm Pl}}(Z).$ Consider first the set ${\mathcal B}_{q}^{(\neq 2)}$. For ${\mathcal H}^{n-2}$ a.e.\ point $Z \in {\mathcal B}_{q}^{(\neq 2)}$, at least one of these tangent functions $\varphi_{Z}$ is \emph{cylindrical}, i.e.\ has the property that 
$S(\varphi_{Z})$ (the linear subspace of $P_{Z}$ along which $\varphi_{Z}$ is invariant under translations) has dimension $n-2$, and thus $\varphi_{Z}$ depends only on 2-variables;  consequently,  $\varphi_{Z}$  takes a very specific and explicit form: namely, identifying $P_{Z}$ with ${\mathbb R}^{n}$ and $P_{Z}^{\perp}$ with ${\mathbb R}^{m}$, in terms of coordinates $X = (x_{1}, \ldots, x_{n}) \in {\mathbb R}^{n}$ chosen so that  
$S(\varphi_{Z}) = \{X \, : \, x_{1} = x_{2} = 0\}$, 
\begin{equation*}
\varphi_{Z}(X) 
	= (q - N_Z q_Z) \llbracket 0 \rrbracket + \sum_{l=1}^{N_Z} \op{Re}(c^{(Z)}_l (x_1+ix_2)^{k_Z/q_Z}) 
\end{equation*}
for some integer $N_{Z} \geq 1;$ some relatively prime (positive) integers $k_{Z}, q_{Z}$ with $k_{Z} > q_{Z} \geq 1$ and $N_{Z}q_{Z} \leq q;$ and some 
$c^{(Z)}_l \in \mathbb{C}^m \setminus \{0\}$. 

The fact that there is a such a simple, explicit classification of $\varphi_{Z}$ is particularly motivating for undertaking a study of  
 higher-order (i.e.\ beyond tangent cone level) analysis of $T$ aimed at proving uniqueness of $\varphi_{Z}$ (i.e.\ that $\varphi_{Z}$ is independent of the sequence of scales used to generate it from $T$). For such an analysis would shed light not just on the local structure of ${\mathcal B}_{q}^{(\neq 2)}$ 
 but also, in view of the explicit form of $\varphi_{Z}$,  
 on the important question of topological and asymptotic nature of $T$ in the vicinity of branch points. 
 In \cite{KrumWicc} we indeed establish higher-order decay estimates valid at ${\mathcal H}^{n-2}$ a.e.\ point $Z \in {\mathcal B}_{q}^{(\neq 2)}.$ These estimates imply, among other things, uniqueness of 
 $\varphi_{Z}$, local $(n-2)$-rectifiability properties of ${\mathcal B}_{q}^{(\neq 2)}$ and, for certain points $Z \in {\mathcal B}_{q}^{(\neq 2)},$ that $T$ near $Z$ is homeomorphic to 
 an $n$-disk admitting a $C^{1, \mu}$ parameterization.
  This last result will be further explained, and put in the context of known results for 2-dimensional area minimizing currents, in the next subsection.  
 
 We note that none of our results in \cite{KrumWicc} nor in the present article or in \cite{KrumWicb} requires the use of center manifolds. 
In \cite{KrumWicd} we obtain for ${\mathcal B}_{q}^{(=2)}$ results analogous to those in \cite{KrumWicc} for ${\mathcal B}_{q}^{(\neq 2)};$ and it is in \cite{KrumWicd} that we make use of the Almgren center manifold,  capitalizing on the special properties of the center manifold and the normal approximation map $N$ valid for planar-frequency 2 points. Indeed, for any point $Z \in {\mathcal B}_{q}^{(=2)}$, there is a \emph{single} center manifold ${\mathcal M}_{Z}$ (which is essentially unique for ${\mathcal H}^{n-2}$-a.e.\ $Z \in {\mathcal B}_{q}^{(=2)}$) containing $Z$ and all nearby points in ${\mathcal B}_{q}^{(=2)}$; the Lipschitz $q$-valued normal map $N_{Z}$ defined on ${\mathcal M}_{Z}$ provides a locally uniformly quantitative approximation for $T$ about all  points $Z^{\prime} \in {\mathcal B}_{q}^{(=2)}$ near $Z$ and at \emph{all} sufficiently small scales (not just on an annular region about $Z$); and moreover, $N_{Z}$ satisfies a locally uniform $L^{2}$-decay estimate about any point in ${\mathcal M}_{Z} \cap {\mathcal B}_{q}^{(= 2)}$ near $Z$. This puts all points in ${\mathcal B}_{q}^{(= 2)}$ that are close to $Z$ in the critical nodal set of $N_{Z}$. These facts about the center manifold and the normal map are of crucial importance for our asymptotic analysis of $T$ near points in ${\mathcal B}_{q}^{(=2)}$, which is essentially an analysis of the critical nodal set, and the behavior near the critical nodal set, of a map satisfying variational requirements and decay conditions. Thus, the reduction of the use of the center manifold to this geometrically canonical case not only allows us to reach the dimension bound for ${\rm sing} \, T$ more efficiently,  
 but also is an indispensable step in our higher order analysis and topological conclusions for the current near branch points.

Omitting the precise statements of various decay estimates, we can summarise as follows (in Theorem~\ref{summary-theorem-1} and Theorem~\ref{summary-theorem-2} below and in Theorem~\ref{topology} in the next subsection) the main structure, uniqueness and topological conclusions resulting from our analysis in \cite{KrumWicb}, \cite{KrumWicc}, \cite{KrumWicd}.

\begin{theorem}[local structure of the singular set]\label{summary-theorem-1} 
If $T$ is an $n$-dimensional locally area-minimizing rectifiable current in an open set $U \subset {\mathbb R}^{n+m}$, then 
for every compact set $K \subset U,$ ${\rm sing} \, T \cap K$ is the union of a finite number of pairwise disjoint sets each of which is locally compact and locally $(n-2)$-rectifiable (and has in particular locally finite $(n-2)$-dimensional Hausdorff measure).
\end{theorem}

This is a direct counterpart in arbitrary codimension to a theorem of L.~Simon (\cite{Sim95}) for codimension 1 area minimizers. Simon's celebrated theorem is applicable to any \emph{multiplicity 1 class} of minimal submanifolds (for which, in particular, branch point singularities are ruled out a priori). For codimension 1 area mininimizers it says that the statement above holds with $n-7$ in place of $n-2$.  Additionally in codimension 1, ${\rm sing} \, T$ has locally finite $(n-7)$-dimensional Hausdorff measure, as established by the more recent work of Naber--Valtorta (\cite{NV15}, which provides an alternative method to that of \cite{Sim95} for proving rectifiability of ${\rm sing} \, T$ in codimension 1. Whether this locally finite measure conclusion (with $n-2$ in place of $n-7$) extends to area minimizers of higher codimension remains an open question. 

\begin{theorem}[asymptotics near branch points]\label{summary-theorem-2} 
If $T$ is an $n$-dimensional locally area-minimizing rectifiable current in an open set $U \subset {\mathbb R}^{n+m}$, then 
\begin{itemize}
\item[{\rm (i)}]  $T$ has a unique tangent cone at ${\mathcal H}^{n-2}$ a.e.\ point; 
\item[{\rm (ii)}] for ${\mathcal H}^{n-2}$ a.e.\ branch point $Z$, there is a unique plane $P_{Z}$ which, taken with multiplicity $q_{Z} = \Theta(T, Z)$ and appropriate orientation, is the unique tangent cone to $T$ at $Z;$ moreover, the following hold:
\begin{itemize}
\item[{\rm (a)}] the planar frequency ${\mathcal N}_{T, {\rm Pl}}(Z) = \lim_{\rho \rightarrow 0^+} N_{\widetilde{T},P_{Z},Z}(\rho)$ exists and is a finite rational number with ${\mathcal N}_{T,{\rm Pl}}(Z) \geq 1+q_{Z}^{-1}$ 
	(where $\widetilde{T} = T \res {\mathbf B}_{\rho_{0}}(Z)$ for appropriately small $\rho_{0}>0$ and $N_{\widetilde{T},P_{Z},Z}(\rho)$ is as in Definition~\ref{freq defn} taken with $\widetilde{T}$ in place of $T$);
\item[{\rm (b)}] the planar frequency ${\mathcal N}_{T, {\rm Pl}}(Z)$ is the order of contact of $T$ with the tangent plane $P_{Z},$ and in particular 
\begin{equation*}
    \lim_{\rho \to 0} \, \rho^{-n-2{\mathcal N}_{T, {\rm Pl}}(Z)} \int_{{\mathbf B}_{\rho}(Z)} {\rm dist}^{2} \, (X, Z + P_{Z}) \, d\|T\|(X)
\end{equation*}
exists and is a finite positive real number; 
\item[{\rm (c)}] there is a unique (possibly zero) single-valued harmonic function $h_{Z} \, : \, P_{Z} \to P_{Z}^{\perp}$ that is homogeneous of degree 
${\mathcal N}_{T, {\rm Pl}}(Z)$;   a (possibly zero) function $\psi_{Z} \, : \, P_{Z} \to P_{Z}^{\perp}$ 
 of class $C^{3, \alpha}$ for some $\alpha \in (0, 1)$ (the ``center manifold'') with 
\begin{equation*}
    \lim_{\rho\rightarrow 0^+} \rho^{-n-2{\mathcal N}_{T, {\rm Pl}}(Z)} \int_{P_Z \cap \mathbf{B}_{\rho}(0)} |\psi_Z - h_Z|^2 = 0 ;
\end{equation*} 
a unique rational number ${\mathcal O}_{T}(Z)$ (the ``branching order'') with ${\mathcal O}_{T}(Z) \geq {\mathcal N}_{T, {\rm Pl}}(Z);$ a unique non-zero $q_{Z}$-valued locally Dirichlet energy minimizing function (tangent function) $$\varphi_{Z} \, : \, P_{Z} \to {\mathcal A}_{q_{Z}}(P_{Z}^{\perp})$$ which is cylindrical (i.e.\ invariant under translation along an $(n-2)$-dimensional subspace of $P_{Z}$), average-free and homogeneous of degree  ${\mathcal O}_{T}(Z)$; and a constant $\gamma >0$  such that 
$$\lim_{\rho \to 0} \, \rho^{-n-2{\mathcal O}_{T}(Z) -\gamma}\int_{{\mathbf B}_{\rho}(Z)} {\rm dist}^{2} \, (X-Z, {\rm graph} \, (\psi_{Z} + \varphi_{Z})) \, d\|T\|(X) = 0;$$ 
moreover: 
\begin{itemize}
\item $\psi_{Z}$ is unique up to an additive term of order $o(|X|^{{\mathcal O}_{T}(Z)})$, that is to say, if the above statement also holds with $\widetilde{\psi}_{Z} \, : \, P_{Z} \to P_{Z}^{\perp}$ in place of $\psi_{Z}$, then 
\begin{equation*}
    \lim_{\rho\rightarrow 0^+} \rho^{-n-2{\mathcal O}_{T}(Z)} \int_{P_Z \cap \mathbf{B}_{\rho}(0)} |\widetilde{\psi}_{Z} - \psi_{Z}|^2 = 0 ; 
    \footnote{It need not be true that $\psi_{Z}(X)  = h_{Z}(X) + o(|X|^{{\mathcal O}_{T}(Z)})$ as $|X| \to 0;$ e.g.\ consider the multiplicity 1 current $T$ associated with the algebraic variety $\{(z, w) \, : \, (z - w^{2} - w^{3})^{2}  = w^{7}\} \subset {\mathbb C} \times {\mathbb C}.$ The origin $Z_{0} = (0, 0)$ is the only singularity of $T$, which is a branch point where the unique tangent plane to $T$ is $P_{Z_{0}} = \{(0, w) \, : \,  w \in {\mathbb C}\}$ with multiplicity $q_{Z_{0}} = 2$, and where ${\mathcal N}_{T, {\rm Pl}}(Z) = 2$. For this example, $h_{Z_{0}}(0, w) = w^{2}$, ${\mathcal O}_{T}(Z_{0}) = 7/2 > {\mathcal N}_{T, {Pl}}(Z)$, $\varphi_{Z_{0}}(0, w) = \llbracket w^{7/2}\rrbracket + \llbracket -w^{7/2}\rrbracket$, and one possible choice of $\psi_{Z_{0}}$ is $\psi_{Z_{0}}(0, w) = h_{Z_{0}}(0,w) + w^{3}.$}
\end{equation*} 
\item ${\mathcal O}_{T}(Z) > {\mathcal N}_{T, {\rm Pl}}(Z)$   $\iff$ ${\mathcal O}_{T}(Z) \geq {\mathcal N}_{T, {\rm Pl}}(Z) + q_{Z}^{-1}$ $\iff$
$$\lim_{\rho \to 0} \, \rho^{-n - 2{\mathcal N}_{T,{\rm Pl}}(Z)  - \gamma^{\prime}}\int_{{\mathbf B}_{\rho}(Z)} {\rm dist}^{2} \, (X-Z, {\rm graph} \, h_{Z}) \, d\|T\|(X) = 0 \;\; \mbox{for some constant $\gamma^{\prime} >0$};$$ 
in particular, if  ${\mathcal N}_{T, {\rm Pl}}(Z)$ is not an integer then ${\mathcal O}_{T}(Z)= {\mathcal N}_{T, {\rm Pl}}(Z)$ (which says, loosely speaking, that the ``sheets'' of the current decay towards each other precisely at the same rate as the current decays towards its tangent plane).
\end{itemize}
\end{itemize}
\end{itemize}
\end{theorem} 

\subsection{Topology of the current near branch points: relation to the work of Chang and Micallef--White} By the work of S.~Chang (\cite{Chang}), singularities of a 2-dimensional area minimizing current $T$ in ${\mathbb R}^{2 +m}$  are isolated, and near any singularity $Z$ where $T$ is locally irreducible (i.e.\ there is no $\rho >0$ such that $T \res {\mathbf B}_{\rho}(Z) = T_{1} + T_{2}$ for non-zero $n$-dimensional rectifiable currents $T_1$ and $T_2$ with $\op{spt} T_1 \cap \op{spt} T_2 = \{Z\}$), the support of $T$ is an immersed $C^{3}$ branched disk. By irreducibility, the multiplicity of $T$ near $Z$ is a constant, which we may assume without loss of generality is 1. 
Let $\Theta \, (T, Z) = q$, and note that $q$ is an integer  $\geq 2$. By combining Chang's theorem with the work of Micallef and White (\cite{MicWhi95}), and assuming further, without loss of generality, that $Z=0$, we have the following: after mapping by an ambient rotation, ${\rm spt} \, T$ near $Z$ can be 
 parameterized by a map of the form $$z \mapsto (z^{q}, f(z)), \;\; z \in B_{1}(0) \subset {\mathbb C} \approx {\mathbb R}^{2} \times \{0\},$$ where $f$ is a $C^{2}$ function taking values in  ${\mathbb R}^{m}$ with $D^{j}f(z) = O(|z|^{q+1-j})$  for $j=0, 1, 2,$ and, for any fixed $q^{\rm th}$ root of unity $\nu \neq 1$,  
 $$f(\nu z) - f(z) = {\rm Re} (az^{p}) + e(z)$$ where 
 \begin{itemize}
 \item $e$ (depending on $\nu$) satisfies $D^{j}e(z) = o(|z|^{p + 1-j})$  for $j=0, 1$; 
 \item $a \in {\mathbb C}^{m} \setminus \{0\}$ is a constant (depending on $\nu$) satisfying $a \cdot a  = 0$ and 
 \item $p$ is an integer (independent of $\nu$) with $p> q$ and $p$ not divisible by $q$ (\cite[Theorems~A~\&~B]{MicWhi95}). 
 \end{itemize}
 Note that in the notation and terminology of Theorem~\ref{summary-theorem-2}, this in particular says that
 $${\mathcal O}_{T}(Z) = p/q$$
 and that the unique ${\mathcal A}_{q}({\mathbb R}^{m})$ valued tangent function $\varphi_{Z}$ of $T$ at $Z$ is given by 
 $$\varphi_{Z}(z) = \sum_{j=0}^{q - 1}\llbracket {\rm Re} \, (a  \, \nu^{j} z^{p/q}) \rrbracket.$$ 
 These striking results say that the current near $Z$, modulo its average height relative to the tangent plane at $Z$, is a $C^{1}$ perturbation of its unique, non-zero tangent function $\varphi_{Z}$, thus providing a complete picture of the behaviour of a 2-dimensional area minimzing current near a branch point.   
 
The extent to which these topological conclusions near branch points may extend to area minimizing currents of general dimension $n \geq 3$ is an important question that has remained largely open. 
Theorem~\ref{summary-theorem-2} provides an inroad into this question via its uniqueness-of-tangent-functions provision. 

To explain this, recall that in general dimensions, by Theorem~\ref{summary-theorem-2}, at ${\mathcal H}^{n-2}$ a.e.\ branch point
 $Z$, $T$ has a unique tangent plane $P_{Z}$ of some integer multiplicity $q \geq 2$ and a unique non-zero, Dirichlet energy minimizing, average-free 
 $q$-valued cylindrical tangent function $\varphi_{Z} \, : \, P_{Z} \to {\mathcal A}_{q}(P_{Z}^{\perp})$ which is homogeneous of 
 degree a rational number ${\mathcal O}_{T}(Z) \geq 1 + 1/q$.  Such $Z$ may or may not be an isolated singularity. In any event, note first that without further hypotheses, 
 the direct analogue of the combined Chang and Micallef--White conclusion (giving that the support of $T$ near $Z$ is homeomorphic to an $n$-disk) 
need not hold;  this is the case even when ${\mathcal O}_{T}(Z) = p/q$ for $p$ co-prime to $q.$ 
 This is demonstrated by  multiplicity 1 currents associated with certain algebraic varieties, such as the current $T_{V_{1}}$ associated with   
 $$V_{1} = \{(x, y,w) \, : \, w^{2} = x(x - y^{2})(x+y^{2})\} \subset {\mathbb C}^{3}.$$ Note that $T_{V_{1}}$ has an isolated singularity at $Z_{0} = (0, 0, 0)$ which is a branch point where 
 the tangent plane is $2\llbracket \{w = 0\} \rrbracket$ and ${\mathcal O}_{T_{V_{1}}}(Z_{0}) = {\mathcal N}_{T_{V_{1}}, {\rm Pl}}(Z_{0}) = 3/2;$ the (unique) tangent function of $T_{V_{1}}$ at $Z_{0}$ is cylindrical, and is (up to a non-zero constant) given by $\varphi_{Z_{0}}(x, y) = \llbracket {\rm Re} \, (x^{3/2}) \rrbracket + \llbracket -{\rm Re} \, (x^{3/2}) \rrbracket.$ An example with similar behaviour near a non-isolated branch point is provided by the current $T_{V_{2}}$ associated with 
 $$V_{2} = \{(x, y, w) \, : \, w^{2} = x^{3}(x^{2} - y^{3}) \} \subset {\mathbb C}^{3}.$$ 
This has ${\rm sing} \, T_{V_{2}} = \{(0, a, 0) \, : \, a \in {\mathbb C}\}$, and each $Z \in {\rm sing} \, T_{V_{2}}$ is a density 2 branch point with unique tangent plane 
$2\llbracket\{w=0\}\rrbracket$; moreover, the unique tangent function of $T_{V_{2}}$ at $Z$ is cylindrical, and is (up to a non-zero constant) given by $\varphi_{Z}(x, y) = \llbracket {\rm Re} \, (x^{5/2}) \rrbracket + \llbracket -{\rm Re} \, (x^{5/2}) \rrbracket$ for $Z = (0, 0, 0)$ and  $\varphi_{Z}(x, y) = \llbracket {\rm Re} \, (ia^{3/2}x^{3/2}) \rrbracket + \llbracket -{\rm Re} \, (ia^{3/2}x^{3/2}) \rrbracket$ for $Z  = (0, a, 0)$, $a \neq 0$. In particular, 
${\mathcal O}_{T_{V_{2}}}(0) = {\mathcal N}_{T_{V_{2}}, {\rm Pl}}(0)= 5/2$ and ${\mathcal O}_{T_{V_{2}}}(Z) = 3/2$ if $Z = (0, a, 0)$, $a \neq 0$. 
Neither of these examples is homeomorphic to a 4-disk near the origin, as can be checked by verifying that the fundamental group of either of them in a small punctured ball about the origin is non-trivial.

 However, a direct analogue of the Chang and Micallef--White results holds  in arbitrary dimension (Theorem~\ref{topology} below) at a density $q$ branch point $Z_{0} \in {\mathcal B}_{q}$ (where ${\mathcal B}_{q}$ is as described at the beginning of Section~\ref{higher-order-asym}) if ${\mathcal O}_{T}(Z_{0}) = p/q$ for some $p$ co-prime to $q$ and, additionally, $Z_{0}$ is non-isolated within the set of points $Z \in {\mathcal B}_{q}$ with ${\mathcal O}_{T}(Z) \geq p/q$ in a  precise strong sense specified in the definition below.  This non-isolatedness requirement cannot be relaxed in view of the example $T_{V_{2}}$.

\begin{definition}[$(n-2)$-strong non-isolatedness of a branch point]\label{strongly-non-isolated} For $Z_{0} \in {\mathcal B}_{q}$ and $\rho > 0$, we say that $Z_{0}$ is \emph{$(n-2)$-strongly non-isolated in ${\mathbf B}_{\rho}(Z_{0})$} if there exists an $(n-2)$-dimensional linear subspace $L$ of ${\mathbb R}^{n+m}$ such that 
\begin{equation}\label{topology hyp2}
 L \cap \mathbf{B}_{\rho/2}(0) 
		\subseteq \pi_L(\{ Z - Z_0 : Z \in \mathcal{B}_q \cap \mathbf{B}_{\rho}(Z_0) \text{ and } {\mathcal O}_{T}(Z) \geq  {\mathcal O}_{T}(Z_{0})\}) , 
\end{equation}
where $\pi_L : \mathbb{R}^{n+m} \rightarrow L$ is the orthogonal projection map onto $L.$ 

We say that $Z_{0} \in {\mathcal B}_{q}$ is \emph{$(n-2)$-strongly non-isolated} if there is a sequence $\rho_{j} \to 0^{+}$ such that $Z_{0}$ is 
$(n-2)$-strongly non-isolated in 
${\mathbf B}_{\rho_{j}}(Z_{0})$ for each $j$.
\end{definition}

\begin{theorem}[\cite{KrumWicc}, \cite{KrumWicd}]\label{topology}
Let $T$ be an $n$-dimensional locally area-minimizing rectifiable current in an open subset of ${\mathbb R}^{n+m}$, $q$ an integer $\geq 2$ and let ${\mathcal B}_{q}$ be as described at the beginning of Section~\ref{higher-order-asym}, noting that ${\mathcal B}_{q}$ contains ${\mathcal H}^{n-2}$ almost all density-$q$ branch points. Let $Z_{0} \in {\mathcal B}_{q}.$ If   
\begin{itemize}
\item [{\rm (a)}] 
${\mathcal O}_{T}(Z_{0}) = p/q$ for some integer $p$ co-prime to $q$ and 
\item[{\rm (b)}] $Z_{0}$ is $(n-2)$-strongly non-isolated (cf.\ Definition~\ref{strongly-non-isolated}) 
\end{itemize}
then ${\rm spt} \, T$ near $Z_{0}$ is homeomorphic to an $n$-dimensional disk and admits a $C^{1, \mu}$ parameterization for some fixed $\mu = \mu(n, m, q) \in (0, 1)$.

More precisely (and more generally), the following holds: Suppose that $Z_{0} \in {\mathcal B}_{q}$ satisfies hypothesis {\rm (a)}. There exists $\rho_{0} \in (0, 1)$ (depending on $Z_{0}$) such that if additionally, in place of hypothesis {\rm (b)}, we have that: 
\begin{itemize} 
\item [{\rm (}${\rm b}^{\prime}${\rm )}] $Z_{0}$ is \emph{$(n-2)$-strongly non-isolated in ${\mathbf B}_{\rho}(Z_{0})$} for some $\rho \in (0, \rho_{0}]$ (cf.\ Definition~\ref{strongly-non-isolated}),  
\end{itemize}
then there exist a number $\mu = \mu(n, m, q) \in (0, 1)$; a rotation $\Gamma$ of $\mathbb{R}^{n+m}$ with $\Gamma(P_{Z_0}) = \mathbb{R}^n \times \{0\}$ (where 
$P_{Z_{0}}$ is the unique tangent plane to $T$ at $Z_{0}$); functions $g \in C^{1,\mu/2}(B^{n-2}_8(0);\mathbb{R}^2)$, $h \in C^{1,\mu/2}(B^{n-2}_8(0);\mathbb{R}^m)$, $\psi \in C^{1,\mu}(B^2_1(0) \times B^{n-2}_1(0);\mathbb{R}^m)$ and $f \in C^{1,\mu/2}(B^2_1(0) \times B^{n-2}_1(0);\mathbb{R}^m)$ such that: 
\begin{enumerate}[itemsep=2mm,topsep=0mm]
	\item[{\rm (i)}]  $\op{sing} \, T \cap {\mathbf B}_{\rho/8}(Z_{0}) = \{Z \in {\mathcal B}_{q} \, : \, {\mathcal O}_{T}(Z) = p/q\} \cap {\mathbf B}_{\rho/8}(Z_{0})$ and 
    \begin{equation*}
		\Gamma(\eta_{Z_0,\rho/64}(\op{sing} T)) \cap \mathbf{C}_8(0) 
		= \{ (g(y), y, h(y)) : y \in B^{n-2}_8(0) \} \cap \mathbf{C}_8(0) ; 
	\end{equation*}
	\item[{\rm (ii)}]  the support of $T$ near $Z_{0}$ is parameterized by
     \begin{align*}
		&\op{spt} \Gamma_{\#} \eta_{Z_0,\rho/64\#} T \cap \mathcal{D} \times \mathbb{R}^m 
		\\=\,& \{ (g(y) + z^q, y, f(z,y)) : z \in B^2_1(0),\, y \in B^{n-2}_1(0) \} , 
	\end{align*}
	where we identify $\mathbb{R}^2 \cong \mathbb{C}$ and $\mathcal{D} = \{ (g(y)+x,y) : x \in B^2_1(0),\, y \in B^{n-2}_1(0) \}$; the mapping $(z,y) \mapsto (g(y)+z^q,y,f(z,y))$ is injective, and hence $\op{spt} \Gamma_{\#} \eta_{Z_0,\rho/64\#} T \cap \mathcal{D} \times \mathbb{R}^m$ is homeomorphic to an $n$-disk; 

	\item[{\rm (iii)}]  there exists a unique tangent function $\varphi$ to $\Gamma_{\#} \eta_{Z_0,\rho/64\#} T$ at the origin which is given by $\varphi(x_1,x_2,y) = \op{Re}(c \,(x_1+ix_2)^{p/q})$ for each $(x_1,x_2,y) \in \mathbb{R}^n$ and for some $c \in \mathbb{C}^m \setminus \{0\}$ with $c \cdot c = 0$ and 
	\begin{equation*}
		f(z,y) = \psi(z^q, y) + \op{Re}(cz^p) + e(z,y) 
	\end{equation*} 
	for $z \in B^2_1(0)$ and $y \in B^{n-2}_1(0)$, where 
	\begin{gather*}
		|e(z,y)| \leq C \,|(z^q,y)|^{p/q+\mu/(2nq)} , \quad 
		|D_z e(z,y)| \leq C \,|(z^q,y)|^{p/q-1+\mu/(2nq)} |z|^{q-1} , \\
		|D_y e(z,y)| \leq C \,|(z^q,y)|^{p/q-1+\mu/(2nq)} , 
	\end{gather*}
	for some constant  
    $C = C(n,m,q,p) \in (0,\infty)$; 
	
	\item[{\rm (iv)}]  for any fixed $q^{\rm th}$ root of unity $\nu \neq 1$
	\begin{equation*}
		f(z,y) - f(\nu z,y) = \op{Re}(az^p) + \widetilde{e}(z,y) 
	\end{equation*} 
	for $z \in B^2_1(0)$ and $y \in B^{n-2}_1(0)$, where $a \in \mathbb{C}^m \setminus \{0\}$ with $a \cdot a = 0$ and 
	\begin{gather*}
		|\widetilde{e}(z,y)| \leq C \,|(z^q,y)|^{p/q+\mu/(2nq)} , \quad 
		|D_z \widetilde{e}(z,y)| \leq C \,|(z^q,y)|^{p/q-1+\mu/(2nq)} |z|^{q-1} , \\
		|D_y \widetilde{e}(z,y)| \leq C \,|(z^q,y)|^{p/q-1+\mu/(2nq)} , 
	\end{gather*}
	where $C = C(n,m,q,p) \in (0,\infty)$ is a constant.
\end{enumerate}
\end{theorem}

\begin{remark} {\rm Suppose that $Z_{0}$ is a density $q$ branch point with ${\mathcal N}_{T, {\rm Pl}}(Z_{0}) = 1 + 1/q$. 
Then hypothesis (a) is automatically satisfied. In this case hypothesis (b) takes the weaker form: \emph{near $Z_{0}$ the set of all density $q$ branch points $Z$ projects fully onto an open ball of an $(n-2)$-dimensional affine subspace}. Example $T_{V_{2}}$ discussed above shows that in general (i.e.\ when ${\mathcal N}_{T, {\rm Pl}}(Z_{0}) > 1 + 1/q$ and hypothesis (a) holds), hypothesis (b) cannot be replaced by this weaker statement. Example $T_{V_{1}}$ shows that this weaker assumption 
cannot be dropped when ${\mathcal N}_{T, {\rm Pl}}(Z_{0}) = 1 + 1/q$.}
\end{remark}
 
A natural extension of local irreducibility of two-dimensional area-minimizing currents \cite[Definition~3.1]{Chang} to  area-minimizing currents of general dimension $n \geq 2$ is as follows:

\begin{definition}\label{irreducible defn}
We say that $T$ is locally irreducible near $Z_{0} \in \op{spt} T$ if there is no radius $\rho >0$ such that $T \res {\mathbf B}_{\rho}(Z_{0}) = T_{1} + T_{2}$ for some non-zero $n$-dimensional rectifiable currents $T_1$ and $T_2$ with $Z_0 \in \op{spt} T_1 \cap \op{spt} T_2$ and $\op{spt} T_1 \cap \op{spt} T_2 \subseteq \op{sing} T$.
\end{definition}

When $n = 2$, $Z_{0}$ is an isolated singular point of $T$ and thus the requirements $Z_0 \in \op{spt} T_1 \cap \op{spt} T_2$ and $\op{spt} T_1 \cap \op{spt} T_2 \subseteq \op{sing} T$ is equivalent to $\op{spt} T_1 \cap \op{spt} T_2 = \{Z_{0}\}$ as in \cite{Chang}.  In higher dimensions $n > 2$, by~\cite{Almgren} the Hausdorff dimension of $\op{sing} T \leq n-2$, and thus if $T \res {\mathbf B}_{\rho}(Z_{0}) = T_{1} + T_{2}$ for $T_{1}$ and $T_{2}$ as in Definition~\ref{irreducible defn}, then $T_{1}$ and $T_{2}$ are locally area-minimizing rectifiable currents of ${\mathbf B}_{\rho}(Z_{0})$ such that $(\partial T_{1}) \res {\mathbf B}_{\rho}(Z_{0}) = 0$, $(\partial T_{2}) \res {\mathbf B}_{\rho}(Z_{0}) = 0$, and $\|T\|\res {\mathbf B}_{\rho}(Z_{0}) = \|T_{1}\| + \|T_{2}\|$.  Definition~\ref{irreducible defn} allows the supports of $T_{1}$ and $T_{2}$ to intersect at a point of the regular set of either $T_{1}$ or $T_{2}$ (for example, consider $T = \{(z,w) : zw = 0 \}$); however, any such point must be in $\op{sing} T$.  If for some radius $\rho > 0$, $Z_{0}$ is contained in the closure of two or more connected components of $\op{reg} T \cap {\mathbf B}_{\rho}(Z_{0})$, then $T$ is not irreducible near $Z_{0}$ as the current $T$ decomposes into locally area-minimizing currents supported on the connected components $\op{reg} T \cap {\mathbf B}_{\rho}(Z_{0})$ with constant integer multiplicity.  On the other hand, if for every $\rho > 0$, $Z_{0}$ is contained in the closure of exactly one connected component $M$ of $\op{reg} T \cap {\mathbf B}_{\rho}(Z_{0})$ (where $M$ depends on $\rho$), then $T$ is irreducible near $Z_{0}$ even if $T$ has integer mulipliticity $\geq 2$ on $M$.  Indeed, if we did have $T \res {\mathbf B}_{\rho}(Z_{0}) = T_{1} + T_{2}$ as in Definition~\ref{irreducible defn}, then the support of $T$ is the closure of $M$ near $Z_{0}$ and thus by the constancy theorem the support of each $T_{i}$ is either empty or the closure of $M$ near $Z_{0}$, contradicting the assumptions on $T_{1}$ and $T_{2}$.

\begin{remark} {\rm  
If $q$ is a prime and $T$ is locally irreducible near $Z_{0}$, then again  Theorem~\ref{topology} holds without hypothesis~(a).}
\end{remark}

\label{prog-summary} The asymptotic analysis leading to Theorems~\ref{summary-theorem-1}, \ref{summary-theorem-2} and \ref{topology}, carried out in \cite{KrumWicb}, \cite{KrumWicc} and  \cite{KrumWicd},  is in part based on the techniques developed in \cite{Wic14} for the analysis of stable codimension 1  integral varifolds, and those  developed in \cite{KrumWic2}, \cite{KrumWic1} for studying branch sets of Dirichlet energy minimizing multi-valued functions and two-valued $C^{1, \alpha}$ minimal graphs respectively. These earlier works were in turn inspired by the fundamental work of Simon (\cite{Sim93}) on the structure of singular sets of minimal submanifolds in a multiplicity 1 class. Among the key new ingredients needed for adaptation of these ideas for our purposes  
are several interior height bounds for area minimising currents. The first of these results (proved in \cite[Theorem~3.11]{KrumWicb}) gives a uniform interior upper bound for the height of an area minimising rectifiable current relative to a union ${\mathbf P}$ of non-intersecting oriented planes in terms of a \emph{linear expression} in the height excess of the current relative to ${\mathbf P}$, under appropriate smallness assumptions on the tilt-excess of the current and of ${\mathbf P}$. The key implication of this estimate for our purposes is the obvious one: the current must separate into a sum of disjoint pieces whenever its height excess relative to ${\mathbf P}$ is much smaller than the smallest distance between any pair of the planes  making up ${\mathbf P}$. 
Two substantial generalizations of this are needed for the higher order analysis and uniqueness-of-tangent-functions results in \cite{KrumWicc}, \cite{KrumWicd}, 
where ${\mathbf P}$ is assumed to constitute graphs of harmonic functions (instead of planes); see \cite[Theorem~4.19]{KrumWicc} and 
\cite[Theorem~4.19]{KrumWicd}. 

\subsection{Program summary}\label{prog-summary}
Our overall program develops a unified framework for the analysis of singularities  of area-minimizing currents by synthesizing several basic questions: tangent cone uniqueness; singular set dimension and structure; higher-order asymptotics for the current; and conditions guaranteeing local topological simplicity of the current near branch points. The methods used within this framework establish uniform a priori estimates as the primary organizing principle. In doing so, this framework brings the higher-codimension theory closer in spirit to the well-established codimension 1 regularity theory (for area minimizers, or more generally, for stable varifolds) and to classical PDE theory. 

A byproduct of this work is a geometrically more direct and technically more efficient proof of Almgren's renowned optimal Hausdorff dimension bound for the singular set. Efficiency is gained by establishing uniform decay estimates at branch points \emph{prior} to constructing the center manifold. 
Once the decay estimates are in place, the center manifold is invoked in our program only to study branch points of quadratic decay.  This is the canonical setting for the center manifold which conforms to PDE principles: here, the center manifold about a given quadratically-decaying branch point contains all nearby quadratically decaying branch points of the same density, while a (multi-valued) normal map on this center manifold which contains all such points in its critical nodal set 
provides a graphical approximation for the current with quantitative, locally uniform accuracy. This is in contrast to the classical approach where for each branch point, an iterative process is employed to construct a sequence of center manifolds which are not required to contain any branch point. In our framework, branch points for which quadratic decay fails are analyzed without center manifolds and through intrinsic geometric arguments based on the approximate monotonicity of the planar frequency function introduced at the outset of the program.

The program yields four principal geometric and topological  conclusions for $n$-dimensional area-minimizing rectifiable currents: 

\begin{itemize}
\item[1.] \emph{Tangent cone uniqueness}: At ${\mathcal H}^{n-2}$-a.e.\ point, the current admits a unique tangent cone supported either on a single plane or a finite union of two or more planes intersecting along a common $(n-2)$-dimensional subspace. 

\item[2.] \emph{Higher-order asymptotics}: At 
${\mathcal H}^{n-2}$-a.e.\ branch point, the current admits a unique expansion of finite order greater than $1 + \epsilon$ with precise remainder decay where  $\epsilon >0$ is a  universal constant, providing an explicit analytic branching model. 

\item[3.] \emph{Structural decomposition of the singular set}: 
The singular set 
locally decomposes into finitely many pairwise disjoint, locally compact, locally $(n-2)$-rectifiable sets (of locally finite measure).

\item[4.] \emph{Topological control via frequency criteria}: Near a branch point satisfying a specific planar frequency/branching order criterion, the support of the current is homeomorphic to an $n$-dimensional disk admitting a 
$C^{1, \mu}$ parameterization. (When $n=2$, this criterion is known to hold at all irreducible points; in higher dimensions, classical algebraic examples show that when the criterion fails, the current may not be locally homeomorphic to a disk).
\end{itemize}

\noindent
\begin{remark}[\emph{Riemannian ambient spaces}] \label{riemannian-ambient-remark}
{\rm All of our main results extend to area minimizing rectifiable currents in sufficiently smooth Riemannian ambient spaces. This is seen by embedding the ambient Riemannian manifold in a Euclidean space and employing fairly straightforward and technical modifications to the arguments 
used for Euclidean ambient space. With regard to the planar frequency function, the key point is that it is possible to establish (via  technical modifications to the 
argument for Euclidean planar frequency function) approximate monotonicity of a certain truncated planar frequency function, which (essentially) agrees with the Euclidean frequency function whenever the current decays sub-quadratically towards the plane.   This then leads to Theorem~1.1, with part~(II)(iii) replaced by the alternatives that for each $Z \in {\mathcal B}_{q}$, either (A) planar frequency at $Z$ exists and has value $\in (1+\alpha,2-\alpha)$, or (B) decay of $T$ about $Z$ to the tangent plane $P_{Z}$ is at a rate $O(|X - Z|^{2-\alpha})$ or faster. Once Theorem~1.1 in this form is in place, all of the other main results,  with appropriate modifications to the statements, are obtained with routine technical adjustments to the arguments in the Euclidean case. We refer the reader to \cite{KrumWice} 
for details. In particular, the key feature of not needing center manifolds to analyze slow-decay branch points and a center manifold is only invoked in its canonical setting (where all relevant branch points are required to be contained in the critical nodal set of a single normal approximation map) continues to hold. As mentioned above this feature is indispensable in our approach to obtaining higher order asymptotics at typical branch points and topological information for the current near branch points.} 
\end{remark} 
 
\section{The work of De~Lellis, Minter and Skorobogatova: a comparison}
In contemporaneous, independent  work \cite{DelSko1}, \cite{DelSko2}, \cite{DelMinSko}
(available on the arXiv at arXiv:2304.11552, arXiv:2304.11555, arXiv:2304.11553 shortly after the initial posting of the present paper and \cite{KrumWicb} on the arXiv at arXiv:2304.10653 and arXiv:2304.10272), 
  De~Lellis, Minter and Skorobogatova  also established two of the main results presented in our work for $n$-dimensional area minimizing rectifiable currents:
\begin{itemize}
\item[(a)] uniqueness of tangent cones  at ${\mathcal H}^{n-2}$ a.e.\ point (Conclusion 1, Section~\ref{prog-summary});  
\item[(b)] countable $(n-2)$-rectifiability of the singular set (part of Conclusion 3, Section~\ref{prog-summary}).
\end{itemize}
Their approach, however, conceptually differs from ours in two fundamental ways:
\begin{itemize}
\item[(i)] in their work, countable rectifiability for the branch set (included in conclusion (b)) is derived from upper Minkowski content bounds for
the branch set; this is achieved without addressing the question of uniqueness of tangent functions at branch points. 
\item [(ii)] a foundational component of their proof of both results (a) and (b) is an iterative 
construction---employed at the outset of the program---of a sequence of auxiliary center manifolds for each branch point; these center manifolds are not required to contain any branch point. 
\end{itemize}
By contrast:
\begin{itemize}
\item[(i)] our primary objective is to determine the asymptotic structure of the current on approach to typical singularities, in particular to establish a higher order asymptotic expansion at typical branch points by proving uniqueness of tangent functions. Both the singular set structure (as in Conclusion~3, Section~\ref{prog-summary}) and local topological information (as in Conclusion~4, Section~\ref{prog-summary}) follow as  corollaries of this; 
\item[(ii)] our methodology relies on planar-frequency-based intrinsic geometric arguments instead of center manifold constructions for arbitrary branch points.
The use of the center manifold is restricted to a canonical situation where it is indispensable: a single center manifold locally containing all relevant branch points and acting as the domain for a single normal approximating map whose critical nodal set contains those branch points. Our structural analysis of the current makes indispensable use of these additional geometric features of the center manifold and the normal map in the setting where they are needed.   
\end{itemize}

In this section, we elaborate on these differences and note (in subsection~\ref{comparison}) a close technical similarity between a specific portion of \cite{DelMinSko} and \cite{KrumWicb}.\footnote{We are grateful to the authors of \cite{DelSko1}, \cite{DelSko2}, \cite{DelMinSko} for an email exchange (in June--July, 2023) on the differences between our work and theirs, and in particular for helpful responses to our pointing out of the relevance of \emph{local uniformity} of the estimates in the conclusions of Theorem~\ref{branch and cones thm-intro}(I)\&(II)  above to local-finiteness-of-measure properties of ${\mathcal S}_{q}$ and ${\mathcal B}_{q}$ (which are respectively the sets ${\mathcal S}_{q, \epsilon, \beta}$ and ${\mathcal B}_{q, \epsilon, \beta}$ as in Corollary~\ref{immediate-cor} once $\epsilon$ and $\beta$ are fixed as described in subsection~\ref{non-branch-structure}). These estimates lead to Theorem~\ref{improved-local-finiteness} below which is an improvement of both \cite{DelSko2} and some of our results; the authors of \cite{DelSko2} helpfully confirmed to us then that with the help of the locally uniform estimate in Theorem~\ref{branch and cones thm-intro}(II)(i), 
their argument in \cite{DelSko2} does indeed carry over to yield the local finiteness of measure conclusions concerning the set of density $q$ branch points, as in Theorem~\ref{improved-local-finiteness}.}

 \subsection{Relation to Almgren's center manifolds}
 Let $Z$ be any branch point of $T$.  The first step of their approach is to select a sequence of intervals of scales at which the rescaling of $T$ about $Z$ has sufficiently small excess, and construct a sequence of center manifolds corresponding to each of these interval of scales.

The use of center manifolds in this generality (i.e.\ for \emph{every} branch point) is a foundational aspect of Almgren's argument.  A technically streamlined, more accessible presentation of the construction was given subsequently by De~Lellis and Spadaro (\cite{DeLSpa2}). 
Heuristically speaking, the purpose of the center manifolds is to facilitate a way to capture an infinitesimal rate (as a fixed positive power of the distance to $Z$) at which $T$ decays, on approach to $Z$, to (the graph of) the ``average height'' of $T$ over a chosen tangent plane $P_{Z}$ of $T$ at $Z$; the center manifolds are sufficiently smooth ($n$-dimensional) graphs, satisfying certain estimates and serving to approximate this average height in appropriate annular regions about $Z$. It should be emphasized that uniqueness of the tangent plane $P_{Z}$ is not known at this stage, which complicates this task. Also, even if $P_{Z}$ is unique, the average height of $T$ makes sense only in an approximate sense, since a parameterization of $T$ near $Z$ as a multi-valued graph over $P_{Z}$ is not available, and there is only an approximate multi-valued Lipschitz graph description at scales where $T$ is weakly close to $P_{Z}$.

For a given branch point $Z$ of some density $q,$ one would ideally like to work with a \emph{single} center manifold about $Z$ that contains $Z$ and all nearby density $q$ branch points $Z^{\prime}$, but the construction in its full generality does not guarantee this. There is an exception to this, however; if at $Z$, the current $T$ decays to a unique tangent plane \emph{quadratically or faster} in the distance to $Z$, then the procedure does yield the \emph{simplest case of a center manifold}: namely, a center manifold ${\mathcal M}$ that automatically contains $Z$ and all nearby density $q$ branch points $Z^{\prime}$ at which uniform quadratic decay towards a tangent plane holds. (We note that this fact, which follows readily from the construction of ${\mathcal M}$, is explicitly observed in the Ph.D.\ thesis work of S.~Chang (on the local structure of 2-dimensional area minimizing currents) supervised by Almgren (see \cite[p.~703]{Chang}); however, there seems to be little indication either in \cite{Almgren} or in \cite{Chang} as to how one might capitalise on this fact. As mentioned in the introduction above, our approach can be seen as providing a way to take full advantage of this fact, as our method reduces the need for center manifolds in the entire program to this simplest case.)

Even if a decay estimate at a branch point is available, the decay rate towards the tangent plane can be less than quadratic. (Consider, e.g.\ the algebraic curve $\{(z, w) \, : \, z^{2} = w^{3}\} \subset {\mathbb C} \times {\mathbb C}$). For this essential reason, the above simplest case of a single center manifold is inadequate for the strategy in  \cite{Almgren} (and by extension, for \cite{DelSko1}, \cite{DelSko2}, \cite{DelMinSko}) since the argument hinges on capturing a priori, for every branch point, 
a rate of decay of $T$ towards its ``average height.''  

Instead, corresponding to each branch point $Z,$ an infinite sequence of center manifolds ${\mathcal M}^{(Z)}_{j}$ (which may not contain any of the branch points) is constructed in \cite{Almgren}, \cite{DelSko1},
followed by  ($q$-valued) normal approximation maps $N^{(Z)}_{j}$ which parameterize most of $T$ 
over the center manifolds. This process requires
the careful selection of a number of parameters and a sequence of disjoint annular regions ${\mathcal A}_{j}^{(Z)}$ about $Z$ in which the excess of $T$ 
remains small. 
The construction ensures that the center manifolds are sufficiently regular so that relevant first variation estimates on the normal maps $N^{(Z)}_{j}$ can be carried out to establish approximate monotonicity of the Almgren frequency function associated with $N_{j}^{(Z)}$. 

The overall argument involves substantial conceptual and technical intricacy, but 
perseverance nonetheless is rewarded when it is ultimately shown that on approach to $Z$, the rate of decay of $T$ towards its average height (that is, the rate of decay of the $N_{j}^{(Z)}$ towards zero) exists as a \emph{finite} positive number. This result---the analytic crown jewel in \cite{Almgren}---follows from the approximate monotonicity of the Almgren frequency function associated with $N_{j}^{(Z)}$. It is a remarkable outcome considering that it is arrived at using very little information about branch points---in particular without the knowledge of uniqueness of the tangent plane at branch points. The technical output needed in the argument is in part a reflection of the absence of such information. Much of these complications can be avoided in the special case of quadratic decay of $T$ towards a unique tangent plane at $Z$ since in that case, as mentioned,  
a single center manifold suffices and the relevant estimates are valid over \emph{all} sufficiently small scales about $Z$.

The finiteness of the rate of decay of $N_{j}^{(Z)}$ to $0$ means that every tangent function (blow-up) of $T$ at $Z$ relative to the chosen sequence of center manifolds is non-zero and hence, since it has zero average (by virtue of the fact that the center manifolds well-approximate the average height of $T$ over $P_{Z}$), such a blow-up is a singular $q$-valued Dirichlet energy minimizing function. In \cite{DelSko1}, the rate of decay of $N_{j}^{(Z)}$ is defined in terms of these (possibly non-unique) blow-ups, and is called the ``singularity degree'' of the branch point $Z$, denoted $I(Z)$; among other things, the properties that $I(Z)  \geq 1$ and that if $I(Z) >1$ then $T$ has a unique (multiplicity $q$) tangent plane at $Z$ then follow. Subsequently, in \cite{DelSko2} it is shown that the set of density $q$ 
branch points $Z$ with $I(Z) >1$ is countably $(n-2)$-rectifiable. Finally, in \cite{DelMinSko} it is established that the set of density $q$ 
branch points $Z$ with $I(Z) =1$ is ${\mathcal H}^{n-2}$ null and that the set of density $q$ singularities that are not branch points is countably $(n-2)$-rectifiable with $T$ having a unique tangent cone at ${\mathcal H}^{n-2}$ a.e.\ such point.

Whichever strategy is being used, having some mechanism by which to produce non-zero blow-ups of $T$ at branch points  is crucial for establishing the Hausdorff dimension upper bound (the main goal of \cite{Almgren}) and local structural properties of the singular set. The work \cite{DelSko1}, \cite{DelSko2}, \cite{DelMinSko} builds on the particular mechanism for this provided by \cite{Almgren}, namely, producing blow-ups relative to sequences of center manifolds corresponding to every branch point $Z$.    
In our approach, blow-ups relative to a tangent plane at ${\mathcal H}^{n-2}$ \emph{almost every} branch point $Z$ are shown to be non-zero and this fact is exploited first to reach a substantial number of conclusions.

\subsection{The present work: uniqueness of tangent cones and other results without center manifolds} 
For the question of uniqueness of tangent planes at branch points, the possibility of sheets of $T$ merging towards each other infinitely fast  on approach to a branch point (i.e.\ $T$ decaying infinitely fast towards its average height)  is not an obstacle but in fact a favourable scenario. Thus, one would expect there to be a simpler argument that avoids center manifolds to address this uniqueness question. With this in mind, and in contrast to \cite{Almgren}, \cite{DelSko1}, \cite{DelSko2}, \cite{DelMinSko}, our method treats this uniqueness question not as secondary to the question of finiteness of the rate of decay of $T$ to its average height, but as an issue tied intimately to the problem of bounding the dimension and analyzing the fine structure of the singular set. 

Unifying, as a strategy, the question of the uniqueness of tangent planes with other questions including the dimension and structure of the singular set, and the structure of the current can be seen, at least in spirit, as an extension to higher codimension of the basic ``codimension 1 philosophy'' in regularity theory: in bounding the size of the singular set of an $n$-dimensional, codimension 1 area minimizing rectifiable current $T$ (giving the well-known result that the singular set has Hausdorff dimension $\leq n-7$)---where it is possible to reduce the problem a priori to the multiplicity 1 setting---the first key step is the De~Giorgi estimate (\cite{DG}) giving decay of $T$ towards a (multiplicity 1) hyperplane when $T$ is weakly close to a hyperplane. Likewise, in extending this size bound to the more general setting of stable codimension 1 integral varifolds with no classical singularities---a setting where reduction to multiplicity 1 is not possible a priori---a fundamental step is to prove that whenever the varifold is weakly close to a multiplicity $q$ hyperplane, it decays, at each interior density $q$ point $Z$, towards a unique multiplicity $q$ tangent hyperplane
(\cite[Lemma~15.1]{Wic14}). (Still more generally, by the more recent work \cite{MW},  
this last result extends to density $q$ branch points $Z$ of stable codimension 1 integral varifolds with no classical singularities of density $< q$, and consequently, also to codimension 1 area minimizing currents mod $2q$.)

In a similar spirit, establishing a decay estimate at branch points (to be precise, at ${\mathcal H}^{n-2}$ a.e.~branch point) is the central first goal in our analysis of area minimizing rectifiable currents of codimension $>1$. However, the way we achieve this estimate is considerably different from the above codimension 1 settings. The basic idea in those settings is to prove and use a decay estimate for the functions (single-valued or multi-valued, depending on whether $q=1$ or $q \geq 2$) arising from an  appropriate linearization scheme. For higher codimension area minimizers, regularity for the linearized setting (that is, for Dirichlet energy minimizing multi-valued functions) is not sufficiently strong to be directly useful in the same way as in codimension 1. We therefore take a different approach and proceed with the help of new ideas.
A main novelty is the introduction of the planar frequency function (in Section~\ref{sec:freq mono sec} below) for the minimizing current, and tying the proof of  the (approximate) monotonicity of the planar frequency function intimately to the proof of the decay of the current towards a plane at branch points, in a strategy that  accomplishes both proofs simultaneously. 

Indeed, the first step in implementing this strategy is to assume that $T$ decays at a point $Z$ towards a plane $P$ (at the rate of a positive power of the scale) over an interval $I$ of scales (with $I$ not necessarily having $0$ as an end point), and then prove approximate monotonicity 
of the planar frequency function (corresponding to $P$ and base point $Z$) over $I$ (Theorem~\ref{mono freq thm} below). 

With this result in hand, the next idea is to focus on the behaviour of $T$ near any singular 
point $Z_{0}$ where $T$ fails to decay (all the way, as the distance to $Z_{0}$ tends to $0$) towards a (unique) plane faster than a fixed power of the distance to $Z_{0}$. As the first main application of the approximate monotonicity of the planar frequency function (Theorem~\ref{mono freq thm}), we establish that near any such point $Z_{0}$, the current satisfies the locally uniform weak approximation property given in Theorem~\ref{branch and cones thm-intro}(I)(i).

Employing this weak approximation property, 
we are then able to prove a number of results concerning the singular set of $T$ without involving center manifolds, reducing thereby the dependence on the center manifold, in the entire program, to its simplest case mentioned above. The first among such results are the following: 
\begin{itemize}
\item[(i)] the tangent cone to $T$ at ${\mathcal H}^{n-2}$ a.e.\ point is unique;
\item[(ii)] for ${\mathcal H}^{n-2}$ a.e.\ branch point $Z$ of a given density $q$: 
\begin{itemize}
\item [(a)] a locally uniform decay estimate holds, giving decay of $T$ towards a unique tangent plane $P_{Z}$ at the rate $o(|X - Z|^{1 + \alpha}),$ $X \in {\rm spt} \, T,$ for some fixed $\alpha = \alpha(n, m, q) \in (0, 1)$, and 
\item[(b)]  the planar frequency function $\rho \mapsto N_{T, Z, P_{Z}}(\rho)$ (see Section~\ref{sec:freq mono sec} for the definition) is approximately monotone (non-decreasing). 
\end{itemize}
  \item[(iii)] the set of non-branch-point singularities of $T$ of any given density $q$ is locally $(n-2)$-rectifiable (and has locally finite ${\mathcal H}^{n-2}$ measure).
\end{itemize}
As mentioned above, conclusions (ii)(a) and (ii)(b) are established simultaneously. By (ii)(b), the planar frequency  ${\mathcal N}_{T, {\rm Pl}}(Z) = \lim_{\rho \to 0^{+}} N_{T, Z, P_{Z}}(\rho)$ exists (which can be regarded as the order of contact of $T$ at $Z$ with the unique tangent plane $P_{Z}$), and is a finite number $\geq 1 + \alpha$, for ${\mathcal H}^{n-2}$ branch point $Z$. This allows us to decompose the set of density $q$ branch points as the disjoint union of the sets 
${\mathcal B}_{q}^{(\neq 2)} = \{Z \, : \,  {\mathcal N}_{T, {\rm Pl}}(Z) \in [1 + \alpha, 2) \cup (2, \infty)\}$,  ${\mathcal B}_{q}^{(=2)} 
= \{Z \, : \, {\mathcal N}_{T, {\rm Pl}}(Z) = 2\}$ and an ${\mathcal H}^{n-2}$ null set. 

As a second application of the approximate monotonicity of the planar frequency function, we obtain that any blow-up (or tangent function) of $T$ at a point $Z \in {\mathcal B}_{q}^{(\neq 2)} \cup {\mathcal B}_{q}^{(=2)}$ relative to the (unique) tangent plane $P_{Z}$ at $Z$ is 
a \emph{non-zero}, homogeneous (of degree $\mu = {\mathcal N}_{T, {\rm Pl}}(Z) \geq 1+ \alpha$), $q$-valued Dirichlet energy minimizing function $\varphi  \, : \, P_{Z} \to {\mathcal A}_{q}(P_{Z}^{\perp})$ (which, in case $\mu$ is an integer $\geq 2$, may well be $q$ copies of a single-valued harmonic function). 
As outlined in the introduction above (in the proof-sketch of Theorem~\ref{Alm big reg thm}), this readily implies  that the Hausdorff dimension of ${\mathcal B}_{q}^{(\neq 2)}$ is $\leq n-2.$

In fact much more can be proved concerning ${\mathcal B}_{q}^{(\neq 2)}$ without involving center manifolds. Specifically, as a third application of the monotonicity formula for the planar frequency function we obtain the following more refined description of the nature of both the current $T$ on approach to ${\mathcal H}^{n-2}$ a.e.\ point in ${\mathcal B}_{q}^{(\neq 2)}$ and the set ${\mathcal B}^{(\neq 2)}_{q}$ itself (\cite{KrumWicc}):
\begin{itemize}
\item[(iv)] At ${\mathcal H}^{n-2}$ a.e.\ point $Z \in {\mathcal B}_{q}^{(\neq 2)},$ the current $T$ has a unique blow-up 
$\varphi_{Z} \, : \, P_{Z} \to {\mathcal A}_{q}(P_{Z}^{\perp})$ relative to the (unique) tangent plane $P_{Z}$ of $T$ at $Z$, and moreover, $\varphi_{Z}$ is a locally Dirichlet energy minimizing $q$-valued function which is homogeneous of some degree $\in ([1+\alpha, 2) \cup (2, \infty)) \cap {\mathbb Q}$ and is invariant under translation along an $(n-2)$-dimensional subspace of $P_{Z}$; hence, in terms of co-ordinates $X = (x_1, x_2, \ldots, x_n) \in {\mathbb R}^{n}$ with respect to an appropriate orthonormal basis for $P_{Z}$, the function $\varphi_{Z}$ has the explicit form as described in the paragraph following Remark~\ref{top-stratum-rem}.
\item[(v)] for every closed ball $B$ of the ambient space, ${\mathcal B}_{q}^{(\neq 2)} \cap B$ is the union of finitely many pairwise disjoint sets, each of which is 
locally $(n-2)$-rectifiable (with locally finite $(n-2)$-dimensional Hausdorff measure). 
\item[(vi)]  Near any point $Z \in {\mathcal B}_{q}^{(\neq 2)}$ such that the planar frequency of $T$ is $\ell_Z/q$ where $\ell_Z > q$ co-prime to $q$ and $Z$ is $(n-2)$-strongly non-isolated, the support of the current $T$ is homeomorphic to an n-dimensional disk admitting a $C^{1,\mu}$ parameterization.
\end{itemize}

Once the results (i)--(vi) are in place, it remains to analyze ${\mathcal B}_{q}^{(=2)}$. Any blow-up of $T$ at a point $Z \in {\mathcal B}_{q}^{(=2)}$ relative to the tangent plane $P_{Z}$ at $Z$ is a non-zero Dirichlet energy minimizing $q$-valued function that is homogeneous of degree 2, and in particular there may be points $Z \in {\mathcal B}_{q}^{(=2)}$ where every blow-up is $q$ copies of a single quadratic harmonic polynomial. 
Thus, in order to analyse ${\mathcal B}_{q}^{(=2)}$, we consider blow-ups of $T$ at an arbitrary point $Z \in {\mathcal B}_{q}^{(=2)}$ relative to a center manifold about $Z$. Since ${\mathcal N}_{T, {\rm Pl}}(Z) = 2$, the scaled current $\eta_{Z, \rho \, \#} \, T$ decays to the tangent plane 
$P_{Z}$ quadratically in $\rho$. Hence, the aforementioned simplest case of the center manifold suffices, which  has the crucial property that it contains $Z$ and all other points in ${\mathcal B}_{q}^{(=2)}$ near $Z$. This facilitates the asymptotic analysis of $T$ on approach to points in ${\mathcal B}_{q}^{(=2)}$ as well as the local structural analysis of ${\mathcal B}_{q}^{(=2)}$, leading to conclusions directly analogous to (iv)--(vi) above. (In fact this analysis carries over, with no additional work, to the larger set $\widetilde{\mathcal B}_{q}^{(\geq 2)} \equiv \{Z \, : \, 
{\mathcal N}_{T, {\rm Pl}}(Z) \geq 2\}$; see \cite{KrumWicd}.) 

Combining (iii), (v) and the analogue of (v) for ${\mathcal B}_{q}^{(=2)}$ then leads to the following conclusion regarding the local structure of ${\rm sing} \, T$ 
(which we note is a more refined conclusion for ${\rm sing} \, T$ than countable $(n-2)$-rectifiability obtained in \cite{DelSko1}, \cite{DelSko2}, \cite{DelMinSko}): \emph{for any 
closed ball $B$ of the ambient space, ${\rm sing} \, T \cap B$ is the union of finitely many pairwise disjoint, locally compact sets, each of which is locally $(n-2)$-rectifiable (with locally finite $(n-2)$-dimensional Hausdorff measure) (\cite{KrumWicd}).}

The uniqueness of blow-ups at ${\mathcal H}^{n-2}$ a.e.\ branch point (as in 
(iv) and its counterpart for ${\mathcal B}_{q}^{(=2)}$) is also an important additional conclusion reached in our approach in contrast to that of  \cite{DelSko1}, \cite{DelSko2}, \cite{DelMinSko}. One significance of this uniqueness result is that the estimates involved in its proof can in some cases reveal the topological nature of the current near branch points. For instance, for certain values $\alpha \, (\geq 1 + q^{-1})$, if near a density $q$ branch point $Z$ of planar frequency $\alpha$ 
the set of all density $q$ branch points of planar frequency $\geq \alpha$ projects fully onto an $(n-2)$-dimensional subspace, then the current near $Z$ is a $C^{1, \alpha}$ $q$-valued graph (\cite[Theorem~1.10]{KrumWicc}). 

Finally, we note the following: 
if one uses 
Theorem~\ref{branch and cones thm-intro} above and \cite[Theorem~6.2]{KrumWicb} in conjunction with the argument in \cite{DelSko2},  
one obtains the following result, which improves both the main conclusion of \cite{DelSko2} and (v) (and its counterpart for ${\mathcal B}_{q}^{\prime\prime}$) above: 

\noindent
\begin{theorem}\label{improved-local-finiteness}
For any integer $q \geq 2$, there is an ambient open set $U_{q}$ with ${\rm sing}_{q} \, T \subset U_{q}$ (where ${\rm sing}_{q} \, T$ is the density $q$ singular points of $T$) such that $\{Z \, : \, \Theta \, (T, Z) \geq q\} \cap U_{q}$ is (in addition to being countably $(n-2)$-rectifiable) the disjoint union of three locally compact sets each having locally finite $(n-2)$-dimensional Hausdorff measure. 
\end{theorem} 

\begin{proof}
Let $\alpha= \alpha (n, m, q) \in (0, 1)$ be as in Theorem~\ref{branch and cones thm-intro} taken with $\epsilon = \beta = \beta_{\star}$, where $\beta_{\star} = \beta_{\star}(n, m, q)$ is as in \cite[Theorem~6.1]{KrumWicb}. Applying conclusion \eqref{branch and cones concl1 intro} of Theorem~\ref{branch and cones thm-intro} (II) (as done in Corollary~\ref{immediate-cor}) provides a \emph{locally uniform} decay estimate for $T$ along  the set ${\mathcal B}_{q}$ of density $q$ branch points of rapid decay, i.e.\ 
for each $Z_{0} \in {\mathcal B}_{q}$, there is a number $\rho_{0} = \rho_{0}(Z_{0})$ such that for each 
$Z \in {\mathcal B}_{q} \cap {\mathbf B}_{\rho_{0}}(Z_{0})$, we have that 
\begin{equation}\label{locally-unif-est}
\rho^{-n-2}\int_{{\mathbf B}_{\rho}(Z)} {\rm dist}^{2} (X, Z +P_{Z}) \, d\|T\|(X) \leq C\rho^{2\alpha} \;\;\;\; \forall \;\rho \in (0, \rho_{0}),
\end{equation} 
where $C = C(n, m, q) \in (0, \infty)$. With this fact in hand, the argument in \cite{DelSko2} gives that the two (disjoint) sets 
$\{Z \in {\mathcal B}_{q} \, : \, {\mathcal N}_{T, {\rm Pl}}(Z) <2\}$ and $\{Z \in {\mathcal B}_{q} \, : \, {\mathcal N}_{T, {\rm Pl}}(Z) \geq 2\}$ have locally finite 
$(n-2)$-dimensional Hausdorff measure (in fact locally finite upper Minkowski content). The fact that each of these two sets is locally compact is 
a consequence of the estimate \eqref{locally-unif-est} and upper semi-continuity of density $\Theta \, (T, \cdot)$ and planar frequency ${\mathcal N}_{T, {\rm Pl}}(\cdot)$. 
Finally, by \cite[Theorem~6.2(iii)]{KrumWicb}, each point $Z \in {\rm sing}_{q} \, T \setminus {\mathcal B}_{q}$ has a neighborhood ${\mathbf B}_{\sigma_{Z}}(Z)$ such that 
$\{Y \, : \, \Theta \, (T, Y) \geq \Theta (T, Z)\} \cap \overline{{\mathbf B}}_{\sigma_{Z}}(Z)$ is compact and has finite $(n-2)$-dimensional Hausdorff measure. 
\end{proof}

\subsection{Comparison between \cite{KrumWicb} and \cite{DelMinSko}}\label{comparison}
Let us now turn to 
a more focused comparison between \cite{KrumWicb} and \cite{DelMinSko}. Between these two papers there is a close similarity, in terms of the techniques vis \`a vis a certain key excess-decay lemma established in either paper, as well as a key difference in the way this decay lemma is used to prove structure results for (part of) the singular set of $T$. Because of the differences in the way a key hypothesis of the decay lemma is verified, in \cite{KrumWicb} the lemma both applies to a larger set of singularities, and yields stronger conclusions for that set, than in \cite{DelMinSko}.

In \cite{KrumWicb} this excess-decay lemma is given as \cite[Lemma~5.14]{KrumWicb} and \cite[Lemma~5.15]{KrumWicb}, and in \cite{DelMinSko} it is 
\cite[Theorem~2.5]{DelMinSko}. It gives improvement of excess for $T$ whenever condition \eqref{branch and cones concl3-intro} of Theorem~\ref{branch and cones thm-intro} above holds (say at scale $\rho = 1$) for some fixed $\beta = \beta(n, m, q) \in (0, 1)$, in addition to other hypotheses including a certain 
``no-large-gaps'' assumption on the set of density $\geq q$ singular points. The close similarity between \cite{KrumWicb} and \cite{DelMinSko} lies in the fact that in either case, the proof of this lemma uses a new pointwise bound for the distance of an area minimizing current to a union of non-intersecting planes 
(\cite[Theorem~3.11]{KrumWicb} and \cite[Theorem~3.2]{DelMinSko}), and proceeds by adapting techniques developed in \cite{Sim93} and \cite[Sections~14]{Wic14}.

To describe how this decay lemma is used in \cite{KrumWicb}, fix $\beta = \beta(n, m, q) \in (0, 1)$ to be the constant as in that lemma 
(i.e.\ \cite[Lemma~5.14]{KrumWicb}), and write 
$${\mathcal S}_{q} = {\rm sing}_{q} \, T \setminus {\mathcal B}_{q}$$ 
where ${\rm sing}_{q} \, T = \{Z \in {\rm sing} \, T \, : \, \Theta \, (\|T\|, Z) = q\}$ and ${\mathcal B}_{q}  = {\mathcal B}_{q, \beta, \beta},$ i.e.\ the set 
${\mathcal B}_{q, \epsilon, \beta}$ as in Corollary~\ref{immediate-cor} (in the introduction above) taken with  $\epsilon = \beta$. 
Thus ${\mathcal S}_{q}$ is the set of density $q$ singularities of $T$ which are not branch points $Z$ where decay of $T$ towards a multiplicity $q$ plane 
holds locally uniformly at a rate $|X - Z|^{1 + \alpha}$ or faster (i.e.\ points $Z$ with planar frequency ${\mathcal N}_{T, {\rm Pl}}(Z) \geq 1 + \alpha$), 
where $\alpha = \alpha(n, m, q, \beta) \in (0, 1)$ (and hence $\alpha = \alpha(n, m, q)$ after we choose $\beta = \beta(n, m, q)$ as in \cite[Lemma~5.14]{KrumWicb}). 
In \cite{KrumWicb} an asymptotic analysis of the behaviour of $T$ on approach to points in ${\mathcal S}_{q}$
is carried out.
The starting point of this analysis is Theorem~\ref{branch and cones thm-intro}(I), which gives that about every point $Z \in {\mathcal S}_{q}$ and at all small scales, either $T$ satisfies condition~\eqref{branch and cones concl3-intro} (which is needed to apply the decay lemma) or condition~\eqref{branch and cones concl4-intro}. Of key significance here is that the validity of these conditions is locally uniform 
in the sense that a scale $\rho_{Z}>0$ can be chosen so that \eqref{branch and cones concl3-intro} or \eqref{branch and cones concl4-intro}  holds about each point of ${\mathcal S}_{q} \cap {\mathbf B}_{\rho_{Z}}(Z)$  and at each scale $\rho \leq \rho_{Z}$. 

In \cite{DelMinSko} on the other hand, the decay lemma is used to analyse the (smaller) set 
$$\Sigma_{q} = \{Z \, : \, Z\;  \mbox{is a branch point with} \, \Theta(T, Z) = q \; \mbox{and} \;  I(Z) = 1\} \cup {\mathcal S}_{q, n-2},$$ 
where ${\mathcal S}_{q, n-2}$ is the set of density $q$ non branch point singularities of $T$ (where each tangent cone has spine dimension at most $(n-2)$). We have 
\emph{a posteriori} that ${\mathcal H}^{n-2}\, ({\mathcal S}_{q} \setminus \Sigma_{q})= 0$, but \emph{a priori} not much can be said about the relationship between the two sets ${\mathcal S}_{q}$, $\Sigma_{q}$ beyond the fact that $\Sigma_{q} \subset {\mathcal S}_{q}$ which follows readily from the definitions of the two sets in view of the fact that $I(Z) \geq 1 + \alpha$ for any 
point $Z \in {\mathcal B}_{q}$. A key difference  in the way the excess decay lemma is used (which influences the conclusions, summarised in the next paragraph) is that in \cite{DelMinSko}, the validity of  \eqref{branch and cones concl3-intro} or \eqref{branch and cones concl4-intro} is checked pointwise, by arguing that for each point $Z \in \Sigma_{q},$ there is a scale $\rho_{Z}$ depending on $Z$ such that for each scale $\rho < \rho_{Z}$ either  \eqref{branch and cones concl3-intro} or \eqref{branch and cones concl4-intro}  holds. This is a fairly direct consequence of either the definition of ${\mathcal S}_{q, n-2}$ (if $Z$ is not a branch point) or the fact that $I(Z) = 1$ (if $Z$ is a branch point), but  because of the pointwise choice of $\rho_{Z}$ (as opposed to a locally uniform choice), the set $\Sigma_{q}$ needs to be decomposed as a union of countably many subsets (in the obvious way by setting the $j$th subset equal to points $Z \in \Sigma_{q}$ for which $\rho_{Z} = j^{-1}$) before the decay lemma is applied (to each of these subsets). The other 
significant contrast to note is that in \cite{DelMinSko}, validity of \eqref{branch and cones concl3-intro} or \eqref{branch and cones concl4-intro} along $\Sigma_{q}$ ultimately relies on center manifold constructions (through the involvement of $I(Z)$ in the definition of $\Sigma_{q}$) whereas in the present work  
the validity of these conditions for ${\mathcal S}_{q}$ is guaranteed by the much more elementary result Theorem~\ref{branch and cones thm-intro}  which is a fairly direct consequence of the approximate monotonicity of the planar frequency function.

The main conclusions in \cite{KrumWicb} are: (i) ${\mathcal S}_{q}$ is countably $(n-2)$-rectifiable; (ii) for ${\mathcal H}^{n-2}$ a.e.\ point in ${\mathcal S}_{q}$ the current $T$ has a unique tangent cone supported on two or more distinct planes meeting along a common $(n-2)$-dimensional subspace, and (iii) there exists an ambient open set 
$V_{q}$ with ${\rm sing}_{q} \, T \subset V_{q}$ such that $V_{q} \cap \{Z \, : \, \Theta(T, Z) \geq q\} \setminus {\mathcal B}_{q}$ is locally compact and has locally finite $(n-2)$-dimensional Hausdorff measure (\cite[Theorem~6.2]{KrumWicb});  in particular ${\mathcal S}_{q}$ has locally finite ${\mathcal H}^{n-2}$ measure (since  ${\mathcal S}_{q} \subset V_{q} \cap \{Z \, : \, \Theta(T, Z) \geq q\} \setminus {\mathcal B}_{q}$). The main conclusions in \cite{DelMinSko} are (i) and (ii) with $\Sigma_{q}$ in place of ${\mathcal S}_{q}$. The additional conclusions (iii) in \cite{KrumWicb} are consequences  
of the fact that ${\mathcal B}_{q}$ 
is relatively closed in $V_{q}$ 
and the fact (emphasized above) that \eqref{branch and cones concl3-intro}  or \eqref{branch and cones concl4-intro} holds in a locally uniform way about every point in ${\mathcal S}_{q}$. 

\section{Preliminaries}

\subsection{General notation}  Throughout the paper, we shall use the following notation:
\begin{enumerate}[itemsep=3mm,topsep=0mm]
	\item[]  $n \geq 2$ and $m \geq 1$ are fixed integers.  $\mathbb{R}^{n+m}$ denotes the $(n+m)$-dimensional Euclidean space.  $X = (x_1,x_2,\ldots,x_{n+m})$ denotes a general point in $\mathbb{R}^{n+m}$.
	\item[]  $e_1,e_2,\ldots,e_{n+m}$ denotes the standard basis of $\mathbb{R}^{n+m}$.
	\item[]  For each $Y \in \mathbb{R}^{n+m}$ and $\rho > 0$, $\mathbf{B}_{\rho}(Y) = \{ X \in \mathbb{R}^{n+m} : |X-Y| < \rho \}$.
	\item[]  For each $Y \in \mathbb{R}^{n+m}$ and $\rho > 0$, $\eta_{Y,\rho} : \mathbb{R}^{n+m} \rightarrow \mathbb{R}^{n+m}$ is defined by $\eta_{Y,\rho}(X) = (X-Y)/\rho$. 
	\item[]  $\mathcal{H}^k$ denotes the $k$-dimensional Hausdorff measure on $\mathbb{R}^{n+m}$.
	\item[]  $\mathcal{L}^k$ denotes the $k$-dimensional Lebesgue measure on $\mathbb{R}^k$.
	\item[]  $\omega_k$ denotes the Lebesgue measure of a ball of radius one in $\mathbb{R}^k$.
	\item[]  For each point $X \in \mathbb{R}^{n+m}$ and set $A \subseteq \mathbb{R}^{n+m}$, $\op{dist}(X,A) = \inf_{Y \in A} |X-Y|$.  
	\item[]  For each pair of sets $A,B \subseteq \mathbb{R}^{n+m}$, $\op{dist}_{\mathcal{H}}(A,B)$ denotes the Hausdorff distance between $A$ and $B$.
	\item[]  For each set $A \subseteq \mathbb{R}^{n+m}$, $\overline{A}$ denotes the closure of $A$ and $\partial A$ denotes the boundary, or frontier, of $A$. 
\end{enumerate}

\subsection{Planes}  We let $\mathcal{P}$ denote the set of all $n$-dimensional planes (containing the origin) in $\mathbb{R}^{n+m}$.  Given an $n$-dimensional plane $P$ in $\mathbb{R}^{n+m}$ we let:
\begin{enumerate}[itemsep=3mm,topsep=0mm]
	\item[]  $P^{\perp}$ denote the orthogonal complement of $P$ in $\mathbb{R}^{n+m}$.
	\item[]  $\pi_P : \mathbb{R}^{n+m} \rightarrow P$ denote the orthogonal projection onto $P$.
	\item[]  $\pi_{P^{\perp}} :  \mathbb{R}^{n+m} \rightarrow P^{\perp}$ denote the orthogonal projection onto $P^{\perp}$.
	\item[]  for each $Z \in \mathbb{R}^{n+m}$ and $\rho > 0$, $B_{\rho}(Z,P) = \{ X \in Z+P : |X-Z| < \rho \}$ and $\mathbf{C}_{\rho}(Z,P) = \{ X \in \mathbb{R}^{n+m} : |\pi_P(X-Z)| < \rho \}$.
\end{enumerate}

We say that $P$ is an oriented $n$-dimensional plane if we equip $P$ with a simple $n$-vector $\vec P = \tau_1 \wedge \tau_2 \wedge\cdots\wedge \tau_n$, called the \emph{orientation} of $P$, where $\{\tau_1,\tau_2,\ldots,\tau_n\}$ is an orthonormal basis for $P$.  We let $\llbracket P \rrbracket$ denote integral current with support $P$, multiplicity one, and orientation $\vec P$.

We will often let $P_0 = \mathbb{R}^n \times \{0\}$, and take $\llbracket P_{0} \rrbracket$ with orientation $\vec P_0 = e_1 \wedge e_2 \wedge\cdots\wedge e_n$.  In this case we will identify $P_0 \cong \mathbb{R}^n$ and $P_0^{\perp} \cong \mathbb{R}^m$.  For each $z \in \mathbb{R}^n$ and $\rho > 0$ we will let $B_{\rho}(z) = \{ x \in \mathbb{R}^n : |x-z| < \rho \}$ and $\mathbf{C}_{\rho}(z) = B_{\rho}(z) \times \mathbb{R}^m$.

\subsection{Rectifiable currents}  Let $U$ be an open subset of $\mathbb{R}^{n+m}$.  An \emph{$n$-dimensional current} $T$ of $U$ is a continuous linear functional on the Fr\'echet space of smooth differential $n$-forms with compact support in $U$.  An \emph{$n$-dimensional integer multiplicity rectifiable current} $T$ of $U$ (abbreviated as \emph{$n$-dimensional rectifiable current} or just \emph{rectifiable current} henceforth) is a current of the form 
\begin{equation}\label{general integral current}
	T(\omega) = \int_M \langle \omega, \vec T \rangle \,\theta \,d\mathcal{H}^n
\end{equation}
for all smooth differential $n$-forms $\omega$ with compact support in $U$, where 
\begin{enumerate}[itemsep=3mm,topsep=0mm]
	\item[(i)]  $M$ is an ${\mathcal H}^{n}$ measurable, countably $n$-rectifiable subset of $U$;
	\item[(ii)]  $\theta : M \rightarrow \mathbb{Z}_+$ is a locally $\mathcal{H}^n$-integrable positive integer-valued function on $M$. The function $\theta$ is called the \emph{multiplicity} of $T$; 
	\item[(iii)]  $\vec T : M \rightarrow \Lambda^n(\mathbb{R}^{n+m})$ is an $\mathcal{H}^n$-measurable function on $M$ such that for $\mathcal{H}^n$-a.e.~$X \in M$, $\vec T(X) = \tau_1 \wedge \tau_2 \wedge\cdots\wedge \tau_n$ for some orthonormal basis $\{\tau_1, \tau_2, \ldots \tau_n\}$ for the approximate tangent plane to $M$ at $X$. The function $\vec T$ is called the \emph{orientation} of $T$. 
\end{enumerate}
We say that an $n$-dimensional current $T$ is an \emph{integral current} if both $T$ and $\partial T$ are integer multiplicity rectifiable.  The general theory of integral currents can be found in~\cite{Fed69} and~\cite{SimonGMT}.  We shall use the following notation associated with integral currents:
\begin{enumerate}[itemsep=3mm,topsep=0mm]
	\item[]  $\|T\| = \theta \,d\mathcal{H}^n$ is the \emph{total variation measure} associated with $T$ (as in~\ref{general integral current}).
	\item[]  $\mathbf{M}_W(T) = \|T\|(W)$ is the \emph{mass} $T$ in an open set $W \subseteq U$.  $\mathbf{M}(T) = \|T\|(U)$ is the \emph{mass} $T$. 
	\item[]  $\Theta(T,X) = \lim_{r \rightarrow 0^+} \omega_n^{-1} r^{-n} \|T\|(\mathbf{B}_r(X))$ is the \emph{density} of $T$ at $X$, whenever it exists.
	\item[]  $\op{spt} T$ is the support of $T$.
	\item[]  $T \llcorner A$ is the restriction of $T$ to a $\|T\|$-measurable set $A \subset U$.
	\item[]  $\partial T$ is the \emph{boundary} of $T$.
	\item[]  $f_{\#} T$ is the \emph{image} or the \emph{push-forward} of $T$ under a Lipschitz map $f : U \rightarrow V$ between two open sets $U,V$ such that $f |_{\op{spt} T}$ is proper, i.e. $f^{-1}(K) \cap \op{spt} T$ is a compact set whenever $K \subset V$ is a compact set.
	\item[]  $S \times T$ is the cartesian product of a pair of integral currents $S$ and $T$. 
	\item[]  $\llbracket M \rrbracket$ is the multiplicity one $k$-dimensional current associated with an oriented $C^{1}$ submanifold $M,$ or an 
	oriented $C^1$-submanifold-with-boundary $M,$ of ${\mathbb R}^{n+m}$.
\end{enumerate}

$G_n(U)$ denotes the \emph{Grassmannian}, which is the fiber bundle consisting of all pairs $(X,S)$ where $X \in U$ and $S$ is an $n$-dimensional plane in $\mathbb{R}^{n+m}$.  An \emph{$n$-dimensional varifold} is a Radon measure on $G_n(U)$.  Varifold convergence is the usual convergence of Radon measures on $G_{n}(U)$. To each pair consisting of an ${\mathcal H}^{n}$ measurable, countably $n$-rectifiable $M \subset U$ and a locally $\mathcal{H}^n$-integrable positive integer-valued function $\theta : M \rightarrow \mathbb{Z}_+$ we associate the \emph{$n$-dimensional integral varifold} $V = \mathbf{v}(M,\theta)$ such that 
\begin{equation*}
\int_{G_{n}(U)} \varphi(X, S) \, dV(X, S) = \int_{M} \varphi(X, \op{Tan}(M, X)) \, \theta(X) d{\mathcal H}^{n}(X)
\end{equation*}
for each $\varphi \in C_{c}(G_{n}(U))$, where $\op{Tan}(M,X)$ is the approximate tangent plane to $M$ at $X$.  In other words, for $V$-a.e.~$(X,S) \in G_n(U)$, $S$ is the approximate tangent plane to $M$ at $X$.  To each $n$-dimensional rectifiable current $T$ we can associate an $n$-dimensional integral varifold $|T| = \mathbf{v}(M,\theta)$, where $M$ and $\theta$ are as in \eqref{general integral current}.  For a further discussion of the theory of integral varifolds and general varifolds, we refer the reader to~\cite[Chapters~4 and 8]{SimonGMT}.

\subsection{Locally area-minimizing rectifiable currents}  Let $U$ be an open subset of $\mathbb{R}^{n+m}$.  

\begin{definition}{\rm 
We say that an $n$-dimensional rectifiable current $T$ in $U$ is \emph{locally area-minimizing in $U$} if 
\begin{equation*}
	\mathbf{M}_W(T) \leq \mathbf{M}_W(S)
\end{equation*}
for every open set $W \subset\subset U$ and every $n$-dimensional rectifiable current $S$ in $U$ such that $\partial S = \partial T$ in $U$ and $\op{spt}(S-T) \subset W.$ (Here $W \subset\subset U$ means the closure of $W$ is a compact subset of $U$.)
 }\end{definition}

Given an $n$-dimensional locally area-minimizing rectifiable current $T$ of $U$, the \emph{regular set} $\op{reg} T$ is the set of all points $Y \in \op{spt} T \cap U \setminus \op{spt} \partial T$ such that for some $\rho > 0$, $\op{spt} T \cap \mathbf{B}_{\rho}(Y)$ is a smooth embedded $n$-dimensional submanifold of $\mathbf{B}_{\rho}(Y)$.  The \emph{singular set} $\op{sing} T = \op{spt} T \cap U \setminus (\op{spt} \partial T \cup \op{reg} T)$. 
 
Whenever $T$ is an $n$-dimensional locally area-minimizing rectifiable current of $U$ with $(\partial T) \llcorner U = 0$, the integral varifold $|T|$ associated with $T$ is \emph{stationary} in the sense that 
\begin{equation}\label{first variation}
	\int_{G_n(U)} \op{div}_S \zeta(X) \,d|T|(X,S) = 0
\end{equation}
for all $\zeta \in C^1_c(U,\mathbb{R}^{n+m})$, where $\op{div}_S$ denotes the divergence with respect to the plane $S$~\cite[Definition~16.3 and Lemma~33.2]{SimonGMT}. 
 
Let $T$ be an $n$-dimensional locally area-minimizing rectifiable current of $U$ with $(\partial T) \llcorner U = 0$ and let $Z \in \op{spt} T \cap U$.  Since $T$ is stationary, $T$ satisfies a well-known \emph{monotonicity formula for area}, which implies
\begin{equation}\label{monotonicity formula for area}
	\frac{\|T\|(\mathbf{B}_{\sigma}(Z))}{\omega_n \sigma^n} \leq \frac{\|T\|(\mathbf{B}_{\rho}(Z))}{\omega_n \rho^n}
\end{equation}
for all $0 < \sigma \leq \rho \leq \op{dist}(Z,\mathbb{R}^{n+m} \setminus U)$~\cite[Theorem~17.6]{SimonGMT}.  If $U = \mathbb{R}^{n+m},$ $Z = 0$ and equality holds true in \eqref{monotonicity formula for area} for all $0 < \sigma \leq \rho < \infty$, then $T$ is a \emph{cone}, i.e.~$\eta_{0,\lambda\#} T = T$ for all $\lambda > 0$ (see \cite[Theorem~19.3]{SimonGMT} and the footnote at the bottom of \cite[p.\ 203]{SimonGMT}).  As a consequence of the monotonicity formula and the compactness theorem for locally area-minimizing rectifiable currents, for each sequence $\rho_k \rightarrow 0^+$ there exists a subsequence $\{k'\} \subset \{k\}$ and $n$-dimensional locally area-minimizing rectifiable current $\mathbf{C}$ which is a cone in $\mathbb{R}^{n+m}$ such that $\eta_{Z,\rho_{k'}\#} T \rightarrow \mathbf{C}$ weakly in $\mathbb{R}^{n+m}$ and $\Theta(\mathbf{C},0) = \Theta(T,Z)$~\cite[Theorem~35.1]{SimonGMT}.  We say that $\mathbf{C}$ is a \emph{tangent cone} of $T$ at $Z$. 

Let $\mathbf{C}$ be any $n$-dimensional locally area-minimizing rectifiable cone in $\mathbb{R}^{n+m}$.  Since $\mathbf{C}$ is a cone, $\Theta(\mathbf{C},Z) = \Theta(\mathbf{C},tZ)$ for all $Z \in \mathbb{R}^{n+m}$ and $t > 0$.  Letting $t \rightarrow 0^+$ using the semi-continuity of density~\cite[Corollary~17.8]{SimonGMT}, $\Theta(\mathbf{C},Z) \leq \Theta(\mathbf{C},0)$ for all $Z \in \mathbb{R}^{n+m}$.  We define the \emph{spine} of $\mathbf{C}$ to be the set 
\begin{equation*}
	\op{spine} \mathbf{C} = \{ Z \in \mathbb{R}^{n+m} : \Theta(\mathbf{C},Z) = \Theta(\mathbf{C},0) \} . 
\end{equation*}
By~\cite[Theorem~2.26]{Almgren}, $\op{spine} \mathbf{C}$ is a subspace of $\mathbb{R}^{n+m}$ and $\eta_{Z,1\, \#} \mathbf{C} = \mathbf{C}$ for all $Z \in \op{spine} \mathbf{C}$.  After an orthogonal change of coordinates, we can assume that $\op{spine} \mathbf{C} = \{0\} \times \mathbb{R}^{n-k}$, in which case $\mathbf{C} = \mathbf{C}_0 \times \mathbb{R}^{n-k}$ for some $k$-dimensional locally area-minimizing rectifiable cone $\mathbf{C}_0$ in $\mathbb{R}^{k+m}$.  If $\dim \op{spine} \mathbf{C} = n$, then $\mathbf{C} = q \llbracket P \rrbracket$ for some integer $q \geq 1$ and some $n$-dimensional oriented plane $P$.  There is no $n$-dimensional locally area-minimizing rectifiable cone $\mathbf{C}$ such that $\dim \op{spine} \mathbf{C} = n-1$.  If $\dim \op{spine} \mathbf{C} = n-2$, then 
\begin{equation*}
	\mathbf{C} = \sum_{i=1}^p q_i \llbracket P_i \rrbracket
\end{equation*}
 for some integers $p \geq 2$ and $q_i \geq 1$ and some $n$-dimensional oriented planes $P_i$ such that $P_i \cap P_j = \op{spine} \mathbf{C}$ whenever $i \neq j.$

Let $T$ be an $n$-dimensional locally area-minimizing rectifiable current in $U$ with $(\partial T) \llcorner U = 0$.  For each $j = 0,1,2,\ldots,n$, define \emph{$j$-th stratum} $S_j$ of the singular set of $T$ to be the set of all points $Z \in \op{sing} T$ such that $\dim \op{spine} \mathbf{C} \leq j$ for every tangent cone $\mathbf{C}$ to $T$ at $Z$.  Observe that 
\begin{equation*}
	S_{0} \subseteq S_{1} \subseteq\cdots\subseteq S_{n-3} \subseteq S_{n-2}  = S_{n-1} \subseteq S_n = \op{sing} T . 
\end{equation*}
$S_{n} \setminus S_{n-2}$ is the set of all branch point singularities, at which $T$ has at least one tangent one which is an integer multiplicity plane.  

\begin{lemma}
Let $T$ be an $n$-dimensional locally area-minimizing rectifiable current of $U$ with $(\partial T) \llcorner U = 0$.  For each $j = 1,2,\ldots,n-2,n$, $S_{j}$ has Hausdorff dimension at most $j$.  For $\alpha > 0$, $\{ Z \in S_{0} : \Theta(T,Z) = \alpha \}$ is discrete.
\end{lemma}

\begin{proof}
See \cite[Theorem~2.26 and Corollary~2.27]{Almgren}.
\end{proof}

\subsection{Multi-valued functions}\label{sec:prelim_multivalued}  For each integer $q \geq 1$, $\mathcal{A}_q(\mathbb{R}^m)$ denotes the space of all sums 
\begin{equation*}
	\sum_{i=1}^q \llbracket a_i \rrbracket
\end{equation*}
of Dirac point masses $\llbracket a_i \rrbracket$ at points $a_i \in \mathbb{R}^m$ (possibly repeating).  We equip $\mathcal{A}_q(\mathbb{R}^m)$ with the metric defined by 
\begin{equation*}
	\mathcal{G}(a,b) = \inf_{\sigma} \, \left( \sum_{i=1}^q |a_i - b_{\sigma(i)}|^2 \right)^{1/2}
\end{equation*}
for each $a = \sum_{i=1}^q \llbracket a_i \rrbracket, b = \sum_{i=1}^q \llbracket b_i \rrbracket \in \mathcal{A}_q(\mathbb{R}^n)$, where the infimum is taken over all permutations $\sigma$ of $\{1,2,\ldots,q\}$.  We write 
\begin{equation*}
	|a| = \mathcal{G}(a,q\llbracket 0\rrbracket) = \left( \sum_{i=1}^q |a_i|^2 \right)^{1/2}
\end{equation*}
for each $a = \sum_{i=1}^q \llbracket a_i \rrbracket \in \mathcal{A}_q(\mathbb{R}^m)$.  

Let $\Omega \subset \mathbb{R}^n$.  A \emph{$q$-valued function} is a map $u : \Omega \rightarrow \mathcal{A}_q(\mathbb{R}^m)$.  For each $x \in \Omega$, we shall express $u(x) = \sum_{i=1}^q \llbracket u_i(x) \rrbracket$ where $u_i(x) \in \mathbb{R}^m$ are the ``$q$ values of $u(x)$''.  We define the \emph{average} $u_a : \Omega \rightarrow \mathbb{R}^m$ of $u$ to be the single-valued function defined by $u_a(x) = \frac{1}{q} \sum_{i=1}^q u_i(x)$ for each $x \in \Omega$.  We say that $u$ is \emph{average-free} if $u_a(x) = 0$ for all $x \in \Omega$.  The \emph{average-free part} of $u$ is the $q$-valued function $u_f : \Omega \rightarrow \mathcal{A}_q(\mathbb{R}^m)$ defined by $u_f(x) = \sum_{i=1}^q \llbracket u_i(x) - u_a(x) \rrbracket$ for each $x \in \Omega$.  Thus for each $x \in \Omega$, we can write $u_f$ and $u$ as 
\begin{equation*}
	u_f(x) = \sum_{i=1}^q \llbracket (u_f)_i(x) \rrbracket, \quad u(x) = \sum_{i=1}^q \llbracket u_a(x) + (u_f)_i(x) \rrbracket
\end{equation*}
where $(u_f)_i(x) = u_i(x) - u_a(x)$.  Observe that for each pair of $q$-valued functions $u,v : \Omega \rightarrow \mathcal{A}_q(\mathbb{R}^m)$, 
\begin{equation}\label{G avg sym}
	\mathcal{G}(u(x),v(x))^2 = q |u_a(x) - v_a(x)|^2 + \mathcal{G}(u_f(x),v_f(x))^2
\end{equation}
for all $x \in \Omega$.  

Let $\Omega \subset \mathbb{R}^n$ be open.  Since $\mathcal{A}_q(\mathbb{R}^m)$ is a metric space, we can define the space $C^0(\Omega,\mathcal{A}_q(\mathbb{R}^m))$ of continuous $q$-valued functions $u : \Omega \rightarrow \mathcal{A}_q(\mathbb{R}^m)$ in the usual way.  For each $\mu \in (0,1]$, we define the space $C^{0,\mu}(\Omega,\mathcal{A}_q(\mathbb{R}^m))$ of H\"{o}lder continuous $q$-valued functions to be the set of all $q$-valued functions $u : \Omega \rightarrow \mathcal{A}_q(\mathbb{R}^m)$ such that 
\begin{equation*}
	[u]_{\mu,\Omega'} = \sup_{x,y \in \Omega', \,x \neq y} \frac{\mathcal{G}(u(x),u(y))}{|x-y|^{\mu}} < \infty  
\end{equation*}
for all $\Omega' \subset\subset \Omega$.  We say that a $q$-valued function $u : \Omega \rightarrow \mathcal{A}_q(\mathbb{R}^m)$ is \emph{Lipschitz} if 
\begin{equation*}
	\op{Lip} u = \sup_{x,y \in \Omega, \,x \neq y} \frac{\mathcal{G}(u(x),u(y))}{|x-y|} < \infty .
\end{equation*}
We say a $q$-valued function $u : \Omega \rightarrow \mathcal{A}_q(\mathbb{R}^m)$ is \emph{differentiable} at $y \in \Omega$ if there exists a $q$-valued function $\ell_y : \mathbb{R}^n \rightarrow \mathcal{A}_q(\mathbb{R}^m)$ of the form $\ell_y(x) = \sum_{i=1}^q \llbracket b_i^y + A_i^y x \rrbracket$ for some $m \times n$ matrices $A_i^y$ and points $b_i^y \in \mathbb{R}^m$ ($1 \leq i \leq q$) such that 
\begin{equation*}
	\lim_{x \rightarrow y} \frac{\mathcal{G}(u(x),\ell_y(x))}{|x-y|} = 0 .
\end{equation*}
If additionally $A^y_i \neq A^y_j$ whenever $b_i^y \neq b_j^y$, then we say that $u$ is \emph{strongly differentiable} at $y$.  Whenever $u$ is differentiable at $y$, the \emph{derivative} of $u$ at $y$ is $Du(y) = \sum_{i=1}^q \llbracket A_i^y \rrbracket$.  We shall use the convention that we write $u(y) = \sum_{i=1}^q \llbracket u_i(y) \rrbracket$ and $Du(y) = \sum_{i=1}^q \llbracket Du_i(y) \rrbracket$ when $u$ is strongly differentiable at $y$ with $\ell_y(x) = \sum_{i=1}^q \llbracket u_i(y) + Du_i(y) \cdot (x-y) \rrbracket$.  By Rademacher's theorem for $q$-valued functions~\cite[Theorem~1.13]{DeLSpaDirMin}, every Lipschitz $q$-valued function is strongly differentiable at $\mathcal{L}^n$-a.e.~$y \in \Omega$.  

Given a Lipschitz $q$-valued function $u : \Omega \rightarrow \mathcal{A}_q(\mathbb{R}^m)$, $T = \op{graph} u$ is an $n$-rectifiable current of $\mathbb{R}^{n+m}$ given by \eqref{general integral current} with 
\begin{gather*}
	M = \{ (x,u_i(x)) : x \in \Omega, \,i \in \{1,2,\ldots,q\}\} , \\
	\theta(x,y) = \#\{ i \in \{1,2,\ldots,q\} : y = u_i(x) \} \text{ for each $(x,y) \in M$,} \\ 
	\vec T(x,u_i(x)) = \frac{\bigwedge_{j=1}^n (e_j,D_j u_i(x))}{\big| \bigwedge_{j=1}^n (e_j,D_j u_i(x)) \big|} 
		\text{ for $\mathcal{L}^n$-a.e.~$x \in \Omega$ and each $i \in \{1,2,\ldots,q\}$,}
\end{gather*} 
where $u(x) = \sum_{i=1}^q \llbracket u_i(x) \rrbracket$ for each $x \in \Omega$ and $Du(x) = \sum_{i=1}^q \llbracket Du_i(x) \rrbracket$ for $\mathcal{L}^n$-a.e.~$x \in \Omega$ following the above convention.  See~\cite[Section~1.5]{Almgren}, in which $\op{graph} u$ is defined via affine approximation, or~\cite[Section~1]{DeLSpa0}, in which $\op{graph} u$ is equivalently defined as the image of the $q$-valued map $F(x) = \sum_{i=1}^q \llbracket (x,u(x)) \rrbracket$ of $x \in \Omega$ using a measurable partition of $\Omega$ and a corresponding measurable selection of $F$.

For each $1 \leq p \leq \infty$, we define the Lebesgue space $L^p(\Omega,\mathcal{A}_q(\mathbb{R}^m))$ to be the set of all Lebesgue measurable $q$-valued functions $u : \Omega \rightarrow \mathcal{A}_q(\mathbb{R}^m)$ such that $\|u\|_{L^p(\Omega)} = \|\mathcal{G}(u,q\llbracket 0 \rrbracket)\|_{L^p(\Omega)} < \infty$.  

The Sobolev space $W^{1,2}(\Omega,\mathcal{A}_q(\mathbb{R}^m))$ of $q$-valued functions is defined in~\cite[Definitions and terminology~2.1]{Almgren} as follows:  Let $N(q,m) \geq 1$ be an integer and $\boldsymbol{\xi} : \mathcal{A}_q(\mathbb{R}^m) \rightarrow \mathbb{R}^N$ be a bi-Lipschitz embedding such that $\op{Lip} \boldsymbol{\xi} \leq 1$ and $\op{Lip} \boldsymbol{\xi} |_{\mathcal{Q}} \leq C(m,q)$, where $\mathcal{Q} = \boldsymbol{\xi}(\mathcal{A}_q(\mathbb{R}^m))$.  Then $W^{1,2}(\Omega,\mathcal{A}_q(\mathbb{R}^m))$ is the set of all Lebesgue measurable $q$-valued functions $u : \Omega \rightarrow \mathcal{A}_q(\mathbb{R}^m)$ such that $\boldsymbol{\xi} \circ u \in W^{1,2}(\Omega,\mathbb{R}^N)$.  $W^{1,2}(\Omega,\mathcal{A}_q(\mathbb{R}^m))$ can be equivalently characterized as the space of Sobolev functions into the metric space $\mathcal{A}_q(\mathbb{R}^m)$, see~\cite[Definition~0.5]{DeLSpaDirMin}.  Every $u \in W^{1,2}(\Omega,\mathbb{R}^n)$ is approximately strongly differentiable at $\mathcal{L}^n$-a.e.~$y \in \Omega$ in the sense that there exists a set $\Omega_y \subset \Omega$ with density one at $y$ such that $u |_{\Omega_y}$ is strongly differentiable at $y$ (\cite[Theorem~2.2]{Almgren}, \cite[Corollary~2.7]{DeLSpaDirMin}).  The derivative of $u$ at $y$ is $Du(y) = D(u |_{\Omega_y})(y)$.  Whenever $u \in W^{1,2}(\Omega,\mathcal{A}_q(\mathbb{R}^m))$, $u \in L^2(\Omega,\mathcal{A}_q(\mathbb{R}^m))$ and $Du \in L^2(\Omega,\mathcal{A}_q(\mathbb{R}^{m\times n}))$.

\subsection{Dirichlet energy minimizing multi-valued functions}\label{Dir-min-fns}  
Let $\Omega$ be an open subset of $\mathbb{R}^n$.  
\begin{definition}{\rm 
We say a $q$-valued function $w \in W^{1,2}(\Omega,\mathcal{A}_q(\mathbb{R}^m))$ is \emph{locally Dirichlet energy minimizing} (or \emph{Dirichlet energy minimizing} for simplicity) if 
\begin{equation*}
	\int_{\Omega'} |Dw|^2 \leq \int_{\Omega'} |Dv|^2 
\end{equation*}
whenever $\Omega' \subset\subset \Omega$ is an open set and $v \in W^{1,2}(\Omega,\mathcal{A}_q(\mathbb{R}^m))$ is a $q$-valued function such that $w(x) = v(x)$ for $\mathcal{L}^n$-a.e.~$x \in \Omega\setminus\Omega'$.
 }\end{definition}
 
The theory of Dirichlet energy minimizing $q$-valued functions was developed by Almgren in~\cite{Almgren} where such functions were used to approximate, in a certain precise sense, locally area-minimizing rectifiable currents weakly close to a multiplicity $q$ plane;  this theory serves as the ``linear theory'' in the study of area-minimizing currents.  For a detailed discussion of this theory 
see~\cite[Chapter~2]{Almgren} or~\cite{DeLSpaDirMin}.  See also the summary of Almgren's existence and regularity theory for Dirichlet energy minimizing multi-valued functions in~\cite[Subsection~3.2 and Section~4]{KrumWic2}.  

Let $\phi : [0,\infty) \rightarrow \mathbb{R}$ be a nonincreasing Lipschitz function such that $\phi = 1$ on $[0,1/2]$ and $\phi = 0$ on $[1,\infty)$.  Let $w \in W^{1,2}_{\rm loc}(\Omega,\mathcal{A}_q(\mathbb{R}^m))$ be a non-zero Dirichlet energy minimizing $q$-valued function, $z \in \Omega$, and $0 < \rho < \op{dist}(z,\partial \Omega)$.  We define the \emph{frequency function} $N_{w,z}(\rho)$ associated with $u$ and $z$ by 
\begin{equation*}
	N_{w,z}(\rho) = \frac{D_{w,z}(\rho)}{H_{w,z}(\rho)}
\end{equation*}
whenever $H_{w,z}(\rho) > 0$, where we let $r = |x-z|$ and 
\begin{equation*}
	H_{w,z}(\rho) = -\rho^{1-n} \int |w|^2 \,\frac{1}{r} \,\phi'(r/\rho) , \quad 
	D_{w,z}(\rho) = \rho^{2-n} \int |Dw|^2 \,\phi(r/\rho) . 
\end{equation*}
By letting $\phi$ increase to the characteristic function on the interval $[0,1)$, we obtain the classical frequency function as defined by Almgren in~\cite[Theorem~2.6]{Almgren}; incorporating a fixed function $\phi$ into the definition provides the convenience, as observed in \cite{DSMV18}, \cite{DeLSpa3},  of avoiding boundary integrals as in Almgren's original definition. As in~\cite{DSMV18}, \cite{DeLSpa3}, here we take $\phi$ to be the Lipschitz function such that $\phi(s) = 2-2s$ for each $s \in (1/2,1)$ (as in \eqref{freq phi defn} below).  One can then show (see \cite[Proposition 3.1]{DSMV18}) that if $\Omega$ is connected and $w$ is not identically $q \llbracket 0 \rrbracket$ on $\Omega$, then $H_{w,z}(\rho) > 0$ for all $0 < \rho < \op{dist}(z,\partial \Omega)$ and 
\begin{equation*}
	N'_{w,z}(\rho) = \frac{-2\rho^{1-n}}{H_{w,z}(\rho)} \int \left(\sum_{i=1}^q |x \cdot Dw_i(x) - N_{w,z}(\rho) \,w_i(x)|^2\right)\frac{1}{r}\phi^{\prime}(r/\rho)\,d\mathcal{L}^n(x) \geq 0 
\end{equation*}
for all $0 < \rho < \op{dist}(z,\partial \Omega)$, where we write $w(x) = \sum_{i=1}^q \llbracket w_i(x) \rrbracket$ for each $x \in \Omega$ and $Dw(x) = \sum_{i=1}^q \llbracket Dw_i(x) \rrbracket$ for $\mathcal{L}^n$-a.e.~$x \in \Omega$ (following the convention from Subsection~\ref{sec:prelim_multivalued}).  In particular, $N_{w,z}(\rho)$ is a monotone nondecreasing function of $\rho \in (0,\op{dist}(z,\partial \Omega))$.  $N'_{w,z}(\tau) = 0$ for each $\tau \in (\sigma,\rho)$ if and only if $w(z+\lambda x) = \sum_{i=1}^q \llbracket \lambda^{\alpha} w_i(z+x) \rrbracket$ whenever $\lambda > 0$ and $x, \lambda x \in B_{\rho}(0) \setminus B_{\sigma}(0)$.  We define the \emph{frequency} $\mathcal{N}_w(z)$ of $w$ at $z$ by 
\begin{equation*}
	\mathcal{N}_w(z) = N_{w, z}(0^{+}) = \lim_{\rho\rightarrow 0^+} N_{w,z}(\rho) .
\end{equation*}
The value of $N_{w,z}(0^+)$ is in fact equal to the classical frequency as defined by Almgren in~\cite[Theorem~2.6]{Almgren}.  Using the monotonicity formula for frequency function and continuity of Dirichlet energy under uniform limits (\cite[Theorem~2.15]{Almgren}, ~\cite[Proposition~3.20]{DeLSpaDirMin}), one can show that frequency is upper semi-continuous in the sense that if $w_k,w \in W^{1,2}(\Omega,\mathcal{A}_q(\mathbb{R}^m))$ are Dirichlet energy minimizing $q$-valued functions such that $w_k \rightarrow w$ uniformly on $\Omega,$ and $z_k,z \in \Omega$ are such that $z_k \rightarrow z$, then 
\begin{equation*}
	\mathcal{N}_w(z) \geq \limsup_{k \rightarrow \infty} \mathcal{N}_{w_k}(z_k) .
\end{equation*}

Suppose that $w \in W^{1,2}_{\rm loc}(\mathbb{R}^n,\mathcal{A}_q(\mathbb{R}^m))$ is a non-zero Dirichlet energy minimizing $q$-valued function and $w$ is homogeneous of some degree $\alpha$, i.e.~$w(\lambda x) = \sum_{i=1}^q \llbracket \lambda^{\alpha} w_i(x) \rrbracket$ for each $x \in \mathbb{R}^n$ and $\lambda > 0$ where we write $w(x) = \sum_{i=1}^q \llbracket w_i(x) \rrbracket$.  Then $N_{w,0}(\rho) = \alpha$ for all $\rho \in (0,\infty)$.  By the homogeneity of $w$ and semi-continuity of frequency, for each $z \in \mathbb{R}^n$ we have that $\mathcal{N}_w(z) = \limsup_{t \rightarrow 0^+} \mathcal{N}_w(tz) \leq \mathcal{N}_w(0) = \alpha$.  We define the \emph{spine} of $w$ by 
\begin{equation*}
	\op{spine} w = \{ z \in \mathbb{R}^n : \mathcal{N}_w(z) = \alpha \} .
\end{equation*}
It follows from~\cite[Theorem~2.14]{Almgren} that $\op{spine} w$ is linear subspace of ${\mathbb R}^{n}$ and that $w(z+x) = w(x)$ for all $z \in \op{spine} w$ and $x \in \mathbb{R}^n$.  If $\dim \op{spine} w= n$, then $w$ is a constant $q$-valued function on $\mathbb{R}^n$.  There is no homogeneous Dirichlet energy minimizing $q$-valued function $w \in W^{1,2}_{\rm loc}(\mathbb{R}^n,\mathcal{A}_q(\mathbb{R}^m))$ with $\dim \op{spine} w = n-1$.  If $\dim \op{spine} w = n-2$, then 
\begin{equation*}
	w(x) = \sum_{i=1}^p q_i \llbracket A_i x \rrbracket
\end{equation*}
where $p$ and $q_i$ are positive integers such that $p \geq 2$ and $\sum_{i=1}^p q_i = q$, and $A_i$ are distinct $m \times n$ matrices such that  $\{ x \in \mathbb{R}^n : A_i x = A_j x \} = \op{spine} w$ if $i \neq j$, and if $w$ is average-free, $\{ x \in \mathbb{R}^n : A_i x = 0 \} = \op{spine} w$.

\subsection{Some elementary estimates}  We have the following well-known ``energy estimate'' which bounds the tilt excess of a stationary integral varifold $V$ relative to a plane $P$ from above in terms of the height-excess, i.e.\ $L^2$-distance of $V$ to $P$.

\begin{lemma}\label{tilt to height estimate lemma}
 If $P$ is an $n$-dimensional plane in $\mathbb{R}^{n+m}$, $V$ is an $n$-dimensional stationary integral varifold on 
$\mathbf{C}_{\rho}(Y, P)$ and $\gamma \in (0,1)$, then 
\begin{equation}\label{tilt to height estimate}
	\frac{1}{\rho^n} \int_{G_n(\mathbf{C}_{\gamma\rho}(Y, P))} \|\pi_S - \pi_P\|^2 \,dV(X,S) 
		\leq \frac{64}{(1-\gamma)^2 \rho^{n+2}} \int_{\mathbf{C}_{\rho}(Y, P)} \op{dist}^2(X,Y+P) \,d\|V\|(X) 
\end{equation}
where $\|\pi_{S} - \pi_{P}\|$ denotes the Frobenius norm of $\pi_{S} - \pi_{P}$. 
\end{lemma}

\begin{proof}
Assume without loss of generality that $Y = 0$, $\rho = 1$, and $P = \mathbb{R}^n \times \{0\}$.  Express points $X \in \mathbb{R}^{n+m}$ as $X = (x,y)$ where $x = (x_1,\ldots,x_n) \in \mathbb{R}^n$ and $y = (x_{n+1},\ldots,x_{n+m}) \in \mathbb{R}^m$.  Inequality \eqref{tilt to height estimate} follows by setting $\zeta(x,y) = \phi^2(x) \,(0,y)$ in (\ref{first variation}), where $\phi \in C^1_c(B_1(0))$ is a cut-off function such that $0 \leq \phi \leq 1$, $\phi = 1$ on $B_{\gamma}(0)$, and $|D\phi| \leq \tfrac{2}{1-\gamma}$.
\end{proof}

 The following well-known result due to Allard (~\cite{Allard}) bounds the $L^{\infty}$-distance of a stationary integral varifold $T$ to a plane $P$ linearly in terms of the $L^2$-distance of $T$ to $P$.  (See \cite{KrumWicb} for a generalisation of this to a union of planes disjoint in a cylinder in place of $P$ in the case that $V$ corresponds to a locally area minimizing rectifiable current).   

\begin{lemma} \label{Allard height lemma}
If $\gamma \in (0,1)$, $V$ is an $n$-dimensional stationary integral varifold on $\mathbf{B}_{\rho}(Y)$ such that $V \llcorner \mathbf{B}_{\gamma\rho}(Y) \neq 0$ and if $P$ is an $n$-dimensional affine plane in $\mathbb{R}^{n+m}$, then 
\begin{equation} \label{Allard height concl}
	\sup_{X \in \op{spt} T \cap \mathbf{B}_{\gamma\rho}(Y)} \op{dist}^2(X,P) \leq \frac{C}{\rho^n} \int_{\mathbf{B}_{\rho}(Y)} \op{dist}^2(X,P) \,d\|T\|(X) 
\end{equation}
for some constant $C = C(n,m,\gamma) \in (0,\infty)$. 
\end{lemma}

\begin{proof}
Without loss of generality assume that $P = \mathbb{R}^n \times \{0\}$.  Since the coordinate functions $x_i$ are a $|T|$-harmonic functions in $\mathbf{B}_{\rho}(Y)$ for all $i = 1,2,\ldots,m$~\cite[7.5(1)(2)]{Allard}, we can apply~\cite[Theorem~7.5(6)]{Allard} to deduce \eqref{Allard height concl}.
\end{proof}

If instead we consider a general closed set $K$ (not necessarily a plane), as a straightforward consequence of the monotonicity formula for area, we have the following more crude  
bound for the $L^{\infty}$-distance of $V$ to $K$ in terms of the $L^2$-distance of $V$ to $K$.

\begin{lemma}\label{coarse L2 distance lemma}
Let $\gamma \in (0,1)$.  Suppose $V$ is an $n$-dimensional stationary integral varifold on  $\mathbf{B}_{\rho}(Y)$ and $K$ is a closed subset of $\mathbb{R}^{n+m}$ such that 
\begin{equation}\label{coarse L2 distance hyp} 
	\frac{1}{\omega_n\rho^{n+2}} \int_{\mathbf{B}_{\rho}(Y)} \op{dist}^2(X,K) \,d\|V\|(X) < \left(\frac{1-\gamma}{2}\right)^{n+2} 
\end{equation}
then 
\begin{equation}\label{coarse L2 distance} 
	\sup_{X \in \op{spt} \|V\| \cap \mathbf{B}_{\gamma\rho}(Y)} \op{dist}(X,K) 	
		\leq 2 \left(\frac{1}{\omega_n} \int_{\mathbf{B}_{\rho}(Y)} \op{dist}^2(X,K) \,d\|V\|(X) \right)^{\frac{1}{n+2}} . 
\end{equation}
\end{lemma}

\begin{proof}
Assume without loss of generality (by traslating and scaling) that $Y = 0$ and $\rho = 1$.  Fix $Z \in \op{spt} T \cap \mathbf{B}_1(0)$ and set $d = \min\{1-\gamma,\op{dist}(Z,K)\}$ so that $\mathbf{B}_d(Z) \subseteq \mathbf{B}_1(0) \setminus K$.  By the monotonicity formula for area and~\cite[Remark~17.9(1)]{SimonGMT}, $\|V\|(\mathbf{B}_{d/2}(Z)) \geq \omega_n (d/2)^n$ and thus 
\begin{equation}\label{coarse L2 distance eqn} 
	\omega_n (d/2)^{n+2} \leq \int_{\mathbf{B}_{d/2}(Z)} \op{dist}^2(X,K) \,d\|V\|(X) \leq \int_{\mathbf{B}_1(0)} \op{dist}^2(X,K) \,d\|V\|(X) . 
\end{equation}
By \eqref{coarse L2 distance eqn} and \eqref{coarse L2 distance hyp}, we must have that $d < 1-\gamma$ and thus $d = \op{dist}(Z,K)$.  Hence \eqref{coarse L2 distance} follows immediately from \eqref{coarse L2 distance eqn}.
\end{proof}

Let $P$ be an $n$-dimensional oriented plane in $\mathbb{R}^{n+m}$ and $T$ be an $n$-dimensional locally area-minimizing rectifiable current in $\mathbf{C}_{\rho}(Y,P)$ with 
\begin{equation*}
	(\partial T) \llcorner \mathbf{C}_{\rho}(Y,P) = 0, \quad \sup_{X \in \op{spt} T} \op{dist}(X,P) < \infty . 
\end{equation*}
By the constancy theorem~\cite[Theorem~26.27]{SimonGMT}, there exists an integer $\theta$ such that 
\begin{equation*}
	\pi_{P\#} (T \llcorner \mathbf{C}_{\rho}(Y,P)) 
	= \theta \llbracket P \rrbracket \llcorner \mathbf{C}_{\rho}(Y,P) . 
\end{equation*}
We define the \emph{excess} $\mathcal{E}(T,P,\mathbf{C}_{\rho}(Y,P))$ of $T$ relative to $P$ in $\mathbf{C}_{\rho}(Y,P)$ by 
\begin{equation}\label{area excess defn}
	\mathcal{E}(T,P,\mathbf{C}_{\rho}(Y,P))^2 = \frac{\|T\|(\mathbf{C}_{\rho}(Y))}{\omega_n\rho^n} - |\theta|,
\end{equation}
whence, 
\begin{equation}\label{area excess eqn1}
	\|T\|(\mathbf{C}_{\rho}(Y,P)) = \big( |\theta| + \mathcal{E}(T,P,\mathbf{C}_{\rho}(Y,P))^2 \big) \,\omega_n \rho^n . 
\end{equation}
By~\cite[5.3.1]{Fed69}, if $\theta \geq 0,$
\begin{align*}
	\mathcal{E}(T,P,\mathbf{C}_{\rho}(Y, P))^2 &= \frac{1}{\omega_n\rho^n} \int_{\mathbf{C}_{\rho}(Y)} (1 - \langle \vec T, \vec P \rangle) \,d\|T\|(X)	
	\\&= \frac{1}{2\omega_n\rho^n} \int_{\mathbf{C}_{\rho}(Y,P)} |\vec T - \vec P|^2 \,d\|T\|(X) . \nonumber 
\end{align*}

\begin{lemma}\label{area excess to height lemma}
Let $\gamma \in (0,1)$.  Let $P$ be an $n$-dimensional oriented plane in $\mathbb{R}^{n+m}$ and $T$ be an $n$-dimensional locally area-minimizing rectifiable current in  $\mathbf{C}_{\rho}(Y,P)$ such that 
\begin{equation*}
	(\partial T) \llcorner \mathbf{C}_{\rho}(Y,P) = 0 , \quad 
	\sup_{X \in \op{spt} T} \op{dist}(X,P) \leq \rho . 
\end{equation*}
Then 
\begin{equation}\label{excess2height concl}
	\mathcal{E}(T,P,\mathbf{C}_{\gamma\rho}(Y,P))^{2} \leq \frac{C}{\rho^{n+2}} \int_{\mathbf{C}_{\rho}(Y, P)} \op{dist}^2(X,P) \,d\|T\|(X) 
\end{equation}
for some constant $C = C(n,m,\gamma) \in (0,\infty)$.  
\end{lemma}

\begin{proof}
The proof is based on~\cite[Lemma~3.2]{HardtSimon} and is included for completion.  Without loss of generality assume that $Y = 0$, $\rho = 1$, and $P = \mathbb{R}^n \times \{0\}$.  Express points $X \in \mathbb{R}^{n+m}$ as $X = (x,y)$ where $x = (x_1,\ldots,x_n) \in \mathbb{R}^n$ and $y = (x_{n+1},\ldots,x_{n+m}) \in \mathbb{R}^m$.  Let $\mu \in C^{\infty}(\mathbf{C}_1)$ be a smooth function such that $0 \leq \mu \leq 1$, $\mu = 1$ on $\mathbf{C}_1 \setminus \mathbf{C}_{(1+\gamma)/2}$, $\mu = 0$ in $\mathbf{C}_{\gamma}$, and $|\nabla \mu| \leq \tfrac{3}{1-\gamma}$.  Define a diffeomorphism $F : \mathbf{C}_1 \rightarrow \mathbf{C}_1$ by $F(X) = (x, \mu(X) \,y)$ and consider the competitor $F_{\#} T$ for $T$ (which fixes points in $\mathbf{C}_1 \setminus \mathbf{C}_{(1+\gamma)/2}$ and projects points in $\mathbf{C}_{\gamma}$ onto $P$).  Since $T$ is area-minimizing, 
\begin{equation*}
	\|T\|(\mathbf{C}_{(1+\gamma)/2}) \leq \|F_{\#} T\|(\mathbf{C}_{(1+\gamma)/2}) . 
\end{equation*}
Thus 
\begin{align}\label{excess2height eqn1}
	\omega_{n}\gamma^{n}\mathcal{E}(T,P,\mathbf{C}_{\gamma})^{2}
	=\,& \|T\|(\mathbf{C}_{\gamma}) - \|\pi_{P \#} T\|(\mathbf{C}_{\gamma})\\ 
	= \,&\|T\|(\mathbf{C}_{\gamma}) - \|F_{\#} T\|(\mathbf{C}_{\gamma}) \nonumber 
	\\=\,& \|T\|(\mathbf{C}_{(1+\gamma)/2}) - \|F_{\#} T\|(\mathbf{C}_{(1+\gamma)/2}) \nonumber 
		\\&+ \|F_{\#} T\|(\mathbf{C}_{(1+\gamma)/2} \setminus \mathbf{C}_{\gamma}) 
		- \|T\|(\mathbf{C}_{(1+\gamma)/2} \setminus \mathbf{C}_{\gamma}) \nonumber
	\\ \leq\,& \|F_{\#} T\|(\mathbf{C}_{(1+\gamma)/2} \setminus \mathbf{C}_{\gamma}) - \|T\|(\mathbf{C}_{(1+\gamma)/2} \setminus \mathbf{C}_{\gamma}) \nonumber
\end{align}

Let $X = (x,y) \in \op{spt} T$ be a point at which the approximate tangent plane $S_X$ to $T$ exists, and let $\{\tau_1,\tau_2,\ldots,\tau_n\}$ be an orthonormal basis for $S_{X}$ such that $\vec T = \tau_1 \wedge \tau_2 \wedge \cdots \wedge \tau_n$.  Noting that $\sup_{X \in \op{spt} T} |y| \leq 1$ and $|\pi_{P^{\perp}} \tau_i| \leq \|\pi_{S_X} - \pi_P\|$ for each $i$, 
\begin{align}\label{excess2height eqn2}
	F_{\#} \vec T &= \bigwedge_{i=1}^n F_{\#} \tau_i = \bigwedge_{i=1}^n (\pi_P \tau_i + \mu(X) \,\pi_{P^{\perp}} \tau_i + \nabla_{\tau_i} \mu(X) \,(0,y)) 
		\\&= \bigwedge_{i=1}^n \pi_P \tau_i + \sum_{i=1}^n (-1)^{i-1} (\mu(X) \,\pi_{P^{\perp}} \tau_i + \nabla_{\tau_i} \mu(X) \,(0,y)) 
			\wedge \bigwedge_{j \neq i} \pi_P \tau_j + \mathcal{R} \nonumber 
\end{align}
where $\mathcal{R}$ is an $n$-vector such that $|\mathcal{R}| \leq C(n,m,\gamma) \,(|y|^2 + \|\pi_{S_X} - \pi_P\|^2)$.  Notice that the first two terms on the second line of \eqref{excess2height eqn2} are mutually orthogonal as $n$-vectors.  Hence, again using the fact that $|\pi_{P^{\perp}} \tau_i| \leq \|\pi_{S_X} - \pi_P\|$, 
\begin{equation*}
	|F_{\#} \vec T|^2 \leq \left| \bigwedge_{i=1}^n \pi_P \tau_i \right|^2 + C (|y|^2 + \|\pi_{S_X} - \pi_P\|^2) 
		= |\pi_{P\#} \vec T|^2 + C (|y|^2 + \|\pi_{S_X} - \pi_P\|^2)
\end{equation*}
for some constant $C = C(n,m,\gamma) \in (0,\infty)$.  Since $\op{Lip} \pi_P \leq 1$, $|\pi_{P\#} \vec T| \leq |\vec T| = 1$ and thus 
\begin{equation}\label{excess2height eqn3}
	|F_{\#} \vec T|^2 \leq 1 + C (|y|^2 + \|\pi_{S_X} - \pi_P\|^2).
\end{equation}
Using \eqref{excess2height eqn3} and the fact that $|y| = \op{dist}(X,P)$, 
\begin{align*}\label{excess2height eqn4}
	&\|F_{\#} T\|(\mathbf{C}_{(1+\gamma)/2} \setminus \mathbf{C}_{\gamma}) - \|T\|(\mathbf{C}_{(1+\gamma)/2} \setminus \mathbf{C}_{\gamma}) 
	= \int_{\mathbf{C}_{(1+\gamma)/2} \setminus \mathbf{C}_{\gamma}} (|F_{\#} \vec T| - 1) \,d\|T\| \nonumber 
	\\ =\,& \int_{\mathbf{C}_{(1+\gamma)/2} \setminus \mathbf{C}_{\gamma}} \frac{|F_{\#} \vec T|^2 - 1}{|F_{\#} \vec T| + 1} \,d\|T\| 
	\leq \int_{\mathbf{C}_{(1+\gamma)/2} \setminus \mathbf{C}_{\gamma}} (|F_{\#} \vec T|^2 - 1) \,d\|T\| \nonumber
	\\ \leq\,& C \int_{\mathbf{C}_1} \op{dist}^2(X,P) \,d\|T\|(X) + C \int_{\mathbf{C}_{(1+\gamma)/2}} \|\pi_{S_X} - \pi_P\|^2 \,d\|T\|(X) 
\end{align*}
for some constant $C = C(n,m,\gamma) \in (0,\infty)$.  Hence using Lemma~\ref{tilt to height estimate lemma}, we conclude that 
\begin{equation}\label{excess2height eqn5}
	\|F_{\#} T\|(\mathbf{C}_{(1+\gamma)/2} \setminus \mathbf{C}_{\gamma}) - \|T\|(\mathbf{C}_{(1+\gamma)/2} \setminus \mathbf{C}_{\gamma})
	\leq C \int_{\mathbf{C}_1} \op{dist}^2(X,P) \,d\|T\|(X) 
\end{equation}
for some constant $C = C(n,m,\gamma) \in (0,\infty)$, which together with \eqref{excess2height eqn1} completes the proof of \eqref{excess2height concl}. 
\end{proof}

\subsection{Lipschitz approximation}  The Strong Lipschitz Approximation Theorem for area-minimizing currents is a fundamental result proved by Almgren in~\cite[Corollaries~3.29 and 3.30]{Almgren} (see also ~\cite[Theorem~2.4]{DeLSpa1}). Here we adapt this result in a straightforward manner to obtain a Lipschitz approximation lemma that can be applied to general domains $\Omega$ and blow-up limits with convergence on compact subsets of $\Omega$.

\begin{theorem}\label{lipschitz approx thm}
Let $q$ be a positive integer and $\Omega$ be a bound open subset of $\mathbb{R}^n$.  For each $\sigma > 0$ let 
\begin{equation*}
	\Omega_{\sigma} = \{ x \in \Omega : \op{dist}(x, \mathbb{R}^n \setminus \Omega) > \sigma \} . 
\end{equation*}
There exists $\varepsilon = \varepsilon(n,m,q) \in (0,1)$ such that if $0 < \delta < \infty$ and $T$ is an $n$-dimensional locally area-minimizing rectifiable current in $\Omega \times \mathbb{R}^m$ such that 
\begin{gather}
	(\partial T) \llcorner (\Omega \times \mathbb{R}^m) = 0, \quad 
	\sup_{X \in \op{spt} T} \op{dist}(X,P_0) < \infty, \quad 
	\pi_{P_0\#} T = q \llbracket \Omega \rrbracket , \nonumber \\ 
	\label{lipschitz approx hyp} \mathcal{E}^2 = \int_{\Omega \times \mathbb{R}^m} |\vec T - \vec P_0|^2 \,d\|T\|(X) < \varepsilon^2 \delta^n ,  
\end{gather}
then there exists a Lipschitz function $u : \Omega_{\delta} \times \mathbb{R}^m \rightarrow \mathcal{A}_q(\mathbb{R}^m)$ and $\mathcal{L}^n$-measurable set $K \subseteq \Omega_{\delta}$ such that 
\begin{equation}\label{lipschitz approx concl1} 
	T \llcorner (K \times \mathbb{R}^m) = (\op{graph} u) \llcorner (K \times \mathbb{R}^m) 
\end{equation}
and for each $\sigma \geq \delta$ 
\begin{gather} 
	\label{lipschitz approx concl2} \sup_{\Omega_{\sigma}} |\nabla u| \leq C \mathcal{E}^{2\gamma} , \\ 
	\label{lipschitz approx concl3} \mathcal{L}^n(\Omega_{\sigma} \setminus K) 
		+ \|T\|((\Omega_{\sigma} \setminus K) \times \mathbb{R}^m) \leq C \mathcal{E}^{2+2\gamma} 
\end{gather}
for some constants $\gamma = \gamma(n,m,q) \in (0,1)$ and $C = C(n,m,q,\sigma,\mathcal{L}^n(\Omega)) \in (0,\infty)$ (independent of $\delta$).
\end{theorem}

\begin{remark}{\rm 
We shall consider blow-up limits of area-minimizing currents $T_k$ in $\Omega \times \mathbb{R}^m$ with $\mathcal{E}_k^2 = \int_{\Omega \times \mathbb{R}^m} |\vec T_k - \vec P_0|^2 \,d\|T_k\|(X) \rightarrow 0$, where $\vec T_k$ denotes the orienting $n$-vector field of $T_k$.  (See for instance Section~\ref{sec:blowup procedure}.)  In this case, we will apply Theorem~\ref{lipschitz approx thm} with $\delta = \delta_k$ and $T = T_k$ where $\delta_k$ is chosen so that $\delta_k^{-n/2} \mathcal{E}_k \rightarrow 0$.  Then \eqref{lipschitz approx concl2} and 
\eqref{lipschitz approx concl3} give us estimates for the Lipschitz approximation of $T = T_k$ on $\Omega_{\sigma}$ with constants $C$ independent of $k$, provide $k$ is large enough that $\sigma \geq \delta_k$.
}\end{remark}

\begin{proof}[Proof of Theorem~\ref{lipschitz approx thm}]
By the Vitali covering lemma, there is a finite collection of points $x_i \in \Omega$ ($1 \leq i \leq N$) such that, setting $d_i =\op{dist}(x_i, \partial \Omega)$, we have that $d_{i+1} \leq d_i$ for all $1 \leq i < N$, $\{B_{d_i/16}(x_i)\}$ covers $\Omega_{\delta/2}$, and $\{B_{d_i/80}(x_i)\}$ is pairwise disjoint.  Note that if $B_{d_i/8}(x_i) \cap B_{d_j/8}(x_j) \neq \emptyset$ then $|d_i - d_j| \leq |x_i - x_j| \leq (d_i+d_j)/8$ and thus $7d_i/9 \leq d_j \leq 9d_i/7$.  For each $\sigma \geq \delta/2$, if $B_{d_i/8}(x_i) \cap \Omega_{\sigma} \neq \emptyset$ then $d_i > 8\sigma/9$.  In particular, setting $\sigma = \delta/2$, $d_i > 4\delta/9$ for all $1 \leq i \leq N$.  For each $\sigma \geq \delta/2$ let $N(\sigma)$ be the largest integer such that $d_{N(\sigma)} > 8\sigma/9$.  Obviously $\{B_{d_i/8}(x_i)\}_{1 \leq i \leq N(\sigma)}$ covers $\Omega_{\sigma}$.  Since $\{B_{d_i/80}(x_i)\}$ is a collection of pairwise disjoint balls in $\Omega$, it follows that 
\begin{equation*}
	N(\sigma) \,\omega_n (\sigma/90)^n 
	\leq \sum_{i=1}^{N(\sigma)} \mathcal{L}^n(B_{\sigma/90}(x_i)) 
	\leq \sum_{i=1}^{N(\sigma)} \mathcal{L}^n(B_{d_i/80}(x_i)) 
	\leq \mathcal{L}^n(\Omega) 
\end{equation*}
and thus $N(\sigma) \leq C(n) \,\sigma^{-n} \mathcal{L}^n(\Omega)$.  Let $\{\psi_i\}$ smooth partition of unity for $\Omega_{\delta}$ subordinate to $\{B_{d_i/8}(x_i)\}$ such that 
\begin{equation}\label{lipschitz approx eqn1}
	0 \leq \psi_i \leq 1, \quad \op{spt} \psi_i \subseteq B_{d_i/8}(x_i), \quad |\nabla \psi_i| \leq \frac{C(n)}{d_i}, \quad \sum_{i=1}^N \psi_i = 1 
\end{equation}
(see for instance~\cite[3.1.13]{Fed69}).  By \eqref{lipschitz approx hyp} and the fact that $d_i > 4\delta/9$, for each $i$ 
\begin{equation*}
	\frac{1}{d_i^n} \int_{B_{d_i}(x_i) \times \mathbb{R}^m} |\vec T - \vec P_0|^2 \,d\|T\|(X) 
		\leq d_i^{-n} \mathcal{E}^2 \leq (4\delta/9)^{-n} \cdot \varepsilon^2 \delta^n = (9/4)^{n} \varepsilon^2 .
\end{equation*}
Hence provided $\varepsilon$ is sufficiently small, by~\cite[Corollary~3.29]{Almgren} (or \cite[Theorem~2.4]{DeLSpa1}), for each $i$ there exists a Lipschitz function $u_i : B_{d_i/4}(x_i) \rightarrow \mathcal{A}_q(\mathbb{R}^m)$ and $\mathcal{L}^n$-measurable set $K_i \subseteq B_{d_i/4}(x_i)$ such that 
\begin{gather}
	\label{lipschitz approx eqn2} T \llcorner (K_i \times \mathbb{R}^m) = (\op{graph} u_i) \llcorner (K_i \times \mathbb{R}^m) , \\
	\label{lipschitz approx eqn3} \sup_{B_{d_i/4}(x_i)} |\nabla u_i| \leq C  d_i^{-n\gamma} \mathcal{E}^{2\gamma} , \\
	\label{lipschitz approx eqn4} \mathcal{L}^n(B_{d_i/4}(x_i) \setminus K_i) + \|T\|((B_{d_i/4}(x_i) \setminus K_i) \times \mathbb{R}^m) 
		\leq C d_i^{-n\gamma} \mathcal{E}^{2+2\gamma} \leq C (9/4)^{n + n\gamma} \varepsilon^{2+2\gamma} d_i^n
\end{gather}
for some constants $\gamma = \gamma(n,m,q) \in (0,1)$ and $C = C(n,m,q) \in (0,\infty)$.  Note that if $B_{d_i/8}(x_i) \cap B_{d_j/8}(x_j) \neq \emptyset$, then using the fact that $d_j \geq 7d_i/9$, there exists a ball of radius $\min\{d_i/16,d_j/16\} \geq 7d_i/144$ contained in $B_{d_i/4}(x_i) \cap B_{d_j/4}(x_j)$.  Thus by \eqref{lipschitz approx eqn4}, 
\begin{equation*}
	\mathcal{L}^n(B_{d_i/4}(x_i) \cap B_{d_j/4}(x_j)) \geq \omega_n (7d_i/144)^n 
		> 2C (9/4)^{n + n\gamma}\varepsilon^{2+2\gamma} d_i^n \geq \mathcal{L}^n(B_{d_i/4}(x_i) \cap B_{d_{j}/4}(x_{j}) \setminus (K_i \cap K_j)) 
\end{equation*}
provided $\varepsilon$ is sufficiently small, where $C$ is as in \eqref{lipschitz approx eqn4}.  Hence $\mathcal{L}^n(K_i \cap K_j) > 0$.  By \eqref{lipschitz approx eqn2}, $u_i(x) = \langle T,\pi_{P_0},x\rangle = u_j(x)$ for $\mathcal{L}^n$-a.e.~$x \in K_i \cap K_j$.  Letting $z \in K_i \cap K_j$ such that $u_i(z) = u_j(z)$ and using \eqref{lipschitz approx eqn3} and the fact that $d_j \leq 9d_i/7$, 
\begin{align}\label{lipschitz approx eqn5}
	\sup_{B_{d_i/8}(x_i) \cap B_{d_j/8}(x_j)} \mathcal{G}(u_i,u_j) 
		&\leq \sup_{B_{d_i/4}(x_i)} \mathcal{G}(u_i,u_i(z)) + \sup_{B_{d_j/4}(x_i)} \mathcal{G}(u_j,u_j(z))
		\\&\leq C d_i^{-\gamma n} \mathcal{E}^{2\gamma} (d_i/2 + d_j/2) \leq \tfrac{8}{7} \,C  d_i^{1-\gamma n} \mathcal{E}^{2\gamma} , \nonumber 
\end{align}
where $C$ is as in \eqref{lipschitz approx eqn3}.  By~\cite[Definition~1.1(6) and Theorem~1.3]{Almgren}, there exists an integer $L(q,m) \geq 1$ and bi-Lipschitz embedding $\boldsymbol{\xi} : \mathcal{A}_q(\mathbb{R}^m) \rightarrow \mathbb{R}^L$ and $\boldsymbol{\rho} : \mathbb{R}^L \rightarrow \mathcal{Q}$ such that $\op{Lip} \boldsymbol{\xi} \leq 1$, $\op{Lip} \boldsymbol{\xi}^{-1} |_{\mathcal{Q}} \leq C(m,q)$, and $\op{Lip} \boldsymbol{\rho} \leq C(m,q)$, where $\mathcal{Q} = \boldsymbol{\xi}(\mathcal{A}_q(\mathbb{R}^m))$.  Define 
\begin{gather}
	\label{lipschitz approx eqn6} K = \Omega_{\delta} \setminus \bigcup_{i=1}^N (B_{d_i/8}(x_i) \setminus K_i) , \\
	\label{lipschitz approx eqn7} u(x) = (\boldsymbol{\xi}^{-1} \circ \boldsymbol{\rho})\left( \sum_{i=1}^N \psi_i(x) \,\boldsymbol{\xi}(u_i(x)) \right) 
		\text{ for each } x \in \Omega_{\delta} .
\end{gather}
By \eqref{lipschitz approx eqn6} and \eqref{lipschitz approx eqn2}, for each $i \in \{1,2,\ldots,n\}$ we have that $K \cap B_{d_i/8}(x_i) \subseteq K_i$ and $u_i(x) = \langle T, \pi_{P_0}, x \rangle$ for $\mathcal{L}^n$-a.e.~$x \in K \cap B_{d_i/8}(x_i)$.  Thus by the fact that $\sum_{i=1}^N \psi_i = 1$ and \eqref{lipschitz approx eqn7}, $u(x) = \langle T, \pi_{P_0}, x \rangle$ for $\mathcal{L}^n$-a.e.~$x \in K$.  Hence $u_i(x) = u(x)$ for all $x \in K \cap B_{d_i/8}(x_i)$ and all $1 \leq i \leq N$.  It follows using \eqref{lipschitz approx eqn2} that \eqref{lipschitz approx concl1} holds true.

Now fix $\sigma \geq \delta$ and let $x \in \Omega_{\sigma}$.  Recalling that $\{B_{d_i/8}(x_i)\}_{1 \leq i \leq N(\sigma)}$ covers $\Omega_{\sigma}$, suppose that $x \in B_{d_i/8}(x_i)$ for some $i \in \{1,2,\ldots,N(\sigma)\}$.  Then by \eqref{lipschitz approx eqn1} and \eqref{lipschitz approx eqn7} 
\begin{equation*}
	u(x) = (\boldsymbol{\xi}^{-1} \circ \boldsymbol{\rho})\left( \boldsymbol{\xi}(u_i(x)) 
		+ \sum_{j=1}^{N(\sigma)} \psi_j(x) \,(\boldsymbol{\xi}(u_j(x)) - \boldsymbol{\xi}(u_i(x))) \right)
\end{equation*}
and thus using \eqref{lipschitz approx eqn1}, \eqref{lipschitz approx eqn3}, and \eqref{lipschitz approx eqn5} 
\begin{equation*}
	|\nabla u(x)| \leq C \sum_{j=1}^{N(\sigma)} \big( |\nabla_j \psi_j(x)| \,\mathcal{G}(u_j(x),u_i(x)) + \psi_j(x) \,|\nabla u_j(x)| \big) 
	\leq C \,N(\sigma) \,\sigma^{-n\gamma} \mathcal{E}^{2\gamma} \leq C_1 \mathcal{E}^{2\gamma} 
\end{equation*}
where $C = C(n,m,q) \in (0,\infty)$ and $C_1 = C_1(n,m,q,\sigma,\mathcal{L}^n(\Omega)) \in (0,\infty)$ are constants, thereby proving \eqref{lipschitz approx concl2}.  Using \eqref{lipschitz approx eqn6} and \eqref{lipschitz approx eqn4}, we obtain \eqref{lipschitz approx concl3}.
\end{proof}

Let $T$, $u$, and $K$ be as in Theorem~\ref{lipschitz approx thm}.  Given an $m \times n$ matrix $p = [p_i^{\kappa}]_{1 \leq i \leq n,\,1 \leq \kappa \leq m}$, let 
\begin{equation*}
	G_{ij}(p) = \delta_{ij} + \sum_{\kappa=1}^m p_i^{\kappa} p_j^{\kappa}
\end{equation*}
for all $i,j \in \{1,2,\ldots,n\}$.  Let $[G^{ij}(p)]$ denote the inverse matrix of $[G_{ij}(p)]$ and $G(p)$ denote the determinant of $[G_{ij}(p)]$.  By Taylor's theorem, 
\begin{equation}\label{lipschitz approx eqn8}
	G^{ij}(p) = \delta_{ij} - \sum_{\kappa=1}^m p_i^{\kappa} p_j^{\kappa} + O(|p|^4), \quad 
	\sqrt{G(p)} = 1 + \frac{1}{2} \,|p|^2 + O(|p|^4) . 
\end{equation}
Let $u(x) = \sum_{l=1}^q \llbracket u_l(x) \rrbracket$ for each $x \in \Omega$ and $Du(x) = \sum_{l=1}^q \llbracket Du_l(x) \rrbracket$ for $\mathcal{L}^n$-a.e.~$x \in \Omega,$ following the conventions from Subsection~\ref{sec:prelim_multivalued}.  Recall from Subsection~\ref{sec:prelim_multivalued} that we can regard $\op{graph} u$ as a rectifiable current in $\Omega_{\delta} \times \mathbb{R}^m$.  Using \eqref{lipschitz approx concl1}, \eqref{lipschitz approx concl2}, and \eqref{lipschitz approx concl3}, 
\begin{equation}\label{lipschitz approx eqn9}
	\int_{\Omega_{\sigma} \times \mathbb{R}^m} \zeta(X) \,d\|T\|(X) 
	= \int_{\Omega_{\sigma}} \sum_{l=1}^q \zeta(x,u_l(x)) \,\sqrt{G(Du_l(x))} \,d\mathcal{L}^n(x) + \mathcal{R} 
\end{equation}
for each $\sigma \in (0,\delta)$ and each bounded $\|\op{graph} u\|+\|T\|$-measurable function $\zeta : \Omega_{\sigma} \times \mathbb{R}^m \rightarrow \mathbb{R}$, where assuming $\varepsilon = \varepsilon(n,m,q,\sigma,|\Omega|)$ is sufficiently small
\begin{equation*}
	|\mathcal{R}| \leq \sup |\zeta| \,\big( (1+C\varepsilon^{2\gamma}) \,\mathcal{L}^n(\Omega_{\sigma} \setminus K) 
		+ \|T\|((\Omega_{\sigma} \setminus K) \times \mathbb{R}^m) \big) 
	\leq C \mathcal{E}^{2+2\gamma} \sup |\zeta| 
\end{equation*}
for some constants $C = C(n,m,q,\sigma,|\Omega|) \in (0,\infty)$.  Note that at $\mathcal{H}^n$-a.e.~point $X = (x,u_l(x)) \in \op{spt}\op{graph}u$, 
\begin{align}
	\label{lipschitz approx eqn10} \frac{1}{2} \,|\vec S_X - \vec P|^2 &= 1 - \frac{1}{\sqrt{G(Du_l(x))}} , \\ 
	\label{lipschitz approx eqn11} \|\pi_{S_X} - \pi_P\|^2 &= \sum_{i,j=1}^n G^{ij}(Du_l(x)) \,D_i u_l(x) \,D_j u_l(x) ,  
\end{align}
where $S_X$ is the approximate tangent plane to $\op{spt}\op{graph}u$ at $X$ and $\vec S_X$ is the orientation $n$-vector $\op{graph}u$ at $X$.  

\section{Planar frequency function and its approximate monotonicity}\label{sec:freq mono sec}

In this section we introduce a frequency function for an area-minimizing current $T$ relative to a plane $P$.

\begin{definition}\label{freq defn}{\rm 
Let $Z \in \mathbb{R}^{n+m}$ and $\rho_0 > 0$.  Let $T$ be an $n$-dimensional locally area-minimizing rectifiable current of $\mathbf{C}_{\rho_0}(Z,P)$ such that 
\begin{equation}\label{freq defn hyp}
	(\partial T) \llcorner \mathbf{C}_{\rho_0}(Z,P) = 0, \quad \sup_{X \in \op{spt} T \cap \mathbf{C}_{\rho_0}(Z,P)} \op{dist}(X,Z + P) < \infty. 
\end{equation}
Let $\phi : [0,\infty) \rightarrow \mathbb{R}$ be the Lipschitz function defined by 
\begin{equation}\label{freq phi defn}
	\phi(s) = \begin{cases} 
		1 &\text{if } 0 \leq s < 1/2 \\
		2 - 2s &\text{if } 1/2 \leq s < 1 \\
		0 &\text{if } 1 \leq s < \infty.
	\end{cases} 
\end{equation}
For $\rho \in (0, \rho_{0}],$ let 
\begin{align}
	\label{H defn0} H_{T,P,Z}(\rho) &= 2\rho^{1-n} \int_{G_n(\mathbf{C}_{\rho}(Z,P) \setminus \overline{\mathbf{C}_{\rho/2}(Z,P)})} 
		\op{dist}^2(X,Z+P) \,|\nabla^S r|^2 \,\frac{1}{r} \,d|T|(X,S) \;\; \mbox{and} \\
	\label{D defn} D_{T,P,Z}(\rho) &= \frac{1}{2} \,\rho^{2-n} \int \|\pi_S - \pi_P\|^2 \,\phi(r/\rho) \,d|T|(X,S) , 
\end{align}
where $|T|$ denotes the $n$-dimensional integral varifold associated with $T$, $r(X) = |\pi_P(X-Z)|$, $\nabla^S$ is the gradient with respect to the plane $S$, and $\|\pi_S - \pi_P\|$ denotes the Frobenius norm of $\pi_S - \pi_P$.  We define the \emph{planar frequency function} $N_{T,P,Z}$ of $T$ at $Z$ relative to the plane $P$ by 
\begin{equation*}
	N_{T,P,Z}(\rho) = \frac{D_{T,P,Z}(\rho)}{H_{T,P,Z}(\rho)} 
\end{equation*}
whenever $H_{T,P,Z}(\rho) > 0$ (see Remark~\ref{H zero rmk} below).}
\end{definition}

\begin{remark}\label{phi deriv rmk}{\rm 
We will often write 
\begin{equation}\label{H defn}
	H_{T,P,Z}(\rho) = -\rho^{1-n} \int \op{dist}^2(X,Z+ P) \,|\nabla^S r|^2 \,\frac{1}{r} \,\phi'(r/\rho) \,d|T|(X,S),
\end{equation}
where we adopt the convention that $\phi'(1/2) = \phi'(1) = 0$.  Note that for $|T|$-a.e.~$(X,S) \in G_n(\partial \mathbf{C}_{\rho/2}(Z,P))$, $S$ is the approximate tangent plane to $\op{spt} T$ at $X$ and $S$ is tangent to $\partial \mathbf{C}_{\rho/2}(Z,P)$, thus $\nabla^S r(X) = 0$.  By similar reasoning, $\nabla^S r(X) = 0$ for $|T|$-a.e.~$(X,S) \in G_n(\partial \mathbf{C}_{\rho}(Z,P))$. 
}\end{remark}

\begin{remark} {\rm (1) Almgren in \cite{Almgren} first introduced and used a frequency function.  The frequency function in \cite{Almgren} is defined ``extrinsically,''  as 
a functional associated with a multi-valued function $f$---specifically, either for a multivalued locally Dirichlet energy minimizing function $f$ over a plane, or for a multi-valued Lipschitz function $f$---the ``normal map''---over a center manifold, whose graph approximates a locally area minimizing rectifiable current lying close to a plane.  
In \cite{DeLSpa3}, the authors observed that in either case, it is technically simpler to incorporate the fixed cut-off function 
$\phi$ into the definition rather than directly using the relevant energy and height terms (i.e.\ $D(\rho)$ and $H(\rho)$ in the notation of \cite{Almgren}), or taking an integral average of these quantities over an interval of scales $\rho$,  as done in \cite{Almgren}.

\noindent
(2) In contrast to the the frequency functions in \cite{Almgren} and \cite{DeLSpa3}, our planar frequency function $N_{T,P,Z}(\rho)$ is intrinsic, and is defined in terms of geometric quantities involving $P$ and $T$ and in terms of integration over $T$.  That way $N_{T,P,Z}(\rho)$ is well-defined for \emph{any} area-minimizing integral current $T$ (provided \eqref{freq defn hyp} holds true) and $N_{T,P,Z}(\rho)$ does not depend on $T$ being close to $P$ or on a choice of Lipschitz approximation of $T$.  This is an important feature of the planar frequency function, which we capitalise on in our analysis. (See for instance Lemma~\ref{freq of cones lemma} below.)}
\end{remark}  

The main result of this section is the following approximate monotonicity of the planar frequency function under the assumption that $T$  is decaying towards the plane $P$ in the $L^{2}$ sense.

\begin{theorem}\label{mono freq thm} 
For each positive integer $q$ there exists $\delta = \delta(n,m,q) \in (0,1)$ and $\eta_0 = \eta_0(n,m,q) \in (0,1)$ such that the following holds true.  Let $Z \in \mathbb{R}^{n+m}$ and $\rho_0 > 0$.  Let $P$ be an $n$-dimensional plane in $\mathbb{R}^{n+m}$ and $T$ be an $n$-dimensional locally area-minimizing rectifiable current in the open cylinder $\mathbf{C}_{7\rho_0/4}(Z,P)$ such that 
\begin{equation}\label{mono freq mass hyp} 
	(\partial T) \llcorner \mathbf{C}_{7\rho_0/4}(Z,P) = 0, \quad
	\Theta(T,Z) \geq q , \quad 
	\|T\|(\mathbf{C}_{7\rho_0/4}(Z,P)) \leq (q + \delta) \,\omega_n (7\rho_0/4)^n .
\end{equation}
Suppose that for some $\eta \in (0,\eta_0]$, $\sigma_0 \in (0,\rho_0)$ and  $\alpha \in (0,1),$
\begin{equation}\label{mono freq decay hyp} 
	\frac{1}{\omega_{n} (7\rho/4)^{n+2}}\int_{{\mathbf C}_{7\rho/4}(Z, P)} {\rm dist}^{2}(X, Z + P) \, d\|T\|(X)  \leq \eta^{2} \Big(\frac{\rho}{\rho_0}\Big)^{2\alpha} 
\end{equation}
for all $\rho \in [\sigma_0,\rho_0]$.  Then:
\begin{eqnarray}\label{mono freq initial}
	&&\exp\left( \frac{C_1 \eta^{2\gamma}}{2\alpha\gamma} \,(\sigma/\rho_{0})^{2\alpha\gamma} 
		+ \frac{C_2}{\gamma} \,\left(\rho_{0}^{-2}D_{T,P,Z}(\sigma)\right)^{\gamma} \right)N_{T,P,Z}(\sigma) \leq \nonumber\\ 
		&&\hspace{1.5in} \exp\left( \frac{C_1 \eta^{2\gamma}}{2\alpha\gamma} \,(\rho/\rho_{0})^{2\alpha\gamma} 
		+ \frac{C_2}{\gamma} \,\left(\rho_{0}^{-2}D_{T,P,Z}(\rho)\right)^{\gamma} \right) N_{T,P,Z}(\rho),  
		\end{eqnarray} 
and hence 
\begin{equation}\label{mono freq concl}
	N_{T,P,Z}(\sigma) \leq e^{\frac{C \eta^{2\gamma}}{2\alpha\gamma} (\rho/\rho_0)^{2\alpha\gamma}} N_{T,P,Z}(\rho)
\end{equation}
for all $\sigma$, $\rho$ with $\sigma_0 \leq \sigma < \rho \leq \rho_0$ provided $H_{T,P,Z}(\tau) > 0$ for all $\tau \in [\sigma,\rho].$ Here 
$C_{1} = C_{1}(n,m , q)$, $C_{2}=  C_{2}(n, m, q)$, $\gamma = \gamma(n,m,q) \in (0,1)$ and $C = C_{1} + 2\alpha \, \overline{C} \in (0,\infty)$  are fixed constants (independent of $\eta$), where $\overline{C} = \overline{C}(n, m, q)$. 
\end{theorem} 

Note that by \eqref{mono freq decay hyp}, Lemma~\ref{Allard height lemma} and the fact that $\left(\partial T \right)\llcorner \mathbf{C}_{7\rho_{0}/4}(Z, P) = 0$, the conditions \eqref{freq defn hyp} hold true and thus $N_{T,P,Z}(\rho)$ is defined for all $\rho \in (0,\rho_0]$.  The remainder of this section will focus on the proofs of Theorem~\ref{mono freq thm}.  Without loss of generality we may let $Z = 0$, $\rho_0 = 1$, and $P = P_0 = \mathbb{R}^n \times \{0\}$.  We will write each point $X = (x_1,x_2,\ldots,x_{n+m}) \in \mathbb{R}^{n+m}$ as $X = (x,y)$ where $x = (x_1,x_2,\ldots,x_n) \in \mathbb{R}^n$ and $y = (x_{n+1},\ldots,x_{n+m}) \in \mathbb{R}^m$.  Note that then $r(X) = |x|$ and $\op{dist}(X,P) = |y|$.

\subsection{Variational formulas} 
In Lemma~\ref{freq identities lemma} below we compute the derivatives of $H_{T,P,0}(\rho)$ and $D_{T,P,0}(\rho)$ and establish identity \eqref{D to I} for $D_{T,P,0}(\rho)$.  
Eventhough we state Lemma~\ref{freq identities lemma} for area-minimizing currents $T$, the only variational property of $T$ used in its proof is stationarity with respect to area (i.e.\ the validity of the first variation formula for area \eqref{first variation}); thus,  the lemma in fact holds true whenever $|T|$ is an stationary integral varifold in $\mathbf{C}_1(0)$ with $(\partial T) \llcorner {\mathbf C}_{1} = 0$. Our proof of Theorem~\ref{mono freq thm} however uses certain results and estimates which are specific for area minimising $T$.

\begin{lemma}\label{freq identities lemma}
Let $T$ be an $n$-dimensional locally area-minimizing rectifiable current in $\mathbf{C}_1(0)$. Suppose that \eqref{freq defn hyp} holds true with $Z = 0$, $\rho_{0} = 1$ and $P = {\mathbb R}^{n} \times \{0\}$.  Then $H_{T,P,0}$ and $D_{T,P,0}$ are absolutely continuous functions on $(0,1],$ and satisfy:
\begin{align}\label{H deriv} 
	H'_{T,P,0}(\rho) =\,& -2\rho^{-n} \int (0,y) \cdot \pi_{P^{\perp}}(\nabla^S r) \,\phi'(r/\rho) \,d|T|(X,S)  
		\\&- \rho^{-n} \int |y|^2 \,\frac{1}{r} \,\phi'(r/\rho) \left( n \,|\nabla^{S^\perp} r|^2 - \frac{1}{2} \,\|\pi_S - \pi_P\|^2 \right) d|T|(X,S) \nonumber
\end{align}
for $\mathcal{L}^1$-a.e.~$\rho \in (0,1]$, 
\begin{equation}\label{D to I} 
	D_{T,P,0}(\rho) = -\rho^{1-n} \int (0,y) \cdot \pi_{P^{\perp}}(\nabla^S r) \,\phi'(r/\rho) \,d|T|(X,S) 
\end{equation}
for all $\rho \in (0,1]$, and 
\begin{align}
	\label{D deriv} D'_{T,P,0}(\rho) =\,& - 2\rho^{-n} \int |\nabla^{S^\perp} r|^2 \,r \,\phi'(r/\rho) \,d|T|(X,S)   
		\\&+ 2 \rho^{1-n} \int \left( 1 - \frac{1}{4} \,\|\pi_S - \pi_P\|^2 - J\pi_P \right) \left( n \,\phi(r/\rho) + \frac{r}{\rho} \,\phi'(r/\rho) \right) d|T|(X,S) 
		\nonumber
\end{align}
for  $\mathcal{L}^1$-a.e.~$\rho \in (0,1]$.  Here, $\nabla^S$ and $\nabla^{S^{\perp}}$ denote the gradient with respect to the linear subspaces $S$ and $S^{\perp}$ respectively.  $J\pi_P$ is the (signed) Jacobian of $\pi_P$ on $T$ defined by $d\pi_{P\#} \vec T(X) = J\pi_P(X) \,\vec P$ for $\|T\|$-a.e.~$X \in \mathbf{C}_1(0)$, where $\vec P = e_1 \wedge e_2 \wedge\cdots\wedge e_n$ is the orientation of $P$.
\end{lemma}

\begin{remark}\label{H zero rmk}{\rm
Let $T$ be as in Lemma~\ref{freq identities lemma}.  We have that $H_{T,P,0}(\rho) = 0$ if and only if $\op{spt} T \cap \mathbf{C}_{\rho} \subset P$.  Clearly by the definition of $H_{T,P,0}(\rho)$ in \eqref{H defn0}, if $\op{spt} T \cap \mathbf{C}_{\rho} \subseteq P$ then $H_{T,P,0}(\rho) = 0$.  To see the converse, suppose $H_{T,P,0}(\rho) = 0$.  Then by \eqref{H defn0}, $X \in P$ or $\nabla^S r(X) = 0$ for $|T|$-a.e. $(X,S) \in G_n(\mathbf{C}_{\rho} \setminus \overline{\mathbf{C}_{\rho/2}})$.  Hence by \eqref{D to I}, $D_{T,P,0}(\rho) = 0$.  Thus by \eqref{D defn}, $P$ must be the approximate tangent plane to $\op{spt} T$ at $\|T\|$-a.e.~$X \in \mathbf{C}_{\rho}$.  This together with $|T|$ being stationary implies that $\op{spt} T \cap {\mathbf C}_{\rho}$ is contained in a finite union of $n$-dimensional planes parallel to $P$.  Also, $|\nabla^S r(X)| = |\nabla r(X)| = 1 \neq 0$ for $|T|$-a.e.~$(X,S) \in G_n(\mathbf{C}_{\rho} \setminus \overline{\mathbf{C}_{\rho/2}})$, so $\op{spt} T \cap \mathbf{C}_{\rho} \setminus \overline{\mathbf{C}_{\rho/2}} \subset P$.  It follows that $\op{spt} T \cap \mathbf{C}_{\rho} \subset P$.
}\end{remark}

\begin{proof}[Proof of Lemma~\ref{freq identities lemma}] 
Let us first consider the case where $H_{T,P,0}$ and $D_{T,P,0}$ are defined by \eqref{H defn} and \eqref{D defn} for $\phi : [0,\infty) \rightarrow \mathbb{R}$ a smooth function such that $\phi(s) = 1$ if $0 \leq s \leq 1/2$ and $\phi(s) = 0$ if $s \geq 1$.  The case where $\phi$ is given by \eqref{freq phi defn} will later follow by approximation.  Clearly when $\phi$ is smooth, $H_{T,P,0}$ and $D_{T,P,0}$ are continuously differentiable on $(0,1]$. 
By direct differentiation,
\begin{align*}
	H'_{T,P,0}(\rho) =\,& (n-1) \,\rho^{-n} \int |y|^2 \,|\nabla^S r|^2 \,\frac{1}{r} \,\phi'(r/\rho) \,d|T|(X,S) 
		\\&+ \rho^{-n-1} \int |y|^2 \,|\nabla^S r|^2 \,\phi''(r/\rho) \,d|T|(X,S). 
	\end{align*}
Setting $\zeta = |y|^2 \,\phi'(r/\rho) \,\nabla r$ in \eqref{first variation}, and noting that $\nabla r = \frac{(x, 0)}{r},$
\begin{align}
	\nabla^S (|y|^2) \cdot \nabla^S r &= \nabla (|y|^2) \cdot \nabla^S r = 2 \,(0,y) \cdot \nabla^S r = 2 \,(0,y) \cdot \pi_{P^{\perp}}(\nabla^S r) \;\; \mbox{and} \nonumber \\
	\label{H deriv eqn1} \op{div}_S(x,0) &= \sum_{i,j=1}^n (\tau_i \cdot e_j)^2 = n - \frac{1}{2} \,\|\pi_S - \pi_P\|^2 
\end{align}
for each $X = (x,y) \in \mathbf{C}_1$ and each $n$-dimensional plane $S \subset \mathbb{R}^{n+m}$, where $\tau_1,\tau_2,\ldots,\tau_n$ is an orthonormal basis for $S$,   gives us 
\begin{align}\label{H deriv eqn2}
	&2 \int (0,y) \cdot \pi_{P^{\perp}}(\nabla^S r) \,\phi'(r/\rho) \,d|T|(X,S) + \frac{1}{\rho} \int |y|^2 \,|\nabla^S r|^2  \,\phi''(r/\rho) \,d|T|(X,S) 
	\\&+ \int |y|^2 \left( n - \frac{1}{2} \,\|\pi_S - \pi_P\|^2 - |\nabla^S r|^2 \right) \frac{1}{r} \,\phi'(r/\rho) \,d|T|(X,S) = 0 \nonumber
\end{align}
for all $\rho \in (0,1]$. Using \eqref{H deriv eqn2} and the above expression for $H^{\prime}_{T, P, 0}$ gives 
\begin{align*}
	H'_{T,P,0}(\rho) =\,& -2\rho^{-n} \int (0,y) \cdot \pi_P^{\perp}(\nabla^S r) \,\phi'(r/\rho) \,d|T|(X,S) 
		\\&- \rho^{-n} \int |y|^2 \left( n \,|\nabla^{S^\perp} r|^2 - \frac{1}{2} \,\|\pi_S - \pi_P\|^2 \right) \frac{1}{r} \,\phi'(r/\rho) \,d|T|(X,S) 
\end{align*}
for all $\rho \in (0,1]$, proving \eqref{H deriv}. 

To see \eqref{D to I}, we set $\zeta = \phi(r/\rho) \,(0,y)$ in \eqref{first variation}, noting that $(0,y) \cdot \nabla^S r = (0,y) \cdot \pi_{P^{\perp}}(\nabla^S r)$ and 
\begin{align*}
	\op{div}_S(0,y) = \sum_{i=1}^n \sum_{j=1}^m (\tau_i \cdot e_{n+j})^2 = \frac{1}{2} \,\|\pi_S - \pi_P\|^2 
\end{align*}
for each $X = (x,y) \in \mathbf{C}_1$ and each $n$-dimensional plane $S \subset \mathbb{R}^{n+m}$.  Thus 
\begin{equation*}
	\int \left( \frac{1}{2} \,\phi(r/\rho) \,\|\pi_S - \pi_P\|^2 + \frac{1}{\rho} \,\phi'(r/\rho) \,\pi_{P^{\perp}}(\nabla^S r) \cdot (0,y) \right) d|T|(X,S) = 0 
\end{equation*}
for all $\rho \in (0,1]$.  Rearranging terms gives us \eqref{D to I}. 

To see \eqref{D deriv}, we set $\zeta = \phi(r/\rho) \,(x,0)$ in \eqref{first variation} and use \eqref{H deriv eqn1} to obtain 
\begin{equation}\label{D deriv eqn1}
	\int \left( n \,\phi(r/\rho) - \frac{1}{2} \,\|\pi_S - \pi_P\|^2 \,\phi(r/\rho) +  \,|\nabla^S r|^2 \,\frac{r}{\rho} \,\phi'(r/\rho) \right) d|T|(X,S) = 0
\end{equation}
for all $\rho \in (0,1]$.  
We can replace $T$ with the multiplicity 1 current $\llbracket P \rrbracket$ in \eqref{D deriv eqn1} to obtain 
\begin{equation*}
	\int_P \left( n \,\phi(r/\rho) + \frac{r}{\rho} \,\phi'(r/\rho) \right) d\mathcal{H}^n = 0 
\end{equation*}
for all $\rho \in (0,1]$.  By the constancy theorem~\cite[Theorem~26.27]{SimonGMT}, $\pi_{P\#} T$ is an integer multiple of $\llbracket P \rrbracket \llcorner \mathbf{C}_1$.  Thus by the area-formula (see~\cite[Remark~27.2(3)]{SimonGMT}) 
\begin{equation}\label{D deriv eqn2}
	\int \left( n \,\phi(r/\rho) + \frac{r}{\rho} \,\phi'(r/\rho) \right) J\pi_P \,d\|T\|(X) = 0 
\end{equation}
for all $\rho \in (0,1]$.  By differentiating $D_{T,P,0}(\rho)$, 
\begin{align}\label{D deriv eqn3}
	D'_{T,P,0}(\rho) =\,& \frac{2-n}{2} \, \,\rho^{1-n} \int \|\pi_S - \pi_P\|^2 \,\phi(r/\rho) \,d|T|(X,S)  
		\\&- \frac{1}{2} \,\rho^{-n} \int \|\pi_S - \pi_P\|^2 \,r \,\phi'(r/\rho) \,d|T|(X,S) \nonumber
\end{align}
for all $\rho \in (0,1]$.  Adding $2\rho^{1-n} \,\cdot\,$\eqref{D deriv eqn1}$\,-2\rho^{1-n} \,\cdot\,$\eqref{D deriv eqn2}$\,+\,$\eqref{D deriv eqn3} and using the fact that $|\nabla^T r|^2 = 1 - |\nabla^{\perp} r|^2$ gives us \eqref{D deriv}.

Now suppose that $H_{T,P,0}$ and $D_{T,P,0}$ are defined by \eqref{H defn} and \eqref{D defn} where $\phi$ is given by \eqref{freq phi defn}.  For each $\varepsilon \in (0,1/2)$ let $\phi_{\varepsilon} : [0,\infty) \rightarrow \mathbb{R}$ be a smooth function such that $\phi_{\varepsilon}(s) = 1$ if $0 \leq s \leq 1/2$, $\phi_{\varepsilon}(s) = 2+3\varepsilon - (2+4\varepsilon) s$ if $1/2+\varepsilon \leq s \leq 1-\varepsilon$, $\phi_{\varepsilon}(s) = 0$ if $s \geq 1$, and $-2-4\varepsilon \leq \phi'_{\varepsilon} \leq 0$.  Let $H^{\varepsilon}_{T,P,0}$ and $D^{\varepsilon}_{T,P,0}$ be given by \eqref{H defn} and \eqref{D defn} with $\phi_{\varepsilon}$ in place of $\phi$.  Then \eqref{H deriv} with $\phi_{\varepsilon}$ in place of $\phi$ gives us  
\begin{align}\label{H deriv eqn3}
	H^{\varepsilon}_{T,P,0}(\rho) &- H^{\varepsilon}_{T,P,0}(\sigma) = -\int_{\sigma}^{\rho} 
		2\tau^{-n} \int (0,y) \cdot \pi_{P^{\perp}}(\nabla^S r) \,\phi'_{\varepsilon}(r/\tau) \,d|T|(X,S) \,d\tau  
		\\&- \int_{\sigma}^{\rho} \tau^{-n} \int |y|^2 \,\frac{1}{r} \,\phi'_{\varepsilon}(r/\tau) 
		\left( n \,|\nabla^{S^\perp} r|^2 - \frac{1}{2} \,\|\pi_S - \pi_P\|^2 \right) d|T|(X,S) \,d\tau \nonumber 
\end{align}
for all $0 < \sigma < \rho \leq 1$.  Notice that $\phi_{\varepsilon} \rightarrow \phi$ uniformly on $[0,\infty)$ and $\phi'_{\varepsilon} \rightarrow \phi'$ pointwise on $[0,\infty)$ (with $\phi'_{\varepsilon}(s) = \phi'(s) = 0$ if $s = 1/2,1$).  Hence letting $\varepsilon \rightarrow 0^+$ in \eqref{H deriv eqn3} using the dominated convergence theorem, we deduce that when $\phi$ is given by \eqref{freq phi defn}, $H_{T,P,0}$ is absolutely continuous and \eqref{H deriv} holds true for $\mathcal{L}^1$-a.e.~$\rho \in (0,1]$.  Similar reasoning shows that  \eqref{D deriv} and \eqref{D to I} hold true with $\phi$ is given by \eqref{freq phi defn}.
\end{proof}

\subsection{Bounding error terms in derivative of $H_{T,P,0}$}  Next in Lemma~\ref{H error lemma} we will estimate certain error terms in \eqref{H deriv}.  We start with a standard consequence of frequency monotonicity of Dirichlet energy minimising functions (Lemma~\ref{H uniq cont lemma0}) and its direct implication to area minimising currents with small excess (Lemma~\ref{H uniq cont lemma}).

\begin{lemma}\label{H uniq cont lemma0}
For every $c > 0$ there exists $\beta = \beta(n,m,q,c) > 0$ such that if $w \in W^{1,2}(B_1(0),\mathcal{A}_q(\mathbb{R}^m))$ is a Dirichlet energy minimizing $q$-valued function such that $\|w\|_{L^2(B_{1/2}(0))} \geq c \|w\|_{L^2(B_1(0))}$ then $\|w\|_{L^2(B_{1/10}(\xi))} \geq \beta \|w\|_{L^2(B_1(0))}$ for each $\xi \in B_{1/5}(0)$.
\end{lemma}

\begin{proof}
Fix $c > 0$ and suppose to the contrary that for $k = 1,2,3,\ldots$ there exists a Dirichlet energy minimizing $q$-valued function $w_k \in W^{1,2}(B_1(0),\mathcal{A}_q(\mathbb{R}^m))$ and point $\xi_k \in B_{1/5}(0)$ such that $\|w_k\|_{L^2(B_1(0))} = 1$, $\|w_k\|_{L^2(B_{1/2}(0))} \geq c$, and $\|w_k\|_{L^2(B_{1/10}(\xi_k))} < 1/k$.  After passing to a subsequence let $\xi_k \rightarrow \xi$ for some $\xi \in \overline{B_{1/5}(0)}$.  By the compactness theorem for Dirichlet energy minimzing $q$-valued functions 
(\cite[Theorem~2.15]{Almgren}, \cite[Proposition~2.11, Theorem~3.20]{DeLSpaDirMin}), after passing to a subsequence there is a Dirichlet energy minimizing $q$-valued function $w \in W^{1,2}(B_1(0),\mathcal{A}_q(\mathbb{R}^m))$ such that $w_k \rightarrow w$ uniformly on compact subsets of $B_1(0)$.  Clearly $\|w\|_{L^2(B_{1/2}(0))} \geq c > 0$ and $w = q \llbracket 0 \rrbracket$ on $B_{1/20}(\xi)$.  On the other hand, it is a standard consequence of the monotonicity of the frequency function associated with $w$ that $w = q \llbracket 0 \rrbracket$ on $B_{1/20}(\xi) \implies w = q \llbracket 0 \rrbracket$ on $B_{1/2}(0).$ This gives the desired contradiction proving the lemma.
\end{proof}

\begin{lemma}\label{H uniq cont lemma} 
For every  $c > 0$ there exists $\varepsilon = \varepsilon(n, m, q, c) \in (0, 1)$, $\delta = \delta(n,m,q,c) \in (0,1)$ and 
$\beta = \beta(n,m,q,c) \in (0, 1)$ such that if $T$ is an $n$-dimensional locally area-minimizing rectifiable current of $\mathbf{C}_8(0)$, $P = {\mathbb R}^{n} \times \{0\}$ and if 
\begin{gather} 
	\label{H uniq cont hyp1} (\partial T) \llcorner \mathbf{C}_8(0) = 0, \quad
    \sup_{X \in \op{spt} T} \op{dist}(X,P) < \infty , \quad 
    \pi_{P\#} T = q \llbracket B_8(0) \rrbracket , \\
	\label{H uniq cont hyp2} \frac{1}{\omega_n 8^{n+2}} \int_{\mathbf{C}_8(0)} \op{dist}^2(X,P) \,d\|T\|(X) < \varepsilon^2 , \\
	\label{H uniq cont hyp3} \int_{\mathbf{C}_{1/2}(0)} \op{dist}^2(X,P) \,d\|T\|(X) \geq c^2 \int_{\mathbf{C}_8(0)} \op{dist}^2(X,P) \,d\|T\|(X) , 
\end{gather}
then for every $\xi \in B_{1/5}(0)$ 
\begin{equation}\label{H uniq cont concl}
	\int_{\mathbf{C}_{1/10}(\xi)} \op{dist}^2(X,P) \,d\|T\|(X) \geq \beta^2 \int_{\mathbf{C}_8(0)} \op{dist}^2(X,P) \,d\|T\|(X) . 
\end{equation}
\end{lemma}

\begin{proof}
By Lemma~\ref{area excess to height lemma}, 
\begin{equation*}
    \mathcal{E}(T,P,\mathbf{C}_6(0))^{2} \leq C\int_{{\mathbf C}_8(0)} {\rm dist}^{2}(X, P) \, d\|T\|(X) < C \varepsilon^2 
\end{equation*}
for some constant $C = C(n, m) \in (0, \infty)$.  Thus by~\cite[Corollaries~3.29 \& 3.30]{Almgren} (or ~\cite[Theorem~2.4]{DeLSpa1}), there exists a Lipschitz function $u : B_1(0) \rightarrow \mathcal{A}_q(\mathbb{R}^m)$ and a set $K \subset B_1(0)$ such that 
\begin{gather}
	\label{H uniq cont eqn1} \sup_{B_1(0)} |u| \leq 1, \quad\quad \sup_{B_1(0)} |\nabla u| \leq C \varepsilon^{\gamma}, \\
	T \llcorner K \times \mathbb{R}^m = (\op{graph} u) \llcorner K \times \mathbb{R}^m, \nonumber \\
	|B_1(0) \setminus K| + \|T\|((B_1(0) \setminus K) \times \mathbb{R}^m) 
		\leq C \varepsilon^{\gamma} \int_{\mathbf{C}_8(0)} \op{dist}^2(X,P) \,d\|T\|(X) , \nonumber 
\end{gather}
where $\gamma = \gamma(n,m,q) \in (0,1)$ and $C = C(n,m,q) \in (0,\infty)$ are constants.  Let $\eta  \in (0,1)$ be a constant to be determined depending only on $n,m,q,c$.  By  \cite[Theorem~2.23]{Almgren} (or \cite[Theorem~2.6]{DeLSpa1}), if $\epsilon = \epsilon(n, m, q, \eta)$ is sufficiently small, there exists a Dirichlet energy minimizing $q$-valued function $w \in W^{1,2}(B_1(0),\mathcal{A}_q(\mathbb{R}^m))$ such that 
\begin{equation}\label{H uniq cont eqn2}
	\int_{B_1(0)} \mathcal{G}(u,w)^2 \leq \eta \int_{\mathbf{C}_8(0)} \op{dist}^2(X,P) \,d\|T\|(X) . 
\end{equation}
Notice that provided $\varepsilon$ and $\eta$ are sufficiently small depending only on $n, m, q$, we have by \eqref{lipschitz approx eqn8}, \eqref{lipschitz approx eqn9}, \eqref{H uniq cont eqn1}, and \eqref{H uniq cont eqn2}, 
\begin{align}\label{H uniq cont eqn3}
	\int_{B_1(0)} |w|^2 &\leq 2 \int_{B_1(0)} |u|^2 + 2\eta \int_{\mathbf{C}_8(0)} \op{dist}^2(X,P) \,d\|T\|(X)  
		\\&\leq 2 \int_{\mathbf{C}_1(0)} \op{dist}^2(X,P) \,d\|T\|(X) + (2\eta + 2C\varepsilon^{2\gamma}) \int_{\mathbf{C}_8(0)} \op{dist}^2(X,P) \,d\|T\|(X) \nonumber 
		\\&\leq 4 \int_{\mathbf{C}_8(0)} \op{dist}^2(X,P) \,d\|T\|(X) . \nonumber 
\end{align}
Similarly, by \eqref{lipschitz approx eqn8}, \eqref{lipschitz approx eqn9}, \eqref{H uniq cont eqn1}, \eqref{H uniq cont eqn2}, and \eqref{H uniq cont hyp3}, 
\begin{align}\label{H uniq cont eqn4}
	\int_{B_{1/2}(0)} |w|^2 &\geq \frac{1}{2} \int_{B_{1/2}(0)} |u|^2 - \eta \int_{\mathbf{C}_8(0)} \op{dist}^2(X,P) \,d\|T\|(X) 
		\\&\geq \frac{1}{2} \int_{\mathbf{C}_{1/2}(0)} \op{dist}^2(X,P) \,d\|T\|(X) 
			- (\eta + C \varepsilon^{2+\gamma}) \int_{\mathbf{C}_8(0)} \op{dist}^2(X,P) \,d\|T\|(X) \nonumber
		\\&\geq \left( \frac{c^2}{2} - \eta - C \varepsilon^{2+\gamma} \right) \int_{\mathbf{C}_8(0)} \op{dist}^2(X,P) \,d\|T\|(X) \nonumber 
		\\&\geq \frac{c^2}{4} \int_{\mathbf{C}_8(0)} \op{dist}^2(X,P) \,d\|T\|(X) , \nonumber
\end{align}
provided $\epsilon, \eta$ are sufficiently small depending only on $n, m, q, c$, where $C = C(n,m,q) \in (0,\infty)$ is a constant.  Hence combining \eqref{H uniq cont eqn3} and \eqref{H uniq cont eqn4}, 
\begin{align*}
	\int_{B_{1/2}(0)} |w|^2 \geq \frac{c^2}{16} \int_{B_1(0)} |w|^2 . 
\end{align*}
Therefore, by Lemma~\ref{H uniq cont lemma0} with $c/4$ in place of $c$, there exists $\beta = \beta(n,m,q,c) > 0$ such that for each $\xi \in B_{1/5}(0)$ 
\begin{equation}\label{H uniq cont eqn5}
	\int_{B_{1/10}(\xi)} |w|^2 \geq \beta^2 \int_{B_1(0)} |w|^2 . 
\end{equation}
Provided $\varepsilon$ and $\eta$ are sufficiently small depending only on $n, m, q, c$, by \eqref{lipschitz approx eqn8}, \eqref{lipschitz approx eqn9}, \eqref{H uniq cont eqn1}, \eqref{H uniq cont eqn2}, \eqref{H uniq cont eqn5}, and \eqref{H uniq cont eqn4}, 
\begin{align*}
	\int_{\mathbf{C}_{1/10}(\xi)} \op{dist}^2(X,P) \,d\|T\|(X) 
		&\geq \int_{B_{1/10}(\xi)} |u|^2 - C \varepsilon^{\gamma} \int_{\mathbf{C}_8(0)} \op{dist}^2(X,P) \,d\|T\|(X)
		\\&\geq \frac{1}{2} \int_{B_{1/10}(\xi)} |w|^2 - (\eta + C \varepsilon^{2+\gamma}) \int_{\mathbf{C}_8(0)} \op{dist}^2(X,P) \,d\|T\|(X)
		\\&\geq \frac{\beta^2}{2} \int_{B_1(0)} |w|^2 - (\eta + C \varepsilon^{2+\gamma}) \int_{\mathbf{C}_8(0)} \op{dist}^2(X,P) \,d\|T\|(X)
		\\&\geq \left( \frac{\beta^2 c^2}{8} - \eta - C \varepsilon^{2+\gamma} \right) \int_{\mathbf{C}_8(0)} \op{dist}^2(X,P) \,d\|T\|(X)
		\\&\geq \frac{\beta^2 c^2}{16} \int_{\mathbf{C}_8(0)} \op{dist}^2(X,P) \,d\|T\|(X) , 
\end{align*}
where $C = C(n,m,q) \in (0,\infty)$ is a constant.  This proves \eqref{H uniq cont concl} with $\beta c/4$ in place of $\beta$. 
\end{proof}

\begin{lemma}\label{H error lemma}
For each positive integer $q$ and $\alpha \in (0,1)$ there exists $\delta = \delta(n,m,q) \in (0,1)$ and $\eta_0 = \eta_0(n,m,q) \in (0,1)$ such that the following holds true.  Let $P = \mathbb{R}^n \times \{0\}$ and let $T$ be $n$-dimensional area-minimizing integral current of $\mathbf{C}_{7/4}(0)$ such that \eqref{mono freq mass hyp} and \eqref{mono freq decay hyp} hold true with $Z = 0$ and $\rho_0 = 1$ for some $\eta \in (0,\eta_0]$ and $\sigma_0 \in (0,1)$.  Then 
\begin{equation}\label{H error concl1}
	-\rho^{1-n}\int \frac{1}{r}\op{dist}^2(X,P) \,\|\pi_S - \pi_P\|^2 \,\phi'(r/\rho) \,d|T|(X,S) \leq C \eta^{2\gamma} \rho^{2\alpha\gamma} H_{T,P,0}(\rho) 
\end{equation}
for all $\rho \in [\sigma_0,1]$ and thus 
\begin{equation}\label{H error concl2}
	\left| H'_{T,P,0}(\rho) + 2 \rho^{-n} \int (0,y) \cdot \pi_{P^{\perp}}(\nabla^S r) \,\phi'(r/\rho) \,d|T|(X,S) \right| 
	\leq C \eta^{2\gamma} \rho^{2\alpha\gamma-1} H_{T,P,0}(\rho) 
\end{equation}
for $\mathcal{L}^1$-a.e.~$\rho \in [\sigma_0,1]$, where $\gamma = \gamma(n,m,q) \in (0,1)$ and $C = C(n,m,q) \in (0,\infty)$ are constants.
\end{lemma} 

\begin{proof}
Fix $\rho \in [\sigma_0,1]$.  With $\eta_{0} = \eta_{0}(n, m, q) \in (0, 1)$ to be chosen, let $\eta \in (0, \eta_{0}]$ and suppose that \eqref{mono freq mass hyp} and 
\eqref{mono freq decay hyp} hold.  Notice that by \eqref{mono freq decay hyp} and Lemma~\ref{Allard height lemma}
\begin{equation*}
    \sup_{X \in \op{spt} T \cap \mathbf{C}_{13/8}(0)} \op{dist}(X,P) \leq C \eta 
\end{equation*}
for some constant $C = C(n,m) \in (0,\infty)$.  By the constancy theorem~\cite[Theorem~26.27]{SimonGMT} 
\begin{equation*}
    \pi_{P\#} (T \llcorner \mathbf{C}_{13/8}(0)) = \theta \llbracket B_{13/8}(0) \rrbracket 
\end{equation*}
for some integer $\theta$.  By Lemma~\ref{area excess to height lemma} 
\begin{equation*}
    \mathcal{E}(T,P,\mathbf{C}_{3/2}(0))^2 \leq C \eta 
\end{equation*}
for some constant $C = C(n,m) \in (0,\infty)$.  Thus by the monotonicity formula for area and \eqref{area excess eqn1}
\begin{equation*}
    |\theta| \leq \frac{\|T\|(\mathbf{C}_{\rho}(0))}{\omega_n \rho^n} 
    \leq \frac{\|T\|(\mathbf{C}_{3/2}(0))}{\omega_n (3/2)^n} \leq |\theta| + C(n,m) \,\eta 
\end{equation*}
which by \eqref{mono freq mass hyp} implies that $|\theta| = q$.  Up to reversing the orientation of $T$, we may assume that $\theta = q$ so that 
\begin{equation}\label{H error eqn1} 
    \pi_{P\#} (T \llcorner \mathbf{C}_{13/8}(0)) = q \llbracket B_{13/8}(0) \rrbracket .
\end{equation}

For each $\xi \in B_{\rho}(0) \setminus \overline{B_{\rho/2}(0)}$ let 
\begin{equation*}
	\overline{\sigma}(\xi) = \sup \{0\} \cup \left\{ \sigma \in (0,\rho/2] : 
		\int_{\mathbf{C}_{\sigma}(\xi)} \op{dist}^2(X,P) \,d\|T\|(X) \geq \eta^2 \omega_n (64\sigma)^{n+2+2\alpha} \right\} . 
\end{equation*}
Using \eqref{mono freq decay hyp}, for each $\xi \in B_{\rho}(0) \setminus \overline{B_{\rho/2}(0)}$ 
\begin{align*}
	\int_{\mathbf{C}_{\sigma}(\xi)} \op{dist}^2(X,P) \,d\|T\|(X)
	&\leq \int_{\mathbf{C}_{3\rho/2}(0)} \op{dist}^2(X,P) \,d\|T\|(X) 
	\\&\leq \eta^2 \omega_n (3\rho/2)^{n+2+2\alpha}
	< \eta^2 \omega_n(64\sigma)^{n+2+2\alpha} 
\end{align*}
for all $\sigma \in [\rho/32,\rho/2]$ and thus  $\overline{\sigma}(\xi) \leq \rho/32$.  Let 
\begin{equation*}
	\Xi = \{ \xi \in B_{\rho}(0) \setminus \overline{B_{\rho/2}(0)} : \overline{\sigma}(\xi) > 0 \} .
\end{equation*}
Notice that by Lemma~\ref{Allard height lemma}, 
\begin{equation}\label{H error eqn3}
	\sup_{X \in \op{spt} T \cap \mathbf{C}_{16\sigma}(\xi)} \op{dist}(X,P) \leq C \eta \sigma^{1+\alpha}
\end{equation}
for all $\sigma \in (\overline{\sigma},\rho/16]$ and some constant $C = C(n,m) \in (0,\infty)$.  In particular, if $\xi \in (B_{\rho}(0) \setminus \overline{B_{\rho/2}(0)}) \setminus \Xi$, we can let $\sigma \rightarrow 0^+$ to obtain 
\begin{equation}\label{H error eqn4}
	\op{spt} T \cap (\{\xi\} \times \mathbb{R}^m) = \{(\xi,0)\} . 
\end{equation}

By the Vitali covering lemma, there is a countable set of points $\xi_i \in \Xi$ such that $\{B_{2\sigma_i/5}(\xi_i)\}$ is a pairwise disjoint collection of balls and $\{B_{2\sigma_i}(\xi_i)\}$ covers of $\Xi$, where $\sigma_i = \overline{\sigma}(\xi_i)$.  By Lemma~\ref{area excess to height lemma}, 
\begin{equation*}
	\sigma_{i}^{-n} \int_{\mathbf{C}_{8\sigma_i}(\xi_i)} |\vec T - \vec P|^2 \,d\|T\|(X) 
	\leq C \sigma^{-n-2} \int_{\mathbf{C}_{16\sigma_i}(\xi_i)} \op{dist}^2(X,P) \,d\|T\|(X) 
	\leq C \,\eta^2 \sigma_i^{2\alpha} , 
\end{equation*}
where $C = C(n,m,q) \in (0,\infty)$.  Hence if $\eta_{0} = \eta_{0}(n, m, q) \in (0, 1)$ is sufficiently small, by~\cite[Corollaries~3.29 \& 3.30]{Almgren} or \cite[Theorem~2.4]{DeLSpa1}, for each $i$ there exists a Lipschitz function $u_i : B_{2\sigma_i}(\xi_i) \rightarrow \mathcal{A}_q(\mathbb{R}^m)$ and a set $K_i \subset B_{2\sigma_i}(\xi_i)$ such that 
\begin{gather}
	\label{H error eqn5} \op{Lip} u_i \leq C \eta^{\gamma} \sigma_i^{\gamma\alpha}, \\
	T \llcorner K_i \times \mathbb{R}^m = (\op{graph} u_i) \llcorner K_i \times \mathbb{R}^m, \nonumber \\
	|B_{2\sigma_i}(\xi_i) \setminus K_i| + \|T\|((B_{2\sigma_i}(\xi_i) \setminus K_i) \times \mathbb{R}^m) 
		\leq C \eta^{2+2\gamma} \sigma_i^{n+2\alpha(1+\gamma)} , \nonumber 
\end{gather}
where $\gamma = \gamma(n,m,q) \in (0,1)$ and $C = C(n,m,q) \in (0,\infty)$ are constants.  By \eqref{lipschitz approx eqn8}, \eqref{lipschitz approx eqn9}, \eqref{lipschitz approx eqn11}, \eqref{H error eqn3}, and \eqref{H error eqn5}, 
\begin{align}\label{H error eqn7}
	&\int_{G_n(\mathbf{C}_{2\sigma_i}(\xi_i))} \op{dist}^2(X,P) \,\|\pi_S - \pi_P\|^2 \,d|T|(X,S) 
	\\ \leq\,& C \int_{K_i \cap B_{2\sigma_i}(\xi_i)} |u_i|^2 \,|Du_i|^2 + C \eta^{4+2\gamma} \sigma_i^{n+2+2\alpha(2+\gamma)} \nonumber 
	\\ \leq\,& C \eta^{2\gamma} \sigma_i^{2\alpha\gamma} \int_{K_i \cap B_{2\sigma_i}(\xi_i)} |u_i|^2 
		+ C \eta^{4+2\gamma} \sigma_i^{n+2+2\alpha(2+\gamma)} \nonumber
	\\ \leq\,& C \eta^{2\gamma} \sigma_i^{2\alpha\gamma} \int_{\mathbf{C}_{2\sigma_i}(\xi_i)} \op{dist}^2(X,P) \,d\|T\|(X) 
		+ C \eta^{4+2\gamma} \sigma_i^{n+2+2\alpha(2+\gamma)} \nonumber
	\\ \leq\,& C \eta^{2+2\gamma} \sigma_i^{n+2+2\alpha(1+\gamma)} , \nonumber 
\end{align}
where $C = C(n,m,q) \in (0,\infty)$ are constants.  Notice that since $\xi_i \in B_{\rho}(0) \setminus \overline{B_{\rho/2}(0)}$ there exists an open ball $B_{\sigma_i/5}(\zeta_i) \subset B_{2\sigma_i/5}(\xi_i) \cap B_{\rho}(0) \setminus \overline{B_{\rho/2}(0)}$.  Also, by the definition of $\sigma_i$ and since $\alpha < 1$, 
\begin{align}\label{H error eqn8}
	\int_{\mathbf{C}_{\sigma_i}(\xi_i)} \op{dist}^2(X,P) \,d\|T\|(X) &\geq \eta^2 \omega_n (64\sigma_i)^{n+2+2\alpha} 
		\\&\geq 16^{-n-4} \int_{\mathbf{C}_{16\sigma_i}(\xi_i)} \op{dist}^2(X,P) \,d\|T\|(X) . \nonumber 
\end{align}
Thus by Lemma~\ref{H uniq cont lemma} there exists a constant $\beta = \beta(n,m,q) > 0$ such that 
\begin{equation*}
	\int_{\mathbf{C}_{\sigma_i/5}(\zeta_i)} \op{dist}^2(X,P) \,d\|T\|(X) \geq \beta^2 \int_{\mathbf{C}_{16\sigma_i}(\xi_i)} \op{dist}^2(X,P) \,d\|T\|(X) . 
\end{equation*}
Using the fact that $\mathbf{C}_{\sigma_i/5}(\zeta_i) \subset \mathbf{C}_{2\sigma_i/5}(\xi_i) \cap \mathbf{C}_{\rho}(0) \setminus \overline{\mathbf{C}_{\rho/2}(0)}$, $\mathbf{C}_{\sigma_i}(\xi_i) \subset \mathbf{C}_{16\sigma_i}(\xi_i)$, and \eqref{H error eqn8}, 
\begin{equation}\label{H error eqn9}
	-\int_{\mathbf{C}_{2\sigma_i/5}(\xi_i)} \op{dist}^2(X,P) \,\frac{1}{r} \,\phi'(r/\rho) \,d\|T\|(X)
	\geq \frac{2}{\rho}\beta^2 \eta^2 \omega_n (64\sigma_i)^{n+2+2\alpha} . 
\end{equation}
Since $\{B_{2\sigma_i}(\xi_i)\}$ covers $\Xi$ and $\{B_{2\sigma_i/5}(\xi_i)\}$ are pairwise disjoint, and by  \eqref{H error eqn4}, ${\rm dist}(X, P) = 0$ 
on $((B_{\rho}(0) \setminus \overline{B_{\rho/2}(0)}) \setminus \Xi) \times {\mathbb R}^{m}$, it follows, using also \eqref{H error eqn7}, and \eqref{H error eqn9}, that 
\begin{align}\label{H error eqn10}
	&-\int \op{dist}^2(X,P) \,\|\pi_S - \pi_P\|^2 \,\frac{1}{r} \,\phi'(r/\rho) \,d|T|(X,S) 
	\\ \leq\,& \frac{4}{\rho} \sum_i \int_{G_n(\mathbf{C}_{2\sigma_i}(\xi_i))} \op{dist}^2(X,P) \,\|\pi_S - \pi_P\|^2 \,d|T|(X,S) \nonumber
	\\ \leq\,& \frac{C}{\rho} \sum_i \eta^{2+2\gamma} \sigma_i^{n+2+2\alpha(1+\gamma)} \nonumber
	\\ \leq\,& -C \eta^{2\gamma} \rho^{2\alpha\gamma} \sum_i \int_{\mathbf{C}_{2\sigma_i/5}(\xi_i)} \op{dist}^2(X,P) \,\frac{1}{r} \,\phi'(r/\rho) \,d\|T\|(X) 
		\nonumber
	\\ \leq\,& -C \eta^{2\gamma} \rho^{2\alpha\gamma} \int \op{dist}^2(X,P) \,\frac{1}{r} \,\phi'(r/\rho) \,d\|T\|(X) , \nonumber
\end{align}
where $C = C(n,m,q) \in (0,\infty)$.  Notice that for each $n$-dimensional plane $S \subset \mathbb{R}^{n+m}$, 
\begin{equation}\label{H error eqn11}
	|\nabla^{S^\perp} r|^2 = |(\pi_S - \pi_P) \,\nabla r|^2 \leq \|\pi_S - \pi_P\|^2 . 
\end{equation}
Hence $|\nabla^S r|^2 = 1 - |\nabla^{S^\perp} r|^2 \geq 1 - \|\pi_S - \pi_P\|^2$, which  by \eqref{H error eqn10} gives us
\begin{equation}\label{H error eqn12}
	-\rho^{1-n} \int \op{dist}^2(X,P) \,\frac{1}{r} \,\phi'(r/\rho) \,d|T|(X,S) \leq 2 H_{T,P,0}(\rho) 
\end{equation}
for all $\rho \in [\sigma_0,1]$ provided $\eta_{0} = \eta_{0}(n, m, q) \in (0, 1)$ is sufficiently small.  It follows from \eqref{H error eqn10} and \eqref{H error eqn12}
that \eqref{H error concl1} holds true.  Using \eqref{H error eqn11} and \eqref{H error concl1} to bound the second integral in \eqref{H deriv}, we obtain \eqref{H error concl2}. 
\end{proof}

\subsection{Bounding error terms in derivative of $D_{T,P,0}$}  Using similar arguments, in Lemma~\ref{D error lemma} we will estimate certain error terms in \eqref{D deriv}. 

\begin{lemma}\label{D uniq cont lemma}
For every  $c > 0$ there exists $\varepsilon = \varepsilon(n, m, q, c) \in (0, 1)$, $\delta = \delta(n,m,q,c) \in (0,1)$ and 
$\beta = \beta(n,m,q,c) \in (0, 1)$ such that if $T$ is an $n$-dimensional locally area-minimizing rectifiable current of $\mathbf{C}_4(0)$, $P = {\mathbb R}^{n} \times \{0\}$ and if  
\begin{gather} 
	\label{D uniq cont hyp1} (\partial T) \llcorner \mathbf{C}_4(0) = 0, \quad 
    \sup_{X \in \op{spt} T} \op{dist}(X,P) < \infty, \quad 
	\pi_{P\#} T = q \llbracket B_8(0) \rrbracket , \\
	\label{D uniq cont hyp2} \int_{\mathbf{C}_4(0)} |\vec T - \vec P|^2 \,d\|T\| < \varepsilon^2 \;\; and\\
	\label{D uniq cont hyp3} \int_{\mathbf{C}_{1/2}(0)} |\vec T - \vec P|^2 \,d\|T\| \geq c^2 \int_{\mathbf{C}_4(0)} |\vec T - \vec P|^2 \,d\|T\|, 
\end{gather}
then for every $\xi \in B_{1/5}(0)$ 
\begin{equation}\label{D uniq cont concl}
	\int_{G_{n}(\mathbf{C}_{1/20}(\xi))} \|\pi_S - \pi_P\|^2 \,d|T|(X,S) \geq \beta^2 \int_{\mathbf{C}_4(0)} |\vec T - \vec P|^2 \,d\|T\| . 
\end{equation}
\end{lemma}

\begin{proof}
By~\cite[Corollaries~3.29 \& 3.30]{Almgren} (or \cite[Theorem~2.4]{DeLSpa1}), if $\epsilon = \epsilon(n, m, q, \eta) \in (0, 1)$ is sufficiently small, there is a Lipschitz function $u : B_1(0) \rightarrow \mathcal{A}_q(\mathbb{R}^m)$ and a set $K \subset B_1(0)$ such that 
\begin{gather}
	\label{D uniq cont eqn1} \quad\quad \sup_{B_1(0)} |\nabla u| \leq C \varepsilon^{\gamma}, \\
	T \llcorner K \times \mathbb{R}^m = (\op{graph} u) \llcorner K \times \mathbb{R}^m, \nonumber \\
	|B_1(0) \setminus K| + \|T\|((B_1(0) \setminus K) \times \mathbb{R}^m) 
		\leq C \varepsilon^{\gamma} \int_{\mathbf{C}_4(0)} |\vec T - \vec P|^2 \,d\|T\|(X) , \nonumber 
\end{gather}
where $\gamma = \gamma(n,m,q) \in (0,1)$ and $C = C(n,m,q) \in (0,\infty)$ are constants.  Let $\eta \in (0,1)$ be a constant to be determined depending only on 
$n, m, q, c$.  By \cite[Theorem~2.23]{Almgren} (or ~\cite[Theorem~2.6]{DeLSpa1}) there exists a Dirichlet energy minimizing $q$-valued function $w \in W^{1,2}(B_1(0),\mathcal{A}_q(\mathbb{R}^m))$ such that 
\begin{equation}\label{D uniq cont eqn2}
	\int_{B_1(0)} (|Du| - |Dw|)^2 \leq \eta \int_{\mathbf{C}_4(0)} |\vec T - \vec P|^2 \,d\|T\|. 
\end{equation}
Provided $\varepsilon$ and $\eta$ are sufficiently small, by \eqref{lipschitz approx eqn8}, \eqref{lipschitz approx eqn9}, \eqref{lipschitz approx eqn10}, \eqref{D uniq cont eqn1}, \eqref{D uniq cont eqn2}, and \eqref{D uniq cont hyp3}
\begin{align}\label{D uniq cont eqn3}
	\int_{B_1(0)} |Dw|^2 &\leq 2 \int_{B_1(0)} |Du|^2 + 2\eta \int_{\mathbf{C}_4(0)} |\vec T - \vec P|^2 \,d\|T\|(X)  
		\\&\leq 2 \int_{\mathbf{C}_1(0)} |\vec T - \vec P|^2 \,d\|T\|(X) 
			+ (2\eta + C \varepsilon^{\gamma}) \int_{\mathbf{C}_4(0)} |\vec T - \vec P|^2 \,d\|T\|(X) \nonumber
		\\&\leq 4 \int_{\mathbf{C}_4(0)} |\vec T - \vec P|^2 \,d\|T\|(X) \nonumber
\end{align}
and  
\begin{align}\label{D uniq cont eqn4}
	\int_{B_{1/2}(0)} |Dw|^2 &\geq \frac{1}{2} \int_{B_{1/2}(0)} |Du|^2 - \eta \int_{\mathbf{C}_4(0)} |\vec T - \vec P|^2 \,d\|T\|(X) 
		\\&\geq \frac{1}{2} \int_{\mathbf{C}_{1/2}(0)} |\vec T - \vec P|^2 \,d\|T\|(X) 
			- (\eta + C \varepsilon^{\gamma}) \int_{\mathbf{C}_4(0)} |\vec T - \vec P|^2 \,d\|T\|(X) \nonumber
		\\&\geq \left( \frac{c^2}{2} - \eta - C \varepsilon^{\gamma} \right) \int_{\mathbf{C}_4(0)} |\vec T - \vec P|^2 \,d\|T\|(X) \nonumber
		\\&\geq \frac{c^2}{4} \int_{\mathbf{C}_4(0)} |\vec T - \vec P|^2 \,d\|T\|(X) ,\nonumber
\end{align}
where $C = C(n,m,q) \in (0,\infty)$ is a constant.  Combining \eqref{D uniq cont eqn3} and \eqref{D uniq cont eqn4}, 
\begin{align}\label{D uniq cont eqn5}
	\int_{B_{1/2}(0)} |Dw|^2 \geq \frac{c^2}{16} \int_{B_1(0)} |Dw|^2 . 
\end{align}
Hence arguing as in the proof of Lemma~\ref{H uniq cont lemma0}, there exists $\beta = \beta(n,m,q,c) > 0$ such that \eqref{D uniq cont eqn5} implies that for each $\xi \in B_{1/5}(0)$ 
\begin{equation}\label{D uniq cont eqn6}
	\int_{B_{1/20}(\xi)} |Dw|^2 \geq \beta^2 \int_{B_1(0)} |Dw|^2 . 
\end{equation}
By \eqref{lipschitz approx eqn8}, \eqref{lipschitz approx eqn9}, \eqref{lipschitz approx eqn11}, \eqref{D uniq cont eqn1}, \eqref{D uniq cont eqn2}, \eqref{D uniq cont eqn6}, and \eqref{D uniq cont eqn4}, 
\begin{align*}
	\int_{\mathbf{C}_{1/20}(\xi)} \|\pi_S - \pi_P\|^2 \,d|T|(X,S) 
		&\geq \int_{B_{1/20}(\xi)} |Du|^2 - C \varepsilon^{\gamma} \int_{\mathbf{C}_4(0)} |\vec T - \vec P|^2 \,d\|T\|(X) 
		\\&\geq \frac{1}{2} \int_{B_{1/20}(\xi)} |Dw|^2 - (\eta + C \varepsilon^{\gamma}) \int_{\mathbf{C}_4(0)} |\vec T - \vec P|^2 \,d\|T\|(X)
		\\&\geq \frac{\beta^2}{2} \int_{B_1(0)} |Dw|^2 - (\eta + C \varepsilon^{\gamma}) \int_{\mathbf{C}_4(0)} |\vec T - \vec P|^2 \,d\|T\|(X)
		\\&\geq \left( \frac{\beta^2 c^2}{8} - \eta - C \varepsilon^{\gamma} \right) \int_{\mathbf{C}_4(0)} |\vec T - \vec P|^2 \,d\|T\|(X)
		\\&\geq \frac{\beta^2 c^2}{16} \int_{\mathbf{C}_4(0)} |\vec T - \vec P|^2 \,d\|T\|(X) , 
\end{align*}
where $C \in (0,\infty)$ is a constant.  This proves \eqref{D uniq cont concl} with $\beta c/4$ in place of $\beta$. 
\end{proof}

\begin{lemma}\label{D error lemma} 
For each positive integer $q$ and $\alpha \in (0,1)$ there exists $\delta = \delta(n,m,q) \in (0,1)$ and $\eta_0 = \eta_0(n,m,q) \in (0,1)$ such that the following holds true.  Let $P = \mathbb{R}^n \times \{0\}$ and let $T$ be $n$-dimensional area-minimizing integral current of $\mathbf{C}_{7/4}(0)$ such that  \eqref{mono freq mass hyp} and \eqref{mono freq decay hyp} hold true with $Z = 0$ and $\rho_0 = 1$ for some $\eta \in (0,\eta_0]$ and $\sigma_0 \in (0,1)$.  Then 
\begin{align}
	\label{D error concl1} &\rho^{2-n} \int_{\mathbf{C}_{\rho}(0)} |\vec T - \vec P|^4 \,d\|T\|(X) 
		\leq C D(\rho)^{\gamma_1} \,((n-1) \,D_{T,P,0}(\rho) + \rho \,D'_{T,P,0}(\rho)) , 
\end{align}
for all $\rho \in [\sigma_0,1]$, 
where $\gamma_1 = \gamma_1(n,m,q) \in (0,1)$ and $C = C(n,m,q) \in (0,\infty)$ are constants.  Thus 
\begin{align}\label{D error concl3} 
	&\left| D'_{T,P,0}(\rho) + 2\rho^{-n} \int_{G_n(\mathbf{B}_2(0)) \cap \{\nabla^S r \neq 0\}} 
		\frac{|\pi_{P^{\perp}}(\nabla^S r)|^2}{|\nabla^S r|^2} \,r \,\phi'(r/\rho) \,d|T|(X,S) \right| 
	\\ \leq\,& C \rho^{-1} D(\rho)^{\gamma_1} \,((n-1) D_{T,P,0}(\rho) + \rho \,D'_{T,P,0}(\rho)) \nonumber 
\end{align}
for $\mathcal{L}^1$-a.e.~$\rho \in [\sigma_0,1]$, where $C = C(n,m,q) \in (0,\infty)$ is a constant.
\end{lemma} 

\begin{proof}
Fix $\rho \in [\sigma_0,1]$.  With $\eta_{0} = \eta_{0}(n, m, q) \in (0, 1)$ to be chosen, let $\eta \in (0, \eta_{0}]$ and suppose that 
\eqref{mono freq mass hyp} and \eqref{mono freq decay hyp} hold. Recall that by \eqref{mono freq mass hyp} and \eqref{mono freq decay hyp}, $T$ has finite height relative to $P$ in $\mathbf{C}_{13/8}(0)$ and thus we may assume that \eqref{H error eqn1} holds true.  By Lemma~\ref{area excess to height lemma}, 
\begin{equation}\label{D error eqn1}
	\int_{\mathbf{C}_{3\rho/2}(0)} |\vec T - \vec P|^2 \,d\|T\|(X) \leq \frac{C_0}{\rho^2} \int_{\mathbf{C}_{7\rho/4}(0)} \op{dist}^2(X,P) \,d\|T\|(X) 
\end{equation}
for some constant $C_0 = C_0(n,m,q) \in (0,\infty)$.  For each $\xi \in B_{\rho}(0)$ let 
\begin{equation*}
	\overline{\sigma}(\xi) = \sup \{0\} \cup \left\{ \sigma \in (0,\rho/2] : 
		\int_{\mathbf{C}_{\sigma}(\xi)} |\vec T - \vec P|^2 \,d\|T\| \geq 4C_0 \eta^2 (64\sigma)^{n+2\alpha} \right\} , 
\end{equation*}
where $C_0$ is as in \eqref{D error eqn1}.  By \eqref{D error eqn1} and \eqref{mono freq decay hyp}, for each $\xi \in B_{\rho}(0)$ 
\begin{align*}
	\int_{\mathbf{C}_{\sigma}(\xi)} |\vec T - \vec P|^2 \,d\|T\| 
	&\leq \int_{\mathbf{C}_{3\rho/2}(0)} |\vec T - \vec P|^2 \,d\|T\| 
	\leq \frac{C_0}{\rho^2} \int_{\mathbf{C}_{7\rho/4}(0)} \op{dist}^2(X,P) \,d\|T\|(X) 
	\\&\leq 4C_0 \eta^2 (7\rho/4)^{n+2\alpha}
	< 4C_0 \eta^2 (64\sigma)^{n+2\alpha} 
\end{align*}
for all $\sigma \in [\rho/32,\rho/2]$ and thus $\overline{\sigma}(\xi) \leq \rho/32$.  Set 
\begin{equation*}
	\Xi = \{ \xi \in B_{\rho}(0) : \overline{\sigma}(\xi) > 0 \}.
\end{equation*}
Since for every $\xi \in B_{\rho}(0) \setminus \Xi$ we have that
\begin{equation*}
	\int_{\mathbf{C}_{\sigma}(\xi)} |\vec T - \vec P|^2 \,d\|T\| < 4C_0 \eta^2 (64\sigma)^{n+2\alpha}
\end{equation*}
for all $\sigma \in (0,\rho/2]$, we must have that $\vec T(X) = \vec P$ for $\|T\|$-a.e.~$X \in (B_{\rho}(0) \setminus \Xi) \times \mathbb{R}^m$,  and hence 
the approximate tangent plane to $\op{spt} T$ is $P$ at $\|T\|$-a.e.~$X \in (B_{\rho}(0) \setminus \Xi) \times \mathbb{R}^m$.  In other words, 
\begin{equation}\label{D error eqn2}
	S = P \text{ for $|T|$-a.e.~$(X,S) \in G_n((B_{\rho}(0) \setminus \Xi) \times \mathbb{R}^m)$.}
\end{equation}

By the Vitali covering lemma, there is a countable set of points $\xi_i \in \Xi$ such that $\{B_{2\sigma_i/5}(\xi_i)\}$ is a pairwise disjoint collection of balls and $\{B_{2\sigma_i}(\xi_i)\}$ covers of $\Xi$, where $\sigma_i = \overline{\sigma}(\xi_i)$.  By \eqref{D error eqn1} and \cite[Corollaries 3.29 \& 3.30]{Almgren} (or~\cite[Theorem~2.4]{DeLSpa1}), for each $i$ there exists a Lipschitz function $u_i : B_{2\sigma_i}(\xi_i) \rightarrow \mathcal{A}_q(\mathbb{R}^m)$ and a set $K_i \subset B_{2\sigma_i}(\xi_i)$ such that \eqref{H error eqn5} holds true.  Hence by \eqref{lipschitz approx eqn8}, \eqref{lipschitz approx eqn9}, \eqref{lipschitz approx eqn10}, \eqref{lipschitz approx eqn11}, and \eqref{H error eqn5}, 
\begin{align}\label{D error eqn3}
	&\int_{\mathbf{C}_{2\sigma_i}(\xi_i)} |\vec T - \vec P|^4 \,d\|T\|(X) 
	\\ \leq\,& C \int_{K_i \cap B_{2\sigma_i}(\xi_i)} |Du_i|^4 + C \eta^{2+2\gamma} \sigma_i^{n+2\alpha(1+\gamma)} \nonumber
	\\ \leq\,& C \eta^{2\gamma} \sigma_i^{2\alpha\gamma} \int_{K_i \cap B_{2\sigma_i}(\xi_i)} |Du_i|^2 
		+ C \eta^{2+2\gamma} \sigma_i^{n+2\alpha(1+\gamma)} \nonumber 
	\\ \leq\,& C \eta^{2\gamma} \sigma_i^{2\alpha\gamma} \int_{\mathbf{C}_{2\sigma_i}(\xi_i)} |\vec T - \vec P|^2 \,d\|T\|(X) 
		+ C \eta^{2+2\gamma} \sigma_i^{n+2\alpha(1+\gamma)} \nonumber 
	\\ \leq\,& C \eta^{2+2\gamma} \sigma_i^{n+2\alpha(1+\gamma)} , \nonumber
\end{align}
$C = C(n,m,q) \in (0,\infty)$.  Since $\xi_i \in B_{\rho}(0)$, there exists a ball $B_{\sigma_i/5}(\zeta_i) \subset B_{2\sigma_i/5}(\xi_i) \cap B_{\rho}(0)$.  By the definition of $\sigma_i$ and $\alpha < 1$, 
\begin{equation}\label{D error eqn4}
	\int_{\mathbf{C}_{\sigma_i}(\xi_i)} |\vec T - \vec P|^2 \,d\|T\| \geq 4C_0 \eta^2 (64\sigma_i)^{n+2\alpha} 
		\geq 8^{-n-2} \int_{\mathbf{C}_{8\sigma_i}(\xi_i)} |\vec T - \vec P|^2 \,d\|T\| . 
\end{equation}
By Lemma~\ref{D uniq cont lemma} there exists a constant $\beta = \beta(n,m,q) > 0$ such that 
\begin{equation}\label{D error eqn5}
	\int_{G_n(\mathbf{C}_{\sigma_i/10}(\zeta_i))} \|\pi_S - \pi_P\|^2 \,d|T|(X,S) 
	\geq \beta \int_{\mathbf{C}_{8\sigma_i}(\xi_i)} |\vec T - \vec P|^2 \,d\|T\| . 
\end{equation}
Using the fact that $\mathbf{C}_{\sigma_i/5}(\zeta_i) \subset \mathbf{C}_{2\sigma_i/5}(\xi_i) \cap \mathbf{C}_{\rho}(0)$, $\overline{\mathbf{C}_{\sigma_i}(\xi_i)} \subset \mathbf{C}_{8\sigma_i}(\xi_i)$, \eqref{D error eqn4}, and \eqref{D error eqn5}, 
\begin{align}\label{D error eqn6}
	\int_{G_n(\mathbf{C}_{2\sigma_i/5}(\xi_i) \cap \mathbf{C}_{\rho}(0))} \|\pi_S - \pi_P\|^2 \,d|T|(X,S) \geq 4\beta C_0 \eta^2 (64\sigma_i)^{n+2\alpha} . 
\end{align}
Since $B_{\sigma_i/5}(\zeta_i) \subset B_{\rho}(0)$, we have that $r \leq \rho-\sigma_i/5$ on $\mathbf{C}_{\sigma_i/10}(\zeta_i)$ and thus $\phi(r/\rho) \geq \sigma_i/(5\rho)$ on $\mathbf{C}_{\sigma_i/10}(\zeta_i)$.  Hence by \eqref{D error eqn4} and \eqref{D error eqn5}
\begin{align}\label{D error eqn7}
	C_0 \beta \eta^2 (64\sigma_i)^{n+1+2\alpha} 
	\leq\,& 16 \sigma_i \int_{G_n(\mathbf{C}_{\sigma_i/10}(\zeta_i))} \|\pi_S - \pi_P\|^2 \,d|T|(X,S) 
	\\ \leq\,& 80\rho \int_{G_n(\mathbf{C}_{\rho}(0))} \|\pi_S - \pi_P\|^2 \,\phi(r/\rho) \,d|T|(X,S) 
	\leq 80 D_{T,P,0}(\rho) \nonumber 
\end{align}
Recalling that $\{B_{2\sigma_i}(\xi_i)\}$ covers $\Xi$ and $\{B_{2\sigma_i/5}(\xi_i)\}$ are pairwise disjoint, and using \eqref{D error eqn2}, \eqref{D error eqn3}, \eqref{D error eqn6}, and \eqref{D error eqn7}, 
\begin{align}\label{D error eqn8}
	&\int_{\mathbf{C}_{\rho}(0)} |\vec T - \vec P|^4 \,d\|T\|(X) 
	\\ \leq\,& \sum_i \int_{\mathbf{C}_{2\sigma_i}(\xi_i)} |\vec T - \vec P|^4 \,d\|T\|(X) \nonumber 
	\\ \leq\,& C \sum_i \eta^{2+2\gamma} \sigma_i^{n+2\alpha(1+\gamma)} \nonumber 
	\\ \leq\,& C \sum_i \eta^{2\gamma} \sigma_i^{2\alpha\gamma} \int_{G_n(\mathbf{C}_{2\sigma_i/5}(\xi_i) \cap \mathbf{C}_{\rho}(0))} \|\pi_S - \pi_P\|^2 
		\,d|T|(X,S) \nonumber
	\\ \leq\,& C D(\rho)^{\gamma_1} \sum_i \int_{G_n(\mathbf{C}_{2\sigma_i/5}(\xi_i) \cap \mathbf{C}_{\rho}(0))} \|\pi_S - \pi_P\|^2 
		\,d|T|(X,S) \nonumber
	\\ \leq\,& C D(\rho)^{\gamma_1} \int_{G_n(\mathbf{C}_{\rho}(0))} \|\pi_S - \pi_P\|^2 \,d|T|(X,S) , \nonumber
\end{align}
where $\gamma_1 = \frac{2\alpha\gamma}{n+1+2\alpha}$ and $C = C(n,m,q) \in (0,\infty)$ are constants.  Notice that by the definition of $\phi$ and by \eqref{D deriv eqn3} (noting also that 
$\phi(r/\rho) - \frac{r}{\rho} \,\phi'(r/\rho) = 1$ on ${\mathbf C}_{\rho/2}$ and 
$\phi(r/\rho) - \frac{r}{\rho} \,\phi'(r/\rho) = 2(1 - \frac{r}{\rho}) + \frac{2r}{\rho}  = 2$ on ${\mathbf C}_{\rho} \setminus {\mathbf C}_{\rho/2}$), we have 
\begin{align}\label{D error eqn9}
	\frac{1}{2} \int_{G_n(\mathbf{C}_{\rho}(0))} \|\pi_S - \pi_P\|^2 \,d|T|(X,S) 
	\leq\,& \frac{1}{2} \int \|\pi_S - \pi_P\|^2 \left( \phi(r/\rho) - \frac{r}{\rho} \,\phi'(r/\rho) \right) d|T|(X,S) 
	\\ =\,& (n-1) \,\rho^{n-2} D_{T,P,0}(\rho) + \rho^{n-1} D'_{T,P,0}(\rho) \nonumber 
\end{align}
which together with \eqref{D error eqn8} gives us \eqref{D error concl1}.  (Note that \eqref{D deriv eqn3} holds true when $\phi$ is given by \eqref{freq phi defn} by the approximation argument in the proof of Lemma~\ref{freq identities lemma}.)

We claim that 
\begin{equation}
	\label{D error eqn10} \left| |\nabla^{S^\perp} r|^2 - \frac{|\pi_{P^{\perp}} \nabla^S r|^2}{|\nabla^S r|^2} \right| \leq C \|\pi_S - \pi_P\|^4 
		\text{ if } \nabla^S r \neq 0
\end{equation}
for each $(X,S) \in G_n(\mathbf{C}_{\rho}(0))$ and that 
\begin{equation}
	\label{D error eqn11} \left| 1 - \frac{1}{4} \,\|\pi_S - \pi_P\|^2 - J\pi_P \right| \leq C |\vec T - \vec P|^4  
\end{equation}
for $|T|$-a.e.~$(X,S) \in G_n(\mathbf{C}_{\rho}(0))$.  
Note that by \eqref{H error eqn11}, if $\nabla^S r = 0$ then in place of \eqref{D error eqn10} we have $|\nabla^{S^\perp} r|^2 = 1 = |\nabla^{S^\perp} r|^4 \leq \|\pi_S - \pi_P\|^4$.  To see \eqref{D error eqn10}, if $|\nabla^S r|^2 \leq 1/2$, then again using \eqref{H error eqn11} we have that $1/2 \leq |\nabla^{S^\perp} r|^2 \leq \|\pi_S - \pi_P\|^2$ and thus 
\begin{equation*}
	\left| |\nabla^{S^\perp} r|^2 - \frac{|\pi_{P^{\perp}} \nabla^S r|^2}{|\nabla^S r|^2} \right| 
		\leq |\nabla^{S^\perp} r|^2 + \frac{|\pi_{P^{\perp}} \nabla^S r|^2}{|\nabla^S r|^2}
		\leq |\nabla^{S^\perp} r|^2 + 1 \leq 2 \leq 8 \|\pi_S - \pi_P\|^4 . 
\end{equation*}
If instead $|\nabla^S r|^2 \geq 1/2$, then since $\nabla r(X) \in P,$ 
\begin{align*}
	|\nabla^{S^\perp} r|^2 - |\pi_{P^{\perp}} \nabla^S r|^2 
	&= |\nabla^{S^\perp} r|^2 - |\pi_{P^{\perp}} \nabla^{S^\perp} r|^2
	= |\pi_P \nabla^{S^\perp} r|^2
	\\&= |\pi_P (\pi_P - \pi_S) \nabla r|^2
	= |(\pi_P - \pi_S)^2 \nabla r|^2, \;\; \mbox{whence} \\
	\left| \frac{|\pi_{P^{\perp}} \nabla^S r|^2}{|\nabla^S r|^2} -  |\nabla^{S^\perp} r|^2 \right|
	&= \frac{\big| |\pi_{P^{\perp}} \nabla^S r|^2 - |\nabla^{S^\perp} r|^2 + |\nabla^{S^\perp} r|^4 \big|}{|\nabla^S r|^2}  
	\\&\leq \frac{|(\pi_P - \pi_S)^2 \nabla r|^2 + |(\pi_P - \pi_S) \nabla r|^4}{|\nabla^S r|^2} \leq 4 \|\pi_S - \pi_P\|^4 . 
\end{align*}
Therefore, \eqref{D error eqn10} holds true.  

To see \eqref{D error eqn11}, without loss of generality suppose that $|\vec T - \vec P| < 1/4$.  Then $|\pi_{P\#} \vec T - \vec P| < 1/4$ and thus 
\begin{equation*}
	\left| \frac{\pi_{P\#} \vec T}{|\pi_{P\#} \vec T|} - \vec P \right| \leq \frac{2 |\pi_{P\#} \vec T - \vec P|}{|\pi_{P\#} \vec T|} \leq 4 |\pi_{P\#} \vec T - \vec P| < 1.
\end{equation*}
Hence (since $\frac{\pi_{P\#} \vec T}{|\pi_{P\#} \vec T|} = \pm \vec{P}$) we must have that $\pi_{P\#} \vec T = |\pi_{P\#} \vec T| \,\vec P$; that is, $J\pi_P = |\pi_{P\#} \vec T| > 0$.  Thus $J\pi_P = \sqrt{\det [\pi_P \tau_i \cdot \pi_P \tau_j]}$ as the square root of the determinant of the $n \times n$ matrix with $(i,j)$-entry $\pi_P \tau_i \cdot \pi_P \tau_j$, where $\tau_1,\tau_2,\ldots,\tau_n$ is an orthonormal basis for the approximate tangent plane $S_X$ of $\op{spt} T$ at $X$.  Since $\pi_P \tau_i \cdot \pi_P \tau_j = \delta_{ij} - \pi_P^{\perp} \tau_i \cdot \pi_P^{\perp} \tau_j$ for all $1 \leq i,j \leq n$, 
\begin{equation*}
	J\pi_P = 1 - \frac{1}{2} \sum_{i=1}^n |\pi_P^{\perp} \tau_i|^2 + \mathcal{R} 
		= 1 - \frac{1}{2} \sum_{i=1}^n \sum_{j=1}^m (\tau_i \cdot e_{n+j})^2 + \mathcal{R} 
		= 1 - \frac{1}{4} \|\pi_{S_X} - \pi_P\|^2 + \mathcal{R} 
\end{equation*}
where $|\mathcal{R}| \leq C(n) \,\|\pi_{S_X} - \pi_P\|^4$.  Noting that $\|\pi_{S_X} - \pi_P\| \leq C(n,m) \,|\vec T(X) - \vec P|$, this completes the proof of \eqref{D error eqn11}.

Again noting that $\|\pi_{S_X} - \pi_P\| \leq C(n,m) \,|\vec T(X) - \vec P|$ and combining \eqref{D deriv}, \eqref{D error eqn10}, \eqref{D error eqn11}, and \eqref{D error concl1}, we conclude that \eqref{D error concl3} holds true.
\end{proof}

\subsection{Proof of Theorem~\ref{mono freq thm}}  

\begin{proof}[Proof Theorem~\ref{mono freq thm}]
Assume without loss of generality that $\rho_0 = 1$ and $Z = 0$.  By differentiating $N_{T,P,0}(\rho)$ using \eqref{D to I}, \eqref{H error concl2} and \eqref{D error concl3}, we obtain  
\begin{align}\label{mono freq eqn1}
	&\frac{d}{d\rho} \log N_{T,P,0}(\rho) = \frac{H_{T,P,0}(\rho) \,D'_{T,P,0}(\rho) - H'_{T,P,0}(\rho) \,D_{T,P,0}(\rho)}{H_{T,P,0}(\rho) \,D_{T,P,0}(\rho)}  
	\\&\geq \frac{2 \rho^{1-2n}}{H_{T,P,0}(\rho) \,D_{T,P,0}(\rho)} \left( \left( \int |y|^2 \,|\nabla^S r|^2 \,\frac{1}{r} \,\phi'(r/\rho) \,d|T|(X, S) \right) \right. \nonumber
		\\&\hspace{5mm} \cdot \left( \int_{\{\nabla^S r \neq 0\}} \frac{|\pi_{P^{\perp}}(\nabla^S r)|^2}{|\nabla^S r|^2} \,r \,\phi'(r/\rho) \,d|T|(X, S) \right) \nonumber 
		\\&\hspace{5mm} \left. - \left( \int (0,y) \cdot \pi_{P^{\perp}}(\nabla^S r) \,\phi'(r/\rho) \,d|T|(X, S) \right)^2 \right) \nonumber
		\\&\hspace{5mm} - C \eta^{2\gamma} \rho^{2\alpha\gamma-1} 
		- \frac{C}{\rho} \,D_{T,P,0}(\rho)^{\gamma-1} \,((n-1) \,D_{T,P,0}(\rho) + \rho \, D'_{T,P,0}(\rho)) \nonumber 
\end{align}
for all $\sigma_0 \leq \rho \leq 1$ and some constants $C = C(n,m,q) \in (0,\infty)$ and $\gamma = \gamma(n,m,q) \in (0,1)$.  By Lemma~\ref{tilt to height estimate lemma} and \eqref{mono freq decay hyp}, 
\begin{align}\label{mono freq eqn2}
	D_{T,P,0}(\rho) &\leq \frac{1}{2} \,\rho^{2-n} \int_{\mathbf{C}_{\rho}(0)} \|\pi_S - \pi_P\|^2 \,d|T|(X,S)  
		\\&\leq \frac{512}{9} \rho^{-n} \int_{\mathbf{C}_{7\rho/4}(0)} \op{dist}^2(X,P) \,d\|T\|(X) \leq C \eta^2 \rho^{2+2\alpha} \nonumber
\end{align}
for each $\sigma_0 \leq \rho \leq 1$ and some constant $C = C(n) \in (0,\infty)$.  By the Cauchy-Schwartz inequality and \eqref{mono freq eqn2}, \eqref{mono freq eqn1} gives us 
\begin{equation}\label{mono freq eqn3}
	\frac{d}{d\rho} \log N_{T,P,0}(\rho) 
	\geq - C_1 \eta^{2\gamma} \rho^{2\alpha\gamma-1} 
	- C_2 D_{T,P,0}(\rho)^{\gamma-1} D'_{T,P,0}(\rho) \nonumber
\end{equation}
for all $\sigma_0 \leq \rho \leq 1$ and some constants $C_1 = C_1(n,m,q) \in (0,\infty)$, $C_2 = C_2(n,m,q) \in (0,\infty)$, and $\gamma = \gamma(n,m,q) \in (0,1)$.  Integrating \eqref{mono freq eqn3} over $[\sigma,\rho]$, 
\begin{equation}\label{mono freq eqn4}
	N_{T,P,0}(\sigma) \leq N_{T,P,0}(\rho) \exp\left( \frac{C_1 \eta^{2\gamma}}{2\alpha\gamma} \,(\rho^{2\alpha\gamma} - \sigma^{2\alpha\gamma}) 
		+ \frac{C_2}{\gamma} \,(D_{T,P,0}(\rho)^{\gamma} - D_{T,P,0}(\sigma)^{\gamma}) \right) . 
\end{equation}
for all $\sigma_0 \leq \sigma < \rho \leq 1$.  Bounding the right-hand side of \eqref{mono freq eqn4} using \eqref{mono freq eqn2} gives us \eqref{mono freq concl}. 
\end{proof}

\section{Preliminary consequences of monotonicity of the planar frequency function}\label{sec:freq mono cor sec}

Here we will draw several preliminary consequences of the monotonicity formula in Theorem~\ref{mono freq thm}, including the existence of planar frequency at ``well-behaved'' branch points (Lemma~\ref{planar frequency lemma} and Definition~\ref{planar frequency defn}), growth estimates for $H_{T,P,Z}$ and for $L^2$-distance to a plane (Corollary~\ref{doubling cor}), and upper semi-continuity of planar frequency with respect to weak convergence of locally area-minimizing rectifiable currents (Corollary~\ref{limit semicont freq cor}).

Here and subsequently we shall use the notation 
$$E(T, P, {\mathbf C}_{\rho}(Z, P)) = \left( \frac{1}{\omega_{n}\rho^{n+2}} \int_{{\mathbf C}_{\rho}(Z, P)} {\rm dist}^{2} \, (X, Z + P) \, d\|T\|(X) \right)^{1/2}.$$

\begin{lemma}\label{planar frequency lemma}
Let $\alpha \in (0,1)$, $Z \in \mathbb{R}^{n+m}$, $\rho_0 > 0$ and $P \subset \mathbb{R}^{n+m}$ be an $n$-dimensional plane.  If 
$T$ is an $n$-dimensional locally area-minimizing rectifiable current of $\mathbf{C}_{7\rho_0/4}(Z, P)$ such that 
\begin{gather}
	\label{planar freq mass hyp} (\partial T) \llcorner \mathbf{C}_{7\rho_0/4}(Z, P) = 0, \quad
	\Theta(T,Z) \geq q , \quad \|T\|(\mathbf{C}_{7\rho_0/4}(Z, P)) < (q+1) \,\omega_n (7\rho_0/4)^n \;\; \mbox{and}\\
	\label{planar freq decay hyp} \sup_{\rho \in (0,\rho_0]} \rho^{-\alpha} E(T,P,\mathbf{C}_{7\rho/4}(Z, P)) < \infty 
\end{gather}
then: 
\begin{enumerate}[itemsep=3mm,topsep=0mm]
	\item[{\rm (i)}]  $P$ can be equipped with an orientation such that $q \llbracket P \rrbracket$ is the unique tangent cone to $T$ at $Z$; 
	\item[{\rm (ii)}]  if additionally $\op{spt} T \cap \mathbf{C}_{\rho}(Z, P) \not\subset P$ for all $\rho \in (0,\rho_0]$, then $N_{T,P,Z}(\rho)$ is well-defined for all sufficiently small $\rho > 0$ and the limit $\lim_{\rho \rightarrow 0^+} N_{T,P,Z}(\rho)$ exists. 
\end{enumerate}
\end{lemma}

\begin{definition}\label{planar frequency defn}{\rm 
Given $T$, $P$, and $Z$ satisfying the hypothesis of Lemma~\ref{planar frequency lemma}, we define the \emph{planar frequency} $\mathcal{N}_{T,{\rm Pl}}(Z)$ of $T$ at $Z$ with respect to $P$ by  
\begin{equation*}
	\mathcal{N}_{T,{\rm Pl}}(Z) = \lim_{\rho \rightarrow 0^+} N_{T,P,Z}(\rho) . 
\end{equation*}
}\end{definition}

\begin{proof}[Proof of Lemma~\ref{planar frequency lemma}]
Without loss of generality, assume that $Z = 0$, $\rho_0 = 1$ and $P = \mathbb{R}^n \times \{0\}$.  For any given $\delta, \eta \in (0, 1),$ 
by \eqref{planar freq decay hyp} and the monotonicity formula for area,  we can rescale and assume that 
\begin{gather}
	\label{planar frequency eqn1} \|T\|(\mathbf{C}_{7/4}(0)) < (\Theta(T,0) + \delta) \,\omega_n (7/4)^n , \\
	\label{planar frequency eqn2} E(T,P,\mathbf{C}_{2\rho}(0)) \leq \eta\rho^{\alpha} \text{ for all } \rho \in (0,1] . 
\end{gather}
Assuming $\eta$ is sufficiently small, by \eqref{planar frequency eqn2} and Lemma~\ref{Allard height lemma}
\begin{equation}\label{planar frequency eqn3}
	\sup_{X \in \op{spt} T \cap \mathbf{B}_{3\rho/2}(0)} \op{dist}(X,P) \leq C(n,m) \,\eta \rho^{1+\alpha} < \rho/2
\end{equation}
for all $\rho \in (0,1]$, where $C = C(n,m) \in (0,\infty)$ is a constant.  Recall that $P = \mathbb{R}^n \times \{0\}$ is oriented by $e_1 \wedge e_2 \wedge\cdots\wedge e_n$.  By \eqref{planar freq mass hyp} and the constancy theorem (cf.\ beginning of the proof of Lemma~\ref{H error lemma}), after possibly reversing the orientation of $P$ 

\begin{equation}\label{planar frequency eqn4}
	(\pi_{P \#} T) \llcorner B_{1/2}(0) = q \llbracket B_{1/2}(0) \rrbracket  .
\end{equation}
Let $\mathbf{C}$ be any tangent cone to $T$ at the origin and $\rho_k \rightarrow 0^+$ such that $\eta_{0,\rho_k\#} T \rightarrow \mathbf{C}$ weakly in $\mathbb{R}^{n+m}$.  By \eqref{planar frequency eqn3}, $\op{spt} \mathbf{C} \subseteq P$.  By continuity of push-fowards of integral currents with respect to weak convergence and \eqref{planar frequency eqn4}, 
\begin{equation*}
	\mathbf{C} = \pi_{P \#} \mathbf{C}  = \pi_{P \#} (\lim_{k \rightarrow \infty} \eta_{0,\rho_k\#} T) = \lim_{k \rightarrow \infty} \pi_{P \#} \eta_{0,\rho_k\#} T
		= \pi_{\#} (q \llbracket P \rrbracket) = q \llbracket P \rrbracket
\end{equation*}
in $\mathbf{B}_1(0)$, where each limit is computed with respect to the weak topology.  
Therefore, $\mathbf{C} = q \llbracket P \rrbracket$ in $\mathbb{R}^{n+m}$.  In light of the arbitrary choice of sequence $(\rho_k)$, 
$\eta_{0,\rho \, \#} T \rightarrow q \llbracket P \rrbracket$ weakly in $\mathbb{R}^{n+m}$ as $\rho \to 0^{+}$; that is, $q \llbracket P \rrbracket$ is the unique tangent cone to $T$ at $Z$.  Note that $\Theta(T,0) = q$.  By \eqref{planar frequency eqn1} and \eqref{planar frequency eqn2}, $H_{T,P,0}(\rho)$ and $D_{T,P,0}(\rho)$ are defined for all $\rho \in (0,\rho_0]$.  By Remark~\ref{H zero rmk}, if we assume that $\op{spt} T \cap \mathbf{B}_{\rho}(Z) \not\subset P$ for all $\rho \in (0,\rho_0]$, then $H_{T,P,0}(\rho) > 0$ and thus $N_{T,P,0}(\rho)$ is well-defined for all $\rho \in (0,1]$.  By \eqref{planar frequency eqn1} and \eqref{planar frequency eqn2}, we can apply Theorem~\ref{mono freq thm} to deduce that $\lim_{\rho \rightarrow 0^+} N_{T,P,Z}(\rho)$ exists.
\end{proof}

\begin{corollary}\label{doubling cor}
For each positive integer $q$  there exists $\delta = \delta(n,m,q) \in (0,1)$ and $\eta_0 = \eta_0(n,m,q) \in (0, 1)$ such that the following holds true.  Let $Z \in \mathbb{R}^{n+m}$ and $\rho_0 > 0$.  Let $P$ be an $n$-dimensional plane in $\mathbb{R}^{n+m}$ and $T$ be $n$-dimensional locally area-minimizing rectifiable current in $\mathbf{C}_{7\rho_0/4}(Z,P)$ such that \eqref{mono freq mass hyp} holds true.  
\begin{enumerate}[itemsep=3mm,topsep=0mm]
	\item[{\rm (i)}]  If there exists $\eta \in (0,\eta_0]$ and $\sigma_0 \in (0,\rho_0)$ such that \eqref{mono freq decay hyp} holds true for some $\alpha \in (0,1)$ and 
	all $\rho \in [\sigma_0,\rho_0],$ and if $H_{T,P,Z}(\rho) > 0$ for some $\rho \in [\sigma_0,\rho_0]$, then $H_{T,P,Z}(\rho) > 0$ for \emph{all} $\rho \in [\sigma_0,\rho_0]$ and 
\begin{align}\label{doubling concl1}
	e^{-C \eta^{2\gamma} (\rho/\rho_0)^{2\alpha\gamma} / (2\alpha\gamma)} 
		&\Big(\frac{\sigma}{\rho}\Big)^{2e^{C\eta^{\gamma} (\rho/\rho_0)^{\alpha\gamma}} N_{T,P,Z}(\rho)} H_{T,P,Z}(\rho) 
	\leq H_{T,P,Z}(\sigma) 
	\\&\leq e^{C (N_{T,P,Z}(\sigma_0) + 1) \eta^{2\gamma} (\rho/\rho_0)^{2\alpha\gamma} / (2\alpha\gamma)} 
		\Big(\frac{\sigma}{\rho}\Big)^{2N_{T,P,Z}(\sigma_0)} H_{T,P,Z}(\rho) & \nonumber 
\end{align}
for all $\sigma_0 \leq \sigma < \rho \leq \rho_0$, where $\gamma = \gamma(n,m,q) \in (0,1)$ is as in Theorem~\ref{mono freq thm} and $C = C(n,m,q,\alpha) \in (0,\infty)$ is a constant (independent of $\eta$ and $\alpha$). 

	\item[{\rm (ii)}]  If there exists $\eta \in (0,\eta_0]$ such that \eqref{mono freq decay hyp} holds true for all $\rho \in (0,\rho_0]$ and if $H_{T,P,Z}(\rho) > 0$ for some $\rho \in (0,\rho_0]$, then $H_{T,P,Z}(\rho) > 0$ for \emph{all} $\rho \in (0,\rho_0]$, $\mathcal{N}_{T,{\rm Pl}}(Z) \geq 1+\alpha$, and 
\begin{align}\label{doubling concl2}
	\frac{1}{4} e^{-C \eta^{2\gamma} (\rho/\rho_0)^{2\alpha\gamma}/(4\alpha\gamma)} 
		\Big(\frac{\sigma}{\rho}\Big)^{e^{C\eta^{\gamma} (\rho/\rho_0)^{\alpha\gamma}} N_{T,P,Z}(\rho) - 1} E(T,P,\mathbf{C}_{\rho}(Z, P)) 
	\leq E(T,P,\mathbf{C}_{\sigma}(Z, P)) &\\
	\leq 4e^{C (\mathcal{N}_{T,{\rm Pl}}(Z) + 2) \,\eta^{2\gamma} (\rho/\rho_0)^{2\alpha\gamma}/(4\alpha\gamma)} 
		\Big(\frac{\sigma}{\rho}\Big)^{\mathcal{N}_{T,{\rm Pl}}(Z)-1} E(T,P,\mathbf{C}_{\rho}(Z, P)) & \nonumber 
\end{align}
for all $0 < \sigma < \rho \leq \rho_0$, where $C = C(n,m,q) \in (0,\infty)$ is a constant (independent of $\eta$ and $\alpha$). 
\end{enumerate}
\end{corollary} 

\begin{proof}
Without loss of generality assume $Z = 0$ and $\rho_0 = 1$.  To see (i), suppose that $\sigma_0 \leq \sigma < \rho \leq 1$ and that $H_{T, P, 0}(\sigma) >0$. 
 Then by Remark~\ref{H zero rmk}, $H_{T, P, 0}(\tau) >0$ for all $\tau \in [\sigma, \rho]$.  By \eqref{H error concl2} and \eqref{D to I}
\begin{equation*}
	\left| \frac{H'_{T,P,0}(\tau)}{H_{T,P,0}(\tau)} - \frac{2 N_{T,P,0}(\tau)}{\tau} \right| \leq C \eta^{2\gamma} \tau^{2\alpha\gamma-1} 
\end{equation*}
for all $\tau \in [\sigma, \rho]$ and some constant $C = C(n,m,q) \in (0,\infty)$.  Hence using Theorem~\ref{mono freq thm} and the fact that $e^{-C\eta^{\gamma}\tau^{\alpha\gamma}} \geq 1 - C\eta^{\gamma}\tau^{\alpha\gamma}$, 
\begin{align}\label{doubling eqn1}
	&\frac{2N_{T,P,0}(\sigma_0)}{\tau} - C (N_{T,P,0}(\sigma_0) + 1) \eta^{2\gamma} \tau^{2\alpha\gamma-1} 
	\\&\hspace{15mm} \leq \frac{H'_{T,P,0}(\tau)}{H_{T,P,0}(\tau)} 
	\leq \frac{2e^{C\eta^{\gamma} \rho^{\alpha\gamma}} N_{T,P,0}(\rho)}{\tau} + C \eta^{2\gamma} \tau^{2\alpha\gamma-1} \nonumber 
\end{align}
for all $\tau \in [\sigma,\rho]$ where $C = C(n,m,q) \in (0,\infty)$.  Integrating \eqref{doubling eqn1} over $\tau \in [\sigma,\rho]$ we obtain \eqref{doubling concl1}, 
subject to the assumption that $H_{T, P, 0}(\sigma) >0$.  Finally, we deduce from this conclusion and Remark~\ref{H zero rmk}  that if 
$H_{T, P, 0}(\tau) >0$ for some $\tau \in [\sigma_{0}, 1]$ then $H_{T, P, 0}(\tau) >0$ for all $\tau \in [\sigma_{0}, 1]$ and hence \eqref{doubling concl1} 
holds for any $\sigma, \rho$ with $\sigma_{0} \leq \sigma \leq \rho \leq 1$. 

To see (ii), fix $0 < \sigma < \rho \leq 1$.  Assume without loss of generality that $P = {\mathbb R}^{n} \times \{0\}$. For $k = 0,1,2,\ldots$ set $\sigma_k = 2^{-k} \sigma$ and $\rho_k = 2^{-k} \rho$.  By Lemma~\ref{H error lemma}, in particular \eqref{H error concl1} together with \eqref{H error eqn11}, 
\begin{equation*}
	\frac{1}{4} \,H_{T,P,0}(\tau) \leq \frac{1}{\tau^{n}} \int_{\mathbf{C}_{\tau}(0) \setminus \mathbf{C}_{\tau/2}(0)} \op{dist}^2(X,P) \,d\|T\|(X) 
		\leq 4(1 + C \eta^{2\gamma} \tau^{2\alpha\gamma}) \,H_{T,P,0}(\tau)
\end{equation*}
for all $\tau \in (0,1]$, where $\gamma = \gamma(n,m,q) \in (0,1)$ and $C = C(n,m,q) \in (0,\infty)$ are constants.  Noting that $1 + C \eta^{2\gamma} \tau^{2\alpha\gamma} \leq e^{C\eta^{2\gamma} \tau^{2\alpha\gamma}}$ and $\sigma_k/\rho_k = \sigma/\rho$ for each $k$, by \eqref{doubling concl1} 
\begin{align}\label{doubling2 eqn1}
	&\frac{1}{16} e^{-C \eta^{2\gamma} \rho^{2\alpha\gamma}/(2\alpha\gamma)} 
		\Big(\frac{\sigma}{\rho}\Big)^{n+2e^{C\eta^{\gamma} \rho^{\alpha\gamma}} N_{T,P_0,0}(\rho)} 
		\int_{\mathbf{C}_{\rho_k}(0) \setminus \mathbf{C}_{\rho_k/2}(0)} \op{dist}^2(X,P_0) \,d\|T\|(X) 
	\\ \leq\,& \int_{\mathbf{C}_{\sigma_k}(0) \setminus \mathbf{C}_{\sigma_k/2}(0)} \op{dist}^2(X,P_0) \,d\|T\|(X) \nonumber
	\\ \leq\,& 16 e^{C ({\mathcal N}_{T,{\rm Pl}}(0) + 2) \,\eta^{2\gamma} \rho^{2\alpha\gamma}/(2\alpha\gamma)} 
		\Big(\frac{\sigma}{\rho}\Big)^{n+2\mathcal{N}_{T,{\rm Pl}}(0)} \int_{\mathbf{C}_{\rho_k}(0) \setminus \mathbf{C}_{\rho_k/2}(0)} 
		\op{dist}^2(X,P_0) \,d\|T\|(X) \nonumber 
\end{align}
for all $k = 0,1,2,\ldots$, where $C = C(n,m,q) \in (0,\infty)$ is a constant.  By summing \eqref{doubling2 eqn1} over $k = 0,1,2,\ldots$ gives us \eqref{doubling concl2}.  By \eqref{mono freq mass hyp}, $E(T,P_{0},\mathbf{C}_{\sigma}(0)) \leq \eta (\sigma/2)^{\alpha}$ for all $\sigma \in (0,1]$, which by in view of \eqref{doubling concl2} gives us $e^{C\eta^{\gamma} \rho^{\alpha\gamma}} N_{T,P_{0},0}(\rho) \geq 1+\alpha$ for all $\rho \in (0,1]$.  Letting $\rho \rightarrow 0^+$ gives us $\mathcal{N}_{T,{\rm Pl}}(0) \geq 1+\alpha$.
\end{proof}

\begin{lemma}\label{limit frequency lemma}
Let $P_k,P$ be $n$-dimensional planes in $\mathbb{R}^{n+m}$ such that 
\begin{equation}\label{limit frequency hyp1}
	\lim_{k \rightarrow \infty} \op{dist}_{\mathcal{H}}(P_k \cap \mathbf{B}_1(0), P \cap \mathbf{B}_1(0)) = 0 . 
\end{equation}
Let $Z_k,Z \in \mathbf{B}_1(0)$, $\rho_k > 0$, and $0 < \rho < \op{dist}(Z,\partial \mathbf{C}_1(0, P))$ such that $Z_k \rightarrow Z$ and $\rho_k \rightarrow \rho$.  Let $T_k,T$ be $n$-dimensional locally area-minimizing rectifiable currents in 
$\mathbf{C}_1(0, P)$ such that 
\begin{equation}\label{limit frequency hyp2}
	\sup_{k \geq 1} \sup_{X \in \op{spt} T_k \cap \mathbf{C}_1(0, P)} \op{dist}(X,Z_k+P_k) < \infty
\end{equation}
and 
\begin{gather}
	\label{limit frequency hyp3} (\partial T_k) \llcorner \mathbf{C}_1(0, P) = 0, \quad 
	\sup_{k \geq 1} \|T_k\|(\mathbf{C}_1(0, P)) < \infty , \\ 
	\op{spt} T \not\subset P, \quad T_k \rightarrow T \text{ weakly in } \mathbf{C}_{1}(0, P) . \nonumber 
\end{gather}
Then 
\begin{equation}\label{limit frequency concl}
	\lim_{k \rightarrow \infty} N_{T_k,P_k,Z_k}(\rho_k) = N_{T,P,Z}(\rho) .  
\end{equation}
\end{lemma}

\begin{proof}
Note that by \eqref{limit frequency hyp1} and \eqref{limit frequency hyp2}, $\op{spt} T_k \cap \overline{\mathbf{C}_{\rho_k}(Z_k,P_k)} \subset \mathbf{C}_1(0,P)$ and thus $(\partial T_k) \llcorner \mathbf{C}_{\rho_k}(Z_k,P_k) = 0$ for all large $k$.  Hence $H_{T_k,P_k,Z_k}(\rho_k)$ and $D_{T_k,P_k,Z_k}(\rho_k)$ are well-defined, and $N_{T_k,P_k,Z_k}(\rho_k)$ is well-defined if $H_{T_k,P_k,Z_k}(\rho_k) > 0$.  Also, by letting $k \rightarrow \infty$ in \eqref{limit frequency hyp2} and \eqref{limit frequency hyp3}, 
\begin{equation*}
	(\partial T) \llcorner \mathbf{C}_1(0, P) = 0, \quad \|T\|(\mathbf{C}_1(0,P)) < \infty , \quad 
	\sup_{X \in \op{spt} T \cap \mathbf{C}_1(0, P)} \op{dist}(X,Z+P) < \infty.
\end{equation*}
Hence $H_{T, P, Z}(\rho)$ and $D_{T, P, Z}(\rho)$ are well-defined, and $N_{T,P,Z}(\rho)$ is well-defined if $H_{T, P, Z}(\rho) >0$.  Next, after rotating and rescaling we may assume that $Z_k = Z = 0$, $\rho_k = \rho$, and $P_k = P = \mathbb{R}^n \times \{0\}$ for all $k$.  Since $T_k$ is area-minimizing in $\mathbf{C}_{1}(0)$ and $(\partial T_k) \llcorner \mathbf{C}_{1}(0) = 0$, $|T_k|$ is stationary in $\mathbf{C}_{1}(0)$.  Thus since $\sup_k \|T_k\|(\mathbf{C}_{1}(0)) < \infty$, by the Allard compactness theorem~\cite[Theorem~42.7]{SimonGMT} after passing to a subsequence $|T_k|$ converges weakly to some stationary integral varifold in $\mathbf{C}_{1}(0)$.  By~\cite[Theorem~34.5]{SimonGMT}, $\|T_k\| \rightarrow \|T\|$ in the sense of Radon measures locally in 
$\mathbf{C}_{1}(0)$, and it follows that $|T_k| \rightarrow |T|$ in the sense of varifolds in $\mathbf{C}_{1}(0)$.  By Remark~\ref{phi deriv rmk}, $|\nabla^S r| = 0$ for $|T|$-a.e.~$(X,S) \in G_n(\mathbf{C}_{1}(0) \cap (\partial \mathbf{C}_{\rho/2}(0) \cup \partial \mathbf{C}_{\rho}(0)))$ and thus 
\begin{equation*}
	\int_{G_n(\partial \mathbf{C}_{\rho/2}(0) \cup \partial \mathbf{C}_{\rho}(0))} 
		\op{dist}^2(X,P) \,|\nabla^S r|^2 \,\frac{1}{r} \,d|T|(X,S) = 0 .
\end{equation*}
Hence using the fact that $|T_k| \rightarrow |T|$ in the sense of varifolds in $\mathbf{C}_{1}(0)$
\begin{align}\label{limit frequency eqn1} 
	\lim_{k \rightarrow \infty} H_{T_k,P,0}(\rho)
	&= \lim_{k \rightarrow \infty} 2\rho^{1-n} \int_{G_n(\mathbf{C}_{\rho}(0) \setminus \mathbf{C}_{\rho/2}(0))} 
		\op{dist}^2(X,P) \,|\nabla^S r|^2 \,\frac{1}{r} \,d|T_k|(X,S) 
	\\&= 2\rho^{1-n} \int_{G_n(\mathbf{C}_{\rho}(0) \setminus \mathbf{C}_{\rho/2}(0))} 
		\op{dist}^2(X,P) \,|\nabla^S r|^2 \,\frac{1}{r} \,d|T|(X,S) \nonumber
	\\&= H_{T,P,0}(\rho) . \nonumber
\end{align}
Similarly, 
\begin{align}\label{limit frequency eqn2} 
	\lim_{k \rightarrow \infty} D_{T_k,P,0}(\rho)
	&= \lim_{k \rightarrow \infty} \rho^{2-n} \int_{G_n(\mathbf{C}_{\rho}(0))} \|\pi_S - \pi_P\|^2 \,\phi(r/\rho) \,d|T_k|(X,S) 
	\\&= \rho^{2-n} \int_{G_n(\mathbf{C}_{\rho}(0))} \|\pi_S - \pi_P\|^2 \,\phi(r/\rho) \,d|T|(X,S) \nonumber
	\\&= D_{T,P,0}(\rho) . \nonumber 
\end{align}
Note that by Remark~\ref{H zero rmk}, since $\op{spt} T \not\subset P$, $H_{T,P,0}(\rho) > 0$.  Hence by dividing \eqref{limit frequency eqn2} by \eqref{limit frequency eqn1} we obtain \eqref{limit frequency concl}.
\end{proof}

\begin{corollary}\label{limit semicont freq cor}
Let $P_k,P$ be $n$-dimensional planes in $\mathbb{R}^{n+m}$ such that 
\begin{equation}\label{limit semicont freq hyp1} 
	\lim_{k \rightarrow \infty} \op{dist}_{\mathcal{H}}(P_k \cap \mathbf{B}_1(0), P \cap \mathbf{B}_1(0)) = 0 . 
\end{equation}
Let $T_k$ and $T$ be $n$-dimensional locally area-minimizing rectifiable currents of $\mathbf{C}_3(0,P)$ such that 
\begin{equation}\label{limit semicont freq hyp2} 
	(\partial T_k) \llcorner \mathbf{C}_3(0,P) = 0, \quad T_k \rightarrow T \text{ weakly in } \mathbf{C}_3(0,P). 
\end{equation}
Let $\alpha \in (0,1)$ be a constant and for $k = 1,2,3,\ldots$ let $Z_k \in \op{spt} T_k \cap \mathbf{C}_1(0)$ such that 
\begin{gather}
	\label{limit semicont freq hyp3} \Theta(T_k,Z_k) \geq q, \quad 
	\limsup_{k \rightarrow \infty} \|T_k\|(\mathbf{C}_{7/4}(Z_k,P_k)) < (q + 1) \,\omega_n (7/4)^n, \\
	\label{limit semicont freq hyp4} \sup_{k \geq 1} \sup_{X \in \op{spt} T_k \cap \mathbf{C}_1(Z_k,P_k)} \op{dist}(X,Z_k+P_k) < \infty , \\
	\label{limit semicont freq hyp5} \sup_{k \geq 1} \sup_{\rho \in (0,1]} \rho^{-\alpha} E(T_k,P_k,\mathbf{C}_{7\rho/4}(Z_k,P_k)) < \infty 
     \end{gather}
and let $Z_k \rightarrow Z$ in $\mathbb{R}^{n+m}$.  
Then 
\begin{gather}
	\label{limit semicont freq concl1} \Theta(T_{k}, Z_{k}) = q \text{ for all } k, \quad \Theta(T,Z) = q, \quad \|T\|(\mathbf{C}_{7/4}(Z,P)) < (q + 1) \,\omega_n (7/4)^n, \\
	\label{limit semicont freq concl2} \sup_{\rho \in (0,1]} \rho^{-\alpha} E(T,P,\mathbf{C}_{7\rho/4}(Z,P)) < \infty, \;\; and\\
	\label{limit semicont freq concl3} \mathcal{N}_{T,{\rm Pl}}(Z) \geq \limsup_{k \rightarrow \infty} \mathcal{N}_{T_k,{\rm Pl}}(Z_k) .
\end{gather}
\end{corollary}

\begin{proof}
The conclusion $\Theta(T_{k}, Z_{k}) = q$ follows from Lemma~\ref{planar frequency lemma}(i). It is clear from \eqref{limit semicont freq hyp2} that $(\partial T) \llcorner \mathbf{C}_3(0,P) = 0$.  By \eqref{limit semicont freq hyp2}, \eqref{limit semicont freq hyp3}, and semi-continuity of density and mass~\cite[Corollary 17.8 and 26.13]{SimonGMT} we have that the rest of the conclusions in \eqref{limit semicont freq concl1} hold true (with equality 
in $\Theta(T, Z) = q$ again following from Lemma~\ref{planar frequency lemma}(i)).  By~\cite[Theorem~34.5]{SimonGMT}, $\|T_k\| \rightarrow \|T\|$ in the sense of Radon measures of $\mathbf{C}_3(0,P)$.  Hence  
\begin{equation*}
	E(T,P,\mathbf{C}_{7\rho/4}(Z,P)) \leq \lim_{k \rightarrow \infty} E(T_k,P_k,\mathbf{C}_{7\rho/4}(Z_k,P_k)) 
\end{equation*}
for each $\rho \in (0,1]$ and thus \eqref{limit semicont freq concl2} holds true.  

By \eqref{limit semicont freq hyp1} and \eqref{limit semicont freq hyp4}, $\op{spt} T_k \cap \overline{\mathbf{C}_1(Z_k,P_k)} \subset \mathbf{C}_1(0)$ and thus $(\partial T_k) \llcorner \mathbf{C}_1(Z_k,P_k) = 0$.  In light of \eqref{limit semicont freq hyp3}, \eqref{limit semicont freq hyp5}, \eqref{limit semicont freq concl1}, and \eqref{limit semicont freq concl2}, by rescaling we may assume that for suitably small constant $\delta = \delta(n,m,q) \in (0,1)$ and $\eta = \eta(n,m,q) \in (0,1)$ and for all sufficiently large $k$ 
\begin{gather*} 
	\|T_k\|(\mathbf{C}_{7/4}(Z_k,P_k)) < (q + \delta) \,\omega_n (7/4)^n, \quad \|T\|(\mathbf{C}_{7/4}(Z,P)) < (q + \delta) \,\omega_n (7/4)^n , \\
	E(T_k,P_k,\mathbf{C}_{7/4}(Z_k,P_k)) \leq \eta\rho^{\alpha}, \quad E(T,P,\mathbf{C}_{7/4}(Z,P)) \leq \eta\rho^{\alpha} 
		\quad\text{for all } \rho \in (0,1] . 
\end{gather*}
Thus by Theorem~\ref{mono freq thm} and Lemma~\ref{limit frequency lemma} 
\begin{equation*}
	\limsup_{k \rightarrow \infty} \mathcal{N}_{T_k,{\rm Pl}}(Z_k) 
	\leq \lim_{k \rightarrow \infty} e^{C \eta^{\gamma} \rho^{\alpha\gamma}} N_{T_k,P_k,Z_k}(\rho)
	= e^{C \eta^{\gamma} \rho^{\alpha\gamma}} N_{T,P,Z}(\rho) 
\end{equation*}
for each $\rho \in (0, 1]$, where $C = C(n,m,q,\alpha) \in (0,\infty)$ is a constant.  Letting $\rho \rightarrow 0^+$ gives us \eqref{limit semicont freq concl3}.\end{proof}

\section{Blow-ups of area-minimizers relative to a plane}\label{sec:blowup}

Let $P$ be an $n$-dimensional oriented plane in $\mathbb{R}^{n+m}$.  For $k = 1,2,3,\ldots,$ let $T_k$ be $n$-dimensional locally area-minimizing rectifiable currents in $\mathbf{C}_{\rho_0}(0)$ which are converging to the multiplicity $q$ plane $q \llbracket P \rrbracket$ as in \eqref{blowup hyp1} and \eqref{blowup hyp2} below.  In Subsection~\ref{sec:blowup procedure}, we discuss, following Almgren's work (\cite{Almgren}), the procedure for constructing a $q$-valued Dirichlet energy minimizing blow-up $w$ of the sequence $(T_k)$ relative to $P$.  
In Lemma~\ref{optimal plane lemma}, we consider the optimal planes 
$\widehat{P}_k$ which minimizes the $L^2$-distance of $T_k$ to a plane in a fixed ball, and show that the blow-up of $\widehat{P}_k$ relative to $P$ is $Dw_a(0) \cdot x$ where $w_{a}$ is the (harmonic) pointwise average of the values of $w$.  
In Lemma~\ref{hardt simon lemma} we establish the Hardt-Simon inequality, which was first introduced in~\cite{HardtSimon} and which we use here to show that points at which $T_k$ has density $\geq q$ blow up to points at which $w$ has Almgren frequency $\geq 1$.  In Lemma~\ref{blowup frequency lemma} and Corollary~\ref{blowup semicont freq cor}, we establish continuity of the planar frequency function associated with $T_{k}$, and upper semi-continuity of planar frequency associated with $T_{k},$ with respect to blowing up of $(T_{k})$, giving in the limit Almgren frequency function and Almgren frequency associated with the blow-up $w$.

\subsection{Blow-up procedure}\label{sec:blowup procedure}  Let $P$ be an $n$-dimensional oriented plane in $\mathbb{R}^{n+m}$.  After an orthogonal change of coordinates, we may assume that $P = P_0 = \mathbb{R}^n \times \{0\}$ and $P_0$ is oriented by $\vec P_0 = e_1 \wedge e_2 \wedge\cdots \wedge e_n$.  We will express points $X \in \mathbb{R}^{n+m}$ as $X = (x,y)$ where $x \in \mathbb{R}^n$ and $y \in \mathbb{R}^m$.  Let $\rho_0 > 0$.  For $k = 1,2,3,\ldots$ let $T_k$ be $n$-dimensional area-minimizing integral currents of $\mathbf{C}_{\rho_0}(0)$ such that 
\begin{gather}
	\label{blowup hyp1} (\partial T_k) \llcorner \mathbf{C}_{\rho_0}(0) = 0, \quad\quad \Theta(T_k,0) \geq q , \quad\quad 
		\lim_{k \rightarrow \infty} \|T_k\|(\mathbf{C}_{\rho_0}(0)) = q \omega_n \rho_0^n , \\
	\label{blowup hyp2} \lim_{k \rightarrow \infty} E(T_k,P,\mathbf{C}_{\rho_0}(0)) = 0 .
\end{gather}
For each $k$ we set $E_k = E(T_k,P,\mathbf{C}_{\rho_0}(0))$.  (More generally, we can take $E_k > 0$ such that $E(T_k,P,\mathbf{C}_{\rho_0}(0)) \leq C E_k$ for a suitable constant $C \in (0,\infty)$.)

By Lemma~\ref{Allard height lemma}, for each $\theta \in (0,1)$ and sufficiently large $k$ 
\begin{equation}\label{blowup eqn1}
	\sup_{X \in \op{spt} T_k \cap \mathbf{C}_{\theta\rho_0}(0)} \op{dist}(X,P) \leq C E_k \rho_0 
\end{equation}
for some constant $C = C(n,m,\theta) \in (0,\infty)$.  Hence arguing as at the beginning of the proof of Lemma~\ref{H error lemma}, it follows using \eqref{blowup hyp1}, \eqref{blowup hyp2}, the constancy theorem, and Lemma~\ref{area excess to height lemma} that for each $\theta \in (0, 1)$ and sufficiently large $k$, 
\begin{equation*}
	\pi_{P \, \#} (T_k \llcorner \mathbf{C}_{\theta\rho_0}(0)) = \pm q \llbracket P \rrbracket \llcorner \mathbf{C}_{\theta\rho_0}(0) . 
\end{equation*}
After reversing the orientation of $T_k$ if necessary, we may assume that for each $\theta \in (0,1)$ and sufficiently large $k$ 
\begin{equation}\label{blowup eqn2}
	\pi_{P \, \#} (T_k \llcorner \mathbf{C}_{\theta\rho_0}(0)) = q \llbracket P \rrbracket \llcorner \mathbf{C}_{\theta\rho_0}(0) . 
\end{equation}
Note that \eqref{blowup hyp1} and the compactness of area-minimizing integral currents (\cite[Theorems~32.2 and 34.5]{SimonGMT}), after passing to a subsequence there exists an $n$-dimensional area-minimizing integral current $T_{\infty}$ of $\mathbf{C}_{\rho_0}(0)$ such that $T_k \rightarrow T_{\infty}$ weakly in $\mathbf{C}_{\rho_0}(0)$.  By \eqref{blowup eqn1}, $\op{spt} T_{\infty} \subset P$ and thus by the constancy theorem $T_{\infty} \llcorner B_{\rho_0}(0)$ is a constant multiple of $\llbracket B_{\rho_0}(0) \rrbracket$.  By the continuity of push-fowards of integral currents with respect to weak convergence 
\begin{equation*}
	T_{\infty} = \pi_{P \#} T_{\infty} = \lim_{k \rightarrow \infty} \pi_{P \#} T_k = \pi_{P \#} (\lim_{k \rightarrow \infty} T_k) 
		= \pi_{P \#} (q \llbracket P \rrbracket) = q \llbracket P \rrbracket
\end{equation*}
in $\mathbf{C}_{\rho_0}(0)$, where each limit is computed with respect to the weak topology.  That is, $T_k \rightarrow q \llbracket P \rrbracket$ weakly in $\mathbf{C}_{\rho_0}(0)$. 

By Theorem~\ref{lipschitz approx thm}, Lemma~\ref{area excess to height lemma}, and \eqref{blowup eqn1}, for each sufficiently large $k$ there exists $\theta_k \in (0,1)$ with $\theta_k \rightarrow 1$, Lipschitz functions $u_k : B_{\theta_k\rho_0}(0) \rightarrow \mathcal{A}_q(\mathbb{R}^m)$, and sets $K_k \subset B_{\theta_k\rho_0}(0)$ such that 
\begin{gather}
	\label{blowup eqn4} T_k \llcorner (K_k \times \mathbb{R}^m) 
		= (\op{graph} u_k) \llcorner (K_k \times \mathbb{R}^m) , \\
	\label{blowup eqn5} \mathcal{L}^n(B_{\theta\rho_0}(0) \setminus K_k) + \|T_k\|((B_{\theta\rho_0}(0) \setminus K_k) \times \mathbb{R}^m) 
		\leq C_{\theta} E_k^{2+2\gamma} \rho_0^n , \\
	\label{blowup eqn6} \sup_{B_{\theta\rho_0}(0)} |u_k| \leq C_{\theta} E_k \rho_0, \quad 
		\sup_{B_{\theta\rho_0}(0)} |\nabla u_k| \leq C_{\theta} E_k^{2\gamma},  
\end{gather}
for each $\theta \in (0,\theta_k]$, where $\gamma = \gamma(n,m,q) \in (0,1)$ and $C_{\theta} = C(n,m,q,\theta) \in (0,\infty)$ are constants.  We will let $u_k(x) = \sum_{i=1}^q \llbracket u_{k,i}(x) \rrbracket$ and $u_{k,a}(x) = \frac{1}{q} \sum_{i=1}^q u_{k,i}(x)$ for each $x \in B_{\theta_k\rho_0}(0)$ and $Du_k(x) = \sum_{i=1}^q \llbracket Du_{k,i}(x) \rrbracket$ for $\mathcal{L}^n$-a.e.~$x \in B_{\theta_k\rho_0}(0)$, where we follow the conventions from Subsection~\ref{sec:prelim_multivalued}.  By \eqref{blowup eqn4}--\eqref{blowup eqn6}, \eqref{lipschitz approx eqn8}, \eqref{lipschitz approx eqn9}, \eqref{lipschitz approx eqn10}, and Lemma~\ref{area excess to height lemma}, for each $\theta \in (0,1)$ and sufficiently large $k$ 
\begin{align}
	\label{blowup eqn8} \int_{B_{\theta\rho_0}(0)} |Du_k|^2 &\leq \int_{\mathbf{C}_{\theta\rho_0}(0)} |\vec T_k - \vec P|^2 \,d\|T_k\|(X) 
		+ C E_k^{2+2\gamma} \rho_0^n
	\\&\leq C \int_{\mathbf{C}_{\rho_0}(0)} \op{dist}^2(X,P) \,d\|T_k\|(X) + C E_k^{2+2\gamma} \rho_0^n \leq C E_k^2 \omega_n \rho_0^n , \nonumber 
\end{align}
where $C = C(n,m,q,\theta) \in (0,\infty)$ are constants.  Hence by~\cite[Proposition~2.11]{DeLSpaDirMin}, after passing to a subsequence there is $q$-valued function $w \in W^{1,2}_{\rm loc}(B_{\rho_0}(0),\mathcal{A}_q(\mathbb{R}^m))$ such that 
\begin{equation}\label{blowup eqn9}
	u_k/E_k \rightarrow w 
\end{equation} 
pointwise $\mathcal{L}^n$-a.e.~on $B_{\rho_0}(0)$ and strongly in $L^2(B_{\theta\rho_0}(0),\mathcal{A}_q(\mathbb{R}^m))$ for all $\theta \in (0,1)$.  Moreover, arguing as in~\cite[Theorem~4.2]{DeLSpa1} (also see~\cite[Theorem~2.19]{Almgren}), $w$ is locally Dirichlet energy minimizing in $B_{\rho_0}(0)$ and 
\begin{equation}\label{blowup eqn10}
	|Du_k|/E_k \rightarrow |Dw|
\end{equation}
in $L^2(B_{\theta\rho_0}(0))$ for all $\theta \in (0,1)$.  We call $w$ a \emph{blow-up} of $T_k$ relative to the plane $P$ in $\mathbf{C}_{\rho_0}(0)$.  We will let $w(x) = \sum_{i=1}^q \llbracket w_i(x) \rrbracket$ and $w_a(x) = \frac{1}{q} \sum_{i=1}^q w_i(x)$ for each $x \in B_{\rho_0}(0)$ for each $x \in B_{\theta_k\rho_0}(0)$ and $Dw(x) = \sum_{i=1}^q \llbracket Dw_i(x) \rrbracket$ for $\mathcal{L}^n$-a.e.~$x \in B_{\rho_0}(0)$, where we follow the conventions from Subsection~\ref{sec:prelim_multivalued}.

\subsection{Blow-up of optimal planes}\label{sec:optimal planes}  

\begin{definition}{\rm 
Let $Z \in \mathbb{R}^{n+m}$, $\rho > 0$, and $T$ be an $n$-dimensional area-minimizing integral current of $\mathbf{B}_{\rho}(Z)$.  We say that an $n$-dimensional plane $\widehat{P}$ is an \emph{optimal plane} for $T$ in $\mathbf{B}_{\rho}(Z)$ if 
\begin{equation}\label{optimal plane eqn1}
	E(T,\widehat{P},\mathbf{B}_{\rho}(Z)) = \inf_{P' \in \mathcal{P}} E(T,P',\mathbf{B}_{\rho}(Z)) , 
\end{equation}
where $\mathcal{P}$ is the set of all $n$-dimensional planes in $\mathbb{R}^{n+m}$. 
}\end{definition}

Notice that since $\mathcal{P}$ is compact with respect to Hausdorff distance between pairs of planes in $\mathbf{B}_1(0)$, at least one optimal plane $\widehat{P}$ as in \eqref{optimal plane eqn1} exists.  However, we do not claim the optimal plane is unique. 

\begin{lemma}\label{optimal plane lemma}
Let $T_k$ be $n$-dimensional locally area-minimizing rectifiable currents in $\mathbf{C}_1(0)$ such that \eqref{blowup hyp1}, \eqref{blowup hyp2}, and \eqref{blowup eqn2} hold true with $P = P_0$ ($= \mathbb{R}^n \times \{0\}$) and $\rho_{0} = 1$.  Let $u_k$ and $E_k$ be as in Subsection~\ref{sec:blowup procedure} with $\rho_0 = 1$ and assume that $w \in W^{1,2}_{\rm loc}(B_1(0),\mathcal{A}_q(\mathbb{R}^m))$ is such that \eqref{blowup eqn9} holds true.  Fix $\rho \in (0,1)$ and let $\widehat{P}_k$ be an optimal plane for $T_k$ in $\mathbf{B}_{\rho}(0)$.  Then there exists a sequence of $m \times n$ matrices $A_k$ such that 
\begin{gather}
	\label{optimal plane concl1} \widehat{P}_k = \{ (x,A_k x) : x \in \mathbb{R}^n \} , \quad \|A_k\| \leq C(n,m,q,\rho) \,E_k, \\
	\label{optimal plane concl2} Dw_a(0) = \lim_{k \rightarrow \infty} A_k/E_k . 
\end{gather}
In particular, in the special case that $\widehat{P}_k = P_0$ for all $k$, we have that $Dw_a(0) = 0$.
\end{lemma}

\begin{proof}
By the triangle inequality, Lemma~\ref{tilt to height estimate lemma}, \eqref{blowup hyp1}, and \eqref{optimal plane eqn1} 
\begin{align*}
	\|\pi_{\widehat{P}_k} - \pi_{P_0}\|^2 
	\leq\,& \frac{2}{q \omega_n (\rho/2)^n} \int_{G_n(\mathbf{B}_{\rho/2}(0))} \|\pi_S - \pi_{\widehat{P}_k}\|^2 \,d|T_k|(X,S) 
		\\&+ \frac{2}{q \omega_n (\rho/2)^n} \int_{G_n(\mathbf{B}_{\rho/2}(0))} \|\pi_S - \pi_{P_0}\|^2 \,d|T_k|(X,S) \nonumber 
	\\ \leq\,& C E(T_k,\widehat{P}_k,\mathbf{B}_{\rho}(0))^2 + C E(T_k,P_0,\mathbf{B}_{\rho}(0))^2  
	\leq 2C E_k^2,  \nonumber 
\end{align*}
where $C = C(n,m,q,\rho) \in (0,\infty)$ is a constant.  Hence for each sufficiently large $k$, there exists an $m \times n$ matrix $A_k$ such that \eqref{optimal plane concl1} holds true. 

By \eqref{optimal plane concl1}, after passing to a subsequence there is an $m \times n$ matrix $\Lambda$ (which a priori depends on $\rho$ since 
$A_{k}$ depend on $\rho$) such that $A_k/E_k \rightarrow \Lambda$.  Since this choice of subsequence is arbitrary, it suffices to show that $\Lambda = Dw_a(0)$.  Let $M$ be any $m \times n$ matrix and let $S_{k}$ be the $n$-dimensional plane defined by $S_k = \{ (x, E_k M x) : x \in \mathbb{R}^n \}$.  Since $\widehat{P}_k$ is an optimal plane of $T_k$ in $\mathbf{B}_{\rho}(0)$, 
\begin{equation}\label{optimal plane eqn3}
	\int_{\mathbf{B}_{\rho}(0)} \op{dist}^2(X, \widehat{P}_k) \,d\|T_k\|(X) \leq \int_{\mathbf{B}_{\rho}(0)} \op{dist}^2(X, S_k) \,d\|T_k\|(X)
\end{equation}
for all $k$.  Take any point $X = (x,y) \in \mathbf{B}_1(0)$.  Since $(x,E_k M x) \in S_k$, clearly 
\begin{equation*}
	\op{dist}(X, S_k) \leq |(x,y) - (x,E_k M x)| = |y - E_k M x|.
\end{equation*}
On the other hand, if $x' \in \mathbb{R}^n$ is such that $|(x,y) - (x',E_k M x')| = \op{dist}(X, S_k)$, then clearly $|x-x'| \leq \op{dist}(X, S_k)$ and $|y-E_k M x'| \leq \op{dist}(X, S_k)$ and thus 
\begin{equation*}
	|y - E_k M x| \leq |y - E_k M x'| + E_k \|M\| \,|x-x'| \leq (1+E_k \|M\|) \op{dist}(X, S_k) .
\end{equation*}
By similar reasoning using \eqref{optimal plane concl1}, 
\begin{equation*}
	\op{dist}(X,  \widehat{P}_k) \leq |y - A_k x| \leq (1+CE_k) \op{dist}(X,  \widehat{P}_k) , 
\end{equation*}
where $C = C(n,m, q, \rho) \in (0,\infty)$ is a constant.  Hence dividing both sides of \eqref{optimal plane eqn3} by $E_k^2$ and letting $k \rightarrow \infty$ using \eqref{blowup eqn4}--\eqref{blowup eqn6}, \eqref{blowup eqn9}, \eqref{lipschitz approx eqn8}, and \eqref{lipschitz approx eqn9}, 
\begin{align*}
	\int_{B_{\rho}(0)} \mathcal{G}(w, q \llbracket \Lambda x \rrbracket)^2 
	&= \lim_{k \rightarrow \infty} \frac{1}{E_k^2} \int_{B_{\rho}(0)} \mathcal{G}(u_k, q \llbracket A_k x \rrbracket)^2  
	= \lim_{k \rightarrow \infty} \frac{1}{E_k^2} \int_{\mathbf{B}_{\rho}(0)} \op{dist}^2(X,\widehat{P}_k) \,d\|T_k\|(X)  
	\\&\leq \lim_{k \rightarrow \infty} \frac{1}{E_k^2} \int_{\mathbf{B}_{\rho}(0)} \op{dist}^2(X,S_k) \,d\|T_k\|(X) 
	= \lim_{k \rightarrow \infty} \frac{1}{E_k^2} \int_{B_{\rho}(0)} \mathcal{G}(u_k, q \llbracket E_k M x \rrbracket)^2  
	\\&= \int_{B_{\rho}(0)} \mathcal{G}(w, q \llbracket Mx \rrbracket)^2 .
\end{align*}
Hence by \eqref{G avg sym}, 
\begin{equation*}
	\int_{B_{\rho}(0)} |w_a - \Lambda x|^2 \leq \int_{B_{\rho}(0)} |w_a - Mx|^2 
\end{equation*}
for all $m \times n$ matrices $M$.  Since $w_a$ is harmonic, it follows (using the series expansion for $w_a$ in terms of $L^{2}$ orthogonal spherical harmonics) that $\Lambda = Dw_a(0)$.  Finally, we note that in the special case $\widehat{P}_k = P_0$ for all $k$, clearly $A_k = 0$ for all $k$ and thus $Dw_a(0) = \Lambda = 0$.
\end{proof}

\subsection{Hardt-Simon inequality}  Here for completion we prove the Hardt-Simon inequality, which was introduced in~\cite{HardtSimon}. 

\begin{lemma}
There exists $\eta = \eta(n) \in (0,1)$ such that if $P$ is an $n$-dimensional plane in $\mathbb{R}^n$ and $T$ is an $n$-dimensional locally area-minimizing rectifiable current in $\mathbf{B}_1(0)$ such that 
\begin{gather}
	\label{hardt simon hyp1} (\partial T) \llcorner \mathbf{B}_1(0) = 0, \quad\quad \Theta(T,0) \geq q , \\
	\label{hardt simon hyp2} E(T,P,\mathbf{B}_1(0)) < \eta , 
\end{gather}
and $Z \in \mathbf{B}_{1/2}(0)$ with $\Theta(T,Z) \geq q$, then 
\begin{gather}
	\label{hardt simon concl1} \op{dist}(Z,P) \leq C E(T,P,\mathbf{B}_1(0)), \;\; \mbox{where} \;\; C = C(n, m, q), \;\; \mbox{and}\\
	\label{hardt simon concl2} \int_{G_n(\mathbf{B}_{1/4}(Z))} \frac{|\pi_{S^{\perp}}(X-Z)|^2}{|X-Z|^{n+2}} \,d|T|(X,S) \leq C E(T,P,\mathbf{B}_1(0))^2
\end{gather}
where $C= C(n,m) \in (0,\infty)$. 
\end{lemma}

\begin{proof}
Without loss of generality assume $P = \mathbb{R}^n \times \{0\}$ and let $Z$ be as in the lemma.    By the monotonicity formula for area (\cite[17.5]{SimonGMT}) and the assumption $\Theta(T,Z) \geq q$, 
\begin{equation}\label{hardt simon eqn1}
	\frac{1}{\omega_n} \int_{G_n(\mathbf{B}_{1/4}(Z))} \frac{|\pi_{S^{\perp}}(X-Z)|^2}{|X-Z|^{n+2}} \,d|T|(X,S) 
	\leq \frac{\|T\|(\mathbf{B}_{1/4}(Z))}{\omega_n (1/4)^n} - q . 
\end{equation}
Recall that by \eqref{hardt simon hyp1}, \eqref{hardt simon hyp2}, and Lemma~\ref{coarse L2 distance lemma}, we have that $\op{dist}(X,P) < 2\eta^{\frac{2}{n+2}}$ for all $X \in \op{spt} T \cap \mathbf{B}_{15/16}(0)$ and we may assume that $\pi_{P\#} (T \llcorner {\mathbf C}_{7/8}(0)) 
= q \llbracket P_0 \rrbracket \llcorner {\mathbf C}_{7/8}(0)$.  By Lemma~\ref{area excess to height lemma} we can bound the right-hand side of \eqref{hardt simon eqn1} 
\begin{align}\label{hardt simon eqn2}
	\frac{\|T\|(\mathbf{B}_{1/4}(Z))}{\omega_n (1/4)^n} - q 
	&= \frac{1}{2\omega_n (1/4)^n} \int_{\mathbf{B}_{1/4}(Z)} |\vec T - \vec P|^2 \,d\|T\|
	\\&\leq C \int_{\mathbf{B}_{3/8}(Z)} \op{dist}^2(X,P) \,d\|T\|(X) 
	\leq C E(T,P,\mathbf{B}_1(0))^2 \nonumber
\end{align}
for some constant $C = C(n,m) \in (0,\infty)$.  Combining \eqref{hardt simon eqn1} and \eqref{hardt simon eqn2} gives us \eqref{hardt simon concl2}.  On the other hand, writing $Z = (\xi,\zeta)$ where $\xi \in \mathbb{R}^n$ and $\zeta \in \mathbb{R}^m$, for each $X = (x,y) \in \mathbb{R}^n$ and each $n$-dimensional plane $S$ 
\begin{align*}
	\pi_{S^{\perp}}(X-Z) &= \pi_{S^{\perp}}(x-\xi,y-\zeta) = \pi_{S^{\perp}}(x-\xi,0) + \pi_{S^{\perp}}(0,y) - \pi_{S^{\perp}}(0,\zeta)
		\\&= (\pi_P - \pi_S)(x-\xi,0) + \pi_{S^{\perp}}(0,y) - (0,\zeta) - (\pi_P - \pi_S)(0,\zeta) 
\end{align*} 
and thus by the triangle inequality 
\begin{equation*}
	|\pi_{S^{\perp}}(X-Z)| \geq |\zeta| - |y| - 2 \|\pi_S - \pi_P\| . 
\end{equation*}
Hence using Lemma~\ref{area excess to height lemma} and the assumption $\Theta(T,Z) \geq q$, we can bound the left-hand side of \eqref{hardt simon eqn1} by 
\begin{align}\label{hardt simon eqn3}
	&\frac{1}{\omega_n} \int_{G_n(\mathbf{B}_{1/4}(Z))} \frac{|\pi_{S^{\perp}}(X-Z)|^2}{|X-Z|^{n+2}} \,d|T|(X,S) 
		\geq \frac{4^{n+2}}{\omega_n} \int_{\mathbf{B}_{1/4}(Z)} |\pi_{S^{\perp}}(X-Z)|^2 \,d\|T\|(X)  
		\\ \geq\,& \frac{4^{n+1}}{\omega_n} \int_{\mathbf{B}_{1/4}(Z)} |\zeta|^2 \,d\|T\|(X) 
			- \frac{4^{n+2}}{\omega_n} \int_{\mathbf{B}_{1/4}(Z)} |y|^2 \,d\|T\|(X) \nonumber
			\\&- \frac{2 \cdot 4^{n+2}}{\omega_n} \int_{G_n(\mathbf{B}_{1/4}(Z))} \|\pi_S - \pi_P\|^2 \,d|T|(X,S) \nonumber
		\\ \geq\,& |\zeta|^2 \frac{\|T\|(\mathbf{B}_{1/4}(Z))}{\omega_n (1/4)^n} - CE(T,P,\mathbf{B}_1(0)) \nonumber
		\\ \geq\,& q \,|\zeta|^2 - CE(T,P,\mathbf{B}_1(0)) , \nonumber
\end{align}
where $C = C(n,m) \in (0,\infty)$ is a constant.  Combining \eqref{hardt simon concl2} and \eqref{hardt simon eqn3} gives us \eqref{hardt simon concl1}. 
\end{proof}

\begin{lemma}\label{hardt simon lemma} 
For $k = 1,2,3,\ldots$ let $T_k$ be $n$-dimensional locally area-minimizing rectifiable currents in $\mathbf{C}_1(0)$ such that \eqref{blowup hyp1}, \eqref{blowup hyp2}, and \eqref{blowup eqn2} hold true with $P = P_0$ and $\rho_{0} = 1$.  Let $E_k = E(T_k,P_0,\mathbf{C}_1(0))$ and $w$ be a blow-up of $T_k$ relative to $P_0$ in $\mathbf{C}_1(0)$ (as in Subsection~\ref{sec:blowup procedure}).  Let $Z_k = (\xi_k,\zeta_k) \in \op{spt} T_k \cap \mathbf{C}_{1/2}(0)$ and $\xi \in \overline{B_{1/2}(0)}$ be such that $\Theta(T_k,Z_k) \geq q$ and $\xi_k \rightarrow \xi$.  Then 
\begin{enumerate}
	\item[{\rm (i)}]  $\zeta_k/E_k \rightarrow w_a(\xi)$ as $k \rightarrow \infty$;
	\item[{\rm (ii)}]  for each $\rho \in (0,1/4]$ 
	\begin{equation}\label{hardt simon concl3} 
		\int_{B_{\rho/4}(\xi)} r^{2-n} \left|\frac{\partial}{\partial r}\bigg(\frac{w - w_a(\xi)}{r}\bigg)\right|^2 \leq C \int_{B_{\rho}(\xi)} |w - w_a(\xi)|^2 , 
	\end{equation}
	where $r = |x-\xi|$, $w(x) - w_a(\xi) = \sum_{i=1}^q \llbracket w_i(x) - w_a(\xi) \rrbracket$ for each $x \in B_1(0)$, and 
	$C = C(n,m,q) \in (0,\infty)$ is a constant;
	\item[{\rm (iii)}]  Assuming $w \not\equiv q \llbracket w_a(\xi) \rrbracket$ in $B_1(0)$, $\mathcal{N}_{w - w_a(\xi)}(\xi) \geq 1$. 
\end{enumerate}
\end{lemma}

\begin{proof}
Applying \eqref{hardt simon concl2} to $\eta_{Z_k,\rho/2 \,\#} T_k$ we obtain 
\begin{equation}\label{hardt simon eqn4}
	\int_{G_n(\mathbf{B}_{\rho/4}(Z_k))} \frac{|\pi_{S^{\perp}}(X-Z_k)|^2}{|X-Z_k|^{n+2}} \,d|T_k|(X,S) 
		\leq C \int_{\mathbf{B}_{\rho}(Z_k)} \op{dist}^2(X,Z_k+P) \,d\|T_k\|(X)
\end{equation}
for each $\rho \in (0,1/4]$.  Let $u_k$ be as in Subsection~\ref{sec:blowup procedure} with $\rho_0 = 1$.  Let $u_k(x) - \zeta_k = \sum_{i=1}^q \llbracket u_{k,i}(x) - \zeta_k \rrbracket$ for each $x \in B_{\theta_k}(0)$.  By~\cite[Theorem~1.13]{DeLSpaDirMin}, $u_k$ is differentiable in the sense of~\cite[Definition~1.9]{DeLSpaDirMin} at $\mathcal{L}^n$-a.e.~$x \in B_{3/4}(0)$.  Hence given $x \in B_{3/4}(0)$ at which $u_k$ is differentiable, $\frac{\partial}{\partial r_k} (x,u_{k,i}(x))$ is tangent to the graph of $u_k$ at $X = (x,u_{k,i}(x))$, where $r_k = |x-\xi_k|$.  It follows that  
\begin{equation*}
	|\pi_{S_X^{\perp}}(X-Z_k)| = |\pi_{S_X^{\perp}}(x-\xi_k, u_{k,i}(x)-\zeta_k)| 
		= r_k^2 \left|\pi_{S_X^{\perp}}\bigg(\frac{\partial}{\partial r_k} \bigg(\frac{(x-\xi_k,u_{k,i}(x)-\zeta_k)}{r_k}\bigg)\bigg)\right| , 
\end{equation*}
where $S_X$ is the approximate tangent plane to the graph of $u_k$ at $X = (x,u_{k,i}(x))$.  Since $\|\pi_{S_X} - \pi_P\| \leq C(n,m) \op{Lip} u_k \leq C(n,m,q) E_k^{2\gamma}$ is small, 
\begin{equation}\label{hardt simon eqn5}
	|\pi_{S_X^{\perp}}(X-Z_k)| \geq \frac{1}{2} \,r_k^2 \left|\frac{\partial}{\partial r_k} \bigg(\frac{(x-\xi_k,u_{k,i}(x)-\zeta_k)}{r_k}\bigg)\right|
		= \frac{1}{2} \,r_k^2 \left|\frac{\partial}{\partial r_k} \bigg(\frac{u_{k,i}(x)-\zeta_k}{r_k}\bigg)\right|. 
\end{equation}
Hence using \eqref{hardt simon eqn5} to bound the left-hand side of \eqref{hardt simon eqn4} and using \eqref{blowup eqn4}--\eqref{blowup eqn6}, 
\begin{align}\label{hardt simon eqn6}
	&\int_{\mathbf{C}_{\rho/4}(\xi_k)} \sum_{i=1}^q \frac{r_k^4}{(r_k^2 + |u_{k,i}(x)-\zeta_k|^2)^{(n+2)/2}} 
		\left| \frac{\partial}{\partial r_k} \bigg(\frac{u_{k,i}(x)-\zeta_k}{r_k}\bigg) \right|^2  
	\\ \leq\,& C \int_{\mathbf{B}_{\rho}(\xi_k)} |u_k-\zeta_k|^2 + C E_k^{2+\gamma} \nonumber
\end{align}
for each $\rho \in (0,1/4]$.  Since $\xi_k \rightarrow \xi$, $r_k \rightarrow r$ uniformly in $B_{3/4}(0)$.  By \eqref{hardt simon concl1}, after passing to a subsequence there exists $\lambda \in \mathbb{R}^m$ such that $\zeta_k/E_k \rightarrow \lambda$.  Let $N(m,q) \geq 1$ be an integer and $\boldsymbol{\xi} : \mathcal{A}_q(\mathbb{R}^m) \rightarrow \mathbb{R}^N$ be the bi-Lipschitz injection as in~\cite[Corollary~2.2]{DeLSpaDirMin}.  By~\cite[Corollary~2.2]{DeLSpaDirMin}, \eqref{blowup eqn6}, \eqref{blowup eqn8}, and \eqref{hardt simon concl1}, for each $\delta \in (0,1/4)$ 
\begin{equation*}
	\left\| \boldsymbol{\xi} \circ \left(\frac{u_{k,i}-\zeta_k}{r_k E_k}\right) \right\|_{W^{1,2}(B_{3/4}(0) \setminus B_{\delta}(\xi))}
	\leq \frac{\|u_{k,i}\|_{W^{1,2}(B_{3/4}(0))} + \omega_{n}^{1/2}(3/4)^{n/2}|\zeta_k|}{\delta E_k} \leq \frac{C(n,m,q)}{\delta} \,. 
\end{equation*}
Thus using the Rellich compactness theorem (Lemma~1 in Section~1.3 of~\cite{Sim96}) and recalling \eqref{blowup eqn9} and $\zeta_k/E_k \rightarrow \lambda$, after passing to a subsequence 
\begin{equation*}
	\boldsymbol{\xi} \circ \left(\frac{u_{k,i}(x)-\zeta_k}{r_k E_k}\right) \rightarrow \boldsymbol{\xi} \circ \left(\frac{w_i(x)-\lambda}{r}\right) 
\end{equation*}
strongly in $L^2(B_{3/4}(0) \setminus B_{\delta}(\xi),\mathbb{R}^N)$ and weakly in $W^{1,2}(B_{3/4}(0) \setminus B_{\delta}(\xi),\mathbb{R}^N)$ for each $\delta \in (0,1/4)$.  It follows that 
\begin{equation*}
	r_k^{(2-n)/2} \frac{\partial}{\partial r_k}\left(\boldsymbol{\xi} \circ \left(\frac{u_{k,i}(x)-\zeta_k}{r_k E_k}\right)\right) 
		\rightarrow r^{(2-n)/2} \frac{\partial}{\partial r}\left(\boldsymbol{\xi} \circ \left(\frac{w_i(x)-\lambda}{r}\right)\right) 
\end{equation*}
weakly in $L^2(B_{3/4}(0) \setminus B_{\delta}(\xi),\mathbb{R}^N)$ for each $\delta \in (0,1/4)$.  Hence for each $\rho \in (0,1]$ and $\delta \in (0,\rho/16]$ 
\begin{equation*}
	\int_{B_{\rho/8}(\xi) \setminus B_{\delta}(\xi)} r^{2-n} 
		\left|\frac{\partial}{\partial r}\bigg( \boldsymbol{\xi} \circ \bigg(\frac{w - \lambda}{r}\bigg)\bigg)\right|^2  
	\leq \liminf_{k \rightarrow \infty} \int_{B_{\rho/8}(\xi_k) \setminus B_{\delta}(\xi)} r^{2-n} 
		\left| \frac{\partial}{\partial r}\bigg( \boldsymbol{\xi} \circ \bigg(\frac{u_k-\zeta_k}{r_k E_k}\bigg)\bigg) \right|^2 , 
\end{equation*}
or equivalently using~\cite[Corollary~2.2]{DeLSpaDirMin} 
\begin{align}\label{hardt simon eqn7}
	\int_{B_{\rho/8}(\xi) \setminus B_{\delta}(\xi)} r^{2-n} \left|\frac{\partial}{\partial r}\bigg(\frac{w - \lambda}{r}\bigg)\right|^2 
	\leq \liminf_{k \rightarrow \infty} \frac{1}{E_k^2} \int_{B_{\rho/8}(\xi_k) \setminus B_{\delta}(\xi)} r_k^{2-n}
		 \left| \frac{\partial}{\partial r_k} \bigg(\frac{u_k-\zeta_k}{r_k}\bigg) \right|^2 .
\end{align}
Dividing both sides of \eqref{hardt simon eqn6} by $E_k$ and letting $k \rightarrow \infty$ using \eqref{blowup eqn9}, \eqref{hardt simon eqn7}, and $\zeta_k/E_k \rightarrow \lambda$, 
\begin{equation*}
	\int_{B_{\rho/8}(\xi) \setminus B_{\delta}(\xi)} r^{2-n} \left|\frac{\partial}{\partial r}\bigg(\frac{w - \lambda}{r}\bigg)\right|^2 
	= C \int_{B_{\rho}(\xi)} |w - \lambda|^2 
\end{equation*}
for each $\rho \in (0,1]$ and $\delta \in (0,\rho/16]$ and for some constant $C = C(n,m,q) \in (0,\infty)$.  Letting $\delta \rightarrow 0^+$, 
\begin{equation}\label{hardt simon eqn8}
	\int_{B_{\rho/8}(\xi)} r^{2-n} \left|\frac{\partial}{\partial r}\bigg(\frac{w - \lambda}{r}\bigg)\right|^2 \leq C \int_{B_{\rho}(\xi)} |w - \lambda|^2 
\end{equation}
for each $\rho \in (0,1]$ and for some constant $C = C(n,m,q) \in (0,\infty)$.  By \eqref{G avg sym},  
\begin{equation}\label{hardt simon eqn9}
	\int_{B_{\rho/8}(\xi)} r^{2-n} \left|\frac{\partial}{\partial r}\bigg(\frac{w_a - \lambda}{r}\bigg)\right|^2 
	\leq C \int_{B_{\rho}(\xi)} |w - \lambda|^2 
\end{equation}
for each $\rho \in (0,1]$ and for some constant $C = C(n,m,q) \in (0,\infty)$.  Since $w_a$ is harmonic, it follows using the series expansion of $w_a(x-\xi)$ that in order for \eqref{hardt simon eqn9} to hold true we must have that $\lambda = w_a(\xi)$.  This together with \eqref{hardt simon eqn8} completes the proof of conclusions~(i) and (ii). 

To show conclusion~(iii), suppose that $w \not\equiv q \llbracket w_a(\xi) \rrbracket$ in $B_1(0).$  Let $\varphi$ be a tangent function of $w - w_a(\xi)$ at $\xi$; that is, for some sequence of radii $\rho_j \rightarrow 0$, 
\begin{equation*} 
	w_j(x) = \frac{w(\xi + \rho_j x) - w_a(\xi)}{\|w - w_a(\xi)\|_{L^2(B_{\rho_j}(\xi))}} \rightarrow \varphi(x)
\end{equation*}
uniformly on compact subsets of $\mathbb{R}^n$ as $j \rightarrow \infty$.  By rescaling \eqref{hardt simon concl3} with $\rho = \rho_j$, 
\begin{equation*} 
	\int_{B_{1/8}(0)} r^{2-n} \left|\frac{\partial}{\partial r}\bigg(\frac{w_j}{r}\bigg)\right|^2 \leq C 
\end{equation*}
for each $j$.  Hence letting $j \rightarrow \infty$, 
\begin{equation}\label{hardt simon eqn11}
	\int_{B_{1/8}(0)} r^{2-n} \left|\frac{\partial}{\partial r}\bigg(\frac{\varphi}{r}\bigg)\right|^2 \leq C . 
\end{equation}
Since $\varphi$ is homogeneous of degree $\mathcal{N}_{w - w_a(\xi)}(\xi)$, in order for \eqref{hardt simon eqn11} to hold true we must have that $\mathcal{N}_{w - w_a(\xi)}(\xi) \geq 1$.
\end{proof}

\subsection{Compatibility of the planar frequency function with blowing up of currents}  

\begin{lemma}\label{blowup frequency lemma}
For $k = 1,2,3,\ldots$, let $T_k$ be $n$-dimensional locally area-minimizing rectifiable currents in $\mathbf{C}_1(0)$ such that \eqref{blowup hyp1}, \eqref{blowup hyp2}, and \eqref{blowup eqn2} hold true with $P = P_0$ and $\rho_{0} = 1$ and 
\begin{equation}\label{blowup frequency hyp1}
	\sup_{k \geq 1} \sup_{X \in \op{spt} T_k} \op{dist}(X,P_0) < \infty . 
\end{equation}
Let $E_k = E(T_k,P_0,\mathbf{C}_1(0))$ and $w$ be a blow-up of $T_k$ relative to $P_0$ in $\mathbf{C}_1(0)$ (as in Subsection~\ref{sec:blowup procedure}).  For $k = 1,2,3,\ldots$, let $A_{k}$ be an $m \times n$ matrix, 
\begin{equation}\label{blowup frequency hyp2}
	P_k = \{(x,A_k x) : x \in \mathbb{R}^n \}, \quad \mbox{and suppose that} \quad A_k/E_k \rightarrow \Lambda 
\end{equation}
for some $m \times n$ matrix $\Lambda$.  Let $Z_k = (\xi_k,\zeta_k) \in \mathbf{C}_1(0)$, $\xi \in B_1(0)$, and $\rho_k,\rho \in (0,1-|\xi|)$ be such that $\xi_k \rightarrow \xi$ and $\rho_k \rightarrow \rho$.  Let $\ell(x) = w_a(\xi) + \Lambda (x-\xi)$, and $(w-\ell)(x) = \sum_{i=1}^q \llbracket w_i(x) - \ell(x) \rrbracket$ for each $x \in B_1(0)$.  Assume that $w \not\equiv q \llbracket \ell \rrbracket$ on $B_1(0)$.  Then $w - \ell$ is locally Dirichlet energy minimizing and 
\begin{equation}\label{blowup frequency concl1}
	\lim_{k \rightarrow \infty} N_{T_k,P_k,Z_k}(\rho_k) = N_{w - \ell,\xi}(\rho)
\end{equation}
where $N_{w-\ell, \xi}$ is the Almgren frequency function (defined in Section~\ref{Dir-min-fns}). 
\end{lemma}

\begin{proof}
The fact that $w - \ell$ is locally Dirichlet energy minimizing follows from \cite[Theorem~2.6(3)]{Almgren}.  By \eqref{blowup frequency hyp1} and \eqref{blowup frequency hyp2}, $\op{spt} T_k \cap \overline{\mathbf{C}_{\rho_k}(Z_k,P_k)} \subset \mathbf{C}_1(0)$ and thus 
\begin{equation*}
	(\partial T_k) \llcorner \mathbf{C}_{\rho_k}(Z_k,P_k) = 0 , \quad\quad
	\sup_{X \in \op{spt} T_k \cap \mathbf{C}_{\rho_k}(Z_k,P_k)} \op{dist}(X,Z_k+P_k) < \infty . 
\end{equation*}
Hence $H_{T_k,P_k,Z_k}(\rho_k)$ and $D_{T_k,P_k,Z_k}(\rho_k)$ are defined and $N_{T_k,P_k,Z_k}(\rho_k)$ is defined if $H_{T_k,P_k,Z_k}(\rho_k) > 0$.  Let $u_k,$ $\theta_{k}$ be as in Subsection~\ref{sec:blowup procedure} with $\rho_0 = 1$.  Let $\ell_k(x) = \zeta_k + A_k (x-\xi_k)$ and $(u_k-\ell_k)(x) = \sum_{i=1}^q \llbracket u_{k,i}(x) - \ell_k(x) \rrbracket$ for each $x \in B_{\theta_k}(0)$.  Let $r_k(X) = |\pi_{P_k}(X-Z_k)|$ for each $X \in \mathbf{C}_1(0)$.  Let $r(x) = |x-\xi|$ for all $x \in B_1(0)$.  Extend $r$ to a function $r(x,y)$ of $(x,y) \in \mathbf{C}_1(0)$ which is independent of $y$.  
Notice that $\xi_k \rightarrow \xi$, $r_k \rightarrow r$ uniformly on compact subsets of $\mathbf{C}_1(0)$, $\phi'(r_k/\rho_k) \rightarrow \phi'(r/\rho)$ pointwise on $B_1(0) \setminus (\partial B_{\rho/2}(\xi) \cup \partial B_{\rho}(\xi))$, and $\nabla r_k \rightarrow \nabla r$ uniformly on compact subsets $\mathbf{C}_1(0) \setminus \mathbf{C}_{\rho/4}(\xi)$.  By Lemma~\ref{hardt simon lemma}, $\zeta_k/E_k \rightarrow w_a(\xi)$ and thus by \eqref{blowup frequency hyp2} we have that $\ell_k/E_k \rightarrow \ell$ uniformly in $B_1(0)$.  By \eqref{blowup eqn6} and \eqref{lipschitz approx eqn8}, for each $\theta \in (0,\theta_k]$ and for $\mathcal{H}^n$-a.e.~$X = (x,u_{k,l}(x)) \in \op{spt} \op{graph} u_k \cap \mathbf{C}_{\theta}(0)$ we have that $1 - C(n,m,q,\theta) \,E_k^{\gamma} \leq |\nabla^{S_X} r_k(X)| \leq 1$, where $S_X$ is the approximate tangent plane to $\op{spt} \op{graph} u_k$ at $X$.  
Hence using \eqref{blowup eqn9}, \eqref{blowup eqn4}--\eqref{blowup eqn6}, \eqref{lipschitz approx eqn8}, and \eqref{lipschitz approx eqn9}, 
\begin{align}\label{blowup frequency eqn1}
	H_{w - \ell,\xi}(\rho) &= -\rho^{1-n} \int |w - \ell|^2 \,\frac{1}{r} \,\phi'(r/\rho)
	\\&= \lim_{k \rightarrow \infty} -\frac{\rho_{k}^{1-n}}{E_k^2} \int \sum_{l=1}^q |u_{k,l} - \ell_k|^2 \,\frac{1}{r_k(x,u_{k,l}(x))} 
		\,\phi'\bigg(\frac{r_k(x,u_{k,l}(x))}{\rho_k}\bigg) \nonumber
	\\&= \lim_{k \rightarrow \infty} -\frac{\rho_{k}^{1-n}}{E_k^2} \int \op{dist}^2(X,Z_k+P_k) \,|\nabla^S r_k|^2 \,\frac{1}{r_k} \,\phi'(r_k/\rho_k) \,d|T_k|(X,S) 
		\nonumber
	\\&= \lim_{k \rightarrow \infty} \frac{H_{T_k,P_k,Z_k}(\rho_k)}{E_k^2} . \nonumber
\end{align}
By \eqref{blowup eqn6}, \eqref{blowup eqn8}, \eqref{blowup eqn9}, \eqref{G avg sym}, and Rellich's compactness theorem, after passing to a subsequence $u_{k,a}/E_k \rightarrow w_a$ weakly in $W^{1,2}(B_{\theta}(0),\mathbb{R}^m)$ as $k \rightarrow \infty$ for each $\theta \in (0,1)$.  Thus using \eqref{blowup eqn10}, \eqref{blowup frequency hyp2}, \eqref{blowup eqn4}--\eqref{blowup eqn8}, \eqref{lipschitz approx eqn8}, \eqref{lipschitz approx eqn9}, and \eqref{lipschitz approx eqn11} 
\begin{align}\label{blowup frequency eqn2}
	D_{w - \ell,\xi}(\rho) &= \rho^{2-n} \int |Dw-\Lambda|^2 \,\phi(r/\rho) 
	\\&= \rho^{2-n} \int (|Dw|^2 - 2q \,Dw_a \cdot \Lambda + q|\Lambda|^2) \,\phi(r/\rho) \nonumber
	\\&= \lim_{k \rightarrow \infty} \frac{\rho_k^{2-n}}{E_k^2} \int \sum_{l=1}^q (|Du_{k,l}|^2 - 2 Du_{k,l} \cdot A_k + |A_k|^2) 
		\,\phi\bigg(\frac{r_k(x,u_{k,l}(x))}{\rho_k}\bigg) \nonumber
	\\&= \lim_{k \rightarrow \infty} \frac{\rho_k^{2-n}}{E_k^2} \int |Du_k - A_k|^2 \,\phi\bigg(\frac{r_k(x,u_{k,l}(x))}{\rho_k}\bigg) \nonumber
	\\&= \lim_{k \rightarrow \infty} \frac{\rho_k^{2-n}}{2E_k^2} \int \|\pi_S - \pi_{P_k}\|^2 \,\phi(r_{k}/\rho_{k}) \,d|T_k|(X,S) \nonumber
	\\&= \lim_{k \rightarrow \infty} \frac{D_{T_k,P_k,Z_k}(\rho_k)}{E_k^2} . \nonumber
\end{align}
(where $Du_k(x) - A_k = \sum_{i=1}^q \llbracket Du_{k,i}(x) - A_k \rrbracket$ for $\mathcal{H}^n$-a.e.~$x \in B_{\rho_k}(Z_k,P_k)$ and $Dw(x) - \Lambda = \sum_{i=1}^q \llbracket Dw_i(x) - \Lambda \rrbracket$ for $\mathcal{L}^n$-a.e.~$x \in B_{\rho}(\xi)$).  By dividing \eqref{blowup frequency eqn2} by \eqref{blowup frequency eqn1} and noting that since $w - \ell$ is non-zero we have $H_{w - \ell,\xi}(\rho) > 0$, we obtain \eqref{blowup frequency concl1}. 
\end{proof}

\begin{corollary}\label{blowup semicont freq cor}
There exists $\eta = \eta(n,m,q) \in (0,1)$ such that the following holds true.  For $k = 1,2,3,\ldots$ let $T_k$ be $n$-dimensional locally area-minimizing rectifiable currents in $\mathbf{C}_3(0)$ such that \eqref{blowup hyp1}, \eqref{blowup hyp2}, and \eqref{blowup eqn2} hold true with $P = P_0$ and $\rho_{0} = 3$ and 
\begin{equation}\label{blowup semicont freq hyp1}
	\sup_{k \geq 1} \sup_{X \in \op{spt} T_k} \op{dist}(X,P_0) < \infty . 
\end{equation}
Let $E_k = E(T_k,P_0,\mathbf{C}_3(0))$ and let $w$ be a blow-up of $T_k$ relative to $P$ in $\mathbf{C}_3(0)$ (as in Subsection~\ref{sec:blowup procedure}).  Let $Z_k = (\xi_k,\zeta_k) \in \op{spt} T_k \cap \mathbf{C}_1(0)$ such that $\Theta(T_k,Z_k) \geq q$ and suppose that for some $\alpha > 0$ and for some $n$-dimensional plane $P_k$ in $\mathbb{R}^{n+m},$ 
\begin{equation}\label{blowup semicont freq hyp2} 
	E(T_k,P_k,\mathbf{C}_{7\rho/4}(Z_k, P_{k})) \leq \eta \rho^{\alpha} \text{ for all } \rho \in (0,1] . 
\end{equation}
Let $\xi_k \rightarrow \xi \in {\mathbb R}^{n}$.  Let $\ell_{\xi}(x) = w_a(\xi)-Dw_a(\xi)\cdot (x-\xi)$ and $w(x)-\ell_{\xi}(x) = \sum_{i=1}^q \llbracket w_i(x)-\ell_{\xi}(x) \rrbracket$ for each $x \in B_3(0)$.  Assume that $w \not\equiv q \llbracket \ell_{\xi} \rrbracket$ on $B_3(0)$.  Then $w - \ell_{\xi}$ is locally Dirichlet energy minimizing and 
\begin{equation}\label{blowup semicont freq concl}
	\mathcal{N}_{w-\ell_{\xi}}(\xi) \geq \lim_{k \rightarrow \infty} \mathcal{N}_{T_k,{\rm Pl}}(Z_k) .
\end{equation}
\end{corollary}

\begin{proof}
The fact that $w - \ell_{\xi}$ is locally Dirichlet energy minimizing follows from \cite[Theorem~2.6(3)]{Almgren}.  
We claim that 
\begin{equation}\label{blowup semicont freq eqn1}
	\lim_{k \rightarrow \infty} \op{dist}_{\mathcal H}(P_k \cap \mathbf{B}_1(0), P_0 \cap \mathbf{B}_1(0)) = 0 .
\end{equation}
Clearly there exists an $n$-dimensional linear plane $P_{\infty}$ of $\mathbb{R}^{n+m}$ such that $\op{dist}_{\mathcal H}(P_k \cap \mathbf{B}_1(0), P_{\infty} \cap \mathbf{B}_1(0)) \rightarrow 0$.  Recall from Subsection~\ref{sec:blowup procedure} that $T_k \rightarrow q \llbracket P_0 \rrbracket$ weakly in $\mathbf{C}_3(0)$.  By~\cite[Theorem~34.5]{SimonGMT}, $\|T_k\| \rightarrow q \,\mathcal{H}^n \llcorner P_0$ in the sense of Radon measures of $\mathbf{C}_3(0)$.  Hence by \eqref{blowup semicont freq hyp2} 
\begin{equation*} 
	E(q \llbracket P_0 \rrbracket,P_{\infty},\mathbf{C}_{7\rho/4}(Z,P_{\infty})) \leq \eta \rho^{\alpha} \text{ for all } \rho \in (0,1] . 
\end{equation*}
It follows that $P_{\infty} = P_0$, completing the proof of \eqref{blowup semicont freq eqn1}. 

Next we claim that there exists an $m \times n$ matrix $A_k$ such that 
\begin{gather}
	\label{blowup semicont freq eqn2} P_k = \{ (x,A_k x) : x \in \mathbb{R}^n \} , \quad \|A_k\| \leq C(n,m,q,\alpha) \,E_k, \;\; \rm{and}\\ 
	\label{blowup semicont freq eqn3} A_k/E_k \rightarrow Dw_a(\xi) . 
\end{gather}
Note that by \eqref{blowup semicont freq hyp1} and \eqref{blowup semicont freq eqn1}, $\op{spt} T_k \cap \overline{\mathbf{C}_{7/4}(Z_{k}, P_{k})} \subset \mathbf{C}_{23/8}(0,P_{0})$ and thus $(\partial T_k) \llcorner \mathbf{C}_{7/4}(Z_{k}, P_{k}) = 0$.  By Lemma~\ref{area excess to height lemma}, $\mathcal{E}(T_k,P_0,{\mathbf C}_{7/4}(Z_{k}, P_{k})) \leq (14/23)^{n/2} \mathcal{E}(T_k,P_0,{\mathbf C}_{23/8}(0,P_{0})) \leq C(n,m) \,E_k$ for all sufficiently large $k$.  Thus by \eqref{area excess eqn1} and \eqref{blowup eqn2}, $\|T_{k}\|({\mathbf C}_{7/4}(Z_{k}, P_{k})) \leq (q + C(n,m) \,E_k^2)\omega_{n}(7/4)^{n}$.  Hence by Corollary~\ref{doubling cor}, $\mathcal{N}_{T_k,{\rm Pl}}(Z_k) \geq 1+\alpha$.  For $\theta = \theta(n,m,q,\alpha) \in (0,1]$ sufficiently small, by Corollary~\ref{doubling cor}
\begin{align}\label{blowup semicont freq eqn4}
	&E(T_k,P_k,\mathbf{C}_{\rho}(Z_k, P_{k})) 
	\\ \leq\,& 4 e^{C (\mathcal{N}_{T_k,{\rm Pl}}(Z_k)+2) \,\eta^{2\gamma}/(4\alpha\gamma)} 
		\rho^{\mathcal{N}_{T_k,{\rm Pl}}(Z_k)-1} E(T_k,P_k,\mathbf{C}_1(Z_k, P_{k})) \nonumber
	\\ \leq\,& 4 e^{C(3+\alpha)\eta^{2\gamma}/(4\alpha\gamma)} (e^{C\eta^{2\gamma}/(4\alpha\gamma)} \theta)^{\mathcal{N}_{T_k,{\rm Pl}}(Z_k)-1-\alpha} 
		\rho^{\alpha} E(T_k,P_k,\mathbf{C}_1(Z_k, P_{k})) \nonumber
	\\ \leq\,& C\rho^{\alpha} E(T_k,P_k,\mathbf{C}_1(Z_k, P_{k})) \nonumber
\end{align}
for all $\rho \in (0,\theta]$, where $C = C(n,m,q,\alpha) \in (0,\infty)$ are constants.  By the triangle inequality, Lemma~\ref{tilt to height estimate lemma}, \eqref{blowup hyp1}, and \eqref{blowup semicont freq eqn4}
\begin{align}\label{blowup semicont freq eqn5}
	\|\pi_{P_k} - \pi_{P_0}\|^2 
	\leq\,& \frac{2}{q \omega_n \rho^n} \int_{G_n(\mathbf{C}_{\rho}(Z_k, P_{k}))} \|\pi_S - \pi_{P_k}\|^2 \,d|T_k|(X,S) 
		\\&+ \frac{2}{q \omega_n \rho^n} \int_{G_n(\mathbf{C}_{\rho}(Z_k, P_{k}))} \|\pi_S - \pi_{P_0}\|^2 \,d|T_k|(X,S) \nonumber 
	\\ \leq\,& C E(T_k,P_k,\mathbf{C}_{2\rho}(Z_k, P_{k}))^2 + C E(T_k,P_0,\mathbf{C}_{2\rho}(Z_k, P_{k}))^2 \nonumber 
	\\ \leq\,& C \rho^{2\alpha} E(T_k,P_k,\mathbf{C}_1(Z_k, P_{k}))^{2} + C\rho^{-n-2} E(T_k,P_0,\mathbf{B}_3(0))^2 \nonumber 
\end{align}
for all $\rho \in (0,\theta/2]$, where $C = C(n,m,q,\alpha) \in (0,\infty)$ are constants.  Noting that $| \op{dist}(X,Z_k+P_k) - \op{dist}(X,Z_k+P_0) | \leq \|\pi_{P_k} - \pi_{P_0}\| \,|X-Z_k|$ for all $X \in \mathbb{R}^{n+m}$ and using the triangle inequality, the fact that $\|T_{k}\|({\mathbf C}_{7/4}(Z_{k}, P_{k})) \leq (q + \delta)\omega_{n}(7/4)^{n}$, and \eqref{blowup semicont freq eqn5} 
\begin{align*}
	E(T_k,P_k,\mathbf{C}_1(Z_k, P_{k}))^2 
	\leq\,& 2 E(T_k,P_0,\mathbf{C}_1(Z_k, P_{k}))^2 + 2(q+1) \,\omega_n \,\|\pi_{P_k} - \pi_{P_0}\|^2 \nonumber 
	\\ \leq\,& C \rho^{2\alpha} E(T_k,P_k,\mathbf{C}_1(Z_k, P_{k})) + C \rho^{-n-2} E(T_k,P_0,\mathbf{B}_3(0))^2
\end{align*}
for all $\rho \in (0,\theta]$, where $C = C(n,m,q,\alpha) \in (0,\infty)$ is a constant.  Choosing $\rho = \rho(n,m,q,\alpha) \in (0,\theta]$ sufficiently small, 
\begin{equation}\label{blowup semicont freq eqn6}
	E(T_k,P_k,\mathbf{C}_1(Z_k, P_{k})) \leq C E_k ,
\end{equation}
where we recall that $E_k = E(T_k,P_0,\mathbf{B}_3(0))$ and where $C = C(n,m,q,\alpha) \in (0,\infty)$ is a constant.  Hence by \eqref{blowup semicont freq eqn5} 
\begin{equation}\label{blowup semicont freq eqn7}
	\|\pi_{P_k} - \pi_{P_0}\| \leq C E_k
\end{equation}
for some constant $C = C(n,m,q,\alpha) \in (0,\infty)$.  It follows that for each $k \in \{1,2,3,\ldots\}$ there exists an $m \times n$ matrix $A_k$ such that \eqref{blowup semicont freq eqn2} holds true.  Thus after passing to a subsequence there exists an $m \times n$ matrix $\Lambda$ such that $A_k/E_k \rightarrow \Lambda$ as $k \rightarrow \infty$.  By \eqref{blowup semicont freq eqn4} and \eqref{blowup semicont freq eqn6}, 
\begin{equation}\label{blowup semicont freq eqn8}
	E(T_k,P_k,\mathbf{C}_{\rho}(Z_k, P_{k})) \leq C \rho^{\alpha} E_k 
\end{equation}
for each $k$, where $C = C(n,m,q,\alpha) \in (0,\infty)$ is a constant.  Dividing both sides of \eqref{blowup semicont freq eqn8} by $E_k$ and letting $k \rightarrow \infty$ using \eqref{blowup eqn4}--\eqref{blowup eqn6}, \eqref{blowup eqn9}, \eqref{lipschitz approx eqn8}, \eqref{lipschitz approx eqn9}, and the fact that $\zeta_k/E_k \rightarrow w_a(\xi)$ and $A_k/E_k \rightarrow \Lambda$, 
\begin{eqnarray*}
	&&\int_{B_{\rho}(\xi)} \mathcal{G}(w, q \llbracket w_a(\xi) + \Lambda \cdot (x-\xi) \rrbracket)^2 \nonumber\\
	&&\hspace{1in}= \lim_{k \to \infty} \frac{1}{E_{k}^{2}} \int_{B_{\rho}(\xi)} {\mathcal G}(u_{k}, q\llbracket \zeta_{k} + A_{k}(x - \xi_{k})\rrbracket)^{2}\nonumber\\
	 &&\hspace{1.5in}\leq \lim_{k \rightarrow \infty} \frac{1}{E_k^2} \int_{{\mathbf C}_{\rho}(\xi_k)} \op{dist}^2(X,P_k) \,d\|T_k\|(X) 
	\leq C \rho^{n+2+2\alpha}
\end{eqnarray*}
for all $\rho \in (0,\theta]$, where $C = C(n,m,q,\alpha) \in (0,\infty)$ is a constant.  Hence by \eqref{G avg sym}, 
\begin{equation}\label{blowup semicont freq eqn9}
	\int_{B_{\rho}(\xi)} |w_a - w_a(\xi) - \Lambda \cdot (x-\xi)|^2 \leq C \rho^{n+2+2\alpha}
\end{equation}
for all $\rho \in (0,\theta]$, where $C = C(n,m,q,\alpha) \in (0,\infty)$ is a constant.  Since $w_a$ is harmonic, it follows using the series expansion of $w_a(x-\xi)$ that in order for \eqref{blowup semicont freq eqn9} to hold true we must have that $\Lambda = Dw_a(\xi)$, completing the proof of \eqref{blowup semicont freq eqn3}.  

Now by Lemma~\ref{blowup frequency lemma} and Theorem~\ref{mono freq thm} 
\begin{equation*}
	\limsup_{k \rightarrow \infty} \mathcal{N}_{T_k,{\rm Pl}}(Z_k) 
	\leq \lim_{k \rightarrow \infty} e^{C \eta^{\gamma} \rho^{\alpha\gamma}} N_{T_k,P_k,Z_k}(\rho)
	= e^{C \eta^{\gamma} \rho^{\alpha\gamma}} N_{w - \ell,\xi}(\rho) 
\end{equation*}
for each $\rho \in (0,1]$.  Letting $\rho \rightarrow 0^+$ gives us \eqref{blowup semicont freq concl}.
\end{proof}

\section{Scales with planar frequency pinching, or without flatness, or without decay}\label{sec:freq pinching sec}

In Lemma~\ref{freq2cone lemma} below, we show that if a locally area-minimizing rectifiable current $T$ in $\mathbf{B}_8(0)$ with no boundary in ${\mathbf B}_{8}(0)$ is sufficiently close to a multiplicity $q$ plane $P$ and if the planar frequency function value $N_{T,P,0}(2) \lesssim 1$, then $T$ is well-approximated by the graph of a homogeneous degree one Dirichlet energy minimizing function in $\mathbf{C}_1(0)$.  As a consequence, we draw the conclusion that 
\begin{itemize}
\item [(i)] either $T$ is significantly closer to a sum of two or more planes intersecting along a common $(n-2)$-dimensional axis than to any plane, or 
\item[(ii)] the points of $T$ of density $\geq q$ concentrate near an $(n-3)$-dimensional linear subspace.  
\end{itemize}
In Section~\ref{sec:nonplanar cones sec} we will use this in combination with monotonicity of the planar frequency function (Theorem~\ref{mono freq thm}) to argue that if $T$ is decaying to a plane $P$ in the sense that \eqref{mono freq decay hyp} holds true for all $\rho \in [\sigma_0,\rho_0]$ and if $N_{T,P,0}(\rho_0) \lesssim 1$, then $N_{T,P,0}(\rho) \lesssim 1$ for all $\rho \in [\sigma_0,\rho_0]$ and thus $T$  satisfies the above conclusion in $\mathbf{C}_{\rho/2}(0)$. 

We shall also need to consider scales where \eqref{mono freq decay hyp} does not hold true, i.e.~where $T$ is not decaying towards a plane.  We look at this in Subsections~\ref{sec:freq pinching subsec2} and \ref{sec:freq pinching subsec3}.  We show in this setting that the frequency of $T$ relative to a plane, not necessarily optimal, is $\lesssim 1$.  We also show that $T$ is close to either an area-minimizing cone or the graph of a homogeneous degree one Dirichlet energy minimizing function in $\mathbf{B}_1(0)$, and thus the conclusion of Lemma~\ref{freq2cone lemma} again holds true.

In what follows, we shall use the notation $\mathcal{P}$ to denote the set of all $n$-dimensional planes in $\mathbb{R}^{n+m}.$ We well also use the 
following:  

\begin{definition}\label{sum of planes defn}{\rm
Let $q$ be a positive integer. We denote by $\mathcal{C}_q$ the set of all $n$-dimensional currents $\mathbf{C}$ of the form $\mathbf{C} = \sum_{i=1}^p q_i \llbracket P_i \rrbracket$ where:
\begin{enumerate}[itemsep=3mm,topsep=0mm]
	\item[(i)] $p \geq 2$, $q_i \geq 1$ are integers such that $\sum_{i=1}^p q_i = q$;
	\item[(ii)]  $P_i$ are $n$-dimensional oriented planes;
	\item[(iii)]  there is an $(n-2)$-dimensional subspace $L = \op{spine} \mathbf{C}$ such that $P_i \cap P_j = L$ for all $i \neq j$.
\end{enumerate} 
}\end{definition}

Note that $\mathcal{C}_q$ contains $n$-dimensional area-minimizing integral cones $\mathbf{C}$ such that $\Theta(\mathbf{C},0) = q$ and $\dim \op{spine} \mathbf{C} = n-2$.  $\mathcal{C}_q$ also contains graphs of all homogeneous degree one, locally Dirichlet energy minimizing $q$-valued functions $w \in W^{1,2}_{\rm loc}(\mathbb{R}^n, \mathcal{A}_q(\mathbb{R}^m))$ with $\dim \op{spine} w = n-2$.  

Given $Z \in \mathbb{R}^{n+m}$, $\rho > 0$, an $n$-dimensional current $T$  of $\mathbf{B}_{\rho}(Z)$, and $\mathbf{C} \in \mathcal{C}_q$, we define 
\begin{align*}
	E(T,\mathbf{C},\mathbf{B}_{\rho}(Z))	= \bigg(& \frac{1}{\rho^{n+2}}\int_{\mathbf{B}_{\rho}(Z)} \op{dist}^2(X, Z+\op{spt} \mathbf{C}) \,d\|T\|(X) \bigg)^{1/2} , \\
	Q(T,\mathbf{C},\mathbf{B}_{\rho}(Z))	= \bigg(& \frac{1}{\rho^{n+2}}\int_{\mathbf{B}_{\rho}(Z)} \op{dist}^2(X, Z+\op{spt} \mathbf{C}) \,d\|T\|(X)
		\\&+ \frac{1}{\rho^{n+2}}\int_{\mathbf{B}_{\rho}(0) \cap \{ X : \op{dist}(X,\op{spine} \mathbf{C}) \geq \rho/16 \}} \op{dist}^2(Z+X, \op{spt} T) 
		\,d\|\mathbf{C}\|(X) \bigg)^{1/2} . 
\end{align*}

\subsection{Planar frequency pinching near one}\label{sec:freq pinching subsec1}
First we show that if $T$ has density $\geq q$ at the origin and $T$ is sufficiently close to a multiplicity $q$ plane, and if the planar frequency function value is not much larger than 1 (as in \eqref{freq2cone hyp4}), then $T$ is significantly closer to the graph of a non-zero homogeneous degree one Dirichlet energy minimizing function than it is to any plane.

\begin{lemma}\label{freq2cone lemma}
For every $\varepsilon \in (0,1)$ and $\beta \in (0,1)$ there exists $\delta = \delta(n,m,q,\varepsilon,\beta) \in (0,1)$, $\eta = \eta(n,m,q,\varepsilon,\beta) \in (0,1)$, and $\nu = \nu(n,m,q,\varepsilon,\beta) \in (0,1)$ such that if $T$ is an $n$-dimensional locally area-minimizing rectifiable current in $\mathbf{B}_8(0)$, $P \in \mathcal{P}$, ${\rm spt} \, T \cap {\mathbf B}_{2}(0) \not\subset P$ and if   
\begin{gather}
	\label{freq2cone hyp1} (\partial T) \llcorner \mathbf{B}_8(0) = 0, \quad\quad \Theta(T,0) \geq q , \quad\quad 
		\|T\|(\mathbf{B}_8(0)) \leq (q + \delta) \,\omega_n 8^n , \\
	\label{freq2cone hyp2} E(T,P,\mathbf{B}_1(0)) = \inf_{P' \in \mathcal{P}} E(T,P',\mathbf{B}_1(0)) , \\
	\label{freq2cone hyp3} E(T,P,\mathbf{B}_4(0)) < \eta , \;\; \mbox{and}\\
	\label{freq2cone hyp4} N_{T \llcorner \mathbf{B}_{15/4}(0),P,0}(2) < 1 + \nu 
\end{gather}
then for each $\rho \in [1/2,1]$ one of the following conclusions {\rm (i)} or {\rm (ii)} holds true: 
\begin{enumerate}[itemsep=3mm,topsep=0mm]
	\item[{\rm (i)}]  there exists $\mathbf{C} \in \mathcal{C}_q$ such that 
	\begin{gather}
		\label{freq2cone concl1} E(T,\mathbf{C},\mathbf{B}_{\rho}(0)) < \varepsilon, \\
		\label{freq2cone concl2} Q(T,\mathbf{C},\mathbf{B}_{\rho}(0)) \leq \beta \inf_{P' \in \mathcal{P}} E(T,P',\mathbf{B}_{\rho}(0)) , \\
		\label{freq2cone concl3} \{ X \in \mathbf{B}_{\rho}(0) : \Theta(T,X) \geq q \} 
			\subset \{ X \in \mathbf{B}_{\rho}(0) : \op{dist}(X,\op{spine} \mathbf{C}) < \varepsilon\rho \} ; 
	\end{gather}
	\item[{\rm (ii)}]  there is an $(n-3)$-dimensional linear subspace $L$ such that 
	\begin{equation}\label{freq2cone concl4}
		\{ X \in \mathbf{B}_{\rho}(0) : \Theta(T,X) \geq q \} \subset \{ X \in \mathbf{B}_{\rho}(0) : \op{dist}(X,L) < \varepsilon\rho \} . 
	\end{equation}
\end{enumerate}
\end{lemma}

\begin{proof}[Proof of Lemma~\ref{freq2cone lemma}]
Without loss of generality suppose that $P = \mathbb{R}^n \times \{0\}$.  Fix $\varepsilon \in (0,1)$ and $\beta \in (0,1)$.  Suppose that for $k = 1,2,3,\ldots$ there exist $\delta_k \rightarrow 0^+$, $\eta_k \rightarrow 0^+$, $\nu_k \rightarrow 0^+$, and an $n$-dimensional locally area-minimizing rectifiable current $T_{k}$ of $\mathbf{B}_8(0)$ such that \eqref{freq2cone hyp1}--\eqref{freq2cone hyp4} hold true with $\delta_k$, $\eta_k$, $\nu_k$, and $T_k$ in place of $\delta$, $\eta$, $\nu$, and $T$.  We want to show that for infinitely many $k$ either conclusion (i) or conclusion (ii) holds true with $T_{k}$ in place of $T$.  Arguing as in Subsection~\ref{sec:blowup procedure}, by \eqref{freq2cone hyp1} and \eqref{freq2cone hyp3} after passing to a subsequence we may assume that $T_k \rightarrow q \llbracket P \rrbracket$ weakly in $\mathbf{B}_4(0)$.  By Lemma~\ref{Allard height lemma} and \eqref{freq2cone hyp3}, 
\begin{equation}\label{freq2cone eqn1}
	\sup_{X \in \op{spt} T_k \cap \mathbf{B}_{15/4}(0)} \op{dist}(X,P) \leq C \eta_k 
\end{equation}
for some constant $C = C(n,m) \in (0,\infty)$.  Set $\widetilde{T}_k = T_k \llcorner \mathbf{B}_{15/4}(0)$ so that by \eqref{freq2cone hyp1} and the monotonicity formula for area, $(\partial \widetilde{T}_k) \llcorner \mathbf{C}_{7/2}(0) = 0$ and 
\begin{equation*}
	\|\widetilde{T}_k\|(\mathbf{C}_{7/2}(0)) \leq \|\widetilde{T}_k\|(\mathbf{B}_{7/2+C \eta_k}(0)) \leq (q + \delta_k) \,\omega_n (7/2+C \eta_k)^n ,
\end{equation*}
where $C$ is as in \eqref{freq2cone eqn1}.  Let $\theta_k \in (0,1)$ with $\theta_k \rightarrow 1^+$, and let $u_k : B_{2\theta_k}(0) \rightarrow \mathcal{A}_q(\mathbb{R}^m)$ and $K_k \subset B_{2\theta_k}(0)$ be such that \eqref{blowup eqn4}--\eqref{blowup eqn6} hold true with $\rho_0 = 2$ and $E_k = E(\widetilde{T}_k,P,\mathbf{C}_2(0)).$ Let $w$ be a blow-up of $\widetilde{T}_k$ relative to $P$ in $\mathbf{C}_2(0)$ (as in Subsection~\ref{sec:blowup procedure}).

We argue using \eqref{freq2cone hyp4} and Corollary~\ref{doubling cor}  that $w$ is non-zero.  Note that by \eqref{freq2cone hyp3} 
\begin{equation}\label{freq2cone eqn2}
	E(\widetilde{T}_k,P,\mathbf{C}_{7\rho/4}(0)) < (32/7)^{(n+2)/2} \eta_k \leq 2(32/7)^{(n+2)/2} \eta_k \Big(\frac{\rho}{2}\Big)^{\alpha} 
\end{equation}
for all $\rho \in [1/4,2]$ and all sufficiently large $k$, where $\alpha = 1/2$.  By \eqref{freq2cone hyp1} and \eqref{freq2cone eqn2} we can apply Lemma~\ref{H error lemma}, in particular \eqref{H error concl1} together with \eqref{H error eqn11}, to obtain 
\begin{equation}\label{freq2cone eqn3}
	\frac{1}{2} \,H_{\widetilde{T}_k,P,0}(\rho) \leq -\frac{1}{\rho^{n}} \int_{\mathbf{C}_{\rho}(0)} \op{dist}^2(X,P) \,\phi'(r/\rho) \,d\|\widetilde{T}_k\|(X) 
		\leq (1 + C \eta_k^{2\gamma}) \,H_{\widetilde{T}_k,P,0}(\rho)
\end{equation}
for all $\rho \in [1/4,2]$, where $\gamma = \gamma(n,m,q) \in (0,1)$ and $C = C(n,m,q) \in (0,\infty)$ are constants.  Again by 
\eqref{freq2cone hyp1} and \eqref{freq2cone eqn2} we can apply Theorem~\ref{mono freq thm} and \eqref{freq2cone hyp4} to see that  
\begin{equation}\label{freq2cone eqn4}
	N_{\widetilde{T}_k,P,0}(\rho) \leq e^{C \eta_k^{2\gamma}} N_{\widetilde{T}_k,P,0}(2) < e^{C \eta_k^{2\gamma}} (1 + \nu_k)
\end{equation} 
for all $\rho \in [1/4,2]$, where $C = C(n,m,q) \in (0,\infty)$ is a constant.  Moreover, by Corollary~\ref{doubling cor}, 
\begin{equation}\label{freq2cone eqn5}
	H_{\widetilde{T}_k,P,0}(2\rho) 
	\leq e^{C\eta_k^{2\gamma}} 2^{2 e^{C\eta_k^{\gamma}} (1+\nu_k)} \,H_{\widetilde{T}_k,P,0}(\rho) 
	\leq 8 \,H_{\widetilde{T}_k,P,0}(\rho)
\end{equation}
for all $\rho \in [1/4,1]$ and all sufficiently large $k$, where $C = C(n,m,q) \in (0,\infty)$ is a constant.  Thus, by the fact that $\phi'(r/\rho) = 2$ if $\rho/2 \leq r \leq \rho$, \eqref{freq2cone eqn3}, and \eqref{freq2cone eqn5}, 
\begin{equation*}
	\frac{1}{2^{n+2}} \int_{\mathbf{C}_{2\rho}(0) \setminus \mathbf{C}_{\rho}(0)} \op{dist}^2(X,P) \,d\|\widetilde{T}_k\|(X) 
	\leq 8 \int_{\mathbf{C}_{\rho}(0) \setminus \mathbf{C}_{\rho/2}(0)} \op{dist}^2(X,P) \,d\|\widetilde{T}_k\|(X) 
\end{equation*}
for all $\rho \in [1/4,1]$.  Thus  
\begin{align*}
	E_k^2 =\,& \int_{\mathbf{C}_2(0)} \op{dist}^2(X,P) \,d\|\widetilde{T}_k\|(X)  
		\\ \leq\,& \int_{\mathbf{C}_{1/4}(0)} \op{dist}^2(X,P) \,d\|\widetilde{T}_k\|(X) 
			+ \int_{\mathbf{C}_{1/2}(0) \setminus \mathbf{C}_{1/4}(0)} \op{dist}^2(X,P) \,d\|\widetilde{T}_k\|(X) \nonumber 
			\\&+ \int_{\mathbf{C}_1(0) \setminus \mathbf{C}_{1/2}(0)} \op{dist}^2(X,P) \,d\|\widetilde{T}_k\|(X) 
			+ \int_{\mathbf{C}_2(0) \setminus \mathbf{C}_1(0)} \op{dist}^2(X,P) \,d\|\widetilde{T}_k\|(X) \nonumber
		\\ \leq\,& (1 + 2^{n+5} + 2^{2(n+5)} + 2^{3(n+5)}) \int_{\mathbf{C}_{1/4}(0)} \op{dist}^2(X,P) \,d\|\widetilde{T}_k\|(X) \nonumber
		\\ \leq\,& \frac{2^{4(n+5)} - 1}{2^{n+5} - 1} \int_{\mathbf{B}_{1/2}(0)} \op{dist}^2(X,P) \,d\|\widetilde{T}_k\|(X) \nonumber 
\end{align*}
Hence using \eqref{blowup eqn4}--\eqref{blowup eqn6}, \eqref{blowup eqn9}, \eqref{lipschitz approx eqn8} and \eqref{lipschitz approx eqn9} 
\begin{align}\label{freq2cone eqn7}
	\frac{2^{n+5} - 1}{2^{4(n+5)} - 1} &\leq \lim_{k \rightarrow \infty} \frac{1}{E_k^2} \int_{\mathbf{C}_{1/2}(0)} \op{dist}^2(X,P) \,d\|\widetilde{T}_k\|(X) 
		\\&= \lim_{k \rightarrow \infty} \frac{1}{E_k^2} \int_{B_{1/2}(0)} |u_k|^2 = \int_{B_{1/2}(0)} |w|^2 . \nonumber 
\end{align}
In particular, $w$ is non-zero. 

Now by Lemma~\ref{hardt simon lemma}(i) we have that $w_{a}(0) = 0$, and by \eqref{freq2cone hyp2} and Lemma~\ref{optimal plane lemma} we have that $Dw_a(0) = 0$. By Lemma~\ref{blowup frequency lemma} and \eqref{freq2cone eqn4}
\begin{equation*} 
	N_{w,0}(2) = \lim_{k \rightarrow \infty} N_{\widetilde{T}_k,P,0}(2) \leq 1.
\end{equation*}
On the other hand, since $\Theta(\widetilde{T}_k,0) = \Theta(T_k,0) \geq q$, it follows from by Lemma~\ref{hardt simon lemma} that $\mathcal{N}_w(0) \geq 1$.  Therefore, by the monotonicity formula for the Almgren frequency function for Dirichlet energy minimizing multi-valued functions, $N_{w,0}(\rho) = 1$ for all $\rho \in (0,2]$ and thus $w$ is homogeneous of degree 1.  In particular, the average $w_a$ of $w$ is a homogeneous degree 1 harmonic function and thus $w_a(x) = Dw_a(0) \cdot x = 0$ for all $x \in B_2(0)$.    Therefore, we have shown that $w$ is a non-zero, average-free, homogeneous degree 1, Dirichlet energy minimizing $q$-valued function. 

If $\dim \op{spine} w = n-2$, then we set $\mathbf{C}_k = \op{graph}(E_k w)$ and note that $\mathbf{C} \in \mathcal{C}_q$.  For each $x \in B_2(0)$ let $u_k(x) = \sum_{i=1}^q \llbracket u_{k,i}(x) \rrbracket$ and $w(x) = \sum_{i=1}^q \llbracket w_i(x) \rrbracket$, where $u_{k,i}(x), w_i(x) \in \mathbb{R}^m$.  Notice that for each $x \in B_2(0)$, at $X = (x,u_{k,i}(x))$ 
\begin{equation*}
	\op{dist}(X, \op{spt} \mathbf{C}_k) \leq \op{dist}( u_{k,i}(x), \op{spt}(E_k w(x)) ) \leq \mathcal{G}(u_k(x), E_k w(x)) ,
\end{equation*}
where $\op{spt}(E_k w(x))$ denotes the set of all values $E_k w_i(x)$ in $\mathbb{R}^m$.  Hence using \eqref{blowup eqn4}--\eqref{blowup eqn6}, \eqref{blowup eqn9}, \eqref{lipschitz approx eqn8} and \eqref{lipschitz approx eqn9}, 
\begin{align}\label{freq2cone eqn8}
	&\lim_{k \rightarrow \infty} \frac{1}{E_k^2} \int_{\mathbf{B}_1(0)} \op{dist}^2(X, \op{spt} \mathbf{C}_k) \,d\|T_k\|(X) 
	\\ \leq\,& \lim_{k \rightarrow \infty} \frac{1}{E_k^2} \int_{B_1(0)} \sum_{i=1}^q \op{dist}^2((x,u_{k,i}(x)), \op{spt} \mathbf{C}_k) \,dx \nonumber 
	\\ \leq\,& \lim_{k \rightarrow \infty} \frac{q}{E_k^2} \int_{B_1(0)} \mathcal{G}(u_k(x), E_k w(x))^2 \,dx = 0 . \nonumber
\end{align}
(Note that here and throughout the remainder of the proof we use the fact that $T_k \llcorner \mathbf{B}_2(0) = \widetilde{T}_k \llcorner \mathbf{B}_2(0)$.)  Similarly, for each $x \in B_2(0)$, at $X = (x,w_i(x))$ we have $\op{dist}(X, \op{spt} T_k) \leq \mathcal{G}(u_k(x), E_k w(x))$ and thus 
\begin{align}\label{freq2cone eqn9}
	&\lim_{k \rightarrow \infty} \frac{1}{E_k^2} \int_{\mathbf{B}_1(0) \cap \{ X : \op{dist}(X,\op{spine} \mathbf{C}_k) \geq 1/32 \}} 
		\op{dist}^2(X, \op{spt} T_k) \,d\|\mathbf{C}_k\|(X) 
	\\ \leq\,& \lim_{k \rightarrow \infty} \frac{q}{E_k^2} \int_{B_1(0)} \mathcal{G}(u_k(x), E_k w(x))^2 \,dx = 0 . \nonumber
\end{align}
Let $\rho \in [1/2,1]$ and $\widetilde{P}_k$ be an optimal plane for $T_k$ in $\mathbf{B}_{\rho}(0)$.  By Lemma~\ref{optimal plane lemma},  $\widetilde{P}_k = \{ (x, \widetilde{A}_k x) : x \in \mathbb{R}^n \}$ for some $m \times n$ matrix $\widetilde{A}_k$ with $\|\widetilde{A}_k\| \leq C(n,m,q) \,E_k$ and $\widetilde{A}_k/E_k \rightarrow Dw_a(0) = 0$.  Recall from the proof of Lemma~\ref{optimal plane lemma} that for each $X = (x,y) \in \mathbf{B}_2(0)$
\begin{equation*}
	\op{dist}(X, \op{spt} \widetilde{P}_k) \leq |y - \widetilde{A}_k x| \leq (1 + CE_k) \op{dist}(X, \op{spt} \widetilde{P}_k) ,
\end{equation*}
where $C = C(n,m) \in (0,\infty)$ is a constant.  Hence using \eqref{blowup eqn4}--\eqref{blowup eqn6}, \eqref{blowup eqn9}, \eqref{lipschitz approx eqn8}, \eqref{lipschitz approx eqn9}, the fact that $\widetilde{A}_k/E_k \rightarrow 0$, and \eqref{freq2cone eqn7}, 
\begin{align}\label{freq2cone eqn10}
	\lim_{k \rightarrow \infty} \frac{1}{E_k^2} \int_{\mathbf{B}_{\rho}(0)} \op{dist}^2(X, \op{spt} \widetilde{P}_k) \,d\|T_k\|(X) 
	\geq\,& \lim_{k \rightarrow \infty} \frac{q}{E_k^2} \int_{B_{\rho}(0)} \mathcal{G}(u_k(x), q \llbracket \widetilde{A}_k x \rrbracket)^2 \,dx 
	\\ =\,& q \int_{B_{\rho}(0)} |w(x)|^2 \,dx \geq \frac{2^{n+7} - 1}{2^{4(n+7)} - 1} . \nonumber
\end{align}
Therefore, combining \eqref{freq2cone eqn8}, \eqref{freq2cone eqn9}, and \eqref{freq2cone eqn10} we conclude that \eqref{freq2cone concl1} and \eqref{freq2cone concl2} hold true for all $\rho \in [1/2,1]$ and all sufficiently large $k$.

To prove \eqref{freq2cone concl3}, suppose to the contrary that for infinitely many $k$ there exists $Z_k \in \mathbf{B}_1(0)$ such that $\Theta(T_k,Z_k) \geq q$ and $\op{dist}(Z_k, \op{spine}\mathbf{C}_k) \geq \varepsilon/2$.  Note that $\op{spine}\mathbf{C}_k = (\op{spine} w) \times \{0\}$.  By Lemma~\ref{hardt simon lemma}, after passing to a subsequence $Z_k \rightarrow (\xi,0)$ for some point $\xi \in \overline{B_1(0)}$ with $\mathcal{N}_w(\xi) \geq 1$ and $\op{dist}(\xi,\op{spine} w) \geq \varepsilon/2$.  But since $w$ is homogeneous degree one and $\mathcal{N}_w(\xi) \geq 1$, $\xi \in \op{spine} w$, contradicting $\op{dist}(\xi,\op{spine} w) \geq \varepsilon/2$.  This completes the proof of conclusion~(i).  By a similar argument, if $\dim \op{spine} w \leq n-3$, then conclusion~(ii) must hold true with $L = \op{spine} w$. 
\end{proof}

\subsection{Distance to a plane not small}\label{sec:freq pinching subsec2}
Let $\eta > 0$.  Let $T$ be an $n$-dimensional locally area-minimizing rectifiable current of $\mathbf{B}_8(0)$ and $P$ be an $n$-dimensional  plane such that 
\begin{gather}
	\label{not flat hyp1} (\partial T) \llcorner \mathbf{B}_8(0) = 0, \quad\quad \Theta(T,0) \geq q , \quad\quad 
		\|T\|(\mathbf{B}_8(0)) \leq (q + \delta) \,\omega_n 8^n , \\
	\label{not flat hyp2} E(T,P,\mathbf{B}_1(0)) = \inf_{P' \in \mathcal{P}} E(T,P',\mathbf{B}_1(0)) .
\end{gather}
We wish to consider the case where 
\begin{equation}\label{not flat hyp3}
	E(T,P,\mathbf{B}_4(0)) \geq \eta . 
\end{equation}
First we show that for any plane $\widetilde{P}$ (not necessarily equal to $P$ nor assumed to be optimal in a ball), the planar frequency function of $T$ relative to $\widetilde{P}$ and any base point close to 0 with density $\geq q$ takes values $\lesssim 1$ at scales $\in [1/2, 1]$.

\begin{lemma}\label{not flat lemma1}
For every $\eta > 0$ and $\nu \in (0,1)$ there exists $\delta = \delta(n,m,q,\eta,\nu) \in (0,1)$ such that if $P, \widetilde{P}$ are $n$-dimensional planes in $\mathbb{R}^{n+m}$ and $T$ is an $n$-dimensional locally area-minimizing rectifiable current in $\mathbf{B}_8(0)$ such that \eqref{not flat hyp1}, \eqref{not flat hyp2}, and \eqref{not flat hyp3} hold true and 
\begin{equation}\label{not flat hyp4} 
	\frac{1}{\omega_n 2^{n+2}} \int_{\mathbf{B}_2(Z)} \op{dist}^2(X,\widetilde{P}) \,d\|T\|(X) < 64^{-n-2} , 
\end{equation}
then for each $Z \in \op{spt} T \cap \mathbf{B}_{\delta}(0)$ with $\Theta(T,Z) \geq q$ we have that ${\rm spt} \, T \cap {\mathbf B}_{1/2}(Z) \not\subset \widetilde{P}$ and that $N_{T \llcorner \mathbf{B}_{15/8}(0),\widetilde{P},Z}(\rho) < 1 + \nu/2$ for each $\rho \in [1/2,1]$.  
\end{lemma}

Lemma~\ref{not flat lemma1} follows from Lemma~\ref{freq of cones lemma}, which states that the frequency function of an area-minimizing cone $\mathbf{C}$ takes values $\leq 1$.  Note that in Lemma~\ref{freq of cones lemma} we do not assume that $\mathbf{C}$ is close to a plane or is represented by the graph of a multi-valued function, showing a key advantage of our intrinsic definition of planar frequency function (Definition~\ref{freq defn}). 

\begin{lemma}\label{freq of cones lemma}
Let $\mathbf{C}$ be a non-zero $n$-dimensional stationary cone $\mathbf{C}$ and $P$ be an $n$-dimensional plane $P$ in $\mathbb{R}^{n+m}$.  Assume that 
\begin{equation}\label{freq of cones hyp}
	\op{spt} \mathbf{C} \cap P^{\perp} = \{0\} , \quad \mbox{and} \quad \op{spt} \mathbf{C} \neq P . 
\end{equation}
Then $N_{\mathbf{C},P,0}(\rho) \leq 1$ for all $\rho > 0$. 
\end{lemma}

\begin{remark}\label{freq of cones rmk}{\rm 
(1)  For an a cone $\mathbf{C}$, the assumption $\op{spt} \mathbf{C} \cap P^{\perp} = \{0\}$ is equivalent to 
\begin{equation}\label{freq of cones rmk eqn1}
	\sup_{X \in \op{spt} \mathbf{C} \cap \mathbf{C}_1(0,P)} \op{dist}(X,P) < \infty .
\end{equation}
To see this, assume that $P = \mathbb{R}^n \times \{0\}$ and express points $X \in \mathbb{R}^{n+m}$ as $X = (x,y)$ where $x \in \mathbb{R}^n$ and $y \in \mathbb{R}^m$.  If \eqref{freq of cones rmk eqn1} holds true then $\op{spt} \mathbf{C} \cap P^{\perp} = \emptyset$ since if there were $(0,y) \in \op{spt} \mathbf{C} \cap P^{\perp}$ for some $y \neq 0$, then $(0,ty) \in \op{spt} \mathbf{C}$ and $\op{dist}((0,ty),P) = t |y|$ for arbitrarily large $t > 0$.  On the other hand, if \eqref{freq of cones rmk eqn1} did not hold true, there would be a sequence of points $X_k = (x_k,y_k) \in \op{spt} \mathbf{C} \cap \mathbf{C}_1(0)$ such that $|y_k| \rightarrow \infty$.  After passing to a subsequence, $X_k/|X_k| \rightarrow (0,\zeta)$ for some $\zeta \in \mathbb{R}^m$ with $|\zeta| = 1$.  Since ${\rm spt} \, {\mathbf C}$ is closed, $(0,\zeta) \in \op{spt} \mathbf{C} \cap P^{\perp}$ and thus $\op{spt} \mathbf{C} \cap P^{\perp} \neq \emptyset$.

\noindent (2)  For an area-minimizing cone $\mathbf{C}$ with $\op{spt} \mathbf{C} \cap P^{\perp} = \emptyset$, it follows from Remark~\ref{H zero rmk} and the constancy theorem~\cite[Theorem~26.27]{SimonGMT} that $\mathbf{C} \neq 0$ and $\op{spt} \mathbf{C} \neq P$ if and only if $H_{\mathbf{C},P,0}(\rho) > 0$ for all $\rho > 0$. 
}\end{remark}

\begin{proof}[Proof of Lemma~\ref{freq of cones lemma}]
Without loss of generality let $P = \mathbb{R}^n \times \{0\}$.  Express points $X \in \mathbb{R}^{n+m}$ as $X = (x,y)$ where $x \in \mathbb{R}^n$ and $y \in \mathbb{R}^m$.  
In light of Remark~\ref{freq of cones rmk}, $H_{\mathbf{C},P,0}(\rho) > 0$ for each $\rho > 0$, so it suffices to show that $D_{\mathbf{C},P,0}(\rho) \leq H_{\mathbf{C},P,0}(\rho)$ for all $\rho > 0$.
Since $\mathbf{C}$ is a cone, for $|\mathbf{C}|$-a.e.~$(X,S) \in G_n(\mathbb{R}^{n+m})$ the radial vector $X$ is tangent to $\op{spt}\mathbf{C}$ at $X$, i.e.~$X \in S$.  Thus 
\begin{equation*}
	0 = \pi_{S^{\perp}}(X) = \pi_{S^{\perp}}(x,y) = \pi_{S^{\perp}}(x,0) + \pi_{S^{\perp}}(0,y) 
\end{equation*}
for $|\mathbf{C}|$-a.e.~$(X,S) \in G_n(\mathbb{R}^{n+m})$.  That is, $\pi_{S^{\perp}}(x,0) = -\pi_{S^{\perp}}(0,y)$.  Hence  
\begin{align*}
	(0,y) \cdot \pi_{P^{\perp}}(\nabla^S r) &= (0,y) \cdot \nabla^S r = (0,y) \cdot \frac{(x,0) - \pi_{S^{\perp}}(x,0)}{r} 
		\\&= (0,y) \cdot \frac{(x,0) + \pi_{S^{\perp}}(0,y)}{r} = \frac{|\pi_{S^{\perp}}(0,y)|^2}{r} 
\end{align*}
for $|\mathbf{C}|$-a.e.~$(X,S) \in G_n(\mathbb{R}^{n+m})$.  Thus by \eqref{D to I} (with $T = \mathbf{C}$) 
\begin{equation}\label{freq of cones eqn2}
	D_{\mathbf{C},P,0}(\rho) = -\rho^{1-n} \int |\pi_{S^{\perp}}(0,y)|^2 \,\frac{1}{r} \,\phi'(r/\rho) \,d|\mathbf{C}|(X,S) 
\end{equation}
for all $\rho > 0$.  For $|\mathbf{C}|$-a.e.~$(X,S) \in G_n(\mathbb{R}^{n+m})$, since $X \in S$, 
\begin{equation*}
	(0,y) = \frac{(0,y) \cdot X}{|X|^2} \,X + V + \pi_{S^{\perp}}(0,y) , 
\end{equation*}
where $V$ is a vector in $S$ and orthogonal to $X$.  Hence 
\begin{equation}\label{freq of cones eqn3}
	|\pi_{S^{\perp}}(0,y)|^2 \leq |y|^2 - \left( (0,y) \cdot \frac{(x,y)}{|X|} \right)^2 = |y|^2 - \frac{|y|^4}{|X|^2} = \frac{|x|^2 |y|^2}{|X|^2} 
\end{equation}
for $|\mathbf{C}|$-a.e.~$(X,S) \in G_n(\mathbb{R}^{n+m})$.  By $\pi_{S^{\perp}}(x,0) = -\pi_{S^{\perp}}(0,y)$ and \eqref{freq of cones eqn3},  
\begin{align}\label{freq of cones eqn4}
	|y|^2 \,|\nabla^S r|^2 &= |y|^2 \left( 1 - \frac{|\pi_{S^{\perp}}(x,0)|^2}{|x|^2} \right) = |y|^2 \left( 1 - \frac{|\pi_{S^{\perp}}(0,y)|^2}{|x|^2} \right)  
		\\&\geq |y|^2 \left( 1 - \frac{|y|^2}{|X|^2} \right) = \frac{|x|^2 |y|^2}{|X|^2} \geq |\pi_{S^{\perp}}(0,y)|^2 \nonumber 
\end{align}
for $|\mathbf{C}|$-a.e.~$(X,S) \in G_n(\mathbb{R}^{n+m})$.  In view of \eqref{H defn} (with $T = \mathbf{C}$), we see from \eqref{freq of cones eqn2} and \eqref{freq of cones eqn4}, that $D_{\mathbf{C},P,0}(\rho) \leq H_{\mathbf{C},P,0}(\rho)$ for all $\rho > 0$.
\end{proof}

\begin{proof}[Proof of Lemma~\ref{not flat lemma1}]
Fix $\eta > 0$ and $\nu \in (0,1)$.  Without loss of generality assume that $\widetilde{P} = \mathbb{R}^n \times \{0\}$.  Suppose to the contrary that for $k = 1,2,3,\ldots,$ there exists $\delta_k \rightarrow 0^+$, an $n$-dimensional plane $P_k$, and an $n$-dimensional locally area-minimizing rectifiable current $T_k$ of $\mathbf{B}_8(0)$ such that \eqref{not flat hyp1}, \eqref{not flat hyp2}, \eqref{not flat hyp3}, and \eqref{not flat hyp4} hold true with $\delta_k$, $P_k$, and $T_k$ in place of $\delta$, $P$, and $T$ but $N_{T_k \llcorner \mathbf{B}_{15/8}(0),\widetilde{P},Z_k}(\rho_k) \geq 1 + \nu/2$ for some $Z_k \in \op{spt} T_k \cap \mathbf{B}_{\delta_k}(0)$ with $\Theta(T_k,Z_k) \geq q$ and some $\rho_k \in [1/2,1]$.  There exit $\rho_{\infty} \in [1/2, 1]$ and an $n$-dimensional plane $P_{\infty}$ such that after passing to a subsequence (without relabeling), $\rho_k \rightarrow \rho_{\infty}$ and 
\begin{equation}\label{not flat eqn1}
	\lim_{k \rightarrow \infty} \op{dist}_{\mathcal{H}}(P_k \cap \mathbf{B}_1(0), P_{\infty} \cap \mathbf{B}_1(0)) = 0 . 
\end{equation}
By \eqref{not flat hyp1} and the compactness of area-minimizing integral currents (\cite[Theorems~32.2 and 34.5]{SimonGMT}), after passing to a further subsequence there exists an $n$-dimensional area-minimizing integral current $\mathbf{C}$ such that $T_k \rightarrow \mathbf{C}$ weakly in $\mathbf{B}_8(0)$.  Since $Z_k \rightarrow 0$, $\Theta(T_k,Z_k) \geq q$, and $\|T_k\|(\mathbf{B}_8(0)) \leq (q + \delta_k) \,\omega_n 8^n$, by the semi-continuity of density and mass (\cite[Corollary 17.8 and 26.13]{SimonGMT}) we have that $\Theta(\mathbf{C},0) \geq q$ and $\|\mathbf{C}\|(\mathbf{B}_8(0)) \leq q \omega_n 8^n$.  By applying the equality case of the monotonicity formula for area (see~\cite[Theorem~19.3]{SimonGMT}), $\mathbf{C}$ is an area-minimizing cone.  Since $\|T_k\| \rightarrow \|\mathbf{C}\|$ in the sense of Radon measures locally in $\mathbf{B}_8(0)$, by \eqref{not flat hyp2} if $\mathbf{C}$ is a multiplicity $q$ plane then 
\begin{align*}
	\int_{\mathbf{B}_1(0)} \op{dist}^2(X,P_{\infty}) \,d\|\mathbf{C}\|(X) 
		&= \lim_{k \rightarrow \infty} \int_{\mathbf{B}_1(0)} \op{dist}^2(X,P_k) \,d\|T_k\|(X) 
		\\&\leq \lim_{k \rightarrow \infty} \int_{\mathbf{B}_1(0)} \op{dist}^2(X, \op{spt} \mathbf{C}) \,d\|T_k\|(X) 
		\\&= \lim_{k \rightarrow \infty} \int_{\mathbf{B}_1(0)} \op{dist}^2(X, \op{spt} \mathbf{C}) \,d\|\mathbf{C}\|(X) = 0 
\end{align*}
and thus after equipping $P_{\infty}$ with an appropriate orientation $\mathbf{C} = q \llbracket P_{\infty} \rrbracket$.  (Note that since $\mathbf{C}$ is a cone, $\|\mathbf{C}\|(\partial \mathbf{B}_{\rho}(0)) = 0$ for all $\rho \in (0,8)$.)  On the other hand, by \eqref{not flat hyp3} 
\begin{equation}\label{not flat eqn2}
	\int_{\mathbf{B}_4(0)} \op{dist}^2(X, P_{\infty}) \,d\|\mathbf{C}\|(X) 
		= \lim_{k \rightarrow \infty} \int_{\mathbf{B}_4(0)} \op{dist}^2(X,P_k) \,d\|T_k\|(X) \geq \eta^2 ,
\end{equation}
so $\mathbf{C} \neq q \llbracket P_{\infty} \rrbracket$.  Therefore $\mathbf{C}$ is not a multiplicity $q$ plane and in particular $\op{spt} \mathbf{C} \not\subseteq \widetilde{P}$.  By \eqref{not flat hyp4} and Lemma~\ref{coarse L2 distance lemma} with $K = \widetilde{P}$,  
\begin{equation*}
	\sup_{X \in \op{spt} T_k \cap \mathbf{B}_{15/8}(0)} \op{dist}(X,\widetilde{P}) \leq 1/16.
\end{equation*}
Hence $(\partial (T_k \llcorner \mathbf{B}_{15/8}(0)) \llcorner \mathbf{C}_1(0) = 0$.  Moreover, 
\begin{equation*}
	\sup_{X \in \op{spt} \mathbf{C} \cap \mathbf{B}_{15/8}(0)} \op{dist}(X,\widetilde{P}) \leq 1/16 
\end{equation*}
and thus by Remark~\ref{freq of cones rmk}(1), $\op{spt} \mathbf{C} \cap \widetilde{P}^{\perp} = \emptyset$.  Since $\op{spt} \mathbf{C} \not\subseteq \widetilde{P}$ and $\op{spt} \mathbf{C} \cap \widetilde{P}^{\perp} = \emptyset$, we can apply Lemma~\ref{freq of cones lemma} to obtain that $N_{\mathbf{C},\widetilde{P},0}(\rho) \leq 1$ for all $\rho \in (0,8)$.  But by Lemma~\ref{limit frequency lemma} and the fact that $N_{T_k \llcorner \mathbf{B}_{15/8}(0),\widetilde{P},Z_k}(\rho_k) \geq 1+\nu/2$, it follows that 
\begin{equation*}
N_{\mathbf{C},\widetilde{P},0}(\rho_{\infty})   = \lim_{k \rightarrow \infty} N_{T_k \llcorner \mathbf{B}_{15/8}(0),\widetilde{P},Z_k}(\rho_k) \geq 1 + \nu/2,
\end{equation*} 
giving us a contradiction.
\end{proof}

Next we show in Lemma~\ref{not flat lemma2} that, if $T$ is not close to the plane $P$ as in \eqref{not flat hyp3}, then $T$ must be $L^2$-close to a non-planar area-minimizing cone in $\mathbf{B}_1(0)$.

\begin{lemma}\label{not flat lemma2}
For every $\varepsilon \in (0,1)$, $\beta \in (0,1)$, and $\eta > 0$ there exists $\delta = \delta(n,m,q,\varepsilon,\beta,\eta) \in (0,1)$ such that if $P$ is an $n$-dimensional plane in $\mathbb{R}^{n+m}$ and $T$ is an $n$-dimensional locally area-minimizing rectifiable current of $\mathbf{B}_8(0)$ such that \eqref{not flat hyp1}, \eqref{not flat hyp2}, and \eqref{not flat hyp3} hold true, then for each $\rho \in [1/2,1]$ one of the following holds true: 
\begin{enumerate}[itemsep=3mm,topsep=0mm]
	\item[{\rm (i)}]  there exists $\mathbf{C} \in \mathcal{C}_q$ such that \eqref{freq2cone concl1}, \eqref{freq2cone concl2}, and \eqref{freq2cone concl3} hold true or 
	\item[{\rm (ii)}]  there is an $(n-3)$-dimensional linear subspace $L$ such that \eqref{freq2cone concl4} holds true. 
\end{enumerate}
\end{lemma}

\begin{proof}
Fix $\varepsilon \in (0,1)$, $\beta \in (0,1)$, and $\eta > 0$.  Without loss of generality assume that $P = \mathbb{R}^n \times \{0\}$.  For $k = 1,2,3,\ldots$ let $\delta_k \rightarrow 0^+$ and $T_k$ be $n$-dimensional area-minimizing integral currents of $\mathbf{B}_8(0)$ such that \eqref{not flat hyp1}, \eqref{not flat hyp2}, and \eqref{not flat hyp3} hold true with $\delta_k$ and $T_k$ in place of $\delta$ and $T$.  We want to show that for infinitely many $k$ either conclusions~(i) or (ii) hold true, which in light of the arbitrary choice of sequence $(T_k)$ proves the theorem.  As in the proof of Lemma~\ref{not flat lemma1}, it follows from \eqref{not flat hyp1} that after passing to a subsequence $T_k \rightarrow \mathbf{C}$ weakly in $\mathbf{B}_8(0)$ for some $n$-dimensional area-minimizing cone $\mathbf{C}$ with $\Theta(\mathbf{C},0) = q$.  Moreover, by \eqref{not flat hyp2} and \eqref{not flat hyp3}, $\mathbf{C}$ is not a multiplicity $q$ plane and 
\begin{equation}\label{not flat eqn3}
	\int_{\mathbf{B}_4(0)} \op{dist}^2(X, P) \,d\|\mathbf{C}\|(X) \geq \eta^2 
\end{equation}
(as in \eqref{not flat eqn2}).  If $\dim \op{spine} \mathbf{C} = n-2$, then using the fact that $T_k \rightarrow \mathbf{C}$ weakly in $\mathbf{B}_8(0)$ and the monotonicity formula for area, for any $\rho \in [1/2,1]$,  
\begin{align}\label{not flat eqn4} 
	\lim_{k \rightarrow \infty} Q(T_k, \mathbf{C}, \mathbf{B}_{\rho}(0))
	\leq \lim_{k \rightarrow \infty} \bigg(& (q+1/2) \omega_n \sup_{X \in \op{spt} T_k \cap \mathbf{B}_1(0)} \op{dist}^2(X,{\rm spt} \, \mathbf{C})
		\\&+ q \omega_n \sup_{X \in \op{spt} \mathbf{C} \cap \mathbf{B}_{1/2}(0)} \op{dist}^2(X,\op{spt} T_k) \bigg) = 0 . \nonumber 
\end{align}
Let $\rho \in [1/2,1]$ and let $\widehat{P}_k$ be an optimal plane for $T_k$ in $\mathbf{B}_{\rho}(0)$.  After passing to a subsequence, let $\widehat{P}_{\infty}$ be an $n$-dimensional linear plane such that 
\begin{equation*}
	\lim_{k \rightarrow \infty} \op{dist}_{\mathcal{H}}(\widehat{P}_k \cap \mathbf{B}_1(0), \widehat{P}_{\infty} \cap \mathbf{B}_1(0)) = 0 . 
\end{equation*}
By the fact that $\|T_k\| \rightarrow \|\mathbf{C}\|$ in the sense of Radon measures locally in $\mathbf{B}_8(0)$, the homogeneity of $\mathbf{C}$, \eqref{not flat hyp2}, and \eqref{not flat eqn3}, 
\begin{align}\label{not flat eqn5} 
	&\lim_{k \rightarrow \infty} \int_{\mathbf{B}_{\rho}(0)} \op{dist}^2(X, \widehat{P}_k) \,d\|T_k\|(X)
	= \int_{\mathbf{B}_{\rho}(0)} \op{dist}^2(X, \widehat{P}_{\infty}) \,d\|\mathbf{C}\|(X)  
	\\ \geq\,& 2^{-n-2} \int_{\mathbf{B}_1(0)} \op{dist}^2(X, \widehat{P}_{\infty}) \,d\|\mathbf{C}\|(X) 
	= \lim_{k \rightarrow \infty} 2^{-n-2} \int_{\mathbf{B}_1(0)} \op{dist}^2(X, \widehat{P}_{\infty}) \,d\|T_k\|(X) \nonumber
	\\ \geq\,& \lim_{k \rightarrow \infty} 2^{-n-2} \int_{\mathbf{B}_1(0)} \op{dist}^2(X, P) \,d\|T_k\|(X) 
	= \lim_{k \rightarrow \infty} 2^{-n-2} \int_{\mathbf{B}_1(0)} \op{dist}^2(X, P) \,d\|\mathbf{C}\|(X) \nonumber
	\\ =\,&  8^{-n-2} \int_{\mathbf{B}_4(0)} \op{dist}^2(X, P) \,d\|\mathbf{C}\|(X) 
	\geq 8^{-n-2} \eta^2 . \nonumber
\end{align}
By \eqref{not flat eqn4} and \eqref{not flat eqn5}, we see that \eqref{freq2cone concl1} and \eqref{freq2cone concl2} hold true.  By upper semi-continuity of density (\cite[Corollary 17.8]{SimonGMT}), \eqref{freq2cone concl3} holds true if $\dim \op{spine} \mathbf{C} = n-2$ and  \eqref{freq2cone concl3} holds true with $L = \op{spine} \mathbf{C}$ if $\dim \op{spine} \mathbf{C} \leq n-3$.  
\end{proof}

\subsection{Distance to a plane not decaying}\label{sec:freq pinching subsec3}
Let $\delta \in (0, 1)$ and let $T$ be an $n$-dimensional locally area-minimizing rectifiable current of $\mathbf{B}_8(0)$ such that and $P,\widehat{P}$ are $n$-dimensional linear planes such that 
\begin{gather}
	\label{not decaying hyp1} (\partial T) \llcorner \mathbf{B}_8(0) = 0, \quad\quad \Theta(T,0) \geq q , \quad\quad 
		\|T\|(\mathbf{B}_8(0)) \leq (q + \delta) \,\omega_n 8^n.
\end{gather}
Let $P$, $\widehat{P}$ be $n$-dimensional planes such that
\begin{gather}
	\label{not decaying hyp2} E(T,P,\mathbf{B}_1(0)) = \inf_{P' \in \mathcal{P}} E(T,P',\mathbf{B}_1(0)) \;\; \mbox{and}\\
	\label{not decaying hyp3} E(T,\widehat{P},\mathbf{B}_2(0)) = \inf_{P' \in \mathcal{P}} E(T,P',\mathbf{B}_2(0)) .  
\end{gather}
For appropriate $\eta \in (0,1)$ and $\alpha \in (0,1)$, we wish to consider the case where 
\begin{gather}
	\label{not decaying hyp4} E(T,P,\mathbf{B}_4(0)) \leq \eta , \\
	\label{not decaying hyp5} E(T,P,\mathbf{B}_2(0)) > 2^{-\alpha} E(T,\widehat{P},\mathbf{B}_4(0)) . 
\end{gather}
First we show that for any plane $\widetilde{P}$ (not necessarily equal to $P$ nor equal to $\widehat{P}$ nor assumed to be optimal in a ball), the frequency of $T$ relative to $\widetilde{P}$ is $\lesssim 1$. 

\begin{lemma}\label{not decaying lemma1}
For every $\nu \in (0,1)$ there exists $\delta = \delta(n,m,q,\nu) \in (0,1)$, $\eta = \eta(n,m,q,\nu) \in (0,1)$, and $\alpha = \alpha(n,m,q,\nu) \in (0,1)$ such that if $T$ is $n$-dimensional locally area-minimizing rectifiable current of $\mathbf{B}_8(0) \subseteq \mathbb{R}^{n+m}$ and $P, \widehat{P}, \widetilde{P}$ are $n$-dimensional planes in $\mathbb{R}^{n+m}$ such that \eqref{not decaying hyp1}--\eqref{not decaying hyp5} and \eqref{not flat hyp4} hold true and ${\rm spt} \, T \cap {\mathbf B}_{4}(0) \not\subset \widehat{P}$, then for each $Z \in \op{spt} T \cap \mathbf{B}_{\delta}(0)$ with $\Theta(T,Z) \geq q$ we have that 
${\rm spt} \, T \cap {\mathbf B}_{1/2}(Z) \not \subset \widetilde{P}$ and that $N_{T \llcorner \mathbf{B}_{15/8}(0),\widetilde{P},Z}(\rho) < 1 + \nu/2$ for each $\rho \in [1/2,1]$.  
\end{lemma}

\begin{proof}
Without loss of generality suppose that $\widetilde{P} = \mathbb{R}^n \times \{0\}$.  Fix $\nu \in (0,1)$.  Suppose to the contrary that for $k = 1,2,3,\ldots$ there exists $\delta_k \rightarrow 0^+$, $\eta_k \rightarrow 0^+$, $\alpha_k \rightarrow 0^+$, an $n$-dimensional area-minimizing integral current $T_k$ of $\mathbf{B}_8(0)$, and $n$-dimensional planes $\widehat{P}_k, P_k$ in $\mathbb{R}^{n+m}$ such that \eqref{not decaying hyp1}--\eqref{not decaying hyp5} hold true with $\delta_k$, $\eta_k$, $\alpha_k$, $T_k$, $P_k$, and $\widehat{P}_k$ in place of $\delta$, $\eta$, $\alpha$, $T$, $P$, and $\widehat{P}$, but $N_{T_k \llcorner \mathbf{B}_{15/8}(0),\widetilde{P},Z_k}(\rho_k) \geq 1+\nu/2$ for some $Z_k \in \op{spt} T_k \cap \mathbf{B}_{\delta_k}(0)$ with $\Theta(T_k,Z_k) \geq q$ and some $\rho_k \in [1/2,1]$.  After passing to a subsequence, let $\rho_k \rightarrow \rho_{\infty}$ in $[1/2,1]$ and let $P_{\infty}$ be an $n$-dimensional plane such that 
\begin{equation}\label{not decaying eqn1} 
	\lim_{k \rightarrow \infty} \op{dist}_{\mathcal{H}}(P_k \cap \mathbf{B}_1(0), P_{\infty} \cap \mathbf{B}_1(0)) = 0 . 
\end{equation}
Arguing as in Subsection~\ref{sec:blowup procedure}, by \eqref{not decaying hyp1} and \eqref{not decaying hyp4} after passing to a subsequence and equipping $P_{\infty}$ with an appropriate orientation we have that $T_k \rightarrow \mathbf{C} = q \llbracket P_{\infty} \rrbracket$ weakly in $\mathbf{B}_4(0)$.  If $P_{\infty} \neq \widetilde{P}$, then we can argue as in the proof of Lemma~\ref{not flat lemma1} that by \eqref{not flat hyp4} we have $\op{spt} \mathbf{C} \cap \widetilde{P}^{\perp} = \emptyset$ and by Lemma~\ref{limit frequency lemma} and Lemma~\ref{freq of cones lemma} we have $N_{T_k \llcorner \mathbf{B}_{15/8}(0),\widetilde{P},Z_k}(\rho_k) \rightarrow N_{\mathbf{C},\widetilde{P},0}(\rho_{\infty}) \leq 1$, contradicting the assumption $N_{T_k \llcorner \mathbf{B}_{15/8}(0),\widetilde{P},Z_k}(\rho_k) \geq 1+\nu/2$.  Thus we may assume that $P_{\infty} = \widetilde{P}$ ($={\mathbb R}^{n} \times \{0\})$).  

Set $E_k = E(T_k,\widetilde{P},\mathbf{B}_4(0))$.  By \eqref{not decaying hyp4}, \eqref{not decaying eqn1} and the fact that $P_{\infty} = \widetilde{P}$, we have that $E_k \rightarrow 0$.  Let $\theta_k \in (0,1)$ such that $\theta_k \rightarrow 1^+$ and $(1-\theta_k)^{-n-2} E_k^2 \rightarrow 0$.  By Lemma~\ref{coarse L2 distance lemma} with $K = \widetilde{P}$,  
\begin{equation*}
	\sup_{X \in \op{spt} T_k \cap \mathbf{B}_{3+\theta_k}(0)} \op{dist}(X,\widetilde{P}) \leq 8 E_k^{\frac{1}{n+2}} < 1-\theta_k .
\end{equation*}
Set $\widetilde{T}_k = T_k \llcorner \mathbf{B}_{3+\theta_k}(0)$ so that by \eqref{not decaying hyp1} and the monotonicity formula for area, $(\partial \widetilde{T}_k) \llcorner \mathbf{C}_{2+2\theta_k}(0) = 0$ and 
\begin{equation}\label{not decaying eqn2}
	\limsup_{k \rightarrow \infty} \|\widetilde{T}_k\|(\mathbf{C}_{2+2\theta_k}(0)) 
		\leq \limsup_{k \rightarrow \infty} \|\widetilde{T}_k\|(\mathbf{B}_{3+\theta_k}(0)) 
		\leq \lim_{k \rightarrow \infty} (q+\delta_k) \,\omega_n (3+\theta_k)^n = q \omega_n 4^n . 
\end{equation}
By Theorem~\ref{lipschitz approx thm} and Lemma~\ref{area excess to height lemma}, there exists $u_k : B_{4\theta_k}(0) \rightarrow \mathcal{A}_q(\mathbb{R}^m)$, and $K_k \subset B_{4\theta_k}(0)$ such that \eqref{blowup eqn4}--\eqref{blowup eqn6} hold true with $\widetilde{T}_k$ in place of $T_k$, $\rho_0 = 4$, and $E_k = E(T_k,\widetilde{P},\mathbf{B}_4(0))$.  After passing to a subsequence, let $w \in W^{1,2}_{\rm loc}(B_4(0),\mathcal{A}_q(\mathbb{R}^m))$ such that \eqref{blowup eqn9} and \eqref{blowup eqn10} hold true with $\rho_0 = 4$ and $E_k = E(T_k,\widetilde{P},\mathbf{B}_4(0))$.  
By \eqref{not decaying hyp2} and Lemma~\ref{optimal plane lemma}, $P_k = \{ (x, A_k x) : x \in \mathbb{R}^n \}$ for some $m \times n$ matrix $A_k$ such that $\|A_k\| \leq C(n,m,q) \,E_k$ and $A_k/E_k \rightarrow Dw_a(0)$.  Similarly, by \eqref{not decaying hyp3} and Lemma~\ref{optimal plane lemma}, $\widehat{P}_k = \{ (x, \widehat{A}_k x) : x \in \mathbb{R}^n \}$ for some $m \times n$ matrix $\widehat{A}_k$ such that $\|\widehat{A}_k\| \leq C(n,m,q) \,E_k$ and $\widehat{A}_k/E_k \rightarrow Dw_a(0)$.  Let $w(x) = \sum_{i=1}^q \llbracket w_i(x) \rrbracket$ and $w(x) - Dw_a(0) \cdot x = \sum_{i=1}^q \llbracket w_i(x) - Dw_a(0) \cdot x \rrbracket$ for each $x \in B_4(0)$.  Note that $w - Dw_a(0) \cdot x$ is Dirichlet energy minimizing (\cite[Theorem~2.6(3)]{Almgren}).   Recall from the proof of Lemma~\ref{optimal plane lemma} that for each $X = (x,y) \in \mathbf{B}_4(0)$
\begin{equation*}
	\op{dist}(X,P_k) \leq |y - A_k x| \leq (1 + CE_k) \op{dist}(X,P_k) , 
\end{equation*}
where $C = C(n,m, q) \in (0,\infty)$ is a constant.  Thus by \eqref{blowup eqn4}--\eqref{blowup eqn6}, \eqref{blowup eqn9}, \eqref{lipschitz approx eqn8}, \eqref{lipschitz approx eqn9}, and the fact that $A_k/E_k \rightarrow Dw_a(0)$, 
\begin{align}\label{not decaying eqn3}
	\int_{B_{\rho}(0)} |w - Dw_a(0) \cdot x|^2 
	&= \lim_{k \rightarrow \infty} \frac{1}{E_k^2} \int_{B_{\rho}(0)} \mathcal{G}(u_k, q \llbracket A_k \rrbracket)^2 
	\\&= \lim_{k \rightarrow \infty} \frac{1}{E_k^2} \int_{\mathbf{B}_{\rho}(0)} \op{dist}^2(X,P_k) \,d\|T_k\|(X) \nonumber 
\end{align}
for all $\rho \in (0,4)$.  (Note that here and throughout the proof we use the fact that $T_k \llcorner \mathbf{B}_{\rho}(0) = \widetilde{T}_k \llcorner \mathbf{B}_{\rho}(0)$ for each $\rho \in (0,2-2\theta_k]$.)  Similarly, using $\widehat{A}_k/E_k \rightarrow Dw_a(0)$, 
\begin{equation}\label{not decaying eqn4}
	\int_{B_{\rho}(0)} |w - Dw_a(0) \cdot x|^2 
	= \lim_{k \rightarrow \infty} \frac{1}{E_k^2} \int_{\mathbf{B}_{\rho}(0)} \op{dist}^2(X,\widehat{P}_k) \,d\|T_k\|(X) 
\end{equation}
for all $\rho \in (0,4)$. 

We claim that $w$ is not identically zero.  We may assume that $Dw_a(0) = 0$ since otherwise $w$ must be non-zero.  By \eqref{not decaying eqn3}, \eqref{not decaying hyp5}, and $\widehat{A}_k/E_k \rightarrow Dw_a(0) = 0$ 
\begin{align}\label{not decaying eqn5}
	\int_{B_2(0)} |w|^2 
	&= \int_{B_2(0)} |w - Dw_a(0) \cdot x|^2 
	\\&= \lim_{k \rightarrow \infty} \frac{1}{E_k^2} \int_{\mathbf{B}_2(0)} \op{dist}^2(X,P_k) \,d\|T_k\|(X) \nonumber 
	\\&\geq \lim_{k \rightarrow \infty} \frac{2^{-n-2-2\alpha_k} }{E_k^2} \int_{\mathbf{B}_4(0)} \op{dist}^2(X,\widehat{P}_k) \,d\|T_k\|(X) \nonumber 
	\\&\geq \lim_{k \rightarrow \infty} \frac{2^{-n-3-2\alpha_k}}{E_k^2} \left( \int_{\mathbf{B}_4(0)} \op{dist}^2(X,\widetilde{P}) \,d\|T_k\|(X) 
		- 2^{2n+1} (q+1) \,\omega_n \|\widehat{A}_k\|^2 \right) \nonumber 
	\\&= 2^{-n-3} \nonumber  
\end{align}
and thus $w$ must be non-zero.

Next we claim that $w$ is homogeneous degree one.  We may assume that $w-Dw_a(0) \cdot x$ is not identically zero since otherwise $w(x) = q \llbracket Dw_a(0) \cdot x \rrbracket$ is certainly homogeneous degree one.   Since $\Theta(\widetilde{T}_k,0) = \Theta(T_k,0) \geq q,$ Lemma~\ref{hardt simon lemma} implies that $\mathcal{N}_w(0) \geq 1$.  Thus by the homogeneity of $Dw_a(0) \cdot x$, $\mathcal{N}_{w - Dw_a(0) \cdot x}(0) \geq 1$.  By $L^2$-growth estimates for Dirichlet energy minimizing multi-valued functions (\cite[Theorem~2.6(8)]{Almgren}), 
\begin{align}\label{not decaying eqn6}
	\int_{B_2(0)} |w - Dw_a(0) \cdot x|^2 
	&\leq 2^{-n-2\mathcal{N}_{w - Dw_{a}(0) \cdot x}(0)} \int_{B_4(0)} |w - Dw_a(0) \cdot x|^2 
	\\&\leq 2^{-n-2} \int_{B_4(0)} |w - Dw_a(0) \cdot x|^2 \nonumber 
\end{align}
with equality if and only if $w$ is homogeneous degree one.  On the other hand, by \eqref{not decaying eqn3}, \eqref{not decaying eqn4}, and \eqref{not decaying hyp5}
\begin{align*}
	\int_{B_2(0)} |w - Dw_a(0) \cdot x|^2 
	&= \lim_{k \rightarrow \infty} \frac{1}{E_k^2} \int_{\mathbf{B}_2(0)} \op{dist}^2(X,P_k) \,d\|T_k\|(X) \nonumber 
	\\&\geq \lim_{k \rightarrow \infty} \frac{2^{-n-2-2\alpha_k}}{E_k^2} \int_{\mathbf{B}_{\rho}(0)} \op{dist}^2(X,\widehat{P}_k) \,d\|T_k\|(X) 
	\\&= 2^{-n-2} \int_{B_{\rho}(0)} |w - Dw_a(0) \cdot x|^2 
\end{align*}
for all $\rho \in (0,4)$.  Letting $\rho \rightarrow 4^+$, 
\begin{equation}\label{not decaying eqn7}
	\int_{B_2(0)} |w - Dw_a(0) \cdot x|^2 \geq 2^{-n-2} \int_{B_4(0)} |w - Dw_a(0) \cdot x|^2 .
\end{equation}
Therefore, equality holds true in \eqref{not decaying eqn6}, proving that $w$ is homogeneous degree one. 

Since $N_{\widetilde{T}_k,\widetilde{P},Z_k}(\rho_k) \geq 1+\nu/2$, by Lemma~\ref{blowup frequency lemma} and the homogeneity of $w$, we have that 
\begin{equation*} 
	1 + \nu/2 \leq \lim_{k \rightarrow \infty} N_{\widetilde{T}_k,\widetilde{P},Z_k}(\rho_k) = N_{w,0}(\rho_{\infty}) = 1
\end{equation*}
giving us a contradiction. 
\end{proof}

Next we show in Lemma~\ref{not decaying lemma2} that if $T$ is not decaying to a plane as in \eqref{not decaying hyp5}, then $T$ must be $L^2$-close to a non-planar area-minimizing cone in $\mathbf{B}_1(0)$.

\begin{lemma}\label{not decaying lemma2}
For every $\varepsilon \in (0,1)$ and $\beta \in (0,1)$ there exists $\delta = \delta(n,m,q,\varepsilon,\beta) \in (0,1)$, $\eta = \eta(n,m,q,\varepsilon,\beta) \in (0,1)$, and $\alpha = \alpha(n,m,q,\varepsilon,\beta) \in (0,1)$ such that if $P,\widehat{P}$ are $n$-dimensional planes in $\mathbb{R}^{n+m}$ and $T$ is $n$-dimensional locally area-minimizing rectifiable current in $\mathbf{B}_8(0)$ such that \eqref{not decaying hyp1}--\eqref{not decaying hyp5} hold true, ${\rm spt} \, T \cap {\mathbf B}_{4}(0) \not\subset \widehat{P}$ and $\Theta(T,0) \geq q$, then for each $\rho \in [1/2,1]$ one of the following holds true: 
\begin{enumerate}[itemsep=3mm,topsep=0mm]
	\item[{\rm (i)}]  there exists $\mathbf{C} \in \mathcal{C}_q$ such that \eqref{freq2cone concl1}, \eqref{freq2cone concl2}, and \eqref{freq2cone concl3} hold true or 
	\item[{\rm (ii)}]  there is an $(n-3)$-dimensional linear subspace $L$ such that \eqref{freq2cone concl4} holds true. 
\end{enumerate}
\end{lemma}

\begin{proof}
Without loss of generality suppose that $P = \mathbb{R}^n \times \{0\}$.  Fix $\varepsilon \in (0,1)$ and $\beta \in (0,1)$.  Suppose that for $k = 1,2,3,\ldots$ there exists $\delta_k \rightarrow 0^+$, $\eta_k \rightarrow 0^+$, $\alpha_k \rightarrow 0^+$, and an $n$-dimensional area-minimizing integral current $T_k$ of $\mathbf{B}_8(0)$ such that \eqref{not decaying hyp1}--\eqref{not decaying hyp5} hold true with $\delta_k$, $\eta_k$, $\alpha_k$, and $T_k$ in place of $\delta$, $\eta$, $\alpha$, and $T$.  We want to show that for infinitely many $k$ either conclusion~(i) or conclusion~(ii) holds true, which in light of arbitrary choice of sequence $(T_k)$ will prove the lemma.  Arguing as in Subsection~\ref{sec:blowup procedure}, by \eqref{not decaying hyp1} and \eqref{not decaying hyp4} after passing to a subsequence we may assume that $T_k \rightarrow q \llbracket P \rrbracket$ weakly in $\mathbf{B}_8(0)$.  Set $E_k = E(T_k,\widetilde{P},\mathbf{B}_4(0))$ and let $\theta_k \in (0,1)$ such that $\theta_k \rightarrow 1^+$ and $(1-\theta_k)^{-n-2} E_k^2 \rightarrow 0$.  By \eqref{not decaying hyp4} and Lemma~\ref{coarse L2 distance lemma} with $K = P$, 
\begin{equation*}
	\sup_{X \in \op{spt} T_k \cap \mathbf{B}_{3+\theta_k}(0)} \op{dist}(X,P) \leq 8 E_k^{\frac{1}{n+2}} < 1-\theta_k .
\end{equation*}
Set $\widetilde{T}_k = T_k \llcorner \mathbf{B}_{3+\theta_k}(0)$ so that $(\partial \widetilde{T}_k) \llcorner \mathbf{C}_{2+2\theta_k}(0) = 0$ and $\limsup_{k \rightarrow \infty} \|\widetilde{T}_k\|(\mathbf{C}_{2+2\theta_k}(0)) \leq q \omega_n 4^n$ (as in \eqref{not decaying eqn2}).  By Theorem~\ref{lipschitz approx thm} and Lemma~\ref{area excess to height lemma}, there exists $u_k : B_{4\theta_k}(0) \rightarrow \mathcal{A}_q(\mathbb{R}^m)$, and $K_k \subset B_{4\theta_k}(0)$ such that \eqref{blowup eqn4}--\eqref{blowup eqn6} hold true with $\widetilde{T}_k$ in place of $T_k$, $\rho_0 = 4$, and $E_k = E(T_k,P,\mathbf{B}_4(0))$.  After passing to a subsequence, let $w \in W^{1,2}_{\rm loc}(B_4(0),\mathcal{A}_q(\mathbb{R}^m))$ such that \eqref{blowup eqn9} and \eqref{blowup eqn10} hold true with $\rho_0 = 4$ and $E_k = E(T_k,\widetilde{P},\mathbf{B}_4(0))$.  
By \eqref{not decaying hyp2} and Lemma~\ref{optimal plane lemma}, $Dw_a(0) = 0$.  By \eqref{not decaying hyp3} and Lemma~\ref{optimal plane lemma}, $\widehat{P}_k = \{ (x, \widehat{A}_k x) : x \in \mathbb{R}^n \}$ for some $m \times n$ matrix $\widehat{A}_k$ such that $\|\widehat{A}_k\| \leq C(n,m,q,\theta) \,E_k$ and $\widehat{A}_k/E_k \rightarrow Dw_a(0) = 0$.  Since $\Theta(T_k,0) \geq q,$ Lemma~\ref{hardt simon lemma} implies that $\mathcal{N}_w(0) \geq 1$.  Thus by $L^2$-growth estimates for Dirichlet energy minimizing multi-valued functions~\cite[Theorem~2.6(8)]{Almgren}, 
\begin{align*}
	\int_{B_2(0)} |w|^2 \leq 2^{-n-2\mathcal{N}_w(0)} \int_{B_4(0)} |w|^2 \leq 2^{-n-2} \int_{B_4(0)} |w|^2  
\end{align*}
Arguing as we did to obtain \eqref{not decaying eqn5} and \eqref{not decaying eqn7} with $P$ in place of $\widetilde{P}$, $A_k = 0$, and $Dw_a(0) = 0$ (since we now blow-up relative to $P$ instead of $\widetilde{P}$), it follows from \eqref{not decaying hyp5} that 
\begin{equation}\label{not decaying eqn8}
	\int_{B_2(0)} |w|^2 \geq 2^{-n-3} , \quad 
	\int_{B_2(0)} |w|^2 \geq 2^{-n-2} \int_{B_4(0)} |w|^2 . 
\end{equation}
Therefore $w$ is non-zero and homogeneous degree one.  Since $w$ is homogeneous degree one and $Dw_a(0) = 0$, $w$ is average-free.  If $\dim \op{spine} w = n-2$, then set $\mathbf{C}_k = \op{graph}(E_k w)$.  Arguing as in the proof of Lemma~\ref{freq2cone lemma}, \eqref{freq2cone eqn8} and \eqref{freq2cone eqn9} hold true.  

Let $\rho \in [1/2,1]$ and $\widetilde{P}_k$ be an optimal plane for $T_k$ in $\mathbf{B}_{\rho}(0)$.  By Lemma~\ref{optimal plane lemma},  
$\widetilde{P}_k = \{ (x, \widetilde{A}_k x) : x \in \mathbb{R}^n \}$ for some $m \times n$ matrix $\widetilde{A}_k$ with $\|\widetilde{A}_k\| \leq C(n,m,q) \,E_k$ and $\widetilde{A}_k/E_k \rightarrow Dw_a(0) = 0$.  Recall from the proof of Lemma~\ref{optimal plane lemma} that for each $X = (x,y) \in \mathbf{B}_2(0)$
\begin{equation*}
	\op{dist}(X, \op{spt} \widetilde{P}_k) \leq |y - \widetilde{A}_k x| \leq (1 + CE_k) \op{dist}(X, \op{spt} \widetilde{P}_k) .
\end{equation*}
Hence using \eqref{blowup eqn4}--\eqref{blowup eqn6}, \eqref{blowup eqn9}, \eqref{lipschitz approx eqn8}, \eqref{lipschitz approx eqn9}, 
the fact that $\widetilde{A}_k/E_k \rightarrow 0$, and \eqref{not decaying eqn8}, 
\begin{align}\label{not decaying eqn9}
	\lim_{k \rightarrow \infty} \frac{1}{E_k^2} \int_{\mathbf{B}_{\rho}(0)} \op{dist}^2(X, \op{spt} \widetilde{P}_k) \,d\|T_k\|(X) 
	\geq\,& \lim_{k \rightarrow \infty} \frac{q}{E_k^2} \int_{B_{\rho}(0)} \mathcal{G}(u_k(x), q \llbracket \widetilde{A}_k x \rrbracket)^2 \,dx 
	\\ =\,& q \int_{B_{\rho}(0)} |w|^2 \geq 4^{-n-2} q\int_{B_2(0)} |w|^2 \,dx \geq 2^{-3n-7}q.  \nonumber
\end{align}
(Recall that $T_k \llcorner \mathbf{B}_{\rho}(0) = \widetilde{T}_k \llcorner \mathbf{B}_{\rho}(0)$.)  Therefore, combining \eqref{freq2cone eqn8}, \eqref{freq2cone eqn9}, and \eqref{not decaying eqn9} we conclude that \eqref{freq2cone concl1} and \eqref{freq2cone concl2} hold true.  Arguing as in the proof of Lemma~\ref{freq2cone lemma} using Lemma~\ref{hardt simon lemma} gives us that \eqref{freq2cone concl3} holds true if $\dim \op{spine} \mathbf{C} = n-2$ and  \eqref{freq2cone concl3} holds true with $L = \op{spine} \mathbf{C}$ if $\dim \op{spine} \mathbf{C} \leq n-3$.
\end{proof}

\section{Proof of the main theorem}\label{sec:nonplanar cones sec} 

In this section we prove our main result, Theorem~\ref{branch and cones thm-intro}. We start with the following:  

\begin{theorem}\label{persistance of cones thm}
For every $\varepsilon \in (0,1)$ and $\beta \in (0,1)$ there exists $R = R(n,m,q,\varepsilon,\beta) \in [2,\infty)$, $\delta = \delta(n,m,q,\varepsilon,\beta) \in (0,1)$, $\eta = \eta(n,m,q,\varepsilon,\beta) \in (0,1)$, and $\alpha = \alpha(n,m,q,\varepsilon,\beta) \in (0,1)$ such that the following holds true.  Let $T$ be $n$-dimensional locally area-minimizing rectifiable current in $\mathbf{B}_{9R}(0)$ such that 
\begin{equation}\label{persist cones hyp1} 
	(\partial T) \llcorner \mathbf{B}_{9R}(0) = 0, \quad\quad \|T\|(\mathbf{B}_{9R}(0)) \leq (q + \delta) \,\omega_n (9R)^n . 
\end{equation}
Let $P$ be an $n$-dimensional plane in $\mathbb{R}^{n+m}$ such that 
\begin{equation}\label{persist cones hyp2} 
	E(T,P,\mathbf{B}_R(0)) = \inf_{P' \in \mathcal{P}} E(T,P',\mathbf{B}_R(0)) 
\end{equation}
and assume that one of the following holds true: 
\begin{enumerate}
	\item[{\rm (A)}] $E(T,P,\mathbf{B}_{4R}(0)) > \eta$; \vspace{3mm}
	\item[{\rm (B)}] $E(T,P,\mathbf{B}_{4R}(0)) \leq \eta$ and $E(T,P,\mathbf{B}_{2R}(0)) > 2^{-\alpha} E(T,\widehat{P},\mathbf{B}_{4R}(0))$, where $\widehat{P}$ is an $n$-dimensional plane in $\mathbb{R}^{n+m}$ such that $E(T,\widehat{P},\mathbf{B}_{2R}(0)) = \inf_{P' \in \mathcal{P}} E(T,P',\mathbf{B}_{2R}(0))$.
\end{enumerate}
Then for each $Z \in \op{spt} T \cap \mathbf{B}_1(0)$ with $\Theta(T,Z) \geq q$ and each $\rho \in (0,R]$ one of the following statements {\rm  (i)} or {\rm (ii)} holds true: 
\begin{enumerate}
	\item[{\rm (i)}]  there exists $\mathbf{C} \in \mathcal{C}_q$ such that 
	\begin{gather}
		\label{persist cones concl1} E(T,\mathbf{C},\mathbf{B}_{\rho}(Z)) < \varepsilon, \\
		\label{persist cones concl2} Q(T,\mathbf{C},\mathbf{B}_{\rho}(Z)) \leq \beta \inf_{P' \in \mathcal{P}} E(T,P',\mathbf{B}_{\rho}(Z)) \;\; and\\
		 \label{persist cones concl3} \{ X \in \mathbf{B}_{\rho}(Z) : \Theta(T,X) \geq q \} \subset 
			\{X \in \mathbf{B}_{\rho}(Z) : \op{dist}(X,Z+\op{spine}\mathbf{C}) < \varepsilon \rho \};
 \end{gather}
	\item[{\rm (ii)}]  there is an $(n-3)$-dimensional linear subspace $L$ of $\mathbb{R}^{n+m}$ such that 
	\begin{equation}\label{persist cones concl3}
		\{ X \in \mathbf{B}_{\rho}(Z) : \Theta(T,X) \geq q \} \subset \{ X \in \mathbf{B}_{\rho}(Z) : \op{dist}(X,Z+L) < \varepsilon \rho \} . 
	\end{equation}
\end{enumerate}
\end{theorem}

\begin{proof}
First note that for any $Z \in \op{spt} T \cap \mathbf{B}_1(0)$ with $\Theta(T,Z) \geq q,$ we have by the monotonicity formula for area and the fact that $\|T\|(\mathbf{B}_{9R}(0)) \leq (q + \delta) \,\omega_n (9R)^n$,  
\begin{equation*}
	q \leq \frac{\|T\|(\mathbf{B}_{\rho}(Z))}{\omega_n \rho^n} 
	\leq \frac{\|T\|(\mathbf{B}_{9R-1}(Z))}{\omega_n (9R-1)^n} \leq \left(\frac{9R}{9R-1}\right)^n \frac{\|T\|(\mathbf{B}_{9R}(0))}{\omega_n (9R)^n} 
	\leq \left(\frac{9R}{9R-1}\right)^n (q + \delta) 
\end{equation*}
for all $\rho \in (0,9R-1]$.  Thus, provided $R = R(n,q,\delta) \geq 2$ is large enough that 
\begin{equation}\label{persist cones R}
	\left(\frac{9R}{9R-1}\right)^n (q+\delta) \leq q + 2\delta,
\end{equation}
we have that 
\begin{equation}\label{persist cones area}
	q \leq \frac{\|T\|(\mathbf{B}_{\rho}(Z))}{\omega_n \rho^n} \leq q + 2\delta 
\end{equation}
for all  $Z \in \op{spt} T \cap \mathbf{B}_1(0)$ with $\Theta(T,Z) \geq q$ and all $\rho \in (0,9R-1]$.

Fix $\varepsilon \in (0,1)$ and $\beta \in (0,1)$.  Choose $\nu$ and $\alpha$ so that 
\begin{align}
	\label{persist parameters1} \nu &= \nu_{\ref{freq2cone lemma}}(n,m,q,\varepsilon,\beta) , \\
	\label{persist parameters2} \alpha &= \min\{\alpha_{\ref{not decaying lemma1}}(n,m,q,\nu), 
		\alpha_{\ref{not decaying lemma2}}(n,m,q,\varepsilon,\beta)\} ,
\end{align}
where $\nu_{\ref{freq2cone lemma}}$ is the constant denoted $\nu$ in Lemma~\ref{freq2cone lemma} and $\alpha_{\ref{not decaying lemma1}}$ and $\alpha_{\ref{not decaying lemma2}}$ are the constants denoted $\alpha$  in Lemma~\ref{not decaying lemma1} and Lemma~\ref{not decaying lemma2} respectively.  Then choose $\eta$ and $\delta$ so that 
\begin{align}
	\label{persist parameters3} \eta &\leq \min\big\{ 64^{-n-2} C_0^{-1}, (7/8)^{(n+2)/2} C_0^{-1} \eta_{\ref{mono freq thm}}(n,m,q), 
		4^{-2\alpha} C_0^{-1} \,\eta_{\ref{freq2cone lemma}}(n,m,q,\varepsilon,\beta), \\&\hspace{15.7mm} \eta_{\ref{not decaying lemma1}}(n,m,q,\nu), 
		\eta_{\ref{not decaying lemma2}}(n,m,q,\varepsilon,\beta) \big\}; \nonumber \\
	\label{persist parameters4} 2\delta &\leq \min\bigg\{ \frac{1}{2}, \delta_{\ref{freq2cone lemma}}(n,m,q,\varepsilon,\beta) , \delta_{\ref{not flat lemma1}}(n,m,q,\eta,\nu), 
		\delta_{\ref{not flat lemma2}}(n,m,q,\eta,\varepsilon,\beta), \\&\hspace{15.7mm} \delta_{\ref{not decaying lemma1}}(n,m,q,\nu), 
		\delta_{\ref{not decaying lemma2}}(n,m,q,\varepsilon,\beta) \bigg\}; \nonumber 
\end{align} 
\begin{gather}
	\label{persist parameters5} (q+2\delta) (1 + C_1 \eta)^n < q + \delta_{\ref{mono freq thm}}(n,m,q) \;\; \mbox{and}\\
	\label{persist parameters6} e^{C\left((8/7)^{(n+2)/2} C_0 \eta\right)^{\gamma}} < \frac{1+\nu}{1+\nu/2} .
\end{gather}
Here $C = C(n, m, q, \alpha)$ is as in Theorem~\ref{mono freq thm}; $C_0 = C_{0}(n, q, \alpha)$ is as in \eqref{persist cones eqn7} below; $C_1 = C_{1}(n, m, q, \alpha)$ is as in \eqref{persist cones eqn8} below; $\eta_{\ref{mono freq thm}}$, $\eta_{\ref{freq2cone lemma}}$, $\eta_{\ref{not decaying lemma1}}$ and $\eta_{\ref{not decaying lemma2}}$ are the constants denoted $\eta$ in Theorem~\ref{mono freq thm}, Lemma~\ref{freq2cone lemma}, Lemma~\ref{not decaying lemma1} and Lemma~\ref{not decaying lemma2} respectively;  $\delta_{\ref{mono freq thm}}$, $\delta_{\ref{freq2cone lemma}}$, $\delta_{\ref{not flat lemma1}}$, $\delta_{\ref{not flat lemma2}}$, $\delta_{\ref{not decaying lemma1}}$ and $\delta_{\ref{not decaying lemma2}}$ are the constants denoted $\delta$ in Theorem~\ref{mono freq thm}, Lemma~\ref{freq2cone lemma}, Lemma~\ref{not flat lemma1}, Lemma~\ref{not flat lemma2}, Lemma~\ref{not decaying lemma1} and Lemma~\ref{not decaying lemma2} respectively.  Finally, choose $R$ so that 
\begin{equation}
	\label{persist parameters7} R \geq 1/\min\{ \delta_{\ref{not flat lemma1}}(n,m,q,\eta,\nu) , \delta_{\ref{not decaying lemma1}}(n,m,q,\nu) \} 
		\text{and with $R$ also satisfying \eqref{persist cones R}} .
\end{equation}

Fix $Z \in \op{spt} T \cap \mathbf{B}_1(0)$ with $\Theta(T,Z) \geq q$.  For each integer $i \geq 0$ set $\rho_i = 2^{-i} R$ and $I_i = [\rho_{i+1},\rho_i]$ so that the interval $(0,R]$ is partitioned into the intervals $I_i$.  Let $P_0 = P$ be such that \eqref{persist cones hyp2} holds true.  (Note that unlike in our notation in earlier sections, $P_{0}$ here is not necessarily ${\mathbb R}^{n} \times \{0\}$). For each integer $i \geq 1$ let $P_i$ be an $n$-dimensional plane in $\mathbb{R}^{n+m}$ such that 
\begin{equation}\label{persist cones eqn1}
	E(T,P_i,\mathbf{B}_{\rho_i}(Z)) = \inf_{P' \in \mathcal{P}} E(T,P',\mathbf{B}_{\rho_i}(Z)).
\end{equation}
For each integer $i \geq 1$ we classify each interval $I_i$ into one of the following three cases: 
\begin{enumerate}[itemsep=3mm,topsep=0mm]
	\item[(a)] $E(T,P_i,\mathbf{B}_{4\rho_i}(Z)) > \eta$; 
	\item[(b)] $E(T,P_i,\mathbf{B}_{4\rho_i}(Z)) \leq \eta$ and $E(T,P_i,\mathbf{B}_{2\rho_i}(Z)) > 2^{-\alpha} E(T,P_{i-1},\mathbf{B}_{2\rho_{i-1}}(Z))$; 
	\item[(c)] $E(T,P_i,\mathbf{B}_{4\rho_i}(Z)) \leq  \eta$ and $E(T,P_i,\mathbf{B}_{2\rho_i}(Z)) \leq  2^{-\alpha} E(T,P_{i-1},\mathbf{B}_{2\rho_{i-1}}(Z))$. 
\end{enumerate}
Recall that $(\partial T) \llcorner \mathbf{B}_{9R}(0) = 0$, $\Theta(T,Z) \geq q$, and \eqref{persist cones area} holds true.  In case~(a) we can apply Lemma~\ref{not flat lemma2} to deduce that either conclusion~(i) or conclusion~(ii) holds true for all $\rho \in I_i$.  In case~(b), noting that 
${\rm spt} \, T \cap {\mathbf B}_{4\rho_{i}}(Z)  = {\rm spt} \, T \cap {\mathbf B}_{2\rho_{i-1}}(Z) \not\subset P_{i-1}$ (for otherwise by \eqref{persist cones eqn1} 
we would have $P_{i} = P_{i-1}$ and hence $E(T,P_i,\mathbf{B}_{2\rho_i}(Z)) = 0$, which is ruled out by the strict inequality in (b)), we can apply Lemma~\ref{not decaying lemma2} to obtain either conclusion~(i) or conclusion~(ii) for all $\rho \in I_i$.  It remains to consider intervals satisfying case~(c).

Let $I = [\rho_{i_1+1},\rho_{i_0}]$ where $1 \leq i_0 \leq i_1$ are integers such that either
\begin{itemize}[itemsep=3mm,topsep=0mm]
	\item $i_0 = 1$ or 
	\item $i_0 > 1$ and either case~(a) or case~(b) holds true when $i = i_0-1$
\end{itemize}	
and
\begin{itemize}[itemsep=3mm,topsep=0mm]
	\item case~(c) holds true whenever $i \in \{i_0,i_0+1,\ldots,i_1\}$.
\end{itemize}
By iteratively applying the inequality $E(T,P_i,\mathbf{B}_{2\rho_i}(Z)) \leq 2^{-\alpha} E(T,P_{i-1},\mathbf{B}_{2\rho_{i-1}}(Z))$ and using the fact that $E(T,P_{i_0},\mathbf{B}_{4\rho_{i_0}}(Z)) \leq \eta$, we obtain for each $i \in \{i_0,i_0+1,\ldots,i_1\},$ 
\begin{equation}\label{persist cones eqn2}
	E(T,P_i,\mathbf{B}_{2\rho_i}(Z)) \leq 2^{-\alpha(i-i_0)} E(T,P_{i_0},\mathbf{B}_{2\rho_{i_0}}(Z)) = \Big(\frac{\rho_i}{\rho_{i_0}}\Big)^{\alpha}E(T,P_{i_0},\mathbf{B}_{2\rho_{i_0}}(Z)) \leq C\eta \Big(\frac{\rho_i}{\rho_{i_0}}\Big)^{\alpha},
\end{equation}
where $C = C(n) \in (0,\infty)$ is a constant.  By the triangle inequality, \eqref{persist cones area}, Lemma~\ref{tilt to height estimate lemma}, and \eqref{persist cones eqn2} 
\begin{align}\label{persist cones eqn3}
	\|\pi_{P_i} - \pi_{P_{i+1}}\|^2 
	\leq\,& \frac{2}{q \omega_n (\rho_i/2)^n} \int_{G_n(\mathbf{B}_{\rho_i}(Z))} \|\pi_S - \pi_{P_i}\|^2 \,d|T|(X,S) 
		\\&+ \frac{2}{q \omega_n \rho_{i+1}^n} \int_{G_n(\mathbf{B}_{\rho_{i+1}}(Z))} \|\pi_S - \pi_{P_{i+1}}\|^2 \,d|T|(X,S) \nonumber 
	\\ \leq\,& C E(T,P_i,\mathbf{B}_{2\rho_i}(Z))^2 + C E(T,P_{i+1},\mathbf{B}_{2\rho_{i+1}}(Z))^2 \nonumber 
	\\ \leq \,& C \Big(\frac{\rho_i}{\rho_{i_0}}\Big)^{2\alpha} E^{2}(T,P_{i_0},\mathbf{B}_{2\rho_{i_0}}(Z))\nonumber 
         \\ \leq\,& C \eta^2 \Big(\frac{\rho_i}{\rho_{i_0}}\Big)^{2\alpha} \nonumber 
\end{align}
for all $i \in \{i_0,i_0+1,\ldots,i_1-1\}$, where $C = C(n,q) \in (0,\infty)$ are constants.  Hence by the triangle inequality, 
\begin{align}\label{persist cones eqn4}
	\|\pi_{P_i} - \pi_{P_j}\| \leq & C E(T,P_{i_0},\mathbf{B}_{2\rho_{i_0}}(Z))\sum_{k=i}^{j-1} \,\Big(\frac{\rho_k}{\rho_{i_0}}\Big)^{\alpha} 
		= C E(T,P_{i_0},\mathbf{B}_{2\rho_{i_0}}(Z))\sum_{k=i}^{j-1} 2^{\alpha(i-k)} \,\Big(\frac{\rho_i}{\rho_{i_0}}\Big)^{\alpha} \nonumber
		\\ \leq \,& \frac{2^{\alpha} C}{2^{\alpha}-1} \,\Big(\frac{\rho_i}{\rho_{i_0}}\Big)^{\alpha}E(T,P_{i_0},\mathbf{B}_{2\rho_{i_0}}(Z)) 
\end{align}
for all $i_0 \leq i \leq j \leq i_1$, where $C = C(n,q) \in (0,\infty)$ is a constant.  For each $i_0 \leq j \leq i \leq i_1$, noting that 
\begin{align*}
	\op{dist}(X,Z+P_i) = |\pi_{P_i^{\perp}} (X-Z)| &\leq |\pi_{P_j^{\perp}} (X-Z)| + |(\pi_{P_i} - \pi_{P_j}) (X-Z)|
		\\&\leq \op{dist}(X,Z+P_j) + \|\pi_{P_i} - \pi_{P_j}\| \,|X-Z|
\end{align*}
for all $X \in \op{spt} T$ and using \eqref{persist cones area}, \eqref{persist cones eqn2}, and \eqref{persist cones eqn4}, we obtain 
\begin{align}\label{persist cones eqn5}
	E^{2}(T,P_i,\mathbf{B}_{2\rho_j}(Z)) = & \frac{1}{\omega_n (2\rho_j)^{n+2}} \int_{\mathbf{B}_{2\rho_j}(Z)} \op{dist}^2(X,Z+P_i) \,d\|T\|(X) \\
	\leq\,& \frac{2}{\omega_n (2\rho_j)^{n+2}} \int_{\mathbf{B}_{2\rho_j}(Z)} \op{dist}^2(X,Z+P_j) \,d\|T\|(X) + 2 (q+1) \,\|\pi_{P_i} - \pi_{P_j}\|^2 \nonumber\\ 
	\leq \, &C \Big(\frac{\rho_j}{\rho_{i_0}}\Big)^{2\alpha} E^{2}(T,P_{i_0},\mathbf{B}_{2\rho_{i_0}}(Z)) \leq C \eta^{2}\Big(\frac{\rho_j}{\rho_{i_0}}\Big)^{2\alpha}, \nonumber 
\end{align}
where $C = C(n,q,\alpha) \in (0,\infty)$ is a constant.  Similarly, arguing using \eqref{persist cones area},  \eqref{persist cones eqn4}, and the fact that $E(T,P_{i_{0}},\mathbf{B}_{4\rho_{i_{0}}}(Z)) \leq \eta$, 
\begin{equation}\label{persist cones eqn6}
	E(T,P_i,\mathbf{B}_{4\rho_{i_0}}(Z)) \leq C \eta
\end{equation}
for all $i$ with $i_0 \leq i \leq i_1$, where $C = C(n,q,\alpha) \in (0,\infty)$ is a constant.  For each $i \in \{i_0,i_0+1,\ldots,i_1\}$ and $\rho \in [\rho_{i+1},2\rho_{i_0}]$ we can find an integer $j$ with $i_0-1 \leq j \leq i$ such that $\rho_{j+1} < \rho \leq \rho_j.$ Thus by \eqref{persist cones eqn5} and \eqref{persist cones eqn6} 
\begin{equation}\label{persist cones eqn7}
	E(T,P_i,\mathbf{B}_{2\rho}(Z)) \leq C_0 \eta \Big(\frac{\rho}{2\rho_{i_0}}\Big)^{\alpha}  
\end{equation}
for $\rho \in [\rho_{i+1},2\rho_{i_0}] = [\rho_{i+1}, \rho_{i_{0}-1}]$, where $C_0 = C_0(n,q,\alpha) \in (0,\infty)$ is a constant.  Note that by Lemma~\ref{Allard height lemma}, \eqref{persist cones eqn7} and \eqref{persist parameters3}
\begin{equation}\label{persist cones eqn8}
	\sup_{X \in \op{spt} T \cap \mathbf{B}_{15\rho/8}(Z)} \op{dist}(X,Z+P_i) \leq C_1 \eta \rho < \rho/4
\end{equation}
for all $\rho \in [\rho_{i+1},2\rho_{i_0}]$, where $C_1 = C_1(n,m,q,\alpha) \in (0,\infty)$ is a constant.  

Now fix any $i \in \{i_0,i_0+1,\ldots,i_1\}$.  Recall that $(\partial T) \llcorner \mathbf{B}_{9R}(0) = 0$, $\Theta(T,Z) \geq q$ and that \eqref{persist cones area} holds true.  Set $\widetilde{T} = T \llcorner {\mathbf B}_{15\rho_{i_{0}}/4}(Z)$.  Note that since $Z \in \mathbf{B}_1(0)$, we have by \eqref{persist parameters7} that 
$Z \in \mathbf{B}_{\min\{ \delta_{\ref{not flat lemma1}}(n,m,q,\eta,\nu) , \delta_{\ref{not decaying lemma1}}(n,m,q,\nu) \} 
R}(0)$.  By \eqref{persist cones eqn7}, $E(T,P_i,\mathbf{B}_{2\rho_{i_0-1}}(Z)) \leq C_0 \eta < 64^{-n-2}$.  Hence when $i_0 = 1$, by virtue of hypothesis~(A) or hypothesis~(B), we can apply respectively Lemma~\ref{not flat lemma1} or Lemma~\ref{not decaying lemma1} (noting also, in the latter case, 
that ${\rm spt} \, T \cap {\mathbf B}_{4R}(0) \not\subset {\widehat P}$ in view of the strict inequality in hypthesis~(B)) to obtain $N_{\widetilde{T},P_i,Z}(\rho) < 1+\nu/2$ for all $\rho \in I_0 = [R/2, R]$.  Arguing in the same way when $i_0 > 1$, we see by virtue of case~(a) or case~(b) that we can again apply Lemma~\ref{not flat lemma1} or Lemma~\ref{not decaying lemma1}  to obtain $N_{\widetilde{T},P_i,Z}(\rho) < 1+\nu/2$ for all $\rho \in I_{i_0-1} = [\rho_{i_{0}}, 2\rho_{i_{0}}]$.  Notice that by \eqref{persist cones eqn8}, \eqref{persist cones area} and \eqref{persist parameters5}, $(\partial\widetilde{T}) \llcorner {\mathbf C}_{7\rho_{i_{0}}/2}(Z,P_i) = 0$ and 
\begin{align*}
	\|\widetilde{T}\|({\mathbf C}_{7\rho_{i_{0}}/2}(Z,P_i)) \leq\,& \|\widetilde{T}\|({\mathbf B}_{(7/2+C_1\eta) \rho_{i_{0}}}(Z)) 
		\leq (q+2\delta) \,\omega_n (7/2+C_1\eta)^n \rho_{i_{0}}^n \\ \leq\,& (q+\delta_{\ref{mono freq thm}}) \,\omega_n (7\rho_{i_{0}}/2)^n . 
\end{align*}
By \eqref{persist cones eqn7} and \eqref{persist cones eqn8}, 
\begin{equation*}
	E(\widetilde{T},P_i,\mathbf{C}_{7\rho/4}(Z,P_i)) \leq (8/7)^{(n+2)/2} E(T,P_i,\mathbf{B}_{2\rho}(Z)) 
		\leq (8/7)^{(n+2)/2} C_0 \eta \Big(\frac{\rho}{2\rho_{i_0}}\Big)^{\alpha} 
\end{equation*}
for all $\rho \in [\rho_{i+1},2\rho_{i_0}]$.  Thus we can apply Theorem~\ref{mono freq thm} together with the fact that $N_{\widetilde{T},P_i,Z}(2\rho_{i_0}) < 1+\nu/2$ and \eqref{persist parameters6} to obtain 
\begin{align}\label{persist cones eqn9}
	N_{\widetilde{T},P_i,Z}(\rho) \leq\,& e^{C\left((8/7)^{(n+2)/2} C_0 \eta\right)^{\gamma} (\rho/(2\rho_{i_0}))^{\alpha\gamma}} N_{\widetilde{T},P_i,Z}(2\rho_{i_0}) 
		\\ \leq\,& e^{C\left((8/7)^{(n+2)/2} C_0 \eta\right)^{\gamma}} (1+\nu/2) \leq 1+\nu \nonumber 
\end{align}
for all $\rho \in [\rho_{i+1},2\rho_{i_0}]$.  Note that by \eqref{persist cones eqn8}, $\widetilde{T} = T \llcorner \mathbf{B}_{15\rho_i/4}(Z)$ in $\mathbf{C}_{2\rho_i}(Z,P_i)$.  Therefore, by \eqref{persist cones eqn1}, the fact that $E(T,P_i,\mathbf{B}_{4\rho_i}(Z)) \leq \eta$ and \eqref{persist cones eqn9}, we can apply Lemma~\ref{freq2cone lemma} to deduce that either conclusion~(i) or conclusion~(ii) holds true for all $\rho \in I_i$. 
\end{proof}

\begin{theorem}[Theorem~\ref{branch and cones thm-intro} in the introduction] \label{branch and cones thm}
For every $\varepsilon \in (0,1)$ and $\beta \in (0,1)$ there exists $R = R(n,m,q,\varepsilon,\beta) > 10$, $\delta = \delta(n,m,q,\varepsilon,\beta) \in (0,1)$ and $\alpha = \alpha(n,m,q,\varepsilon,\beta) \in (0,1)$ such that the following holds true.  Suppose that $T$ is an $n$-dimensional locally area-minimizing rectifiable current of $\mathbf{B}_R(0)$ such that 
\begin{equation}\label{branch and cones hyp1} 
	(\partial T) \llcorner \mathbf{B}_R(0) = 0, \quad\quad \|T\|(\mathbf{B}_R(0)) \leq (q + \delta) \,\omega_n R^n.
\end{equation}
Then 
\begin{equation*}
	\{ X \in \op{spt} T \cap \mathbf{B}_1(0) : \Theta(T,X) \geq q \} = \mathcal{S} \cup \mathcal{B}
\end{equation*}
where $\mathcal{B}$ and $\mathcal{S}$ are disjoint, locally compact sets for which the following hold:
\begin{enumerate}[itemsep=3mm,topsep=0mm]
	\item[{\rm (I)}]  For each $Z_0 \in \mathcal{S}$ there exists $\rho_0 \in (0,1]$ (depending on $T$ and $Z_0$) such that for every $Z \in \op{spt} T \cap \mathbf{B}_{\rho_0}(Z_0)$ with $\Theta(T,Z) \geq q$ and every $\rho \in (0,\rho_0]$ one of the following holds true:  
	\begin{enumerate}[itemsep=3mm,topsep=3mm]
	\item[{\rm (i)}]  there exists $\mathbf{C} \in \mathcal{C}_q$ such that 
	\begin{gather}
		\label{branch and cones concl1} E(T,\mathbf{C},\mathbf{B}_{\rho}(Z)) < \varepsilon, \\
		\label{branch and cones concl2} Q(T,\mathbf{C},\mathbf{B}_{\rho}(Z)) \leq \beta \inf_{P' \in \mathcal{P}} E(T,P',\mathbf{B}_{\rho}(Z)), \\
		\label{branch and cones concl3} \{ X \in \mathbf{B}_{\rho}(Z) : \Theta(T,X) \geq q \} \subset 
			\{ X \in \mathbf{B}_{\rho}(Z) : \op{dist}(X,Z+\op{spine}\mathbf{C}) < \varepsilon \rho \}; 
	\end{gather}
	\item[{\rm (ii)}]  there is an $(n-3)$-dimensional linear subspace $L$ of $\mathbb{R}^{n+m}$ such that 
	\begin{equation}\label{branch and cones concl4}
		\{ X \in \mathbf{B}_{\rho}(Z) : \Theta(T,X) \geq q \} \subset \{ X \in \mathbf{B}_{\rho}(Z) : \op{dist}(X,Z+L) < \varepsilon \rho \} . 
	\end{equation}
	\end{enumerate}
	
	\item[{\rm (II)}]  ${\mathcal B}$ is relatively closed in ${\mathbf B}_{1}(0),$ and:  
	\begin{enumerate}[itemsep=3mm,topsep=3mm]
	\item[{\rm (i)}]  if $Z_0 \in \mathcal{B}$ then $\Theta(T, Z_{0}) = q$ and there is a unique $n$-dimensional plane $P_{Z_0}$ such that 
	\begin{equation}\label{branch and cones concl5} 
		E(T,P_{Z_0},\mathbf{B}_{\sigma}(Z_0)) \leq C_0 \Big(\frac{\sigma}{\rho}\Big)^{\alpha} E(T, P_{Z_{0}}, \mathbf{B}_{\rho}(Z_{0}) ) 
	\end{equation}
	for all $0 < \sigma \leq \rho \leq 2$, where $C_{0} = C_{0}(n, m, q, \varepsilon, \beta) \in (0, \infty)$ is a constant;

	\item[{\rm (ii)}]  if $Z_{0} \in {\mathcal B}$ then $P_{Z_{0}}$ (as in {\rm (i)}) can be equipped with an orientation such that $q \llbracket P_{Z_0} \rrbracket$ is the unique tangent cone to $T$ at $Z_0$; 
	
	\item[{\rm (iii)}]  if $Z_{0} \in {\mathcal B}$ then the planar frequency $\mathcal{N}_{T,{\rm Pl}}(Z_0) = \lim_{\rho \rightarrow 0^+} N_{\widetilde{T},P_{Z_0},Z_0}(\rho)$ exists and $\mathcal{N}_{T,{\rm Pl}}(Z_0) \geq 1+\alpha$, where $\widetilde{T} = T \res {\mathbf B}_{15/8}(Z_{0})$; 
	 \item[{\rm (iv)}] for any two points $Z_{1}, Z_{2} \in {\mathcal B} \, \cap\, {\mathbf B}_{1/2}(0)$, we have that 
	$$\|\pi_{P_{Z_{1}}} - \pi_{P_{Z_{2}}}\| \leq C|Z_{1} - Z_{2}|^{\alpha}\left(E(T, P_{Z_{1}}, {\mathbf B}_{2}(Z_{1})) +E(T, P_{Z_{2}}, {\mathbf B}_{2}(Z_{2}))\right),$$ 
    where $C = C(n, m, q, \epsilon, \beta)$ is a constant.
	\end{enumerate}
\end{enumerate}
\end{theorem}

\begin{proof}[Proof of Theorem~\ref{branch and cones thm}]
Fix $\varepsilon \in (0,1)$ and $\beta \in (0,1)$ and let $\alpha$, $\eta$ be the constants (depending only on $n$, $m$, $q$, $\epsilon$ and $\beta$) 
as in Theorem~\ref{persistance of cones thm}. Denote by $\delta_{1}$ the constant $\delta = \delta(n, m, q, \epsilon, \beta) \in (0, 1)$ and by $R_{1}$ the constant $R = R(n, m, q, \epsilon, \beta)$ both given by 
Theorem~\ref{persistance of cones thm}. Set $\delta = \delta_{1}/2$, and with $R>10$ to be chosen, suppose that the hypotheses 
$(\partial T) \llcorner \mathbf{B}_R(0) = 0$ and $\|T\|(\mathbf{B}_R(0)) \leq (q + \delta) \,\omega_n R^n$ as in the present theorem hold. Then by the monotonicity formula for area (see the argument leading to \eqref{persist cones area}), we may choose $R = R(n,q,\delta) = R(n, m, q, \epsilon, \beta)$ sufficiently large to ensure 
\begin{equation}\label{branch and cones eqn1}
	q \leq \frac{\|T\|(\mathbf{B}_{\rho}(Z))}{\omega_n \rho^n} \leq q + \delta_{1} 
\end{equation}
for all $Z \in \op{spt} T \cap \mathbf{B}_1(0)$ with $\Theta(T,Z) \geq q$ and all $\rho \in (0,R-1].$  

Set $\rho_i = 2^{-i},$ $i=0, 1, 2, \ldots,$ and for each $Z_0 \in \op{spt} T \cap \mathbf{B}_1(0)$ with $\Theta(T,Z_{0}) \geq q$ and each $i = 0,1,2,\ldots$,  let $P_{Z_0,i}$ be an $n$-dimensional plane such that 
\begin{equation}\label{branch and cones eqn2}
	E(T,P_{Z_0,i},\mathbf{B}_{\rho_i}(Z_0)) = \inf_{P' \in \mathcal{P}} E(T,P',\mathbf{B}_{\rho_i}(Z_0))
\end{equation}
Define $\widetilde{\mathcal{S}}$ to be the set of all $Z_0 \in \op{spt} T \cap \mathbf{B}_1(0)$ with $\Theta(T,Z_0) \geq q$ such that  one of the following holds true: 
\begin{enumerate}[itemsep=3mm,topsep=0mm]
	\item[(A)] $E(T,P_{Z_0,i},\mathbf{B}_{4\rho_i}(Z_0)) > \eta$ for for some $i  = i(Z_{0}) \in \{0, 1,2,3, \ldots\}$; 
	\item[(B)] $E(T,P_{Z_0,i},\mathbf{B}_{4\rho_i}(Z_0)) \leq  \eta$ and $E(T,P_{Z_0,i},\mathbf{B}_{2\rho_i}(Z_0)) > 2^{-\alpha} E(T,P_{Z_0,i-1},\mathbf{B}_{4\rho_i}(Z_0))$ for some $i  = i(Z_{0}) \in \{1,2,3, \ldots\}$. 
\end{enumerate}
If $Z_{0} \in \widetilde{{\mathcal S}}$, then  we can set $\rho_0 = \rho_i/R_{1}$ and deduce from Theorem~\ref{persistance of cones thm} that either conclusion (I)(i) 
or conclusion I(ii) of the present theorem holds for each $Z \in \mathbf{B}_{\rho_0}(Z_0)$ with $\Theta(T,Z) \geq q$ and each 
$\rho \in (0,\rho_0]$, i.e.\ for each such $Z$ and $\rho$, either there exists $\mathbf{C} \in \mathcal{C}_q$ such that \eqref{branch and cones concl1}, \eqref{branch and cones concl2} and \eqref{branch and cones concl3} hold true or there is an $(n-3)$-dimensional linear subspace $L$ such that \eqref{branch and cones concl4} holds true. 

Now set 
\begin{eqnarray*}
&{\mathcal S} = \{Z_{0} \in {\mathbf B}_{1}(0) \, : \, \Theta(T, Z_{0}) \geq q, \; \exists \rho_{0} >0 \; \mbox{such that conclusion (I)(i) or conclusion (I)(ii) of}\nonumber\\
&\hspace{.5in}\mbox{the present theorem holds for every $Z \in {\mathbf B}_{\rho_{0}}(Z_{0})$ with $\Theta(T, Z) \geq q$ and every $\rho \in (0, \rho_{0}]$}\}
\end{eqnarray*}
and 
\begin{equation*}
	\mathcal{B} = \{ Z \in \op{spt} T \cap \mathbf{B}_1(0) : \Theta(T,Z) \geq q \} \setminus \mathcal{S}. 
\end{equation*}
It is clear (from the definition of ${\mathcal S}$) that ${\mathcal S}$ is relatively open in $\{Z \, : \, \Theta(T, Z) \geq q\} \cap {\mathbf B}_{1}(0)$, 
${\mathcal B}$ is relatively closed in $\{Z \, : \, \Theta(T, Z) \geq q\} \cap {\mathbf B}_{1}(0)$ and hence (since $\{Z \, : \, \Theta(T, Z) \geq q\}$ is closed) 
both sets are locally compact. Moreover, we have just shown that 
$\widetilde{\mathcal S} \subset {\mathcal S}$ and hence  if $Z_0 \in \mathcal{B}$ then we have that 
$Z_{0} \in {\mathbf B}_{1}(0)$, $\Theta(T, Z_{0}) \geq q,$   
\begin{gather}
	\label{branch and cones eqn3} E(T,P_{Z_0,i},\mathbf{B}_{4\rho_i}(Z_0)) \leq \eta \;\; \mbox{for each $i = 0, 1,2,\ldots$ and}   \\
	\label{branch and cones eqn4} E(T,P_{Z_0,i},\mathbf{B}_{2\rho_i}(Z_0)) \leq 2^{-\alpha} E(T,P_{Z_0,i-1},\mathbf{B}_{4\rho_i}(Z_0)) \;\; 
	\mbox{for each $i = 1,2,\ldots$}. 
\end{gather}

Fix $Z_0 \in \mathcal{B}$ and to simplify notation let $P_i = P_{Z_0,i}$.  Iteratively applying \eqref{branch and cones eqn4} and using \eqref{branch and cones eqn3} with $i=0$, we obtain for each $i \in \{0, 1,2,\ldots\},$ 
\begin{equation}\label{branch and cones eqn5}
	E(T,P_i,\mathbf{B}_{2\rho_i}(Z_0)) \leq 2^{-\alpha i} E(T,P_{Z_0,0},\mathbf{B}_2(Z_0)) \leq C \eta \,2^{-\alpha i},
\end{equation}
where $C = C(n) \in (0,\infty)$ is a constant.  

Arguing as we did to obtain \eqref{persist cones eqn3}, it follows from \eqref{branch and cones eqn5} that 
\begin{equation}\label{branch and cones eqn7}
	\|\pi_{P_i} - \pi_{P_j}\| \leq \,C\, 2^{-\alpha i} E(T,P_{Z_0,0},\mathbf{B}_2(Z_0))
\end{equation}
for all $0 \leq i \leq j$, where $C = C(n,q) \in (0,\infty)$ is a constant.  Hence the sequence $(P_i)$ is a Cauchy sequence and there must be an $n$-dimensional plane $P_{Z_0}$ such that $\|\pi_{P_i} - \pi_{P_{Z_0}}\| \rightarrow 0$ as $i \rightarrow \infty$.  Equivalently, $\op{dist}_{\mathcal{H}}(P_i \cap \mathbf{B}_1(0), P_{Z_0} \cap \mathbf{B}_1(0)) \rightarrow 0$ as $i \rightarrow \infty$.  By letting $j \rightarrow \infty$ in \eqref{branch and cones eqn7}, 
\begin{equation}\label{branch and cones eqn8}
	\|\pi_{P_i} - \pi_{P_{Z_0}}\| \leq  C 2^{-\alpha i}
E(T,P_{Z_0,0},\mathbf{B}_2(Z_0)) \leq C \eta \,2^{-\alpha i} 
\end{equation}
for each $i \geq 0$.  Noting that 
\begin{align*}
	\big| \op{dist}(X,Z_0+P_i) - \op{dist}(X,Z_0+P_{Z_0}) \big| 
	&= \left| |\pi_{P_i^{\perp}} (X-Z_0)| - |\pi_{P_{Z_0}^{\perp}} (X-Z_0)| \right| 
	\\&\leq \|\pi_{P_i} - \pi_{P_{Z_{0}}}\| \,|X-Z_0|
\end{align*}
for all $X \in \mathbb{R}^{n+m}$ and using the triangle inequality, \eqref{branch and cones eqn5}, \eqref{branch and cones eqn8}, and \eqref{branch and cones eqn1}, 
\begin{align}\label{branch and cones eqn9}
	&E^{2}(T,P_{Z_0},\mathbf{B}_{2\rho_i}(Z_0)) 
	= \frac{1}{\omega_n (2\rho_i)^{n+2}} \int_{\mathbf{B}_{2\rho_i}(Z_0)} \op{dist}^2(X,Z_0+P_{Z_0}) \,d\|T\|(X) 
	\\ \leq\,& \frac{2}{\omega_n (2\rho_i)^{n+2}} \int_{\mathbf{B}_{2\rho_i}(Z_0)} \op{dist}^2(X,Z_0+P_i) \,d\|T\|(X) 
		+ 2(q+1) \,\omega_n \,\|\pi_{P_i} - \pi_{P_{Z_0}}\|^2 \nonumber 
	\\ \leq\,& C \rho_{i}^{2\alpha} E^{2}(T,P_{Z_0,0},\mathbf{B}_2(Z_0))\leq C \eta^2 \rho_i^{2\alpha} \nonumber
\end{align}
for each $i \geq 0$, where $C = C(n,q,\alpha) \in (0,\infty)$ are constants.  For each $\rho \in (0,1]$ we can find an integer $i \geq 0$ such that $\rho_{i+1} < \rho \leq \rho_i$ and thus by \eqref{branch and cones eqn9}, we see that 
\begin{equation}\label{branch and cones eqn10}
    E(T,P_{Z_0},\mathbf{B}_{2\rho}(Z_0)) \leq C \eta\rho^{\alpha} 
\end{equation}
for all $\rho \in (0, 1]$, where $C = C(n,q,\alpha) \in (0,\infty)$ is a constant.

Set $\widetilde{T} = T \llcorner \mathbf{B}_{15/4}(Z)$.  Arguing as in the proof of Theorem~\ref{persistance of cones thm}, by Lemma~\ref{Allard height lemma} and \eqref{branch and cones eqn10} 
\begin{equation*}
    \sup_{X \in \op{spt} T \cap \mathbf{B}_{15\rho/8}(Z_0)} \op{dist}(X,Z_0+P_{Z_0}) \leq C \eta \rho < \rho/8 
\end{equation*}
for all $\rho \in (0, 1]$, where $C = C(n,q,\alpha) \in (0,\infty)$ is a constant, and it follows that $\widetilde{T}$ is a locally area-minimizing rectifiable current of $\mathbf{C}_{7/4}(Z_0,P_{Z_0})$ with $(\partial\widetilde{T}) \llcorner \mathbf{C}_{7/4}(Z_0,P_{Z_0}) = 0$, $\|\widetilde{T}\|(\mathbf{C}_{7/4}(Z_0,P_{Z_0})) \leq (q+2\delta) \,\omega_n (7/2+C\eta)^n \rho^n$, and $E(\widetilde{T},P_{Z_0},\mathbf{C}_{7\rho/4}(Z_0,P_{Z_0})) \leq C\eta\rho^{\alpha}$, where $C = C(n,q,\alpha) \in (0,\infty)$ are constants.  Hence as was discussed in Lemma~\ref{planar frequency lemma}, we can equip $P_{Z_0}$ with an orientation such that $\eta_{Z_0,\rho\#} T \rightarrow q \llbracket P_{Z_0} \rrbracket$ weakly in $\mathbb{R}^{n+m}$ as $\rho \rightarrow 0^+$; that is, $q \llbracket P_{Z_0} \rrbracket$ is the unique tangent cone to $T$ at $Z_0$.  By Corollary~\ref{doubling cor}, $N_{T,P_{Z_0},Z_0}(0^+) \geq 1+\alpha$ 
\begin{align}\label{persist cones eqn11}
    E(T,P_{Z_0},\mathbf{B}_{\sigma}(Z_0)) 
    \leq\,& E(\widetilde{T},P_{Z_0},\mathbf{C}_{\sigma}(Z_0)) 
    \leq C \Big(\frac{\sigma}{\rho}\Big)^{\alpha} E(\widetilde{T},P_{Z_0},\mathbf{C}_{\rho}(Z_0))
    \\ \leq\,& C \Big(\frac{\sigma}{\rho}\Big)^{\alpha} E(T,P_{Z_0},\mathbf{B}_{2\rho}(Z_0)) \nonumber 
\end{align}
for all $0 < \sigma \leq \rho \leq 1$, where $C = C(n, m, q, \epsilon, \beta)$ is a constant.  The inequality \eqref{branch and cones concl5} holds true by \eqref{persist cones eqn11} when $\sigma \leq \rho/2$ and is obvious when $\rho/2 < \sigma \leq \rho$.

Finally, given $Z_{1}, Z_{2} \in {\mathcal B} \cap {\mathbf B}_{1/2}(0)$, set $\sigma = |Z_{1} - Z_{2}|$ and argue as in \eqref{persist cones eqn3} to obtain, if $\sigma < 1/2$,  
\begin{align*}
	\|\pi_{P_{Z_{1}}} - \pi_{P_{Z_{2}}}\|^2 
	\leq\,& \frac{2}{q \omega_n \sigma^n} \int_{G_n(\mathbf{B}_{\sigma}(Z_{1}))} \|\pi_S - \pi_{P_{Z_{1}}}\|^2 \,d|T|(X,S) 
		\\&+ \frac{2}{q \omega_n \sigma^n} \int_{G_n(\mathbf{B}_{2\sigma}(Z_{2}))} \|\pi_S - \pi_{P_{Z_{2}}}\|^2 \,d|T|(X,S) \nonumber 
	\\ \leq\,& C E^{2}(T,P_{Z_{1}},\mathbf{B}_{2\sigma}(Z_{1})) + C E^{2}(T,P_{Z_{2}},\mathbf{B}_{4\sigma}(Z_{2})) \nonumber 
	\\ \leq\,& CC_{0}^{2} \sigma^{2\alpha} \left((E^{2}(T, P_{Z_{1}}, \mathbf{B}_{2}(Z_{1})) + E^{2}(T, P_{Z_{2}}, \mathbf{B}_{2}(Z_{2}))\right)\nonumber 
	\\=\,&  CC_{0}^{2} |Z_{1} - Z_{2}|^{2\alpha} \left((E^{2}(T, P_{Z_{1}}, \mathbf{B}_{2}(Z_{1})) + E^{2}(T, P_{Z_{2}}, \mathbf{B}_{2}(Z_{2}))\right) \nonumber
\end{align*}
where $C = C(n, q)$ and we have used \eqref{branch and cones concl5} in the third inequality. If on the other hand $\sigma \geq 1/2$ estimate again as above noting that  $\frac{3}{2}(\frac{1}{4} + \sigma) <2$ (since $\sigma <1$) to obtain 
\begin{align*}
	\|\pi_{P_{Z_{1}}} - \pi_{P_{Z_{2}}}\|^2 
	\leq\,& \frac{2(4^{n})}{q \omega_n} \int_{G_n(\mathbf{B}_{1/4}(Z_{1}))} \|\pi_S - \pi_{P_{Z_{1}}}\|^2 \,d|T|(X,S) 
		\\&+ \frac{2(5^{n})(1/4 + \sigma)^{-n}}{q \omega_n} \int_{G_n(\mathbf{B}_{1/4+ \sigma}(Z_{2}))} \|\pi_S - \pi_{P_{Z_{2}}}\|^2 \,d|T|(X,S) \nonumber 
	\\ \leq\,& C E^{2}(T,P_{Z_{1}},\mathbf{B}_{1/2}(Z_{1})) + C E^{2}(T,P_{Z_{2}},\mathbf{B}_{\frac{3}{2}(\frac{1}{4} + \sigma)}(Z_{2})) \nonumber 
	\\ \leq\,& C 4^{n+2} E^{2}(T, P_{Z_{1}}, \mathbf{B}_{2}(Z_{1})) + C \left(\frac{3}{4}\bigg(\frac{1}{4} + \sigma\bigg)\right)^{-n-2}E^{2}(T, P_{Z_{2}}, \mathbf{B}_{2}(Z_{2}))\nonumber 
	\\\leq\,& \overline{C}\left((E^{2}(T, P_{Z_{1}}, \mathbf{B}_{2}(Z_{1})) + E^{2}(T, P_{Z_{2}}, \mathbf{B}_{2}(Z_{2}))\right) \nonumber
	\\\leq\,& \overline{C}2^{2\alpha}|Z_{1} - Z_{2}|^{2\alpha} \left((E^{2}(T, P_{Z_{1}}, \mathbf{B}_{2}(Z_{1})) + E^{2}(T, P_{Z_{2}}, \mathbf{B}_{2}(Z_{2}))\right) \nonumber
\end{align*}
 where $\overline{C} = \overline{C}(n, m, q, \epsilon, \beta)$. This completes the proof of the theorem.
\end{proof}

 \bigskip
\hskip-.2in\vbox{\hsize3in\obeylines\parskip -1pt 
  \small 
Brian Krummel
School of Mathematics \& Statistics 
University of Melbourne
Parkville,VIC  3010, Australia
\vspace{4pt}
{\tt brian.krummel@unimelb.edu.au}} 
\vbox{\hsize3in
\obeylines 
\parskip-1pt 
\small 
Neshan Wickramasekera
DPMMS 
University of Cambridge 
Cambridge CB3 0WB, United Kingdom
\vspace{4pt}
{\tt N.Wickramasekera@dpmms.cam.ac.uk}
}


\begin{thebibliography}{DeLMinSko23}	
	\bibitem[All72]{Allard} 
		Allard, W.~K.~\textit{On the first variation of a varifold.} Ann.\ of Math.\ 95.3 (1972): 417-491.
		\bibitem[Alm66]{A1} Almgren, Jr.~F.~J.~\textit{Some interior regularity theorems for minimal surfaces and an extension of Bernstein's theorem.}
Ann.\ of Math.\ {\bf 84} (1966), 277--292.
	\bibitem[Alm83]{Almgren}
		Almgren, Jr.~F.~J.~\textit{Almgren's big regularity paper: Q-valued functions minimizing Dirichlet's integral and the regularity of area minimizing 
		rectifiable currents up to codimension two.}  World Scientific Monograph Series in Mathematics.  {\bf 1} 
		World Scientific Publishing Co.~Inc., River Edge, NJ (2000).
\bibitem[Cha88]{Chang} Chang, S.~X.~\textit{Two dimensional area minimizing integral currents are classical minimal surfaces.} Jour.~Amer.~Math.~Soc.~1(4), (1988): 699--778.
\bibitem[DeG61]{DG} De~Giorgi,~E.  {\it Frontiere orientate di misura minima} Sem. Mat. Scuola Norm. Sup. Pisa (1961), 1-56.
\bibitem[DMS23]{DelMinSko} De~Lellis,~C., Minter,~P.~and~Skorobogatova,~A. \textit{The fine structure of the singular set of area-minimizing integral currents III: frequency 1 flat singular points and ${\mathcal H}^{m-2}$ a.e.\ uniqueness of tangent cones.} arXiv:2304.11553 (2023).
\bibitem[DeLSko23-I]{DelSko1} De~Lellis,~C.~and~Skorobogatova,~A. \textit{The fine structure of the singular set of area-minimizing integral currents I: the singularity 
degree of flat singular points.} arXiv:2304.11552 (2023).
\bibitem[DeLSko23-II]{DelSko2}  De~Lellis,~C.~and~Skorabogatova,~A. \textit{The fine structure of the singular set of area-minimizing integral currents II: rectifiability of flat 
singular points of singularity degree larger than 1.} arXiv:2304.11555 (2023)
	\bibitem[DeLSpa11]{DeLSpaDirMin}
		De~Lellis,~C.~and~Spadaro,~E.  \textit{Almgren's Q-valued functions revisited.}  Mem.~Amer.~Math.~Soc.~211 (2011) no. 991, vi+79 pp.
	\bibitem[DeLSpa14]{DeLSpa1}
		De~Lellis,~C.~and~Spadaro,~E.  \textit{Regularity of area minimizing currents I: gradient L p estimates.}  GAFA 24.6 (2014): 1831--1884.
	\bibitem[DeLSpa15]{DeLSpa0}
		De~Lellis,~C.~and~Spadaro,~E.  \textit{Multiple valued functions and integral currents.}  Ann.~Sc.~Norm.~Super.~Pisa Cl.~Sci. (5) {\bf 14} (2015), 
	\bibitem[DeLSpa16-I]{DeLSpa2}
		De~Lellis,~C.~and~Spadaro,~E.  \textit{Regularity of area minimizing currents II: center manifold.}  Ann.\ of Math. (2016): 499--575.
	\bibitem[DeLSpa16-II]{DeLSpa3}
		De~Lellis,~C.~and~Spadaro,~E.  \textit{Regularity of area minimizing currents III: blow-up.}  Ann.\ of Math. (2016): 577--617.
	\bibitem[DMSV18]{DSMV18}
		De~Lellis,~C., Marchese,~A., Spadaro,~E.~and~Valtorta,~E.  \textit{Rectifiability and upper Minkowski bounds for singularities of harmonic
		$Q$-valued maps.}  Commentarii Mathematici Helvetici {\bf 93.4} (2018):737--779.
	\bibitem[Fed69]{Fed69}	
		Federer,~H.  \textit{Geometric measure theory.}  Springer, 2014.
		\bibitem[Fed70]{Fed70} Federer,~H. {\it The singular sets of area minimizing rectifiable currents with codimension one and area minimizing flat chains modulo two with arbitrary codimension} Bull. Amer. Math. Soc. 76 (1970), 767-771.
		\bibitem[Fle62]{FW} Fleming~W.~H. \textit{On the oriented Plateau problem.}
Rend.\ Circ.\ Mat.\ Palermo (2) {\bf 11} (1962), 69--90.
	\bibitem[HarSim79]{HardtSimon}
		Hardt,~R.~and~Simon,~L.  \textit{Boundary regularity and embedded solutions for the oriented Plateau problem.}  
		Annals of Math. {\bf 110} (1979): 439--486.
	\bibitem[KruWic-b]{KrumWicb}
		Krummel,~B.~and~Wickramasekera,~N. \emph{Analysis of singularities of area minimizing currents, Part~II: a uniform height bound, estimates away 
		from branch points of rapid decay, and uniqueness of tangent cones.} preprint. 
		\bibitem[KruWic-c]{KrumWicc}
		Krummel,~B.~and Wickramasekera,~N. \textit{Analysis of singularities of area minimizing currents, Part~III: branch points of planar frequency $\neq 2$, higher order asymptotics, and the local topology.} preprint. 
	\bibitem[KruWic-d]{KrumWicd}
		Krummel,~B.~and Wickramasekera,~N. \textit{Analysis of singularities of area minimizing currents, Part~IV: branch points of planar frequency $\geq 2,$ higher order asymptotics, and the local topology.} preprint.
    \bibitem[KruWic-e]{KrumWice} Krummel,~B.~and Wickramasekera,~N. \textit{Analysis of singularities of area minimizing currents, Part~V: general Riemannian ambient spaces.} preprint.  
	\bibitem[KruWic21]{KrumWic1}
		Krummel,~B.~and~Wickramasekera,~N. \emph{Fine proporties of branch point singularities: stationary graphs and stable minimal 
		hypersurfaces near points of density $< 3$.} arXiv:2111.12246v1 (2021). 
	\bibitem[KruWic17]{KrumWic2}
		Krummel,~B.~and~Wickramasekera,~N. \emph{Fine properties of branch point singularities: Dirichlet energy minimizing multi-valued 
		functions.}  arXiv:1711.06222 (2017).
	\bibitem[MicWhi95]{MicWhi95}
		Micallef,~M.~and White, B.  \emph{The structure of branch points in minimal surfaces and in pseudoholomorphic curves.}  
		Ann.\ of Math., {\bf 141} (1995), 35--85.	
	\bibitem[MinWic24]{MW}
	Minter,~P.~and~Wickramasekera,~N. \emph{A structure theory for stable codimension 1 integral varifolds with applications to area minimising hypersurfaces mod $p$}. Jour.~Amer.~Math.~Soc., 37(2024), 861--927. 	
	\bibitem[NabVal15]{NV15}
	Naber,~A.~and~Valtorta,~D. \emph{The singular structure and regularity of stationary and minimizing varifolds.} Jour.~Euro.~Math.~Soc.
	22(10) (2015), 3305-3382.
		\bibitem[Sim83]{SimonGMT}
		Simon,~L.  \textit{Lecture Notes on Geometric Measure Theory.} Proceedings of the Centre for Mathematical Analysis, 
		Australian National University, 3 (1983). 
	\bibitem[Sim93]{Sim93}
		Simon,~L.  \textit{Cylindrical tangent cones and the singular set of minimal submanifolds.}  J.~Diff.~Geo. (1993)  38:585--652.
	\bibitem[Sim95]{Sim95}
		Simon,~L.  \textit{Rectifiability of the singular sets of multiplicity 1 minimal surfaces and energy minimizing maps.}  Surveys in Diff.\ Geom.\ (1995), 246--305.
\bibitem[Sim96]{Sim96}
		Simon,~L.  \emph{Theorems on regularity and singularity of energy minimizing map.}  Springer Science \& Business Media (1996).
		\bibitem[JSim68]{Simons} Simons,~J.  \textit{Minimal varieties in Riemannian manifolds} Ann. of Math. 88 (1968), 62-105.
	\bibitem[Whi83]{Whi83} White,~B.
\textit{Tangent cones to 2-dimensional area-minimizing integral currents are unique.} Duke Math.\ J.\ {\bf 50} (1983), 143--160. 
\bibitem[Wic08]{Wic08} Wickramasekera,~N. \emph{A regularity and compactness theory for immersed stable minimal hypersurfaces of multiplicity at most 2.} J.~Diff.~Geo.~80.1 (2008): 79--173.
	\bibitem[Wic14]{Wic14}
		Wickramasekera,~N. \emph{A general regularity theory for stable codimension 1 integral varifolds.} Ann.\ of Math. (2014): 843--1007.
\end{thebibliography}
\end{document}